\documentclass[openany]{memo-l}

\usepackage{amsthm,amscd,amssymb}
\usepackage[all]{xy}

\newtheorem*{Beil}{Beilinson's theorem}
\newtheorem*{CatSch}{Theorem (Catanese--Schreyer)}
\newtheorem{theo}{Theorem}[section]
\newtheorem{lemm}[theo]{Lemma}
\newtheorem{coro}[theo]{Corollary}
\newtheorem{prop}[theo]{Proposition}

\theoremstyle{definition}

\newtheorem{defi}[theo]{Definition}
\newtheorem{exem}[theo]{Example}

\theoremstyle{remark}

\newtheorem{rema}[theo]{Remark}

\numberwithin{section}{chapter}

\numberwithin{equation}{section}

\newcommand{\good}{good}

\newcommand{\triv}{trivial}
\newcommand{\qtriv}{quasi--trivial}
\newcommand{\std}{standard}
\newcommand{\qstd}{quasi--standard}
\newcommand{\gbwcp}{good birational weighted canonical projection}
\newcommand{\nzdiv}{non--zerodivisor}

\newcommand{\iso}{\cong}

\newcommand{\all}{\forall}
\newcommand{\exi}{\exists}
\newcommand{\nexi}{\nexists}
\newcommand{\exiun}{\exists!}
\newcommand{\id}{{\rm id}}
\newcommand{\minus}{\setminus}
\newcommand{\st}{\,|\,} 
\newcommand{\dual}{^{\vee}} 
\newcommand{\tdual}{^*} 
\newcommand{\dtdual}{^{**}} 
\newcommand{\depth}{{\rm depth}}
\newcommand{\rk}[1]{{\rm rk}(#1)}
\newcommand{\djun}{\sqcup} 
\newcommand{\mono}{\hookrightarrow} 
\newcommand{\epi}{\twoheadrightarrow} 
\newcommand{\mor}[1]{\xrightarrow{#1}} 
\newcommand{\isomor}{\mor{\sim}} 
\newcommand{\rat}{\dashrightarrow} 
\newcommand{\comp}{\circ} 
\newcommand{\card}[1]{\##1} 
\newcommand{\abs}[1]{\lvert#1\rvert} 
\newcommand{\sw}[1]{\lvert#1\rvert} 
\newcommand{\rest}[1]{|_{#1}} 
\newcommand{\linsys}[1]{\lvert#1\rvert} 
\newcommand{\cp}[1]{#1^{\bullet}} 
\newcommand{\ring}[1]{\mathbb{#1}}
\newcommand{\K}{\ring{K}} 
\newcommand{\C}{\ring{C}} 
\newcommand{\N}{\ring{N}} 
\newcommand{\Z}{\ring{Z}} 
\newcommand{\I}[1]{\mathtt{#1}} 
\newcommand{\cat}[1]{\begin{bf}#1\end{bf}} 
\newcommand{\opp}{^{\circ}} 
\newcommand{\Pic}{{\rm Pic}} 
\newcommand{\mrs}{^{\#}} 
\newcommand{\s}[1]{\mathcal{#1}} 
\newcommand{\g}[1]{\overline{#1}} 
\newcommand{\ung}[1]{\hat{#1}} 
\renewcommand{\hom}{^{\rm hom}} 
\newcommand{\Hom}{{\rm Hom}} 
\newcommand{\Ext}{{\rm Ext}} 
\newcommand{\gHom}{\g{\Hom}} 
\newcommand{\gExt}{\g{\Ext}} 
\newcommand{\sHom}{\s{H}om} 
\newcommand{\sExt}{\s{E}xt} 
\newcommand{\gsHom}{\g{\sHom}} 
\newcommand{\gsExt}{\g{\sExt}} 
\newcommand{\gsec}{\g{\Gamma}} 
\DeclareMathOperator{\im}{im} 
\DeclareMathOperator{\cok}{coker} 
\DeclareMathOperator{\spec}{Spec} 
\DeclareMathOperator{\Spec}{\bf{Spec}} 
\DeclareMathOperator{\proj}{Proj} 
\DeclareMathOperator{\gspec}{\g{Spec}} 
\DeclareMathOperator{\gSpec}{\g{\bf{Spec}}} 
\DeclareMathOperator{\gproj}{\g{Proj}} 
\newcommand{\mult}{\mathbb{G}_m} 
\newcommand{\roots}[1]{\mu_{#1}} 
\newcommand{\p}{{\rm P}} 
\newcommand{\pd}{{\rm p}} 
\renewcommand{\P}{\mathbb{P}} 
\newcommand{\gP}{\g{\P}} 
\newcommand{\sh}[1]{^{\sim_{#1}}} 
\newcommand{\gsh}[1]{^{\g{\sim}_{#1}}} 
\newcommand{\mo}[1]{\cat{Mod}(#1)} 
\newcommand{\fmo}[1]{\cat{mod}(#1)} 
\newcommand{\smo}[1]{\cat{Mod}(#1)} 
\newcommand{\gmo}[1]{\g{\cat{Mod}}(#1)} 
\newcommand{\gfmo}[1]{\g{\cat{mod}}(#1)} 
\newcommand{\gsmo}[1]{\g{\cat{Mod}}(#1)} 
\newcommand{\coh}[1]{\cat{Coh}(#1)} 
\newcommand{\qco}[1]{\cat{Qcoh}(#1)} 
\newcommand{\gcoh}[1]{\g{\cat{Coh}}(#1)} 
\newcommand{\gqco}[1]{\g{\cat{Qcoh}}(#1)} 
\newcommand{\prng}{\g{\cat{PRing}}} 
\newcommand{\gsch}{\g{\cat{Sch}}} 
\newcommand{\stdsch}[2]{\g{#1(#2)}} 
\newcommand{\umf}[1]{-_0} 
\newcommand{\omf}[1]{#1\otimes_{#1_0}-} 
\newcommand{\hmf}[1]{\gHom_{#1_0}(#1,-)} 
\newcommand{\usmf}[1]{-_0} 
\newcommand{\osmf}[1]{\s{O}_{#1}\otimes_{\s{O}_{#1_0}}-} 
\newcommand{\hsmf}[1]{\gsHom_{#1_0}(\s{O}_{#1},-)} 
\newcommand{\V}[1]{V(#1)} 
\newcommand{\D}[1]{D(#1)} 
\newcommand{\QS}[1]{QS(#1)} 
\newcommand{\QT}[1]{QT(#1)} 
\newcommand{\compl}[1]{#1^c} 
\newcommand{\w}{{\rm w}} 
\newcommand{\lotimes}{\overset{L}{\otimes}} 
\newcommand{\md}{\Omega} 
\newcommand{\smd}{Syz} 
\newcommand{\so}{\s{O}} 
\newcommand{\sd}{\Omega} 
\newcommand{\diff}[1]{d_{#1}} 
\newcommand{\cdiff}[1]{d_{\cp{#1}}} 
\newcommand{\ko}{K} 
\newcommand{\sko}{\s{\ko}} 
\newcommand{\esd}[1]{{\mathtt S}_{#1}} 
\newcommand{\omd}[1]{\alpha^{#1}} 
\newcommand{\dom}[1]{\beta^{#1}} 
\newcommand{\ind}[1]{_{(#1)}} 
\newcommand{\bo}[1]{\s{M}\ind{#1}} 
\newcommand{\bd}[1]{\s{N}\ind{#1}} 
\newcommand{\abo}[1]{\tilde{\s{M}}_{(#1)}} 
\newcommand{\abd}[1]{\tilde{\s{N}}_{(#1)}} 
\newcommand{\bro}{\s{X}} 
\newcommand{\brd}{\s{Y}} 
\newcommand{\gbro}{\g{\bro}} 
\newcommand{\gbrd}{\g{\brd}} 
\newcommand{\oset}{O} 
\newcommand{\dset}{D} 
\newcommand{\goset}{\g{\oset}} 
\newcommand{\gdset}{\g{\dset}} 
\newcommand{\cbrd}{y} 
\newcommand{\mco}[1]{\cat{M}_{#1}} 
\newcommand{\mcom}[1]{\cat{\widehat{M}}_{#1}} 
\newcommand{\co}[1]{\cat{O}_{#1}} 
\newcommand{\com}[1]{\cat{\widehat{O}}_{#1}} 
\newcommand{\ohom}[1]{\Hom^{\co{#1}}} 
\newcommand{\coq}[1]{\cat{Q}_{#1}} 
\newcommand{\into}{]-\sw{\w},0]} 
\newcommand{\intd}{]n-\sw{\w},0[} 
\newcommand{\fo}{F_{\p}} 
\newcommand{\fd}{G_{\p}} 
\newcommand{\fde}{F_{\Lambda}} 
\newcommand{\gfo}{\g{F}_{\p}} 
\newcommand{\gfd}{\g{G}_{\p}} 
\newcommand{\tri}[6]{#1\mor{#4}#2\mor{#5}#3\mor{#6}#1[1]} 
\newcommand{\MC}[1]{MC(#1)} 
\newcommand{\mc}[1]{MC(\cp{#1})} 
\newcommand{\inc}[1]{j(\cp{#1})} 
\newcommand{\pro}[1]{p(\cp{#1})} 
\newcommand{\tbo}[2]{\delta_{#1,#2}} 
\newcommand{\R}{\mathcal{R}} 
\renewcommand{\S}{{\rm S}_0} 
\newcommand{\gS}{{\rm S}} 
\newcommand{\X}{{\rm X}_0} 
\newcommand{\gX}{{\rm X}} 
\newcommand{\Y}{{\rm Y}_0} 
\newcommand{\gY}{{\rm Y}} 
\newcommand{\ds}{\omega} 
\newcommand{\ph}{\phi_0} 
\newcommand{\gph}{\phi} 
\newcommand{\ps}{\psi_0} 
\newcommand{\gps}{\psi} 
\newcommand{\pr}{\pi_0} 
\newcommand{\gpr}{\pi} 
\newcommand{\rh}{\varrho} 
\newcommand{\f}{{\rm f}} 
\newcommand{\cre}{c} 
\newcommand{\ep}{\epsilon}
\newcommand{\rel}[1]{rel'_{#1}} 
\newcommand{\crel}[1]{rel_{#1}} 
\newcommand{\cov}{\mathbf{U}} 

\begin{document}

\frontmatter
\title{The Beilinson complex and canonical rings of irregular surfaces}

\author{Alberto Canonaco}

\address{Dipartimento di Matematica ``F. Casorati''\\
Universit\`a di Pavia\\
Via Ferrata 1\\
27100 Pavia\\
Italy}


\email{alberto.canonaco@unipv.it}


\date{February 24, 2003}

\subjclass[2000]{Primary $14$A$20$, $14$M$99$, $14$F$05$, $14$J$29$; Secondary
$13$A$02$, $18$E$30$, $14$K$05$}


\begin{abstract}
An important theorem by Beilinson describes the bounded derived category of
coherent sheaves on $\mathbb{P}^n$, yielding in particular a resolution of
every coherent sheaf on $\mathbb{P}^n$ in terms of the vector bundles
$\Omega_{\mathbb{P}^n}^j(j)$ for $0\le j\le n$. This theorem is here extended
to weighted projective spaces. To this purpose we consider, instead of the
usual category of coherent sheaves on $\mathbb{P}({\rm w})$ (the weighted
projective space of weights $\rm w=({\rm w}_0,\dots,{\rm w}_n)$), a suitable
category of graded coherent sheaves (the two categories are equivalent if and
only if ${\rm w}_0=\cdots={\rm w}_n=1$, i.e. $\mathbb{P}({\rm w})=
\mathbb{P}^n$), obtained by endowing $\mathbb{P}({\rm w})$ with a natural
graded structure sheaf. The resulting graded ringed space
$\overline{\mathbb{P}}({\rm w})$ is an example of {\em graded scheme} (in
chapter $1$ graded schemes are defined and studied in some greater generality
than is needed in the rest of the work). Then in chapter $2$ we prove for
graded coherent sheaves on $\overline{\mathbb{P}}({\rm w})$ a result which is
very similar to Beilinson's theorem on $\mathbb{P}^n$, with the main difference
that the resolution involves, besides $\Omega_{\overline{\mathbb{P}}({\rm
w})}^j(j)$ for $0\le j\le n$, also $\mathcal{O}_{\overline{\mathbb{P}}({\rm
w})}(l)$ for $n-\sum_{i=0}^n{\rm w}_i<l<0$.

This weighted version of Beilinson's theorem is then applied in chapter $3$ to
prove a structure theorem for {\em good birational weighted canonical
projections} of surfaces of general type (i.e., for morphisms, which are
birational onto the image, from a minimal surface of general type $S$ into a
$3$--dimensional $\mathbb{P}({\rm w})$, induced by $4$ sections $\sigma_i\in
H^0(S,\mathcal{O}_S({\rm w}_iK_S))$). This is a generalization of a theorem by
Catanese and Schreyer (who treated the case of projections into
$\mathbb{P}^3$), and is mainly interesting for irregular surfaces, since in the
regular case a similar but simpler result (due to Catanese) was already known.
The theorem essentially states that giving a good birational weighted canonical
projection is equivalent to giving a symmetric morphism of (graded) vector
bundles on $\overline{\mathbb{P}}({\rm w})$, satisfying some suitable
conditions. Such a morphism is then explicitly determined in chapter $4$ for a
family of surfaces with numerical invariants $p_g=q=2$, $K^2=4$, projected into
$\mathbb{P}(1,1,2,3)$.
\end{abstract}

\maketitle

\setcounter{page}{4}

\tableofcontents

\mainmatter

            \chapter*{Introduction}

A classical method to study algebraic surfaces of general type is to consider
the behaviour of pluricanonical maps (see e.g. the important paper \cite{Bo}).
If $\S$ is a (smooth, projective) minimal surface of general type (over an
algebraically closed field $\K$), denoting by $K_{\S}$ the canonical divisor on
$\S$ and by $\ds_{\S}=\so_{\S}(K_{\S})$ the canonical sheaf, the $m^{th}$
canonical map is the rational map $\varphi_{\linsys{mK_{\S}}}:\S\rat\P^N$
determined by the linear system $\linsys{mK_{\S}}$. More generally, one can
consider a map into a weighted projective space, induced by sections of
different multiples of $K_{\S}$. More precisely, given positive numbers
$\w=(\w_0,\dots,\w_n)$, we denote by $\P(\w)$ (or simply by $\P$) the weighted
projective space of weights $\w$, which is defined as $\proj\p$, where
$\p=\p(\w)$ is the polynomial ring $\K[x_0,\dots, x_n]$ graded by
$\deg(x_i)=\w_i$. Then the choice of sections $\sigma_i\in
H^0(\S,\ds_{\S}^{\w_i})$ determines a rational map $\ph:\S\rat\P$, called a
{\em weighted canonical projection}. $\ph$ is called {\em good birational} if
it is a morphism which is birational onto its image $\Y\subset\P$. Such a map
always factors through the canonical model $\X$ of $\S$ (which is defined as
$\proj\R$, where $\R:=\bigoplus_{d\ge0}H^0(\S,\ds_{\S}^d)$ is the canonical
ring of $\S$):
\[\xymatrix{\S \ar[rr]^{\ph}\ar[dr]_{\pr} & & \Y\subset\P \\
& \X \ar[ur]_{\ps}}\]
$\X$ is a surface with only rational double points as singularities (the
birational morphism $\pr$ contracts $(-2)$--curves of $\S$), which essentially
behaves like $\S$. The advantage of considering $\ps$ instead of $\ph$ is that
it is a finite morphism, induced by the morphism of graded rings
\[\begin{split}
\rh:\p(\w) & \to \R \\
x_i & \mapsto \sigma_i
\end{split}\]
It follows from known results on pluricanonical maps that every surface of
general type admits a {\gbwcp} to a $3$--dimensional weighted projective space
with weights $\w_0,\w_1\le2$, $\w_2\le3$ and $\w_3\le5$.

So let's assume that $\ph:\S\to\P$ is a {\gbwcp} with $\dim\P=3$. If $\S$ is
regular (i.e. $q(\S):=h^1(\S,\so_{\S})=0$) good results can be obtained by
looking for (graded) free resolutions of the canonical ring $\R$, regarded as a
graded $\p$--module via the morphism $\rh$ defined before. In \cite{C1} it is
proved that $\R$ always admits a length $1$ resolution of the form
\[0\to\bigoplus_{j=1}^h\p(l_j-\sw{\w}-1)\mor{\alpha}\bigoplus_{i=1}^h\p(-l_i)
\to\R\to0\]
(with $\sw{\w}:=\sum_{i=0}^3\w_i$) where $\alpha$ is a minimal symmetric matrix
of homogeneous polynomials $\alpha_{i,j}\in\p_{\sw{\w}+1-l_j-l_i}$ for some
integers $0=l_1<l_2\le\cdots\le l_h$; $\alpha$ being minimal means that
$\alpha_{i,j}=0$ if $l_i+l_j>\sw{\w}$. Moreover, $\alpha$ satisfies the
following properties:
\begin{enumerate}

\item $\f:=\det(\alpha)$ is an irreducible homogeneous polynomial defining
$\Y=\proj\p/(\f)\subset\P$;

\item (rank condition) denoting by $\beta$ the matrix of cofactors of $\alpha$,
\[(\beta_{1,i}\st i=1,\dots,h)=(\beta_{i,j}\st i,j=1,\dots,h)\]
as ideals of $\p$.
\end{enumerate}

This last condition is seen to be equivalent to the fact that $\cok\alpha$
carries a natural structure of $\p/(\f)$--algebra (``rank condition =
ring condition''). Then the important fact is that, conversely, if $\alpha$ is
a minimal symmetric matrix satisfying the two conditions above (so that
$\cok\alpha$ is a ring) and if moreover $\proj(\cok\alpha)$ is a surface
with only rational double points as singularities, then $\cok\alpha$ is the
canonical ring of a surface of general type.

Using this result Catanese could prove interesting
applications to moduli problems: for instance, he showed that the surfaces with
numerical invariants $K^2=7$, $p_g=4$ (and $q=0$, of course) such that the
canonical map is base point free, have a unirational irreducible moduli space.

The reason why this method works well only for regular surfaces is that only
in this case is $\R$ a Gorenstein ring (hence it has a resolution of length $1$
as graded $\p$--module). On the other hand, a closely related
approach is possible also in the case of irregular surfaces: the idea is to
look for locally free resolutions of the coherent sheaf $\R\sh{}\iso
{\ph}_*\so_{\S}$ on $\P$. In particular, for projections into $\P^3$ the
following theorem is proved in \cite{CS} (see also \cite{C2}).

        \begin{CatSch}
Let $\S$ be a minimal surface of general type with $p_g(\S)=p_g$, $q(\S)=q$ and
$K_{\S}^2=K^2$, let $\ph:\S\to\Y\subset\P^3$ be  a good birational canonical
projection, and consider the vector bundle
\[\s{E}:=\so(-2)^{p_g-q+K^2-9}\oplus\sd^1(-1)^q\oplus(\sd^2)^{p_g-4}\]
on $\P^3$. Then there exists an exact sequence of coherent sheaves on $\P^3$
\[0\to(\so\oplus\s{E})\dual(-5)\mor{\alpha=\begin{pmatrix}
\alpha^{(1)} \\
\alpha'
\end{pmatrix}}
\so\oplus\s{E}\to{\ph}_*\so_{\S}\to0,\]
such that $\alpha$ is a minimal morphism which satisfies the following
properties:
\begin{enumerate}

\item $\alpha$ is symmetric, i.e. $\alpha=\alpha\dual(-5)$;

\item $\f=\det(\alpha)$ is an irreducible polynomial (defining $\Y=
\proj\p/(\f)$);

\item $\s{I}_r(\alpha)=\s{I}_r(\alpha')$ (where $r:=\rk{\s{E}}$ and, if $\beta:
\s{E}_1\to\s{E}_2$ is a morphism of vector bundles, $\s{I}_k(\beta)$ denotes
the sheaf of ideals defined by the image of the natural map
$\Lambda^k(\s{E}_1)\otimes\Lambda^k(\s{E}_2)\dual\to\so_{\P^3}$ induced by
$\Lambda^k(\beta)$);

\item $\X:=\Spec(\cok{\alpha})$ is a surface with only rational double points
as singularities.
\end{enumerate}

Conversely, if $\alpha:(\so\oplus\s{E})\dual(-5)\to\so\oplus\s{E}$ is minimal
and satisfies properties $1$, $2$, $3$, then $\cok\alpha$ is in a natural way a
coherent sheaf of commutative $\so_{\Y}$--algebras, and if $4$ also holds,
then, denoting by $\pr:\S\to\X$ a minimal resolution of singularities of $\X$,
$\S$ is a minimal surface of general type with $p_g(\S)=p_g$, $q(\S)=q$ and
$K_{\S}^2=K^2$, $\X$ is its canonical model and $\ph:\S\mor{\pr}\X\to\Y$ is a
good birational canonical projection.
        \end{CatSch}

The main purpose of this work is to extend the above result to the weighted
case. Before we can say how this is done, we have to discuss briefly the main
elements in the proof of the above theorem. The essential ingredient is a
famous theorem by Beilinson (\cite{B1}), which (in its more concrete version)
gives, for every coherent sheaf $\s{F}$ on $\P^n$, a two--sided resolution
(i.e., a complex $\cp{\brd}=\cp{\brd}(\s{F})$ whose cohomology is $\s{F}$ in
position $0$ and $0$ elsewhere) in terms of the vector bundles of twisted
differential forms $\sd_{\P^n}^j(j)$. One can always take a minimal such
resolution, and then it is explicitly given by
\[\brd^i=\bigoplus_{0\le j\le n}\sd^j(j)\otimes_{\K}
H^{i+j}(\P^n,\s{F}(-j)).\] In our case Beilinson's theorem is applied to the
coherent sheaf $\s{F}:= {\ph}_*\so_{\S}(2)$ on $\P^3$: the complex
$\cp{\brd}(\s{F})(-2)$ is then a resolution of ${\ph}_*\so_{\S}$, which is
easily seen to be quasi--isomorphic to one of the form $0\to
(\so\oplus\s{E})\dual(-5)\mor{\alpha}\so\oplus\s{E}\to0$ (with $\alpha$ not
necessarily symmetric). The fact that $\alpha$ can be actually chosen to be
symmetric is a consequence of the birationality of $\ph$. Then the rest of the
proof is rather straightforward, using a result which says that (under quite
general assumptions) the condition $\s{I}_r(\alpha)= \s{I}_r(\alpha')$ (which
is just the sheafified version of the rank condition already discussed for
modules) is equivalent to the fact that $\cok\alpha$ has a natural structure of
a sheaf of commutative $\so_{\Y}$--algebras .

The principal difficulties in trying to extend both theorems to the weighted
case reside in the fact that weighted projective spaces are singular varieties,
presenting several pathologies: in particular, the sheaves of the form $\so(i)$
or $\sd^j(i)$ are not locally free in general. The technique used to overcome
such problems is to consider, instead of the usual category of sheaves of
modules on $\P$, a category of graded sheaves, which we are now going to
describe, and which turns out to be much better behaved (a similar method is
used in \cite{B} and \cite{GL}, where however a different graduation is used).

Let $R=\bigoplus_{d\ge0}R_d$ be a noetherian positively graded ring ($R=\p$ and
$R=\R$ are the examples we are interested in): one can consider a ``graded
scheme'' $\gproj R$ (we will write $\gP$ for $\gproj\p$), which is equal to
$\proj R$ as a topological space, but which has a $\Z$--graded structure sheaf
\[\so_{\gproj R}:=\bigoplus_{d\in\Z}\so_{\proj R}(d).\]
In other words, $\gproj R$ is defined as $\proj R$, except that in the
definition of the structure sheaf one takes as stalk at a point $\I{p}$
($\I{p}$ is a homogeneous, non maximal prime ideal of $R$) the whole local
$\Z$--graded ring $R_{\I{p}}:=S^{-1}R$ (where $S$ denotes the set of
homogeneous elements of $R\minus\I{p}$), instead of its degree $0$ part
$R_{(\I{p})}$. It is clear that there is a natural notion of graded sheaf of
modules on $\gproj R$; in particular, if $M$ is a graded $R$--module, one
can consider the associated graded sheaf $M\gsh{}:=\bigoplus_{d\ge0}M(d)\sh{}$
(as before, the stalk at $\I{p}$ of $M\gsh{}$ is $M_{\I{p}}$ instead of
$M_{(\I{p})}$). Denoting $\gproj R$ by $Z$ and $\proj R$ by $Z_0$, we will
use the (abelian) category $\gsmo{Z}$ of graded sheaves of $\so_Z$--modules,
instead of the usual category $\smo{Z_0}$ of sheaves of $\so_{Z_0}$--modules
(similar considerations can be made for the subcategories of quasi--coherent
and coherent sheaves, of which natural definitions can be given also in the
graded case). The two categories are in any case closely related: sending a
graded sheaf $\s{F}$ to $\s{F}_0$ clearly defines an exact functor
\[\usmf{Z}:\gsmo{Z}\to\smo{Z_0}.\]
Moreover, it is easy to see that there exists (not unique) a functor
\[F:\smo{Z_0}\to\gsmo{Z}\]
such that $\usmf{Z}\comp F\iso\id$. This is important for us, because if
$\s{G}$ is a sheaf on $Z_0$ and $\cp{\s{C}}$ is a resolution of $F(\s{G})$,
then $\cp{\s{C}}_0$ will be a resolution of $\s{G}$. One can also prove that
$\usmf{Z}$ is an equivalence if and only if $\all\I{p}\in\gproj R$ the graded
ring $R_{\I{p}}$ contains an invertible element of degree $1$. This is always
the case if $R=\R$ (the canonical ring of a surface of general type), whereas
if $R=\p(\w)$ this condition is satisfied if and only if $\w=(1,\dots,1)$.

As we mentioned before, the category $\gsmo{Z}$ is usually nicer than
$\smo{Z_0}$, mainly because in the graded case there is a natural way to define
twist functors: if $\s{F}$ is a graded sheaf and $d\in\Z$, one sets, as usual,
$\s{F}(d)_i:=\s{F}_{d+i}$. These functors satisfy all the good properties one
expects: they are exact, $(d)\comp(d')=(d+d')$, $(d)$ is isomorphic to
$-\otimes_{\so_Z}\so_Z(d)$, if $M$ is a graded $R$--module then $M(d)\gsh{}\iso
M\gsh{}(d)$, and so on. On the contrary, there is not a good notion of twist
functors on $Z_0$ in general (for instance when $Z_0=\P(\w)$ is not isomorphic
to a $\P^n$): one can define $(d)$ as $-\otimes_{\so_{Z_0}}\so_{Z_0}(d)$ or as
$\sHom_{\so_{Z_0}}(\so_{Z_0}(-d),-)$ (which is usually better), but in either
case they are not well behaved. Moreover, the sheaves $\so_Z(j)$ are always
locally free by definition.

It was pointed out by the referee that (at least when $R_0$ is a
field) our category of graded (quasi)coherent sheaves on $Z$ is
equivalent to the category considered in \cite{AKO} (which appeared on
the web after the submission of the present paper) of (quasi)coherent
sheaves on the quotient algebraic stack of $\spec R\minus\{R_+\}$ by the action
of the multiplicative group $\mult$ induced by the graduation of $R$.

The use of graded sheaves allows us to prove a version of Beilinson's theorem
on a (graded) weighted projective space $\gP=\gP(\w)$, which is much better
than the one of \cite{Ca1}. As in the case of $\P^n$, we prove it first in an
abstract form (as an equivalence between the bounded derived category of graded
coherent sheaves and a suitable homotopy category of graded modules, see
\ref{wbeilthm}; similar results are also proved in \cite{B} and
\cite{AKO}) and then give the explicit form of the minimal resolution of
every (graded coherent) sheaf, or, more generally, of every bounded complex of
sheaves (\ref{wbeilfor}). Such a resolution involves, besides the sheaves
$\sd_{\gP}^j(j)$ (which can be defined in a natural way, and which are locally
free), also the sheaves $\so_{\gP}(j)$ for $n-\sw{\w}<j<0$. While the
coefficients of the $\sd^j(j)$ in the resolution of a (graded coherent) sheaf
$\s{F}$ are again given by the dimensions of cohomology groups of twists of
$\s{F}$ (the same as in the case of $\P^n$), the coefficients of the $\so(j)$
are more complicated, since they are expressed as the dimensions of  cohomology
groups of $\s{F}$ tensored with suitably defined complexes of sheaves.

This weighted version of Beilinson's theorem can  be usefully applied to extend
the theorem by Catanese and Schreyer to the case of weighted projections.
Before one can do that, however, it is necessary to generalize somehow the
graded theory we have considered so far, mainly because we need an analogue, in
the graded case, of the definition of $\Spec$ of a sheaf of algebras on $\P^3$.
Without entering into details, we will just say that we give a general
definition of graded scheme (which includes both the $\gproj R$ defined before
and usual schemes as particular cases), such that every graded scheme is
locally isomorphic to what we denote by $\gspec$ of a $\Z$--graded ring $A$ (as
a set, it consists of homogeneous prime ideals of $A$, and it is just $\spec$
if $A=A_0$). Then one can define, as usual, a sheafified version $\gSpec$ of
$\gspec$ for graded sheaves of algebras on a graded scheme. We need moreover to
extend also the surfaces $\X$, $\Y$ and $\S$ to graded schemes $\gX$, $\gY$ and
$\gS$: this is done by setting (obviously) $\gX:=\gproj\R$,
$\gY:=\gproj\p/(\f)$, whereas $\gS$ is equal to $\S$ as a topological space and
has structure sheaf $\so_{\gS}:=\bigoplus_{d\in\Z}\ds_{\S}^d$. Also the maps
$\pr$, $\ph$ and $\ps$ extend naturally to morphisms of graded schemes $\gpr$,
$\gph$ and $\gps$. As we have already said, $\gX$ has the property that the
categories of graded sheaves and of usual sheaves are equivalent, and the same
is true for $\gS$. Our main theorem (\ref{mthm}) is then quite similar to that
of Catanese and Schreyer, except that everything is graded (of course, if $\w=
(1,1,1,1)$, our version coincides with theirs, modulo the equivalences of
categories mentioned before). The general case is however more complicated,
since the vector bundle $\s{E}$ now contains some extra terms of the form
$\so_{\gP}(j)$ (coming from the Beilinson resolution), whose coefficients are
no more uniquely determined by the numerical invariants of the surface, but
depend on the projection $\gph$. Apart from that, everything is very similar to
the not weighted case, also in the proof (where, again, the point which
requires more care is the explicit determination of $\s{E}=\s{E}(\gph)$).

It must be said that our result, though interesting from a theoretical point of
view (since it provides, in principle, a description as locally closed
algebraic sets of the moduli spaces of surfaces of general type together with a
{\gbwcp), doesn't seem to be of great practical utility. The problem is that
for irregular surfaces (as we are interested in) the vector bundle $\s{E}$ has
always a high rank, and then the conditions to impose on the morphism $\alpha$
(in particular the rank condition) become almost intractable. However, to
illustrate our result we explicitly compute the resolution for a class of
surfaces with $p_g=q=2$ and $K^2=4$, namely the double covers of a principally
polarized abelian surface $(A,\Theta)$, branched on a divisor of
$\linsys{2\Theta}$. In order to do that, we need (among other things) to
determine the structure of the ring of theta functions of an abelian surface.

As for the contents of the single chapters (hopefully clear from the table of
contents and from what we said) we just want to say that chapter $1$ contains
also some material on graded schemes which is not needed in the sequel, and
which has been included for the sake of completeness; a short
description of the relation between graded schemes and algebraic
stacks is also present. The appendix mainly contains the essential
facts about derived categories and derived functors; the notation
introduced there is freely used in the text, especially from chapter $2$.

{\bf Acknowledgements.} This work originates from my Ph. D. thesis. It is a
pleasure for me to thank here my thesis advisor Fabrizio Catanese, for posing
the problem and helping me in several ways to solve it in the course of the
years. Thanks also to Miles Reid and Frank--Olaf Schreyer for fruitful
discussions I had with them and to Sorin Popescu for pointing out the existence
of \cite{B}.

                        \chapter{Graded schemes}

By the word graded we will always mean $\Z$--graded.\index{graded} If $X$ is a
graded object (ring, module, sheaf,\dots) we will write it as
$\bigoplus_{d\in\Z}X_d$, where $X_d$ is the component of degree $d$. We will
denote by $X\hom$\index{$\hom$} the set of homogeneous elements of $X$.

As we said in the introduction, the main motivation for this chapter came from
the necessity of defining, given a noetherian positively graded ring $R$, the
``graded scheme'' $\gproj R$ (by which we mean a topological space equal to
$\proj R$ and a graded structure sheaf given by $\so_{\gproj R}=
\bigoplus_{d\in\Z}\so_{\proj R}(d)$). On the other hand, since such an object
is not a scheme in the usual meaning of the word, we need to give a definition
of graded scheme. This can be done in a quite natural way (at least with
noetherian hypotheses, which we put, for simplicity, in the definition of
graded scheme): first of all, a graded scheme will be a locally graded ringed
space $(X,\so_X)$, i.e. a (noetherian) topological space $X$, together with a
structure sheaf of graded rings $\so_X$ such that every stalk $\so_{X,x}$ is a
local graded ring (a graded ring with a unique maximal homogeneous ideal).
Moreover, we require $X$ to be locally isomorphic to $\gspec A$ for some
noetherian graded ring $A$, where we denote by $\gspec A$ the following locally
graded ringed space: as a topological space, $\gspec A$ is given by the
homogeneous prime ideals of $A$ (with the obvious topology), and the graded
structure sheaf $\so_{\gspec A}$ is defined in such a way that the stalk at a
point $\I{p}\in\gspec A$ is given by the localization $A_{\I{p}}:=
(A\hom\minus\I{p})^{-1}A$.

With this definition $\gproj R$ is indeed a graded scheme: it is locally
isomorphic to $\gspec R_r$ for some $r\in R_+\hom$ ($R_r:=\{r^n\st
n\in\N\}^{-1}R$). Notice also that graded schemes include usual (noetherian)
schemes. It is clear that one can extend many aspects of the theory of schemes
to graded schemes, but here we are not much interested in that. Again
influenced by the example of $\gproj R$, where our aim is to compare the
category of graded sheaves on $\gproj R$ with the category of sheaves on $\proj
R$, we are rather interested in those graded schemes where a similar comparison
is possible, i.e. graded schemes $X$ such that $X_0:=(X,(\so_X)_0)$ is a
scheme. This is certainly the case if $X$ is locally isomorphic to $\gspec A$,
with $A$ a graded ring such that $\gspec A$ and $\spec A_0$ are homeomorphic:
we call such schemes (and such rings) \good. Examples of {\good} rings are
graded rings containing an invertible element of positive degree (such as $R_r$
above), which we call {\qstd} (we call \std, instead, those containing an
invertible element of degree $1$), and graded rings such that every element of
degree $\ne0$ is nilpotent, which we call {\qtriv} (we call \triv, instead,
those concentrated in degree $0$). Then we prove that a {\good} noetherian
local graded ring is either {\qstd} or \qtriv, so that a connected {\good}
graded scheme is also either {\qstd} or \qtriv.

As we said before, given a {\good} graded scheme $X$ we can compare the
category $\gsmo{X}$ of graded sheaves of $\so_X$--modules and the category
$\smo{X_0}$ of sheaves of $\so_{X_0}$--modules. There is an obvious exact
functor $\usmf{X}:\gsmo{X}\to\smo{X_0}$, and we prove that it admits a right
and a left adjoint, which are both right inverse of it, and that it
is an equivalence of categories if and only if $X$ is \std. Similar
results hold also for the subcategories of (quasi)coherent graded
sheaves, which can be defined in the usual way.

When $X$ is of the form $\gproj R$ (with $R_0$ a field) it turns out that the
category of (quasi)coherent graded sheaves on $X$ is equivalent to the category
(considered in \cite{AKO}) of (quasi)coherent sheaves on the quotient algebraic
stack naturally associated to $\proj R$, and we show that a similar result
holds for more general graded schemes.

                \section{Graded rings and modules}

        \begin{defi}
Let $A$ be a graded ring.\footnote{Unless otherwise stated, all rings are
assumed to be (strictly) commutative with unit.} The topological space $\gspec
A$\index{Spec@$\gspec$} is the set of homogeneous prime ideals of $A$ with
topology defined in the following way: a subset of $\gspec A$ is closed if and
only if it is of the form
\[\V{\I{a}}:=
\{\I{p}\in\gspec A\st\I{a}\subseteq\I{p}\}\]\index{V(a)@$\V{\I{a}}$}
for some homogeneous ideal $\I{a}$ of $A$.
        \end{defi}

As usual this is well defined, and the open subsets of the form $\D{a}:= \gspec
A\minus\V{(a)}$\index{D(a)@$\D{a}$} (for $a\in A\hom$) form a base of the
topology of $\gspec A$. A morphism of graded rings (i.e. a morphism of rings
preserving degrees) $\varphi:A\to B$ induces a continuous map $f:\gspec B\to
\gspec A$ (defined by $f(\I{p}):=\varphi^{-1}(\I{p})$), so that $\gspec$ is a
contravariant functor from the category of graded rings to the category of
topological spaces.

Of course, $\gspec A=\spec A$ if $A_d=0$ for $d\ne0$.

        \begin{lemm}\label{specsurj}
Let $A$ be a graded ring. The natural map
\[\begin{split}
f:\gspec A & \to \spec A_0 \\
\I{p} & \mapsto \I{p}_0
\end{split}\]
is continuous, closed and surjective.

Moreover, $\all\I{q}\in\spec A_0$ $\exiun\g{\I{q}}\in f^{-1}(\I{q})$
such that $\I{p}\subseteq\g{\I{q}}$ $\all \I{p}\in f^{-1}(\I{q})$.

If $\I{p}\in\gspec A$ then $\I{p}=\g{\I{p}_0}$ if and only if the following
condition is satisfied $\all d\in\Z$: $\I{p}_d=A_d\iff\I{p}_{-d}=A_{-d}$.
        \end{lemm}

        \begin{proof}
$f$ is clearly continuous (it is the map induced by the morphism of graded
rings $A_0\mono A$).

$\all\I{q}\in\spec A_0$ the homogeneous ideal $\g{\I{q}}$ of $A$ defined by
\begin{align*}
& \g{\I{q}}_d:=\{a\in A_d\st(a)_0\subseteq\I{q}\} & \all d\in\Z
\end{align*}
is prime: indeed, if $a,b\in A\hom\minus\g{\I{q}}$, by definition of
$\g{\I{q}}$ there exist $c,d\in A\hom$ such that $ac,bd\in A_0\minus\I{q}$.
Since $\I{q}$ is prime this implies that $abcd\in A_0\minus\I{q}=
(A\minus\g{\I{q}})_0$, whence $ab\notin\g{\I{q}}$. Now it is clear that
$\g{\I{q}}$ satisfies the required property, and so, in particular, $f$ is
surjective.

In order to prove that $f$ is closed, it is enough to show that for every
homogeneous ideal $\I{a}$ of $A$, $f(\V{\I{a}})=\V{\I{a}_0}$. Obviously
$f(\V{\I{a}})\subseteq\V{\I{a}_0}$, and conversely $\all \I{q}\in
\V{\I{a}_0}$ we have $\I{q}=f(\g{\I{q}})$ and $\g{\I{q}}\in\V{\I{a}}$.

As for the last statement, if $\I{p}\in\gspec A$ satisfies $\I{p}=\g{\I{p}_0}$
and $\I{p}_d=A_d$ for some $d$, then $\all a\in A_{-d}$ we have $(a)_0\subseteq
\I{p}_0$, which implies that $\I{p}_{-d}=A_{-d}$. Conversely, if $\I{p}
\subsetneq\g{\I{p}_0}$ then $\exi a\in(\g{\I{p}_0}\minus\I{p})_d$ for some
$d\in\Z$, whence $\I{p}_d\subsetneq A_d$ and $\I{p}_{-d}=A_{-d}$ ($\all b
\in A_{-d}$ we have $ab\in(\g{\I{p}_0})_0=\I{p}_0$).
        \end{proof}

        \begin{defi}
A graded ring $A$ is {\em \good}\index{good!graded ring} if the natural map
$\gspec A\to\spec A_0$ is injective (and hence a homeomorphism by
\ref{specsurj}).
        \end{defi}

        \begin{rema}\label{goodcrit}
It follows from \ref{specsurj} that a graded ring $A$ is {\good} if and only if
$\I{p}=\g{\I{p}_0}$ $\all \I{p}\in\gspec A$, which is true if and only if the
following condition is satisfied $\all d\in\Z$: $\I{p}_d=A_d\iff\I{p}_{-d}=
A_{-d}$.
        \end{rema}

        \begin{defi}
A graded ring $A$ is {\em \triv}\index{trivial!graded ring} (respectively {\em
\qtriv})\index{quasi--trivial!graded ring} if every element of $A_{\ne0}$ is
zero (respectively nilpotent).
        \end{defi}

        \begin{defi}
A graded ring is {\em \std}\index{standard!graded ring} (respectively {\em
\qstd})\index{quasi--standard!graded ring} if it contains an invertible
homogeneous element of degree $1$ (respectively $>0$).
        \end{defi}

        \begin{lemm}\label{qtrivcrit}
A graded ring $A$ is {\qtriv} if and only if $A_{\ne0}\subset\I{p}$ $\all\I{p}
\in\gspec A$.
        \end{lemm}

        \begin{proof}
It is an immediate consequence of the following algebraic fact: if $\I{a}$ is
a homogeneous ideal of a graded ring $A$, then $\sqrt{\I{a}}=\bigcap_{\I{p}\in
\V{\I{a}}}\I{p}$.
        \end{proof}

        \begin{lemm}\label{stdcrit}
A graded ring $A$ is {\std} if and only if it is isomorphic to
$A_0[t,t^{-1}]$, where $\deg(t)=1$.
        \end{lemm}

        \begin{proof}
One implication is obvious, since $t$ is invertible and of degree $1$ in
$A_0[t,t^{-1}]$. Conversely, if $A$ is a graded ring and $a\in A_1$ is
invertible, then the morphism of graded rings $\varphi:A_0[t,t^{-1}]\to A$
defined by $\varphi(t)=a$ is easily seen to be an isomorphism.
        \end{proof}

        \begin{lemm}\label{qstdgood}
For a graded ring $A$ the following implications are true:
\begin{gather*}
A\text{ \std}\implies A\text{ \qstd}\implies A\text{ \good,}\\
A\text{ \triv}\implies A\text{ \qtriv}\implies A\text{ \good.}
\end{gather*}
        \end{lemm}

        \begin{proof}
The implications on the left are obvious.

If $A$ is \qstd, by \ref{goodcrit} we have to prove that if $\I{p}\in\gspec A$
is such that $\I{p}_d=A_d$ for some $d\in\Z$, then also $\I{p}_{-d}=
A_{-d}$. Now, if $s\in A_n$ ($n>0$) is invertible and $a\in A_{-d}$, then
$s^da^{n-1}\in A_d\subset\I{p}$, whence $n>1$ and $a\in\I{p}$.

If $A$ is \qtriv, then $A_{\ne0}\subset\I{p}$ $\all\I{p}\in\gspec A$ by
\ref{qtrivcrit}, whence $A$ is {\good} by \ref{goodcrit}.
        \end{proof}

If $A$ is a graded ring, $S\subset A\hom$ will be called a {\em homogeneous
multiplicative system} if $0\notin S$ and $a,b\in S\implies ab\in S$. As usual,
one can define the localized graded ring $S^{-1}A$ and, for every graded
$A$--module $M$, the graded $S^{-1}A$--module $S^{-1}M$. If $S=\{a^n\st n\in
\N\}$ for some $a\in A\hom$ not nilpotent and $T=A\hom\minus\I{p}$ for
some $\I{p}\in\gspec A$, we will write $A_a$, $M_a$ and $A_{\I{p}}$,
$M_{\I{p}}$ instead of $S^{-1}A$, $S^{-1}M$ and $T^{-1}A$,
$T^{-1}M$.

        \begin{exem}\label{Rrqstd}
If $R$ is a positively graded ring (i.e. $R_d=0$ if $d<0$) and $r\in R_n$
($n>0$) is not nilpotent, then clearly the localization $R_r$ is a {\qstd} (and
also {\std} if $n=1$) graded ring.

On the other hand, $R$ itself is never \qstd, and it is {\good} if and only
if it is {\qtriv} (this follows easily from \ref{goodcrit} and
\ref{qtrivcrit}).
        \end{exem}

        \begin{lemm}\label{0loc=loc0}
Let $A$ be a graded ring, $M$ a graded $A$--module and $S\subset A$ a
homogeneous multiplicative system such that $\all d\in\Z$ $S_d=
\emptyset\iff S_{-d}=\emptyset$. Then the natural morphism of
$S_0^{-1}A_0$--modules $f:S_0^{-1}M_0\to(S^{-1}M)_0$ is an isomorphism.
        \end{lemm}

        \begin{proof}
$f$ is injective: if $m\in M_0$ and $s\in S_0$ are such that $m/s=0$ in
$(S^{-1}M)_0$, there exist $d\in\Z$ and $t\in S_d$ such that $tm=0$. By
hypothesis there exist $t'\in S_{-d}$, and so $t'tm=0$, which implies that
$m/s=0$ in $S_0^{-1}M_0$ because $t't\in S_0$.

$f$ is surjective: given $m/s\in(S^{-1}M)_0$ (with $m\in M_d$ and $s\in S_d$
for some $d\in\Z$) by hypothesis there exists $s'\in S_{-d}$, whence $m/s=
f((s'm)/(s's))$.
        \end{proof}

        \begin{coro}\label{0p=p0}
If $A$ is a {\good} graded ring and $M$ a graded $A$--module, then $\all
\I{p}\in\gspec A$ the natural morphism of $(A_0)_{\I{p}_0}$--modules
$(M_0)_{\I{p}_0}\to(M_{\I{p}})_0$ is an isomorphism.
        \end{coro}

        \begin{proof}
By \ref{goodcrit} the multiplicative system $A\hom\minus\I{p}$ satisfies the
condition of \ref{0loc=loc0}.
        \end{proof}

        \begin{prop}\label{goodlocgood}
Let $A$ be a graded ring and $S\subset A$ a homogeneous multiplicative system.
If $A$ is {\good} or {\qstd} or {\std} or {\qtriv} or {\triv}, then so is
$S^{-1}A$.
        \end{prop}

        \begin{proof}
Obvious if $A$ is {\qstd} or {\std} or \triv; if $A$ is \qtriv, the statement
follows immediately from the fact that $S\subset A_0$.

If $A$ is \good, given $\I{q},\I{q}'\in\gspec(S^{-1}A)$ such that $\I{q}\ne
\I{q}'$, we have to show that $\I{q}_0\ne\I{q}'_0$. Since there is a one to one
correspondence between $\{\I{p}\in\gspec A\st\I{p}\cap S=\emptyset\}$ and
$\gspec(S^{-1}A)$ (given by $\I{p}\mapsto S^{-1}\I{p}$), there exist (unique)
$\I{p},\I{p}'\in\gspec A$ such that $\I{p}\cap S=\I{p}'\cap S=\emptyset$,
$\I{q}=S^{-1}\I{p}$, $\I{q}'=S^{-1}\I{p}'$ and $\I{p}\ne\I{p}'$. This implies
$\I{p}_0\ne\I{p}'_0$ because $A$ is \good, and we can assume that $\I{p}_0
\subsetneq\I{p}'_0$. So let $a\in\I{p}_0\minus\I{p}'_0$: then clearly $a\in
\I{q}_0$ and $a\notin\I{q}'_0$.
        \end{proof}

        \begin{lemm}\label{Apqstd}
Given a graded ring $A$ and $\I{p}\in\gspec A$, consider the ideal of $\Z$
\[I(\I{p}):=(i\in\Z\st\I{p}_i\ne A_i).\]
Then $A_{\I{p}}$ is {\std} (respectively \qstd) if and only if $I(\I{p})=(1)$
(respectively $I(\I{p})\ne(0)$, i.e. $A_{\ne0}\nsubseteq\I{p}$).
        \end{lemm}

        \begin{proof}
If $A_{\I{p}}$ is {\std} (respectively \qstd), then there exist $x\in
(A_{\I{p}})_d$ invertible with $d=1$ (respectively $d>0$). $x$ can be written
in the form $x=a/s$ with $a,s\in A\hom\minus\I{p}$ and $\deg(a)=\deg(s)+d$.
Since $\deg(a),\deg(s)\in I(\I{p})$, also $d\in I(\I{p})$.

Conversely, if $I(\I{p})=(d)$ with $d=1$ (respectively $d>0$), then there exist
elements $a_i\in A_{d_i}\minus\I{p}$ ($i=1,\dots,n$) such that $(d_1,\dots,d_n)
=(d)$. Therefore there exist $m_1,\dots,m_n\in\Z$ such that $d=\sum_{i=1}^n
m_id_i$, whence $x:=\prod_{i=1}^na_i^{m_i}\in A_{\I{p}}$ is an invertible
element of degree $d$, i.e. $A_{\I{p}}$ is {\std} (respectively \qstd).
        \end{proof}

        \begin{prop}\label{locgoodgood}
If $A$ is a graded ring such that $A_{\I{p}}$ is {\good} $\all \I{p}\in
\gspec A$, then $A$ is \good, too. The same statement holds with {\qtriv} or
{\triv} in place of \good.
        \end{prop}

        \begin{proof}
By \ref{goodcrit} given $\I{p}\in\gspec A$ we have to prove that $\I{p}=
\g{\I{p}_0}$, or, equivalently, that $\I{p}A_{\g{\I{p}_0}}=\g{\I{p}_0}
A_{\g{\I{p}_0}}$. Since $A_{\g{\I{p}_0}}$ is \good, this is true if and only if
$(\I{p}A_{\g{\I{p}_0}})_0=(\g{\I{p}_0}A_{\g{\I{p}_0}})_0$. Now, given $a/s\in
(\g{\I{p}_0}A_{\g{\I{p}_0}})_0$ (where $a\in(\g{\I{p}_0})_d$ and $s\in
(A\minus\g{\I{p}_0})_d$), by \ref{specsurj} there exists $t\in
(A\minus\g{\I{p}_0})_{-d}$, whence $a/s=(at)/(st)\in(\I{p}A_{\g{\I{p}_0}})_0$
because $at\in(\g{\I{p}_0})_0=\I{p}_0$.

If $\I{p}\in\gspec A$ is such that $A_{\I{p}}$ is \qtriv, then $A_{\ne0}
\subset\I{p}$ by \ref{Apqstd}. Therefore if $A_{\I{p}}$
is {\qtriv} $\all \I{p}\in\gspec A$, then $A$ is {\qtriv} by
\ref{qtrivcrit}.

If $A_{\I{p}}$ is {\triv} $\all \I{p}\in\gspec A$, let $a\in A_{\ne0}$:
$\all \I{p}\in\gspec A$ $\exi s\in A\minus\I{p}=A_0\minus\I{p}_0$ such
that $sa=0$, which implies that $a=0$, whence $A$ is \triv.
        \end{proof}

On the other hand, if $A_{\I{p}}$ is {\qstd} (respectively \std) $\all \I{p}\in
\gspec A$, then $A$ need not be {\qstd} (respectively \std), too, as the
following example shows.

        \begin{exem}\label{lstdnotqstd}
There exists a graded ring $A$ which is not \qstd, but such that $A_{\I{p}}$ is
{\std} $\all \I{p}\in\gspec A$. To see this, let $B$ be a noetherian ring and
$M$ a $B$--module such that $M_{\I{q}}$ is a free $B_{\I{q}}$--module of rank
one $\all\I{q}\in\spec B$ and $\all n>0$ $M^{\otimes n}$ is not a free
$B$--module (this is equivalent to require that $M\sh{}$ is an invertible sheaf
on $\spec B$ whose equivalence class in $\Pic(\spec B)$ is not a torsion
element). Setting $M^{\otimes n}:=\Hom_B(M^{\otimes(-n)},B)$ for $n<0$, it is
then clear that $A:=\bigoplus_{d\in\Z}M^{\otimes d}$ is a graded ring (with
$A_0=B$) in a natural way. By construction each $A_{\I{p}}$ is \std, but $A$ is
not \qstd: assume on the contrary that there is $n>0$ and $a\in A_n$, $b\in
A_{-n}$ such that $ab=1$. Regarding $b$ as a morphism $b:M^{\otimes n}\to B$ of
$B$--modules, $b$ is surjective (since $b(a)=1$) and hence it is an isomorphism
(because $M^{\otimes n}$ and $B$ are both projective $B$--modules of rank one),
which contradicts the hypothesis on $M$.
        \end{exem}

        \begin{lemm}\label{lqstdlocqstd}
Let $A$ be a graded ring. If $\I{p}\in\gspec A$ is such that $A_{\I{p}}$ is
{\qstd} (respectively \std), then there exists $t\in A\hom\minus\I{p}$ such
that $A_t$ is also {\qstd} (respectively \std).
If moreover $A$ is \good, then $t$ can be taken in $A_0\minus\I{p}_0$.
        \end{lemm}

        \begin{proof}
If $a/s\in(A_{\I{p}})_n$ (with $s\in A_d\minus\I{p}$, $a\in A_{d+n}$ and $n>0$)
is an invertible element, it is enough to prove that $(A_t)_n$ contains an
invertible element for some $t\in A\hom\minus\I{p}$. By hypothesis there exist
$s'\in A_{d'}\minus\I{p}$, $a'\in A_{d'-n}$ and $s''\in A_{d''}\minus\I{p}$
such that $aa's''=ss's''$. Then, setting $t:=ss's''\in A_{d+d'+d''}\minus
\I{p}$, $as's''/t\in(A_t)_n$ is an invertible element (its inverse is
$a'ss''/t$).

If moreover $A$ is \good, then (by \ref{goodcrit}) $\all u\in A\hom\minus
\I{p}$ there exists $u'\in A\hom\minus\I{p}$ such that $uu'\in A_0$, and this
implies that in the above argument we can assume $d=d'=d''=0$, whence $t\in A_0
\minus\I{p}_0$.
        \end{proof}

        \begin{defi}
$A$ is a {\em local graded ring}\index{local graded ring} if $A$ is a graded
ring with a unique maximal homogeneous ideal.
        \end{defi}

        \begin{prop}\label{gspecAgood}
Let $\varphi:A\to B$ be a morphism of graded rings.
\begin{enumerate}

\item If $A$ is {\qstd} (respectively \std), then $B$ is {\qstd}
(respectively \std), too.

\item If $A$ is local graded and $\varphi$ is surjective (whence $B$ is also
local graded), then $A$ is {\qstd} (respectively \std) if and only if $B$ is
{\qstd} (respectively \std).

\end{enumerate}
        \end{prop}

        \begin{proof}
\begin{enumerate}

\item Clear, since $\varphi$ preserves degrees and $\varphi(a)$ is
invertible if $a$ is invertible.

\item Just observe that, denoting by $\I{m}$ and $\I{n}$ the maximal
homogeneous ideals of $A$ and $B$, $\all d\in\Z$ we have $\I{m}_d=A_d$ if and
only if $\I{n}_d=B_d$.
\end{enumerate}
        \end{proof}

        \begin{prop}\label{Anoeth}
For a graded ring $A$ the following conditions are equivalent:
\begin{enumerate}

\item $A_0$ is a noetherian ring and $A$ is a finitely generated
$A_0$--algebra;

\item $A$ is a noetherian ring;

\item $A$ satisfies the ascending chain condition for homogeneous ideals, or,
equivalently, every homogeneous ideal of $A$ is finitely generated.

\end{enumerate}
If the above conditions are satisfied and $M$ is a finitely generated graded
$A$--module, then $M_d$ is a finitely generated $A_0$--module $\all d\in\Z$.
        \end{prop}

        \begin{proof}
(A proof of $2\iff3$ can also be found in \cite{NO}). First we show that $3$
implies the last statement. Clearly $M$ satisfies the ascending chain condition
for graded $A$--submodules. Since moreover for every $A_0$--submodule $N$ of
$M_d$ we have $(AN)_d=N$, it follows that $M_d$ is a noetherian (hence finitely
generated) $A_0$--module. Notice that, in particular, $A_0$ is a noetherian
ring.

The implication $1\implies2$ follows from Hilbert's base theorem, whereas
$2\implies3$ is obvious. So let's prove that $3\implies 1$. We have already
proved that $A_0$ is noetherian. Let $A_+$ and $A_-$ be the ideals of $A$
generated by $A_{>0}$ and $A_{<0}$ respectively: the hypothesis clearly implies
that there exist $a_1,\dots,a_n\in A\hom_{>0}$ and $b_1,\dots,b_m\in
A\hom_{<0}$ such that $A_+=(a_1,\dots,a_n)$ and $A_-=(b_1,\dots,b_m)$. Let
$r,s\in\Z$ be such that $r>\deg(a_i)$ for $i=1,\dots,n$ and $s<\deg(b_i)$ for
$i=1,\dots,m$. Since each $A_d$ is a finitely generated $A_0$--module, there
exist homogeneous elements $c_1,\dots,c_l$  which generate
$\bigoplus_{s<d<r}A_d$ as an $A_0$--module. Then it is easy to see that
$a_1,\dots,a_n,b_1,\dots,b_m,c_1,\dots,c_l$ generate $A$ as an $A_0$--algebra.
        \end{proof}

        \begin{prop}\label{locgqtqs}
Let $A$ be a noetherian local graded ring. If $A$ is \good, then $A$ is either
{\qstd} or \qtriv.
        \end{prop}

        \begin{proof}
Assume on the contrary that $A$ is {\good} and neither {\qstd} nor \qtriv:
then, denoting by $\I{m}$ the maximal homogeneous ideal of $A$, $A_{\ne0}
\subset\I{m}$ and $\exi \I{p}\in\gspec A$ such that $A_{\ne0}\nsubseteq \I{p}$
(by \ref{qtrivcrit}). As $A$ is noetherian, we can also choose $\I{p}$ in such
a way that $\nexi\I{q}\in\gspec A$ such that $\I{p}\subsetneq \I{q}$ and
$A_{\ne0}\nsubseteq\I{q}$. Let $B:=A/\I{p}$ and $\I{n}:=\I{m}/ \I{p}$: clearly
$B$ is again a {\good} noetherian local graded ring with maximal homogeneous
ideal $\I{n}\ne(0)$. Let $0\ne b\in\I{n}_0$ (notice that $\I{n}_0\ne(0)$
because $\I{n}\ne(0)$ and $B$ is a {\good} domain): by the choice of $\I{p}$ we
have $\sqrt{(b)}\supset B_{\ne0}$, and this easily implies (since $B$ is a
finitely generated $B_0$--algebra by \ref{Anoeth}) that there exists $l\in\N$
such that $B_d\subset(b)$ for $\abs{d}>l$. Hence for $\abs{d}>l$ we have $B_d=
\I{n}_0B_d$, so that $B_d=(0)$ by Nakayama's lemma ($B_d$ is a finitely
generated $B_0$--module again by \ref{Anoeth}). As $B$ is a domain, this
actually implies that $B_{\ne0}=(0)$, whence $A_{\ne0}\subset\I{p}$, a
contradiction.
        \end{proof}

        \begin{exem}
If $A$ is a {\good} noetherian graded ring, then it is not true in general that
$A$ is {\qtriv} or $A_{\I{p}}$ is {\qstd} $\all \I{p}\in\gspec A$. Take for
instance
\begin{align*}
& A=\Z/6\Z[x,y]/(2x,2y,xy-3) & \deg(x)=-\deg(y)=1.
\end{align*}
It is easy to see that $A_0=\Z/6\Z$ and $\gspec A=\{(2),(x,y)\}$
(so that $A$ is \good), but $A_{(2)}$ is {\std} and $A_{(x,y)}$ is \triv.
        \end{exem}

        \begin{prop}\label{clopen}
Let $A$ be a {\good} noetherian graded ring. Then
\begin{align*}
& \QT{A}:=\{\I{p}\in\gspec A\st A_{\I{p}}\text{ \qtriv}\}, \\
&\QS{A}:=\{\I{p}\in\gspec A\st A_{\I{p}}\text{ \qstd}\}
\end{align*}
are closed and open in $\gspec A$ and $\gspec A=\QT{A}\djun\QS{A}$.
        \end{prop}

        \begin{proof}
Notice first that the last statement follows from \ref{locgqtqs}. Then, since
$\QT{A}=\gspec A\minus\QS{A}=\V{(A_{\ne0})}$ by \ref{Apqstd}, it is enough to
prove that $\QT{A}$ is open. Now, if $\I{p}\in\QT{A}$, then $\all a\in
A_{\ne0}$ there exist $n\in\N$ and $t\in A\minus\I{p}\subset A_0$ such that
$ta^n=0$. Since $A$ is a finitely generated $A_0$--algebra by \ref{Anoeth},
this implies that there exists $s\in A\minus\I{p}$ such that $\all a\in
A_{\ne0}$ $\exi n\in\N$ with $sa^n=0$. Therefore, $\I{q}\in\D{s}$ and $a\in
A_{\ne0}$ imply $a\in\I{q}$, so that $\I{q}\in\QT{A}$. This proves that $\D{s}
\subset\QT{A}$, and so $\QT{A}$ is open.
        \end{proof}

With a proof completely analogous to that of the not graded case (see
\cite[prop. II.3.2]{H1}), one can prove the following result.

        \begin{lemm}\label{locnoethnoeth}
Let $A$ be a graded ring and $M$ a graded $A$--module. If there exist $a_1,
\dots,a_n\in A\hom$ such that $(a_1,\dots,a_n)=(1)$ and each $M_{a_i}$ is a
noetherian $A_{a_i}$--module, then $M$ is a noetherian $A$--module.
        \end{lemm}

                \section{Modules versus graded modules}

If $A$ is a graded ring, we will denote by $\gmo{A}$\index{Mod(A) @$\gmo{A}$}
the (abelian) category of graded $A$--modules (with morphisms preserving
degrees) and by $\gfmo{A}$\index{mod(A) @$\gfmo{A}$} its full subcategory of
finitely generated graded $A$--modules. $\mo{A}$\index{Mod(A)@$\mo{A}$} and
$\fmo{A}$\index{mod(A)@$\fmo{A}$} will denote instead the categories of
$A$--modules and of finitely generated $A$--modules respectively. Notice that
if $A=A_0$ then $\mo{A}$ can be identified with a full subcategory of $\gmo{A}$
(by regarding a module as concentrated in degree $0$).

$\all n,m\in\Z$ the natural functors $(n):\gmo{A}\to\gmo{A}$ (defined
by $M(n)_d=M_{n+d}$) are exact autoequivalences and $(n)\comp(m)=(n+m)$.

$\gmo{A}$ has both enough projectives ($\all n\in\Z$ the free graded
$A$--module $A(n)$ is obviously projective) and enough injectives (this can be
easily seen slightly adapting the standard proof that $\mo{A}$ has enough
injectives).

If $M$ and $N$ are (finitely generated) graded $A$--modules, so is $M\otimes_A
N$. As usual $M\otimes_A-:\gmo{A}\to\gmo{A}$ is right exact and $\all n\in
\Z$ the functors $(n)$ and $A(n)\otimes_A-$ are isomorphic.

We will write $\Hom_A(M,N)$\index{Hom_A@$\Hom_A$} instead of
$\Hom_{\gmo{A}}(M,N)$, and we will denote by
$\gHom_A(M,N)$\index{Hom_A@$\gHom_A$} the graded $A$-module defined by
\begin{align*}
& \gHom_A(M,N)_d:=\Hom_A(M,N(d))=\Hom_A(M(-d),N) & \all d\in\Z.
\end{align*}
Clearly $\gHom_A(M,N)_0=\Hom_A(M,N)$ and the functor $\gHom_A(A(-n),-)$ is
isomorphic to $(n)$. The functors
\begin{gather*}
\Hom_A(-,-):\gmo{A}\opp\times\gmo{A}\to\mo{A_0},\\
\gHom_A(-,-):\gmo{A}\opp\times\gmo{A}\to\gmo{A}
\end{gather*}
are left exact in both arguments. Notice that $\gHom_A(M,N)=\Hom_{\mo{A}}(M,N)$
(as $A$--modules) if $M$ is finitely generated, and that $\gHom_A(M,N)\in
\gfmo{A}$ and $\Hom_A(M,N)\in\fmo{A_0}$ if $M,N\in\gfmo{A}$ and $A$ is
noetherian.

A morphism of graded rings $\varphi:A\to B$ induces functors
\begin{align*}
& B\otimes_A-:\gmo{A}\to\gmo{B}, & \gHom_A(B,-):\gmo{A}\to\gmo{B}.
\end{align*}
If $M\in\gfmo{A}$ then $B\otimes_A M\in\gfmo{B}$ and $\gHom_A(B,M)\in\gfmo{B}$
if $A$ is noetherian and $\varphi$ is finite. As in the not graded case (see
\cite[p. 693]{E}) it is straightforward to prove the following result.

        \begin{lemm}\label{adj}
$B\otimes_A-$ and $\gHom_A(B,-)$ are respectively left and right adjoint of
the natural forgetful functor $\gmo{B}\to\gmo{A}$.
        \end{lemm}

Now we are going to study the relation between $\gmo{A}$ and $\mo{A_0}$, where
$A$ is a graded ring. Clearly taking the degree $0$ component defines an exact
functor
\[\umf{A}:\gmo{A}\to\mo{A_0},\]
which sends $\gfmo{A}$ into $\fmo{A_0}$ if $A$ is noetherian (by \ref{Anoeth}).

        \begin{prop}\label{adj0}
The functors
\begin{align*}
& \omf{A}:\mo{A_0}\to\gmo{A}, & \hmf{A}:\mo{A_0}\to\gmo{A}
\end{align*}
are respectively left and right adjoint of $\umf{A}$, and they are both right
inverse of $\umf{A}$, i.e.
\[\umf{A}\comp\omf{A}\iso\umf{A}\comp\hmf{A}\iso\id.\]
$\omf{A}$ sends $\fmo{A_0}$ into $\gfmo{A}$, and the same is true for
$\hmf{A}$ if $A$ is {\good} and noetherian.
        \end{prop}

        \begin{proof}
Taking into account that given $M\in\gmo{A}$ and $N\in\mo{A_0}\subset\gmo{A_0}$
\begin{align*}
& \Hom_{A_0}(M_0,N)=\Hom_{A_0}(M,N) & &\text{and} &
\Hom_{A_0}(N,M_0)=\Hom_{A_0}(N,M),
\end{align*}
by \ref{adj} everything is clear except the last statement.

If $A$ is {\good} and noetherian, then there exist $a_1,\dots,a_n\in A_0$ such
that each $A_{a_i}$ is either {\qstd} or {\qtriv} and with
$(a_1,\dots,a_n)=(1)$: taking into account that $\gspec A$ is quasi--compact
because $A$ is noetherian, this follows immediately from \ref{clopen},
\ref{lqstdlocqstd} and \ref{locgoodgood} (just notice that if $a\in A\hom$ is
such that $A_a$ is \qtriv, then $a\in A_0$). By \ref{locnoethnoeth} it is
enough to prove that $\gHom_{A_0}(A,N)_{a_i}$ is a finitely generated
$A_{a_i}$--module $\all N\in\fmo{A_0}$ and $\all i=1,\dots,n$. Since
$\gHom_{A_0}(A,N)_{a_i}\iso\gHom_{(A_0)_{a_i}}(A_{a_i},N_{a_i})$ (by
\cite[prop. 2.10]{E}) and since $(A_0)_{a_i}\iso(A_{a_i})_0$ (by
\ref{0loc=loc0}), this means that we can assume that $A$ is either {\qstd} or
\qtriv. If $M:=\gHom_{A_0}(A,N)$ for some $N\in\fmo{A_0}$, $\all d\in \Z$ we
have $M_d=\Hom_{A_0}(A_{-d},N)\in\fmo{A_0}$ (because $A_{-d}\in \fmo{A_0}$ by
\ref{Anoeth}), whence it is enough to show that $\exi m>0$ such that $M$ is
generated by $\bigoplus_{-m\le d\le m}M_d$ as an $A$--module. Now, if $A$ is
{\qstd} we can take $m$ such that there exists $s\in A_{2m}$ invertible
(because then $\all d\in\Z$ multiplication by $s$ is bijective between $M_d$
and $M_{d+2m}$), whereas if $A$ is {\qtriv} we can take $m$ such that $A_d=(0)$
for $\abs{d}>m$ (such an $m$ exists because $A$ is a finitely generated
$A_0$--algebra by \ref{Anoeth}).
        \end{proof}

        \begin{prop}\label{modequivcrit}
$\umf{A}:\gmo{A}\to\mo{A_0}$ is an equivalence of categories if and only if
$A_{\I{p}}$ is {\std} $\all \I{p}\in\gspec A$.

If this is the case, then $\omf{A}\iso\hmf{A}$ and they are quasi--inverse of
$\umf{A}$ and exact. Moreover, the restrictions of these functors give an
equivalence of categories between $\gfmo{A}$ and $\fmo{A_0}$ if $A$ is
noetherian.
        \end{prop}

        \begin{proof}
By \ref{adj0} it is clear that $\umf{A}$ is an equivalence of categories if and
only if $\omf{A}$ is quasi--inverse of $\umf{A}$ (and in this case $\omf{A}\iso
\hmf{A}$ and they are exact since $\umf{A}$ is). Therefore $\umf{A}$ is an
equivalence if and only if $\all M\in\gmo{A}$ the natural morphism (in
$\gmo{A}$) $\alpha(M):A\otimes_{A_0}M_0\to M$ is an isomorphism,
i.e. if and only if the morphism of $A_0$--modules
\[\alpha(M)_d:A_d\otimes_{A_0}M_0\to M_d\]
is an isomorphism $\all d\in\Z$. First we notice that if $\umf{A}$ is
an equivalence, then $A$ is \good: $\all d\in\Z$ we have
$(A(d)/(A_d)A(d))_0=0$, whence $(A_d)=A$, which implies that $\I{p}_d\ne A_d$
$\all \I{p}\in\gspec A$, so that $A$ is {\good} by \ref{goodcrit}.
Therefore we can assume that $A$ is \good, and so (taking into account
\ref{0p=p0}) $\alpha(M)_d$ is an isomorphism if and only if
\[\alpha(M_{\I{p}})_d:(A_{\I{p}})_d\otimes_{(A_{\I{p}})_0}(M_{\I{p}})_0\to
(M_{\I{p}})_d\]
is an isomorphism $\all \I{p}\in\gspec A$. Now, if $A_{\I{p}}$ is \std, then
clearly $\alpha(M_{\I{p}})_d$ is an isomorphism, since if $t\in(A_{\I{p}})_1$
is invertible, then multiplication by $t^d$ induces an isomorphism between
$(M_{\I{p}})_0$ and $(M_{\I{p}})_d$ (and between $(A_{\I{p}})_0$ and
$(A_{\I{p}})_d$). Conversely, if $\umf{A}$ is an equivalence, then, in
particular, the natural multiplication map
\[\alpha(A_{\I{p}}(1))_{-1}:(A_{\I{p}})_{-1}\otimes_{(A_{\I{p}})_0}
(A_{\I{p}})_1\to(A_{\I{p}})_0\]
is an isomorphism $\all \I{p}\in\gspec A$. So there exist
$a\in(A_{\I{p}})_1$ and $b\in(A_{\I{p}})_{-1}$ such that $ab\notin
\I{p}A_{\I{p}}$ (otherwise $\im\alpha(A_{\I{p}}(1))_{-1}\subseteq\I{p}
A_{\I{p}}$), whence $a$ is invertible. Therefore $A_{\I{p}}$ is {\std} $\all
\I{p}\in\gspec A$.

The last statement follows from the fact that if $A$ is noetherian, then
$\umf{A}$ sends $\gfmo{A}$ into $\fmo{A_0}$ and $\omf{A}$ and $\hmf{A}$ send
$\fmo{A_0}$ into $\gfmo{A}$.
        \end{proof}

                \section{Graded schemes}

If $X$ is a topological space and $\s{F}$ a sheaf of graded abelian groups on
$X$, we will denote by $\gsec(U,\s{F})$\index{Gamma(U,F)@$\gsec(U,\s{F})$} the
graded abelian group of sections of $\s{F}$ over an open subset $U$ of $X$.
Notice that if $\s{F}_d$ ($d\in\Z$) are sheaves of abelian groups on $X$, then
the presheaf $\s{F}:= \bigoplus_{d\in\Z}\s{F}_d$ is not a sheaf in general.
However, $\s{F}$ is certainly a sheaf if every open subset of $X$ is
quasi--compact (which is true if and only if $X$ is noetherian).

        \begin{defi}
A {\em graded ringed space}\index{graded!ringed space} $X=(X,\so_X)$ is a
topological space $X$ with a sheaf of graded rings
$\so_X=\bigoplus_{d\in\Z}(\so_X)_d$ on $X$.

A morphism of graded ringed spaces\index{morphism of!graded ringed spaces}
$f=(f,f\mrs):X\to Y$ is a continuous map $f:X\to Y$ together with a morphism
(of degree $0$) $f\mrs:\so_Y\to f_*\so_X$ of sheaves of graded rings on $Y$.
        \end{defi}

If $X$ is a graded ringed space, we will denote by $\gsmo{\so_X}$, or simply by
$\gsmo{X}$,\index{Mod(X)@$\gsmo{X}$} the (abelian) category of sheaves of
graded $\so_X$--modules on $X$. With a proof completely analogous to that of
the not graded case, one can show that $\gsmo{X}$ has enough injectives.

$\all n\in\Z$ the natural {\em twist functors} $(n):\gsmo{X}
\to\gsmo{X}$ (defined by $\s{F}(n)_d:=\s{F}_{n+d}$) are exact autoequivalences
and $(n)\comp(m)=(n+m)$.

        \begin{defi}
A sheaf of graded $\so_X$--modules $\s{F}$ is {\em locally free} \index{locally
free graded sheaf} if it is locally isomorphic to a sheaf of the form
$\bigoplus_{i\in I}\so_X(n_i)$ (where $n_i\in\Z$). $\card{I}$ (if finite and
constant) is the {\em rank} of $\s{F}$.
        \end{defi}

As in the case of graded modules over a graded ring, one can define the
following functors (the first right exact and the other two left exact in both
arguments):\index{Hom_X@$\Hom_X$}\index{Hom_X @$\gHom_X$}
\begin{gather*}
-\otimes_X-:=-\otimes_{\so_X}-:\gsmo{X}\times\gsmo{X}\to\gsmo{X},\\
\Hom_X(-,-):=\Hom_{\gsmo{X}}(-,-):\gsmo{X}\opp\times\gsmo{X}\to
\mo{\gsec(X,\so_X)_0},\\
\gHom_X(-,-):=\bigoplus_{d\in\Z}\Hom_X(-,-(d)):\gsmo{X}\opp\times
\gsmo{X}\to\gmo{\gsec(X,\so_X)}.
\end{gather*}
One can consider moreover the functors\index{Hom_X@$\sHom_X$}\index{Hom_X
@$\gsHom_X$}
\begin{gather*}
\sHom_X(-,-):\gsmo{X}\opp\times\gsmo{X}\to\smo{(\so_X)_0},\\
\gsHom_X(-,-)=\bigoplus_{d\in\Z}\sHom_X(-,-(d)):\gsmo{X}\opp\times
\gsmo{X}\to\gmo{X},
\end{gather*}
where, for every open subset $U$ of $X$, $\sHom_X(-,-)(U):=
\Hom_U(-\rest{U},-\rest{U})$ and $\gsHom_X(-,-)(U):=
\gHom_U(-\rest{U},-\rest{U})$. By definition we have
$\gHom_X(-,-)_0=\Hom_X(-,-)$ and $\gsHom_X(-,-)_0=\sHom_X(-,-)$, and clearly
$\all n\in\Z$ the functors $(n)$, $\so_X(n)\otimes_X-$ and
$\gsHom_X(\so_X(-n),-)$ are isomorphic. If $\s{L}$ is a locally free graded
sheaf of finite rank, we will write $\s{L}\dual$\index{$\dual$} for
$\gsHom_X(\s{L},\so_X)$. More generally, if $\cp{\s{L}}$ is a complex of
locally free graded sheaves of finite rank, $(\cp{\s{L}})\dual$ will denote the
complex $\gsHom_X(\cp{\s{L}},\so_X)$ (which in position $i$ has
$(\s{L}^{-i})\dual$).

Given $\s{F}\in\gsmo{X}$ and $i\in\N$, the $i^{\rm th}$ right derived functor
of $\Hom_X(\s{F},-)$, $\gHom_X(\s{F},-)$, $\sHom_X(\s{F},-)$ and
$\gsHom_X(\s{F},-)$ will be denoted by
$\Ext^i_X(\s{F},-)$,\index{Ext^i_X@$\Ext^i_X$}
$\gExt^i_X(\s{F},-)$,\index{Ext^i_X @$\gExt^i_X$}
$\sExt^i_X(\s{F},-)$\index{Ext^i_X@$\sExt^i_X$} and
$\gsExt^i_X(\s{F},-)$\index{Ext^i_X @$\gsExt^i_X$}, respectively. Notice that
$\gExt^i_X(\s{F},-)\iso \bigoplus_{d\in\Z}\Ext^i_X(\s{F},-(d))$ and that $\all
n\in\Z$ there are natural isomorphisms
\[\gExt^i_X(\s{F},-)(n)\iso\gExt^i_X(\s{F},-(n))\iso
\gExt^i_X(\s{F}(-n),-)\]
(whence, in particular, also $\Ext^i_X(\s{F},-(n))\iso
\Ext^i_X(\s{F}(-n),-)$). Analogous results hold, of course, also for
$\gsExt^i_X(\s{F},-)$ and $\sExt^i_X(\s{F},-)$. It is also clear that
$\gsExt^i_X(\so_X(n),-)=0$ if $i>0$. Moreover, if $\s{L}$ is a locally free
graded sheaf of finite rank, $\all i\in\N$ there are natural isomorphisms
\begin{gather*}
\gExt^i_X(\s{F}\otimes\s{L},-)\iso\gExt^i_X(\s{F},-\otimes\s{L}\dual),\\
\gsExt^i_X(\s{F}\otimes\s{L},-)\iso\gsExt^i_X(\s{F},-\otimes\s{L}\dual)\iso
\gsExt^i_X(\s{F},-)\otimes\s{L}\dual.
\end{gather*}

$H^i(X,-)$\index{H^i(X,-)@$H^i(X,-)$} ($i\in\N$) will denote the $i^{\rm th}$
right derived functor of the degree zero global section functor
$\gsec(X,-)_0:\gsmo{X}\to\mo{\gsec(X,\so_X)_0}$. Since clearly
$\Hom_X(\so_X(-n),-)\iso\gsec(X,-)_n$, there are natural isomorphisms
$H^i(X,-(n))\iso\Ext^i_X(\so_X(-n),-)$ $\all n\in\Z$. We define more generally
$\all i\in\Z$
\[H^i(X,-):=R^i\gsec(X,-)_0:D^+(\gsmo{X})\to\mo{\gsec(X,\so_X)_0}.\]
As usual, a morphism $f:X\to Y$ of graded ringed spaces induces a left exact
functor $f_*:\gsmo{X}\to\gsmo{Y}$ and a right exact functor $f^*:\gsmo{Y}\to
\gsmo{X}$, which is left adjoint of $f_*$.

        \begin{defi}
A {\em locally graded ringed space}\index{locally graded ringed space} is a
graded ringed space $X$ such that $\all x\in X$ the stalk $\so_{X,x}$ is a
local graded ring.

A morphism of locally graded ringed spaces\index{morphism of!locally graded
ringed spaces} $f:X\to Y$ is a morphism of graded ringed spaces such that $\all
x\in X$ the induced map on stalks $f\mrs_x: \so_{Y,f(x)}\to\so_{X,x}$ is a
local graded homomorphism of local graded rings.
        \end{defi}

Given a noetherian graded ring $A$ one can consider the graded ringed space $X$
such that $X=\gspec A$ as a topological space and with structure sheaf $\so_X$
that we are going to define. In general, for every graded $A$--module $M$ one
can consider the graded sheaf $M\gsh{A}$ (usually denoted simply by
$M\gsh{}$)\index{$\gsh{A}$} such that if $U\subseteq X$ is open, $M\gsh{}(U)$
is the set of functions $s:U \to\djun_{\I{p}\in U}M_{\I{p}}$ satisfying the
following conditions $\all\I{p} \in U$:
\begin{enumerate}

\item $s(\I{p})\in M_{\I{p}}$;

\item there exist $V\subseteq U$ neighbourhood of $\I{p}$ and $a\in A\hom$,
$m\in M$ such that $\all \I{q}\in V$, $a\notin\I{q}$ and $s(\I{q})=m/a$ in
$M_{\I{q}}$.

\end{enumerate}
Notice that since $A$ is noetherian, $U$ is quasi--compact, and this implies
that $M\gsh{}(U)=\bigoplus_{d\in\Z}M\gsh{}(U)_d$, where $M\gsh{}(U)_d$ is the
subset of $M\gsh{}(U)$ consisting of those sections which are locally of the
form $m/a$ with $m\in M_{d+\deg(a)}$.

Then $\so_X:=A\gsh{}$ is a sheaf of graded rings, so that $X$ is a graded
ringed space, and $M\gsh{}$ is a sheaf of graded $\so_X$--modules for every
graded $A$--module $M$. Clearly $A(n)\gsh{}\iso\so_X(n)$ $\all n\in\Z$.

As in the not graded case, it is not difficult to prove the following facts.
\begin{enumerate}

\item $\all \I{p}\in X$ the stalk $\so_{X,\I{p}}$ is isomorphic to $A_{\I{p}}$,
whence $X$ is a locally graded ringed space; more generally, $(M\gsh{})_{\I{p}}
\iso M_{\I{p}}$ if $M$ is a graded $A$--module.

\item Denoting (for $a\in A\hom$) by $\D{a}$ the open subset $X\minus\V{(a)}$
of $X$, $(\D{a},\so_X\rest{\D{a}})$ is isomorphic (as locally graded ringed
space) to $\gspec A_a$, and under this identification $(M\gsh{A})\rest{\D{a}}
\iso(M_a)\gsh{A_a}$.

\item A morphism $\varphi:A\to B$ of noetherian graded rings induces a morphism
$f:\gspec B\to\gspec A$ of locally graded ringed spaces.

\end{enumerate}

        \begin{defi}
An {\em affine graded scheme}\index{affine!graded scheme} is a locally graded
ringed space which is isomorphic to $\gspec A$ for some noetherian graded ring
$A$.

A {\em graded scheme}\index{graded!scheme} is a locally graded ringed space $X$
which is noetherian as a topological space and in which every point has an open
neighbourhood $U$ such that $(U,\so_X\rest{U})$ is an affine graded scheme.

A morphism of graded schemes\index{morphism of!graded schemes} is just a
morphism of locally graded ringed spaces.
        \end{defi}

        \begin{rema}
Of course the category $\gsch$\index{Sch@$\gsch$} of graded schemes includes as
a full subcategory the category of noetherian schemes (whose structure sheaves
are considered concentrated in degree $0$).
        \end{rema}

        \begin{rema}
If $X$ is a graded scheme, one can also define in the obvious way the category
$\gsch_X$ of graded schemes over $X$. If $A$ is a noetherian graded ring,
$\gsch_{\gspec A}$ will be usually denoted by $\gsch_A$.\index{Sch_A@$\gsch_A$}
        \end{rema}

Much of the theory of (noetherian) schemes extends naturally (with some obvious
modifications) to graded schemes. Here we just give the definitions and state
the results that we will need. Unless otherwise stated, the proofs are
completely analogous to those in the not graded case.

        \begin{defi}
Let $X$ be a graded scheme. A sheaf of graded $\so_X$--modules $\s{F}$ is {\em
quasi--coherent}\index{quasi--coherent graded sheaf} (respectively {\em
coherent})\index{coherent graded sheaf} if for some (and hence any) open cover
of $X$ by affine graded schemes $U_i\iso\gspec A^i$, there exist graded
$A^i$--modules (respectively finitely generated graded $A^i$--modules) $M^i$
such that $\s{F}\rest{U_i}\iso(M^i)\gsh{}$.
        \end{defi}

If $X$ is a graded scheme, $\gqco{X}$\index{Qcoh(X)@$\gqco{X}$} (respectively
$\gcoh{X}$)\index{Coh(X)@$\gcoh{X}$} will be the full subcategory of $\gsmo{X}$
whose objects are quasi--coherent (respectively coherent) graded sheaves. As
usual, kernels, cokernels, images and extensions of quasi--coherent
(respectively coherent) graded sheaves are also quasi--coherent (respectively
coherent), so that $\gqco{X}$ and $\gcoh{X}$ are abelian categories. $\gqco{X}$
has enough injectives.

The restrictions of the twist functors, of $-\otimes_X-$ and of $\gsHom_X(-,-)$
to $\gqco{X}$ (respectively $\gcoh{X}$) have image contained in $\gqco{X}$
(respectively $\gcoh{X}$).

If $f:X\to Y$ is a morphism of graded schemes, then $f^*(\gqco{Y})\subseteq
\gqco{X}$, $f^*(\gcoh{Y})\subseteq\gcoh{X}$ and $f_*(\gqco{X})\subseteq
\gqco{Y}$ (but, of course, it is false in general that $f_*(\gcoh{X})\subseteq
\gcoh{Y}$).

        \begin{rema}
A locally free graded sheaf is quasi--coherent, and also coherent if it is of
finite rank. In the latter case it will be often called {\em vector
bundle}.\index{vector bundle}

It is easy to prove that a sheaf of graded $\so_X$--modules $\s{F}$ is
quasi--coherent (respectively coherent) if and only if it is locally isomorphic
to the cokernel of a morphism of (locally) free sheaves (respectively of vector
bundles).
        \end{rema}

        \begin{rema}
A sheaf of graded $\so_X$--modules $\s{F}$ is locally free if and only if
$\all x\in X$ the stalk $\s{F}_x$ is a free graded $\so_{X,x}$--module.
        \end{rema}

        \begin{prop}\label{gspecequiv}
The contravariant functor $\gspec$ induces an antiequivalence between the
category of noetherian graded rings and the category of affine graded schemes.
Its quasi--inverse is defined on objects by $X\mapsto \gsec(X,\so_X)$.
Moreover, if $X$ is a graded scheme and $A$ is a noetherian graded ring, the
natural map
\[\begin{split}
\Hom_{\gsch}(X,\gspec A) & \to \Hom_{\g{\cat{Ring}}}(A,\gsec(X,\so_X)) \\
f & \mapsto f\mrs(\gspec A)
\end{split}\]
(where $\g{\cat{Ring}}$ denotes the category of graded rings) is bijective.
        \end{prop}

        \begin{prop}
If $A$ is a noetherian graded ring, the associated graded sheaf functor
$\gsh{A}:\gmo{A}\to\gsmo{\gspec A}$ is exact and its right adjoint is the
global sections functor $\gsec(\gspec A,-):\gsmo{\gspec A}\to\gmo{A}$.
These two functors induce quasi--inverse equivalences of categories between
$\gmo{A}$ (respectively $\gfmo{A}$) and $\gqco{\gspec A}$ (respectively
$\gcoh{\gspec A}$).
        \end{prop}

        \begin{prop}\label{modaff}
Let $\varphi:A\to B$ be a morphism of noetherian graded rings and let $f:Y:=
\gspec B\to X:=\gspec A$ be the induced morphism of graded schemes.
\begin{enumerate}

\item If $M,N\in\gmo{A}$, then $(M\otimes_A N)\gsh{}\iso M\gsh{}\otimes_X
N\gsh{}$.

\item If $M\in\gmo{A}$ and $N\in\gmo{B}$, then $f^*(M\gsh{A})\iso
(M\otimes_AB)\gsh{B}$ and $f_*(N\gsh{B})\iso N\gsh{A}$.
\end{enumerate}
        \end{prop}

        \begin{lemm}\label{secext}
Let $X$ be a graded scheme, $f\in\gsec(X,\so_X)_d$ ($d\in\Z$) and let $X_f$ be
the open subset of $X$ defined by $X_f:=\{x\in X\st f_x\notin\I{m}_x\}$ (where
$\I{m}_x\subset\so_{X,x}$ is the maximal homogeneous ideal). Then $\all
\s{F}\in\gqco{X}$ we have:
\begin{enumerate}

\item $\all s\in\gsec(X,\s{F})$ such that $s\rest{X_f}=0$, $\exi n>0$ such that
$f^ns=0$;

\item $\all t\in\gsec(X_f,\s{F})$ $\exi n>0$ and $\exi s\in\gsec(X,\s{F})$ such
that $s\rest{X_f}=f^nt$.
\end{enumerate}
        \end{lemm}

        \begin{proof}
The proof is similar (but easier, since the graduation is already in the
structure sheaf, and there is no need of an extra line bundle) to that of
\cite[II, lemma 5.14]{H1}.\footnote{Actually this result can be  viewed as a
generalization of the cited lemma, which just deals with the particular case of
{\std} graded schemes (see \ref{stdchar}).}
        \end{proof}

        \begin{lemm}\label{injflasque}
Let $X$ be a graded scheme and $\s{I}\in\gsmo{X}$ an injective object. Then
$\s{I}$ is flasque (i.e., the restriction map $\gsec(U,\s{I})\to\gsec(V,\s{I})$
is surjective if $V\subseteq U$ are open subsets of $X$).
        \end{lemm}

        \begin{defi}
A morphism $f:X\to Y$ of graded schemes is {\em affine}\index{affine!morphism}
if for some (and then for any) open affine cover $\{U_i\st i\in I\}$ of $Y$,
$f^{-1}(U_i)$ is affine $\all i\in I$. $f$ is {\em finite}\index{finite
morphism} if moreover the induced maps of graded rings
$\so_Y(U_i)\to\so_X(f^{-1}(U_i))$ are finite. $f$ is a {\em closed
immersion}\index{closed immersion} if, as a map of topological spaces, it gives
a homeomorphism between $X$ and a closed subset of $Y$, and if $f\mrs$ is
surjective.
        \end{defi}

Given a graded scheme $Y$ and a quasi--coherent sheaf of graded
$\so_Y$--algebras $\s{A}$, there exist unique (up to isomorphism) a graded
scheme $X$ and a morphism $f:X\to Y$ such that $f^{-1}(U)\iso\gspec\s{A}(U)$
and $f\rest{f^{-1}(U)}$ is the morphism corresponding to the natural map
$\so_Y(U)\to\s{A}(U)$ for every open affine subset $U$ of $Y$, and such that
the inclusion morphism $f^{-1}(V)\mono f^{-1}(U)$ corresponds to the
restriction map $\s{A}(U)\to\s{A}(V)$ if $V\subseteq U$ are open affine subsets
of $Y$. The graded scheme $X$ will be denoted by $\gSpec\s{A}$.\index{Spec
@$\gSpec$}

        \begin{prop}\label{gSpecA}
\begin{enumerate}

\item If $\s{A}$ is a quasi--coherent (respectively coherent) sheaf of graded
$\so_Y$--algebras, then the natural morphism of graded schemes $f:X:=
\gSpec\s{A}\to Y$ is affine (respectively finite) and $\s{A}\iso f_*\so_X$.
Conversely, if $f:X\to Y$ is an affine (respectively finite) morphism of graded
schemes, then $\s{A}:=f_*\so_X$ is a quasi--coherent (respectively coherent)
sheaf of graded $\so_Y$--algebras and $X\iso\gSpec\s{A}$.

\item Let $f:X\to Y$ be an affine morphism of graded schemes and let $\s{A}:=
f_*\so_X$. Then
\[f_*:\gqco{X}\to\gqco{Y}\cap\gsmo{\s{A}}\]
is an exact equivalence of abelian categories. If moreover $f$ is finite, then
$f_*$ restricts to an equivalence $f_*:\gcoh{X}\to\gcoh{Y}\cap\gsmo{\s{A}}$.

\end{enumerate}
        \end{prop}

        \begin{lemm}\label{gcgqgc}
Let $X$ be a graded scheme and let $\phi:\s{G}\to\s{F}$ be a surjective
morphism in $\gqco{X}$ with $\s{F}\in\gcoh{X}$. Then there exists a morphism
$\psi:\s{H} \to\s{G}$ in $\gqco{X}$ such that $\s{H}\in\gcoh{X}$ and
$\phi\comp\psi$ is surjective.
        \end{lemm}

        \begin{lemm}\label{gsExtgcoh}
Let $X$ be a graded scheme and let $\s{F},\s{G}\in\gcoh{X}$. Then
$\gsExt^i_X(\s{F},\s{G})\in\gcoh{X}$ $\all i\in\N$.
        \end{lemm}

                \section{Good graded schemes}

        \begin{defi}\label{defgoodsch}
Let $X$ be a graded scheme. A point $x\in X$ is called
{\good}\index{good!point} or {\qstd}\index{quasi--standard!point} or
{\std}\index{standard!point} or {\qtriv}\index{quasi--trivial!point} or
{\triv}\index{trivial!point} if the noetherian local graded ring $\so_{X,x}$
has the same property.

$X$ is called {\good}\index{good!graded scheme} or
{\qstd}\index{quasi--standard!graded scheme} or {\std}\index{standard!graded
scheme} or {\qtriv}\index{quasi--trivial!graded scheme} or
{\triv}\index{trivial!graded scheme} if every point of $X$ has the same
property.
        \end{defi}

        \begin{rema}\label{locgsch}
It follows from \ref{goodlocgood}, \ref{locgoodgood} and \ref{lqstdlocqstd}
that a graded scheme $X$ has one of the properties of \ref{defgoodsch} if and
only if $\all x\in X$ there is an open neighbourhood $U$ of $x$ which is
isomorphic to $\gspec A$ for some noetherian graded ring $A$ with the same
property. Moreover, if $X$ is {\good} or {\qtriv} or \triv, and if $U$ is any
open subset of $X$ such that $U\iso\gspec A$ for some noetherian graded ring
$A$, then $A$ has the same property.
        \end{rema}

From now on we will consider only {\good} graded schemes.

        \begin{rema}
Since a {\triv} graded scheme is just a noetherian scheme, we see that (usual)
noetherian schemes are {\good} graded schemes.
        \end{rema}

        \begin{prop}
If $X$ is a {\good} graded scheme, then $X_0:=(X,(\so_X)_0)$ is a noetherian
scheme. If $f:X\to Y$ is a morphism of {\good} graded schemes, then
$f_0:=(f,f\mrs_0): X_0\to Y_0$ is a morphism of noetherian schemes.
        \end{prop}

        \begin{proof}
If $U\subseteq X$ is an open subset such that $U\iso\gspec A$ for some {\good}
noetherian graded ring $A$, then $U_0\iso\spec A_0$: this follows immediately
from \ref{0p=p0}. The last statement is obvious.
        \end{proof}

        \begin{rema}
If $A$ is a {\good} noetherian graded ring, then $(\gspec A)_0\iso\spec A_0$.
        \end{rema}

        \begin{prop}\label{YgoodXgood}
Let $f:X\to Y$ be a morphism of graded schemes.
\begin{enumerate}

\item If $x\in X$ is such that $f(x)$ is {\qstd} (respectively \std), then $x$
is also {\qstd} (respectively \std). In particular, if $Y$ is {\qstd} or \std,
then so is $X$.

\item If $f$ is a closed immersion, then $x\in X$ is {\qstd} (respectively
\std) if and only if $f(x)$ is {\qstd} (respectively \std).
\end{enumerate}
        \end{prop}

        \begin{proof}
It follows immediately from \ref{gspecAgood}.
        \end{proof}

        \begin{prop}
Let $X$ be a {\good} graded scheme. Then
\begin{align*}
& \QT{X}:=\{x\in X\st x\text{ \qtriv}\}, & \QS{X}:=\{x\in X\st x\text{ \qstd}\}
\end{align*}
are closed and open in $X$ and $X=\QT{X}\djun\QS{X}$. In particular, a
connected {\good} graded scheme is either {\qstd} or \qtriv.
        \end{prop}

        \begin{proof}
It follows immediately from \ref{clopen}, taking into account \ref{locgsch}.
        \end{proof}

        \begin{rema}
We have already noticed that a {\triv} graded scheme is just a usual noetherian
scheme, and it is very easy to see that a graded scheme $X$ is {\qtriv} if and
only if $X_0$ is a noetherian scheme and $(\so_X)_d=0$ for $\abs{d}>>0$.
        \end{rema}

Also {\std} graded schemes are easy to characterize, as the following result
shows.

        \begin{prop}\label{stdchar}
If $X$ is a noetherian scheme and $\s{L}$ an invertible sheaf on $X$, then
$(X,\bigoplus_{d\in\Z}\s{L}^{\otimes d})$ is a {\std} graded scheme, which will
be denoted by $\stdsch{X}{\s{L}}$.\index{X(L)@$\stdsch{X}{\s{L}}$}

Conversely, if $Z$ is a {\std} graded scheme, then there is an invertible sheaf
$\s{L}$ on $Z_0$ such that $Z\iso\stdsch{Z_0}{\s{L}}$.

If $\stdsch{X}{\s{L}}$ and $\stdsch{Y}{\s{M}}$ are {\std} graded schemes,
then to give a morphism of graded schemes from $\stdsch{X}{\s{L}}$ to
$\stdsch{Y}{\s{M}}$ is equivalent to giving a morphism of schemes $f:X\to Y$
together with an isomorphism $\varphi:f^*\s{M}\isomor\s{L}$ of invertible
sheaves on $X$.
        \end{prop}

        \begin{proof}
If $X$ is a noetherian scheme and $\s{L}$ is an invertible sheaf on $X$, then
$\all x\in X$ there is an affine open neighbourhood $U\iso\spec A$ of $x$
where $\s{L}$ is trivial. Then $\stdsch{X}{\s{L}}\rest{U}\iso
\gspec A[t,t^{-1}]$ (with $\deg(t)=1$), whence $\stdsch{X}{\s{L}}$ is \std.

Conversely, if $Z$ is \std, then $\s{L}:=(\so_Z)_1$ is an invertible sheaf on
$Z_0$ and the natural morphism
$(\so_Z)_{-1}\to\s{L}^{\otimes-1}=\sHom_{Z_0}(\s{L},\so_{Z_0})$ is an
isomorphism (because it is an isomorphism on each stalk). It follows that there
is a natural morphism of sheaves of graded rings
$\bigoplus_{d\in\Z}\s{L}^{\otimes d}\to\so_Z$, which is an isomorphism (again,
because it is an isomorphism on each stalk), so that $Z\iso
\stdsch{Z_0}{\s{L}}$ as graded schemes.

If $g:\stdsch{X}{\s{L}}\to\stdsch{Y}{\s{M}}$ is a morphism of {\std} graded
schemes, then $f:=(g,g\mrs_0):X\to Y$ is a morphism of schemes and the morphism
of $\so_Y$--modules $g\mrs_1:\s{M}\to f_*\s{L}$ corresponds by adjunction to a
morphism $\varphi:f^*\s{M}\to\s{L}$ of $\so_X$--modules, which is easily seen
to be an isomorphism. Conversely, given $f:X\to Y$ and an isomorphism
$\varphi:f^*\s{M}\isomor\s{L}$, we can naturally extend $\varphi$ to an
isomorphism $\tilde{\varphi}:\bigoplus_{d\in\Z}f^*(\s{M}^{\otimes d})\iso
f^*(\bigoplus_{d\in\Z}\s{M}^{\otimes d})\isomor
\bigoplus_{d\in\Z}\s{L}^{\otimes d}$ of sheaves of graded rings on $X$, which,
by adjunction, corresponds to a morphism
$\alpha:\bigoplus_{d\in\Z}\s{M}^{\otimes d}\to
f_*\bigoplus_{d\in\Z}\s{L}^{\otimes d}$ of sheaves of graded rings on $Y$, thus
yielding a morphism of graded schemes
$g=(f,\alpha):\stdsch{X}{\s{L}}\to\stdsch{Y}{\s{M}}$.
        \end{proof}

In the following we will be mainly interested in {\qstd} graded schemes. Even
when we have an ordinary noetherian scheme, we will often regard it as a {\std}
graded scheme by choosing a line bundle on it (in the problems we are going to
consider, there will be a natural way of making such a choice).

For the rest of this section $X$ will be a {\good} graded scheme: we are going
to study the relation between graded sheaves on $X$ and sheaves on $X_0$. We
will denote by $\smo{X_0}$\index{Mod(X_0)@$\smo{X_0}$} the category of sheaves
of $\so_{X_0}$--modules and by $\qco{X_0}$\index{Qcoh(X_0)@$\qco{X_0}$}
(respectively $\coh{X_0}$)\index{Coh(X_0)@$\coh{X_0}$} its full subcategory of
quasi--coherent (respectively coherent) sheaves.

Sending a graded sheaf of $\so_X$--modules $\s{F}$ to $\s{F}_0$ defines an
exact functor
\[\usmf{X}:\gsmo{X}\to\smo{X_0}.\]
In analogy with the case of graded modules over a graded ring, we can prove the
following result.

        \begin{prop}\label{adjsheaf}
The functors
\begin{align*}
& \osmf{X}:\smo{X_0}\to\gsmo{X}, & \hsmf{X}:\smo{X_0}\to\gsmo{X}
\end{align*}
are respectively left and right adjoint of $\usmf{X}$, and they are both right
inverses of $\usmf{X}$.

Moreover, $\usmf{X}$ sends $\gqco{X}$ (respectively $\gcoh{X}$) into
$\qco{X_0}$ (respectively $\coh{X_0}$). $\osmf{X}$ and $\hsmf{X}$ send
$\qco{X_0}$ (respectively $\coh{X_0}$) into $\gqco{X}$ (respectively
$\gcoh{X}$).
        \end{prop}

        \begin{proof}
It follows easily from \ref{adj0}.
        \end{proof}

        \begin{rema}
By \ref{locsubcat} the above result implies that $\smo{X_0}$ is equivalent to
$\gsmo{X}/\cat{C}$, where $\cat{C}$ is the full subcategory of $\gsmo{X}$ whose
objects are the graded sheaves $\s{F}$ such that $\s{F}_0=0$. Notice that
$\cat{C}$ is both left and right localizing. The same result holds, of course,
also for quasi--coherent and coherent sheaves.
        \end{rema}

        \begin{coro}
$\osmf{X}$ and $\hsmf{X}$ identify $\smo{X_0}$ (respectively $\qco{X_0}$,
respectively $\coh{X_0}$) with an additive (but not abelian in general)
subcategory of $\gsmo{X}$ (respectively $\gqco{X}$, respectively $\gcoh{X}$).
        \end{coro}

        \begin{coro}\label{esssurj}
$\usmf{X}:\gsmo{X}\to\smo{X_0}$ is essentially surjective. The same is true for
its restriction to quasi--coherent and coherent sheaves.
        \end{coro}

        \begin{prop}\label{shequivcrit}
$\usmf{X}:\gsmo{X}\to\smo{X_0}$ is an equivalence of categories if and only if
$X$ is \std.

If this is the case, then the restriction of $\usmf{X}$ also gives an
equivalence between $\gqco{X}$ (respectively $\gcoh{X}$) and $\qco{X_0}$
(respectively $\coh{X_0}$). Moreover, $\osmf{X}$ and $\hsmf{X}$ are isomorphic,
exact and quasi--inverse of $\usmf{X}$.
        \end{prop}

        \begin{proof}
By \ref{adjsheaf} $\usmf{X}$ is an equivalence if and only if $\all \s{F}
\in\gsmo{X}$ the natural morphism of $\so_X$--modules $\s{F}_0
\otimes_{\so_{X_0}}\so_X\to\s{F}$
is an isomorphism. Then everything follows easily from \ref{modequivcrit}.
        \end{proof}

        \begin{prop}\label{goodcohom}
$\all\cp{\s{F}}\in D^+(\gsmo{X})$ and $\all i\in\Z$ there is a natural
isomorphism $H^i(X,\cp{\s{F}})\iso H^i(X_0,\cp{\s{F}}_0)$.
        \end{prop}

        \begin{proof}
Let $\cp{\s{I}}\in D^+(\gsmo{X})$ be an injective resolution of $\cp{\s{F}}$,
so that (by definition)
\[H^i(X,\cp{\s{F}})\iso H^i(\gsec(X,\cp{\s{I}})_0)\iso
H^i(\Gamma(X_0,\cp{\s{I}_0})).\]
On the other hand, since each $\s{I}^j$ is flasque by \ref{injflasque} and
since $\usmf{X}$ is an exact functor which sends flasque sheaves into flasque
sheaves, $\cp{\s{I}_0}$ is a flasque resolution of $\cp{\s{F}}_0$ on $X_0$,
which can be used to compute cohomology. Therefore
$H^i(\Gamma(X_0,\cp{\s{I}_0}))\iso H^i(X_0,\s{F}_0)$.
        \end{proof}

                \section{$\gproj$ of a noetherian positively graded ring}

We will denote by $\prng$\index{PRing@$\prng$} the full subcategory of the
category of graded rings whose objects are noetherian positively graded rings
(i.e., noetherian graded rings $R$ such that $R_d=0$ for $d<0$).

        \begin{rema}
It follows from \ref{Anoeth} that $R\in\prng$ is isomorphic to a quotient of a
polynomial ring $R_0[x_0,\dots,x_n]$, where $\deg(x_i)>0$ and $R_0$ is a
noetherian ring.
        \end{rema}

If $A$ is a noetherian ring and $\w=(\w_0,\dots,\w_n)\in\N_+^{n+1}$, the
noetherian positively graded ring $A[x_0,\dots,x_n]$ where $\deg(x_i)=w_i$ will
be denoted by $\p_A(\w)$ (or simply by $\p(\w)$\index{P(w)@$\p$, $\p(\w)$} if
there can be no doubt about $A$).

        \begin{defi}
Given $R\in\prng$, $\gproj R$\index{Proj@$\gproj$} is the open graded subscheme
$\gspec R\minus\V{R_+}$ of $\gspec R$.
        \end{defi}

        \begin{rema}
Since $\gproj R=\bigcup_{r\in R_+\hom}\D{r}$ and (if $r$ is not nilpotent)
$\D{r}\iso\gspec R_r$, we see that $\gproj R$ is a {\qstd} graded scheme
(actually it is just the open graded subscheme of $\gspec R$ consisting of
{\qstd} points). It is also clear that $(\gproj R)_0=\proj R$ and $\so_{\gproj
R}=\bigoplus_{d\in\Z}\so_{\proj R}(d)$.
        \end{rema}

        \begin{exem}
If $R=\p_A(\w)$ for some noetherian ring $A$ and some $\w=(\w_0,\dots,\w_n)\in
\N_+^{n+1}$, $\gproj R$ will be denoted by $\gP_A(\w)$ and $\proj R$ by
$\P_A(\w)$ (again, we will simply write $\gP(\w)$\index{P(w) @$\gP$, $\gP(\w)$}
and $\P(\w)$\index{P(w)@$\P$, $\P(\w)$} if $A$ is clear from the context). It
will be called (graded) weighted projective space (over $A$) of
weights $\w$. Of course, if $\w=(1,\dots,1)$, it will be denoted by
$\gP^n$ or $\P^n$.
        \end{exem}

        \begin{rema}
Clearly $\gproj R$ is a {\std} graded scheme if $R$ is generated by $R_1$ as an
$R_0$--algebra. The converse is false in general (examples are given by the
canonical rings of surfaces of general type, see \ref{canmodstd}), but it is
true if $R$ is of the form $\p(\w)$, as the following result shows.
        \end{rema}

        \begin{prop}\label{stdmori}
Let $A$ be a noetherian ring and $\w=(\w_0,\dots,\w_n)\in\N_+^{n+1}$. Then
$\I{p}\in\gP_A(\w)$ is {\std} if and only if $\gcd(\w_i\st x_i\notin\I{p})=1$.
In particular, $\gP_A(\w)$ is {\std} if and only if $\w=(1,\dots,1)$.
        \end{prop}

        \begin{proof}
By \ref{Apqstd} it is enough to prove that $I(\I{p})=(\w_i\st x_i\notin
\I{p})$. Obviously $\w_i\in I(\I{p})$ if $x_i\notin\I{p}$. Conversely, if $d>0$
is such that $\I{p}_d\ne\p_A(\w)_d$, then there exists a monomial $x_0^{m_0}
\cdots x_n^{m_n}\in\p_A(\w)_d\minus\I{p}$. It follows that
$d=\sum_{i=0}^{n}m_i\w_i\in I(\I{p})$, since $x_i\notin\I{p}$ if $m_i>0$.

The last statement follows at once, taking into account that $\all i=0,\dots,n$
there exists $\I{p}\in\gP_A(\w)$ such that $(x_0,\dots,\hat{x_i},\dots,x_n)
\subseteq\I{p}$.
        \end{proof}

    \begin{rema}
The above result says that the (open) subset of {\std} points in a weighted
projective space coincides with the regular locus which was first considered in
\cite{M}: it is the locus where ``everything behaves well''(see also
\cite{BR}), and it is contained in the smooth locus (in general strictly, e.g.
it is empty if $\gcd(\w_0,\dots,\w_n)\ne1$).
    \end{rema}

        \begin{prop}\label{projmor}
Let $A$ be a noetherian ring, let $X\in\gsch_A$ and, for $i=0,\dots,n$, let
$s_i\in\gsec(X,\so_X)_{\w_i}$ for some $\w=(\w_0,\dots,\w_n)\in\N_+^{n+1}$.
Then the $s_i$ induce a rational map $f:X\rat\gP_A(\w)$ in $\gsch_A$, defined
on the open subset
\[U:=\{x\in X\st ((s_0)_x,\dots,(s_n)_x)=\so_{X,x}\}\subseteq X,\]
such that $f^*(x_i)=s_i\rest{U}$ for $i=0,\dots,n$. In particular, $f$ is a
morphism (i.e. $U=X$) if and only if $\{s_0,\dots,s_n\}$ generate $\so_X$ (and
conversely, of course, given a morphism $g:X\to\gP_A(\w)$ in $\gsch_A$, the
sections $\{g^*(x_0),\dots,g^*(x_n)\}$ generate $\so_X$).
        \end{prop}

        \begin{proof}
The $s_i$ clearly determine a morphism of graded $A$-algebras
$\varphi:\p(\w)=\p_A(\w) \to\gsec(X,\so_X)$ (defined by $\varphi(x_i):=s_i$),
which (by \ref{gspecequiv}) corresponds to a morphism $f':X\to\gspec\p(\w)$ in
$\gsch_A$. Since $f'$, as a map of topological spaces, is defined by
$f'(x):=\{a\in\p(\w)\st \varphi(a)_x\in\I{m}_x\}$, it is easy to see that
${f'}^{-1}(\gP(\w))=U$, and then everything follows.
        \end{proof}

        \begin{rema}
$U$ is {\qstd} (even {\std} if $\w=(1,\dots,1)$) by \ref{YgoodXgood}, and in
this case the sections $s^i$ also determine a rational map $f_0:X_0\rat
\P_A(\w)$ (over $A$), defined on $U_0$. Notice that if $X$ is {\std} (say, by
\ref{stdchar}, $X\iso\stdsch{Y}{\s{L}}$ for some noetherian scheme $Y$ and some
invertible sheaf $\s{L}$ on $Y$), then $s_i\in\Gamma(Y,\s{L}^{\w_i})$.
Therefore \ref{projmor} implies in particular the well known result
that to give a morphism of (noetherian) schemes over $A$ from $Y$ to
$\P^n_A$ is equivalent to giving an invertible sheaf on $Y$ and $n+1$
global sections which generate it.
        \end{rema}

If $R\in\prng$ and $M\in\gmo{R}$, by abuse of notation we will denote
$M\gsh{R}\rest{\gproj R}$ again by $M\gsh{R}$ (or simply by
$M\gsh{}$,\index{$\gsh{R}$} if there is no doubt about the ring). Also, if
$\s{F}\in\gsmo{\gproj R}$, we will usually write
$\gsec(\s{F})$\index{Gammma(F)@$\gsec(\s{F})$} instead of $\gsec(\gproj
R,\s{F})$.

A morphism $\varphi:R\to S$ in $\prng$ induces a rational map of graded schemes
$f:\gproj S\rat\gproj R$, defined on the open subset $U:= \{\I{p}\in\gproj
S\st\varphi(R_+)\nsubseteq\I{p}\}$ of $\gproj S$ ($f$ is just the restriction
of the induced morphism $f':\gspec S\to\gspec R$). Notice that $f$ is affine,
since $f'$ is. Observe also that $f$ is a morphism (i.e. $U=\gproj S$) if and
only if $\sqrt{(\varphi(R_+))}\supseteq S_+$. The following result is an
immediate consequence of \ref{modaff}.

        \begin{prop}\label{affproj}
Let $\varphi:R\to S$ be a morphism in $\prng$ and let $f:U\subseteq\gproj S\to
X:=\gproj R$ be the induced morphism of graded schemes.
\begin{enumerate}

\item If $M,N\in\gmo{R}$, then $(M\otimes_R N)\gsh{}\iso M\gsh{}\otimes_X
N\gsh{}$.

\item If $M\in\gmo{R}$ and $N\in\gmo{S}$, then $f^*(M\gsh{R})\iso
(M\otimes_RS)\gsh{S}\rest{U}$ and $f_*(N\gsh{S}\rest{U})\iso N\gsh{R}$.
\end{enumerate}
        \end{prop}

        \begin{rema}
The above result (except that $f_*(N\sh{S}\rest{U})\iso N\sh{R}$ still holds)
is not true in general with $\proj$ instead of $\gproj$.
        \end{rema}

        \begin{prop}\label{prngmor}
Let $\varphi:R\to S$ be a morphism in $\prng$ and let $f:\gproj S\rat\gproj R$
be the induced rational map. Then the following conditions are equivalent:
\begin{enumerate}

\item $\varphi:R\to S$ is finite;

\item $f$ is a morphism and $\varphi_0:R_0\to S_0$ is finite.

\end{enumerate}
If these conditions are satisfied, then $f:\gproj S\to\gproj R$ (and also
$f_0:\proj S\to\proj R$) is finite.
        \end{prop}

        \begin{proof}
\begin{description}

\item [$1\implies2$] By hypothesis there exist $s_1,\dots,s_n\in S\hom$ which
generate $S$ as an $R$--module. Now $\all s\in S\hom_+$ we can choose $m>0$
such that $m\deg(s)>\deg(s_i)$ for $i=1,\dots,n$. Then there exist $r_1,
\dots,r_n\in R\hom_+$ such that $s^m=\sum_{i=1}^nr_is_i$, whence $s\in
\sqrt{(\varphi(R_+))}$, which proves that $f$ is a morphism. Moreover, the
$s_i$ of degree $0$ generate $S_0$ as an $R_0$--module.

\item [$2\implies1$] Since $S$ is a finitely generated $S_0$--algebra and
$\sqrt{(\varphi(R_+))}\supseteq S_+$, it is easy to see that there exists $m>0$
such that $(\varphi(R_+))_{>m}=S_{>m}$. By \ref{Anoeth} and since $\varphi_0$
is finite, there exist $s_1,\dots,s_n\in S\hom$ which generate
$\bigoplus_{0\le d\le m}S_d$ as an $R_0$--module, and then it is not difficult
to prove that $s_1,\dots,s_n$ generate $S$ as an $R$--module.

\end{description}
As for the last statement, it is enough to observe that $\all r\in R\hom_+$ non
nilpotent
\[f\rest{f^{-1}(\D{r})}:f^{-1}(\D{r})=\D{\varphi(r)}\iso
\gspec S_{\varphi(r)}\to\D{r}\iso\gspec R_r\] is the morphism of affine graded
schemes induced by the finite morphism of graded rings $\varphi_r:R_r\to
S_{\varphi(r)}$.
        \end{proof}

For the rest of this section $R$ will be a noetherian positively graded ring,
$\gproj R$ will be denoted by $X$ and $R_0$ by $A$. We are going to see that a
generalized version of Serre's correspondence holds between (finitely
generated) graded $R$--modules and (quasi)coherent sheaves on $X$ (see also
\cite{B} and \cite{AKO}).

        \begin{prop}\label{gmogqco}
\begin{enumerate}

\item $\gsec(-):\gsmo{X}\to\gmo{R}$ is right adjoint of the exact functor
$\gsh{}:\gmo{R}\to\gsmo{X}$.

\item If $M\in\gmo{R}$ (respectively $\gfmo{R}$), then $M\gsh{}$ is
quasi--coherent (respectively coherent).

\item If $\s{F}\in\gqco{X}$, then the natural map $\gsec(\s{F})\gsh{}\to\s{F}$
is an isomorphism.

\end{enumerate}
        \end{prop}

        \begin{proof}
As for $1$, denoting by $j:X\mono\gspec R:=Y$ the inclusion, $\all M\in\gmo{R}$
and $\all\s{F}\in\gsmo{X}$ there are natural isomorphisms
\[\Hom_R(M,\gsec(\s{F}))\iso\Hom_R(M,\gsec(Y,j_*\s{F}))\iso
\Hom_Y(M\gsh{},j_*\s{F})\iso\Hom_X(M\gsh{},\s{F})\] (using the fact that
$\gsec(Y,-)$ is right adjoint of $\gsh{}$ and $j^*$ is left adjoint of $j_*$).
$2$ is clear, since the same is true for $M\gsh{}$ on $Y$, and the proof of $3$
is similar to that of \cite[II, prop. 5.15]{H1} (using \ref{secext} instead of
\cite[II, lemma 5.14]{H1}).
        \end{proof}

If $\s{F}\in\gcoh{X}$, it is not true in general that $\gsec(\s{F})\in\gfmo{R}$
(it can happen that $\gsec(\s{F})_d\ne0$ $\all d\in\Z$, for instance if $\s{F}$
is supported at a point), but, as we are going to see, the truncations
$\gsec(\s{F})_{>d}$ are finitely generated.

        \begin{lemm}\label{lbcohom}
Given $\w\in\N_+^{n+1}$, let $\p(\w)=\p_A(\w)$ and $\gP(\w)=\gP_A(\w)$.
\begin{enumerate}

\item The natural map $\p(\w)\to\gsec(\so_{\gP(\w)})$ is an isomorphism in
$\gmo{\p(\w)}$; in particular, $H^0(\gP(\w),\so_{\gP(\w)}(m))\iso\p(\w)_m$
$\all m\in\Z$.

\item $H^i(\gP(\w),\so_{\gP(\w)}(m))=0$ for $0<i<n$ and $\all m\in\Z$.

\item $H^n(\gP(\w),\so_{\gP(\w)}(m))\iso\p(\w)_{-m-\sum_{i=0}^n\w_i}$.
\end{enumerate}
        \end{lemm}

        \begin{proof}
Since $H^i(\gP(\w),\so_{\gP(\w)}(m))\iso H^i(\P(\w),\so_{\P(\w)}(m))$ by
\ref{goodcohom}, one can just consider ordinary weighted projective spaces, and
then the result is proved, for instance, in \cite[thm. 3.1]{D}.
        \end{proof}

Using \ref{lbcohom} one can then prove the following result, with a proof
similar to those of \cite[II, prop. 5.17 and III, thm. 5.2]{H1}.

        \begin{prop}\label{cohcohom}
Every $\s{F}\in\gcoh{X}$ satisfies the following properties:

\begin{enumerate}

\item $\s{F}$ is generated by a finite number of homogeneous global sections;

\item $H^i(X,\s{F})$ is a finitely generated $A$--module $\all i\in\N$;

\item $\exi m_0\in\Z$ such that $H^i(X,\s{F}(m))=0$ for $m\ge m_0$ and
$\all i>0$.
\end{enumerate}
        \end{prop}

        \begin{coro}\label{gcohfg}
If $\s{F}\in\gcoh{X}$, then $\gsec(\s{F})_{>d}\in\gfmo{R}$ $\all d\in\Z$.
        \end{coro}

        \begin{proof}
$R=\p(\w)/I$ for some $n\in\N$, some $\w\in\N_+^{n+1}$ and some ideal $I\subset
\p(\w):=\p_A(\w)$, and (by \ref{prngmor}) the projection $R\to\p(\w)$ induces a
finite morphism $f:X\to\gP(\w)$ (which is actually a closed immersion), whence
$f_*\s{F}\in\gcoh{\gP(\w)}$. Moreover, since $\gsec(X,\s{F})\iso\gsec(\gP(\w),
f_*\s{F})$ and since an $R$--module is finitely generated if and only if it is
finitely generated as a $\p(\w)$--module, we can assume $R=\p(\w)$ and $X=
\gP(\w)$.

By part 1 of \ref{cohcohom} there is an exact sequence of the form
\[0\to\s{R}\to\s{E}\to\s{F}\to0\]
with $\s{E}=\bigoplus_{j=1}^r\so_{\gP(\w)}(-m_j)$ for some $r\in\N$ and some
$m_j\in\Z$. Since $\s{R}$ is also coherent, by part 3 of \ref{cohcohom} $\exi
m_0\in\Z$ such that $H^1(\gP(\w),\s{R}(m))=0$ for $m\ge m_0$, so that the map
$H^0(\gP(\w),\s{E}(m))\to H^0(\gP(\w),\s{F}(m))$ is surjective. From this and
from part 1 of \ref{lbcohom} it follows that $M:=\gsec(\s{F})$ is generated by
$M_{<l}$, where $l>m_j$ for $j=0\dots,r$. Now, $\all d\in\Z$ take $d'\in\Z$
such that $d'\ge l$ and $d'>d+\w_j$ for $j=0\dots,n$. Then it is clear that
$M_{>d}$ is generated by $\bigoplus_{d<j<d'}M_j$, hence is finitely generated
(because every $M_j\in\fmo{A}$ by part 2 of \ref{cohcohom}).
        \end{proof}

We include here also a result which will be needed later, and whose proof is
completely similar to that of \cite[III, prop. 6.9]{H1}.

        \begin{lemm}\label{secgsExt}
$\all\s{F},\s{G}\in\gcoh{X}$ and $\all i\in\N$ there exists $m\in\Z$ such that
\begin{align*}
& \Ext^i_X(\s{F},\s{G}(d))\iso H^0(X,\gsExt^i_X(\s{F},\s{G}(d))) & \all d>m.
\end{align*}
        \end{lemm}

We are now ready to prove the generalized version of Serre's correspondence;
see \ref{abelquot} for the notion of localizing subcategory and of quotient
of an abelian category by a localizing subcategory.

        \begin{prop}\label{Serre}
Let $\cat{Z}(R)$ (respectively $\cat{z}(R)$) be the full right localizing
subcategory of $\gmo{R}$ (respectively $\gfmo{R}$) whose objects are the
modules $M$ such that $M\gsh{}=0$, and let $T:\gmo{R}\to\gmo{R}/\cat{Z}(R)$
(respectively $T':\gfmo{R}\to\gfmo{R}/\cat{z}(R)$) be the natural functor. Then
the functor $\gsh{}:\gmo{R}\to\gsmo{X}$ induces an exact equivalence of
categories $F:\gmo{R}/\cat{Z}(R)\to\gqco{X}$ (respectively $F':
\gfmo{R}/\cat{z}(R)\to\gcoh{X}$), whose quasi--inverse is the functor $T\comp
\gsec(-)$ (respectively $T'\comp\gsec(-)_{>0}$).
        \end{prop}

        \begin{proof}
The statement about $\gmo{R}$ and $\gqco{X}$ is an immediate consequence of
\ref{locsubcat} and \ref{gmogqco}. The other statement follows from
\ref{gcohfg}, if one observes that (by definition of quotient category)
$\gfmo{R}/\cat{z}(R)$ can be identified with a full subcategory of
$\gmo{R}/\cat{Z}(R)$ and that $\all M\in\gmo{R}$ the inclusion $M_{>0}\mono M$
is an isomorphism in $\gmo{R}/\cat{Z}(R)$.
        \end{proof}

        \begin{rema}
It is not difficult to prove (using the fact that $R$ is a noetherian ring)
that $\all M,N\in\gmo{R}$
\[\Hom_{\gmo{R}/\cat{Z}(R)}(M,N)\iso\lim_{M'\in I(M)}\Hom_{\gmo{R}}(M',N),\]
where $I(M):=\{M'\subseteq M\st M/M'\in\cat{Z}(R)\}$. It follows (since if $L
\in\gfmo{R}$, $L\in\cat{z}(R)$ if and only if $L_d=0$ for $d>>0$) that if
$M\in\gfmo{R}$ then
\[\Hom_{\gmo{R}/\cat{Z}(R)}(M,N)\iso\lim_{n\in\N}\Hom_{\gmo{R}}(M_{>n},N).\]
Therefore, in the particular case $X=\gP^n$ (and more generally when $X$ is
\std), considering the equivalence (given by \ref{shequivcrit}) between
$\gcoh{X}$ and $\coh{X_0}$, we obtain the usual description of $\coh{X_0}$, as
proved in \cite{S}.
        \end{rema}

Now we will give also a graded version of Serre duality. Still denoting
$\gproj R$ by $X$ (and $\proj R$ by $X_0$), we assume now that $R_0$ is a
field $\K$ and that $\dim X_0=n$. We recall that on $X_0$ there exists unique
(up to isomorphism) a dualizing sheaf $\ds_{X_0}$ (see
\cite[III, prop.7.5]{H1}).

        \begin{defi}
A dualizing sheaf on $X$ is a coherent graded sheaf
$\ds_X$\index{omega_X@$\ds_X$} together with a map $t:H^n(X,\ds_X)\to\K$ such
that $\all\s{F}\in\gcoh{X}$ the bilinear map
\[\Hom_X(\s{F},\ds_X)\times H^n(X,\s{F})\to H^n(X,\ds_X)\mor{t}\K\]
is a perfect pairing of $\K$--vector spaces (i.e., the induced map
$\Hom_X(\s{F},\ds_X)\to H^n(X,\s{F})\dual$ is an isomorphism).
        \end{defi}

        \begin{prop}\label{Serredual}
A dualizing sheaf on $X$ exists unique (up to isomorphism), and it is given by
$\ds_X=\gsHom_{X_0}(\so_X,\ds_{X_0})$. If moreover $X_0$ is Cohen--Macaulay and
equidimensional and $\sExt^i_{X_0}(\so_{X_0}(m),\ds_{X_0})=0$ for $i>0$ and
$\all m\in\Z$, then $\all i\ge0$ there are natural isomorphisms
$\Ext^i_X(\s{F},\ds_X)\iso H^{n-i}(X,\s{F})\dual$.
        \end{prop}

        \begin{proof}
Uniqueness follows as usual from the uniqueness of an object representing a
functor. Since, given $\s{F}\in\gcoh{X}$, we have $H^n(X,\s{F})\iso
H^n(X_0,\s{F}_0)$ by \ref{goodcohom} and $\Hom_X(\s{F},\ds_X)\iso
\Hom_{X_0}(\s{F}_0,\ds_{X_0})$ by \ref{adjsheaf}, it is clear that $\ds_X$ is a
dualizing sheaf.

As for the last statement, remember that (by \cite[III, thm. 7.6]{H1}) $X_0$ is
Cohen--Macaulay and equidimensional if and only if $\all\s{G}\in\coh{X_0}$
there are natural isomorphisms $\Ext^i_{X_0}(\s{G},\ds_{X_0})\iso
H^{n-i}(X_0,\s{G})\dual$. As before, if $\s{F}\in\gcoh{X}$ then
$H^{n-i}(X,\s{F})\iso H^{n-i}(X_0,\s{F}_0)$, so that it is enough to prove that
$\Ext^i_X(\s{F},\ds_X)\iso\Ext^i_{X_0}(\s{F}_0,\ds_{X_0})$. Let $\cp{\s{I}}$ be
an injective resolution of $\ds_{X_0}$, and observe that then $\cp{\s{J}}:=
\gsHom_{X_0}(\so_X,\cp{\s{I}})$ is an injective resolution of $\ds_X$ (it is
injective because $\hsmf{X}$ preserves injective objects, as it is easy to
check, and it is a resolution because $\gsExt^i_{X_0}(\so_X,\ds_{X_0})\iso
\bigoplus_{m\in\Z}\sExt^i_{X_0}(\so_{X_0}(m),\ds_{X_0})=0$ for $i>0$). Then,
again by \ref{adjsheaf}, we have
\[\Ext^i_X(\s{F},\ds_X)\iso H^i(\Hom_X(\s{F},\cp{\s{J}}))\iso
H^i(\Hom_{X_0}(\s{F}_0,\cp{\s{I}}))\iso\Ext^i_{X_0}(\s{F}_0,\ds_{X_0}).\]
        \end{proof}

                \section{Graded schemes and algebraic stacks}\label{stacks}

We assume the reader is familiar with algebraic stacks (see e.g.
\cite{LM}); here we just recall that, if a (smooth) algebraic group
$G$ acts on a scheme $Z$, the quotient stack $[Z/G]$ is algebraic (the
projection $Z\to[Z/G]$ is a smooth presentation) and $\coh{[Z/G]}$
(respectively $\qco{[Z/G]}$) is equivalent to $G$--$\coh{Z}$
(respectively $G$--$\qco{Z}$), the category of (quasi)coherent
$G$--equivariant sheaves on $Z$.

If $R\in\prng$ and $R_0=\K$ is a field, \ref{Serre} and \cite[prop. 2.3]{AKO}
imply that the category $\gcoh{\gproj R}$ (respectively $\gqco{\gproj R}$) is
equivalent to $\coh{[(\spec R\minus\{R_+\})/\mult]}$ (respectively
$\qco{[(\spec R\minus\{R_+\})/\mult]}$), where the action of
$\mult:=\spec\K[x,x^{-1}]$ is the natural one induced by the graduation of $R$
(such that the quotient scheme is the usual $\proj R$). We are going to see
briefly how this result can be generalized; a more accurate study of the
relation between graded schemes (considered in a more general sense, namely
allowing the grading group to be different from $\Z$) and algebraic stacks will
appear in a future paper.

Let $X$ be a graded scheme, which we assume for simplicity of finite type over
a field $\K$: we are going to define a scheme $\ung{X}$ together with an action
of $\mult$ on it. Assume first that $X$ is affine, say $X=\gspec A$ for some
graded finitely generated $\K$--algebra $A$. To avoid confusion, we denote by
$\ung{A}$ the ring obtained from $A$ by forgetting the graduation (i.e.,
$\ung{A}=\ung{A}_0=A$). The graduation of $A$ yields an action
$\mult\times\spec\ung{A}\iso\spec\ung{A}[x,x^{-1}]\to\spec\ung{A}$ of $\mult$
on $\ung{X}:=\spec\ung{A}$ (induced by the morphism of rings
$\ung{A}\to\ung{A}[x,x^{-1}]$ defined by $a\mapsto ax^{\deg(a)}$ for $a\in
A\hom$). It is then easy to see that for general $X$ we can take an open cover
of $X$ by affine graded schemes $U_i$ and glue the schemes $\ung{U_i}$ to a
scheme $\ung{X}$ with an action $\rho_0:\mult\times\ung{X}\to\ung{X}$: this is
possible (and the result does not depend, up to isomorphism, on the chosen
cover) because the above construction is compatible with localization, in the
sense that, if $a\in A\hom$ (so that $\gspec A_a$ can be identified with an
open graded subscheme of $\gspec A$), then $\spec\ung{A_a}$ can be identified
with an open subscheme of $\spec\ung{A}$ and the action of $\mult$ on
$\spec\ung{A_a}$ is the restriction of the action on $\spec\ung{A}$. It is also
clear that, if $X=\gproj R$, then $\ung{X}=\spec\ung{R}\setminus\V{R_+}$ and
$\rho_0$ is the restriction of the action on $\spec\ung{R}$.

        \begin{prop}
There is an equivalence of categories between $\gcoh{X}$ (respectively
$\gqco{X}$) and $\coh{[\ung{X}/\mult]}$ (respectively $\qco{[\ung{X}/\mult]}$).
    \end{prop}

    \begin{proof}
We claim that $\rho_0:\mult\times\ung{X}\to\ung{X}$ can be extended to
an action $\rho:\mult\times\tilde{X}\to\tilde{X}$ of $\mult$ on the
{\std} graded scheme $\tilde{X}:=\stdsch{\ung{X}}{\so_{\ung{X}}}$
and that there is a morphism of graded schemes $p:\tilde{X}\to X$ such
that the diagram
\[\begin{CD}
\mult\times\tilde{X} @>{\rho}>> \tilde{X} \\
@V{\pi}VV @VVpV \\
\tilde{X} @>>p> X
\end{CD}\]
(where $\pi$ is the projection) is cartesian in $\gsch_{\K}$ (and also
cocartesian, so that $X$ is the quotient $\tilde{X}/\mult$ as graded scheme).
If $X=\gspec A$ is affine, $\tilde{X}\iso\gspec\ung{A}[t,t^{-1}]$,
$\mult\times\tilde{X}\iso\gspec\ung{A}[t,t^{-1},x,x^{-1}]$ (with $\deg(t)=1$
and $\deg(x)=0$) and $\pi$ and $\rho$ are induced by the morphisms of graded
rings $\iota,\varphi:\ung{A}[t,t^{-1}]\to \ung{A}[t,t^{-1},x,x^{-1}]$ defined
by $\iota(at^i):=at^i$ and $\varphi(at^i):=at^ix^{\deg(a)-i}$ for $a\in A\hom$
and $i\in\Z$. We take as $p$ the morphism induced by the morphism of graded
rings $\alpha:A\to\ung{A}[t,t^{-1}]$ defined by $a\mapsto at^{\deg(a)}$ for
$a\in A\hom$. Indeed, it is straightforward to check that the diagram of graded
rings
\[\begin{CD}
\ung{A}[t,t^{-1},x,x^{-1}] @<{\varphi}<< \ung{A}[t,t^{-1}] \\
@A{\iota}AA @AA{\alpha}A \\
\ung{A}[t,t^{-1}] @<<{\alpha}< A
\end{CD}\]
is cocartesian (and also cartesian), and the general claim follows easily from
this.

Then, extending classical results of faithfully flat descent theory to the
graded setting, the fact that $p$ is faithfully flat (and even smooth) because
the same is true for $\alpha$ (since $\ung{A}[t,t^{-1}]$, as a graded
$A$--module via $\alpha$, is isomorphic to $\bigoplus_{i\in\Z}A(i)$) implies
that $p^*$ gives an equivalence of categories between  $\gcoh{X}$ (respectively
$\gqco{X}$) and $\mult$--$\gcoh{\tilde{X}}$ (respectively
$\mult$--$\gqco{\tilde{X}}$). As $\tilde{X}$ is \std, the natural
generalization of \ref{shequivcrit} to the case of $\mult$--equivariant sheaves
allows to conclude that $\mult$--$\gcoh{\tilde{X}}$ (respectively
$\mult$--$\gqco{\tilde{X}}$) is equivalent to $\mult$--$\coh{\ung{X}}$
(respectively $\mult$--$\qco{\ung{X}}$), hence to $\coh{[\ung{X}/\mult]}$
(respectively $\qco{[\ung{X}/\mult]}$).
    \end{proof}

    \begin{rema}
If $X$ is {\qstd} (in particular, if $X$ is of the form $\gproj R$) then
$[\ung{X}/\mult]$ is a Deligne--Mumford stack (and even a scheme if $X$ is
\std), at least if $\K$ is algebraically closed of characteristic $0$. In fact,
we can reduce to the case $X=\gspec A$ with $A$ \qstd, say $s\in A_d$ is
invertible for some $d>0$. Denoting by
$\roots{d}:=\spec\K[x]/(x^d-1)\iso\Z/d\Z$ the subgroup of $\mult$ of $d^{\rm
th}$ roots of unit, the above described action of $\mult$ induces an action of
$\roots{d}$ on $\spec(\ung{A}/(s-1))$ (which can be again extended to an action
of $\roots{d}$ on $\gspec(\ung{A}/(s-1)[t,t^{-1}])$), and it can be proved that
$[\spec\ung{A}/\mult]$ is isomorphic to $[\spec(\ung{A}/(s-1))/\roots{d}]$ (for
which the projection $\spec(\ung{A}/(s-1))\to[\spec(\ung{A}/(s-1))/\roots{d}]$
is an \'etale presentation, and obviously an isomorphism if $d=1$). Observe
that the morphism of graded rings (induced by $\alpha$)
$A\to\ung{A}/(s-1)[t,t^{-1}]\iso A[t]/(t^d-s)$ is \'etale because $s$ is
invertible.

In particular, in the case of weighted projective space $\gP(\w)$ one obtains a
smooth Deligne--Mumford stack locally isomorphic to
$[\mathbb{A}^n/\roots{\w_i}]$ (this is the point of view of \cite{K}, where
important results about bounded derived categories of coherent sheaves are
extended from smooth varieties to smooth Deligne--Mumford stacks).
    \end{rema}

                        \chapter{Beilinson's theorem on $\gP(\w)$}

Throughout this chapter we will consider (graded) weighted projective spaces
over a fixed field $\K$.\footnote{Probably much of what follows is still true
if $\K$ is replaced by an arbitrary noetherian ring.} Therefore, given $n\in\N$
and $\w\in\N_+^{n+1}$, we will write $\p(\w)$,\index{P(w)@$\p$, $\p(\w)$}
$\P(\w)$\index{P(w)@$\P$, $\P(\w)$} and $\gP(\w)$\index{P(w) @$\gP$, $\gP(\w)$}
(or simply $\p$, $\P$ and $\gP$ if no confusion is possible) instead of
$\p_{\K}(\w)=\K[x_0,\dots,x_n]$ (with $\deg(x_i)=\w_i$),
$\P_{\K}(\w)=\proj\p(\w)$ and $\gP_{\K}(\w)=\gproj\p(\w)$. $\all m\in\Z$ we set
moreover $\pd_m=\pd_m(\w):=\dim_{\K}\p(\w)_m$.\index{p_m(w)@$\pd_m$,
$\pd_m(\w)$} $\w$ is called {\em normalized}\index{normalized} if
$\gcd(\w_0,\dots,\hat{\w}_i,\dots,\w_n)=1$ for $i=0,\dots,n$. We recall that
$\P(\w)$ is always isomorphic to $\P(\w')$ for some $\w'$ normalized (of
course, this is not true for $\gP(\w)$). If $I \subseteq\{0,\dots,n\}$, we set
$\sw{\w_I}:=\sum_{i\in I}\w_i$;\index{w_I@$\sw{\w_I}$}
$\sw{\w_{\{0,\dots,n\}}}$ will be usually denoted by
$\sw{\w}$.\index{w@$\sw{\w}$} We will write $\compl{I}$\index{I^c@$\compl{I}$}
for $\{0,\dots,n\}\minus I$.

We will extend to weighted projective spaces a classical theorem by Beilinson
(see \cite{B1}), which gives equivalences between $D^b(\coh{\P^n})$ (the
bounded derived category of coherent sheaves on $\P^n$) and two homotopy
categories of modules. More precisely, let $\p:=\p(1,\dots,1)$ and let
$\Lambda=\bigoplus_{i\in\Z}\Lambda_i$\index{Lambda@$\Lambda$} be the graded
(skew--commutative) ring defined by $\Lambda_i:=\Lambda^i(\p_1\dual)$.
Moreover, if $A$ is a graded ring, we will denote by $\mco{[-n,0]}(A)$ the full
subcategory of $\gfmo{A}$ whose objects are finite direct sums of graded
$A$--modules of the form $A(j)$ for $-n\le j\le0$. Then the theorem can be
stated in the following way.

        \begin{Beil}
The natural additive functors $\mco{[-n,0]}(\p)\to\coh{\P^n}$ and
$\mco{[-n,0]}(\Lambda)\to\coh{\P^n}$ (defined on objects by $\p(j)\mapsto
\so(j)$ and $\Lambda(j)\mapsto\sd^{-j}(-j)$ for $-n\le j\le0$) extend to exact
functors between triangulated categories
\begin{align*}
& \fo:K^b(\mco{[-n,0]}(\p))\to D^b(\coh{\P^n}), &
\fde:K^b(\mco{[-n,0]}(\Lambda))\to D^b(\coh{\P^n}).
\end{align*}
$\fo$ and $\fde$ are equivalences of categories.
        \end{Beil}

In particular, $D^b(\coh{\P^n})$ is generated (as a triangulated category) by
each of the following two sets of vector bundles:
\begin{align*}
& \oset:=\{\so_{\P^n}(j)\st-n\le j\le0\},
& \dset:=\{\sd_{\P^n}^j(j)\st0\le j\le n\}.
\end{align*}
This means that $\all\cp{\s{F}}\in D^b(\coh{\P^n})$ there exist bounded
complexes $\cp{\bro}=\cp{\bro}(\cp{\s{F}})$ and $\cp{\brd}=
\cp{\brd}(\cp{\s{F}})$ isomorphic to $\cp{\s{F}}$ (in $D^b(\coh{\P^n})$) and
with each $\bro^i$ (respectively $\brd^i$) isomorphic to a finite direct sum of
sheaves of $\oset$ (respectively $\dset$). $\cp{\bro}$ and $\cp{\brd}$ are also
unique up to isomorphism in $C^b(\coh{\P^n})$ if we require them to be minimal,
and in this case their terms are given explicitly by
\begin{align*}
\bro^i & =\bigoplus_{0\le j\le n}\so(-j)\otimes_{\K}
H^{i+j}(\P^n,\cp{\s{F}}\otimes\sd^j(j)),\\
\brd^i & =\bigoplus_{0\le j\le n}\sd^j(j)\otimes_{\K}
H^{i+j}(\P^n,\cp{\s{F}}(-j)).
\end{align*}

As we said in the introduction, it turns out that, in order to get a good
version of this theorem on weighted projective spaces, it is better to consider
$\gcoh{\gP(\w)}$ instead of $\coh{\P(\w)}$ (remember that, by \ref{shequivcrit}
and \ref{stdmori}, these two categories are equivalent if and only if
$\w=(1,\dots,1)$). The reason of this is that working with graded sheaves on
$\gP$ one avoids all the pathologies of $\P$ (basically because the sheaves of
the form $\so_{\gP}(j)$ are always invertible, whereas the $\so_{\P}(j)$ are
not, in general).

In the weighted case the two sets of vector bundles defined above will be
replaced by the following two sets\index{O@$\goset$}\index{D@$\gdset$}
\begin{align*}
\goset & :=\{\so_{\gP}(j)\st-\sw{\w}<j\le0\}, \\
\gdset & :=\{\so_{\gP}(j)\st n-\sw{\w}<j<0\}\cup \{\sd_{\gP}^j(j)\st0\le j\le
n\}
\end{align*}
(notice that they are the same as before in the particular case $\w=
(1,\dots,1)$), where $\sd_{\gP}^j$ are defined as the graded sheaves associated
to the $\p$--modules obtained as syzygies of the Koszul complex $\cp{\ko}$ of
the regular sequence $(x_0,\dots,x_n)$. Then one can define similar functors to
$\fo$ and $\fde$, which one must prove to be equivalences. As in the case of
$\P^n$, fully faithfulness follows from the fact that
$\Ext_{\gP}^i(\s{A},\s{B})=0$ for $i>0$, $\all\s{A},\s{B}\in\goset$ and
$\all\s{A},\s{B}\in\gdset$ (here it is essential to use $\gcoh{\gP(\w)}$
instead of $\coh{\P(\w)}$), whereas essential surjectivity corresponds to the
fact that $\goset$ and $\gdset$ generate $D^b(\gcoh{\gP})$ as a triangulated
category. Now, this last fact follows easily from Hilbert's syzygy theorem, but
it must be said that Beilinson proved it using a resolution of the diagonal in
$\P^n\times\P^n$, and we don't know if a similar resolution exists also in the
weighted case. The drawback of this fact is that we cannot apply the technique
used in \cite{AO} to obtain the explicit form of the minimal resolutions. On
the other hand, once one knows that they exist unique, they can be determined
in another way (in the case of $\P^n$, this was already observed by Beilinson
in \cite{B2}). For instance, the fact that the coefficient of $\so(-j)$ ($0\le
j\le n$) in $\bro^i(\cp{\s{F}})$ is given by
$H^{i+j}(\P^n,\cp{\s{F}}\otimes\sd^j(j))$ is a formal consequence of the fact
that $\all\s{A}\in O$ and $\all k\in\Z$ we have $H^k(\P^n,\s{A}\otimes\sd^j(j))
=0$, except for $\s{A}=\so(-j)$ and $k=j$, when it is $\K$. Therefore, also in
the weighted case it is enough to find objects of $D^b(\gcoh{\gP})$ with
similar properties of cohomology vanishing as $\sd^j(j)$ in the example above.
It turns out that such objects are given by suitable complexes obtained from
$\cp{\ko}$, except for those which determine the coefficients of the
$\sd^j(j)$, which are again $\so(-j)$. The presence of these complexes makes it
less easy to compute the minimal resolution in general. As for the
problem of determining the differentials of the minimal resolutions,
we can prove very little: in the case of $\P^n$, a method to construct
them has been found recently (see \cite{EFS}), but it doesn't seem easy to
extend it to the weighted case.

Our approach is close to those of \cite{B} (see also \cite{GL}) and of
\cite{AKO}, where $\gP(\w)$ is replaced, respectively, by a different
graded variety (it is obtained grading $\p(\w)$ not over $\Z$, but
over an abelian group depending on $\w$), and by its associated
algebraic stack (having an equivalent category of coherent sheaves, as
we have seen in \ref{stacks}). In both papers, however, only the case
of $\goset$ (and not of $\gdset$) is treated, and the explicit form of
the minimal resolution is not given. In \cite{AKO} (where
noncommutative deformations of weighted projective spaces are also
considered) another generating set for $D^b(\gcoh{\gP})$ (``dual'' to
$\goset$) is found, and we prove that it coincides (up to twists
and shifts) with the above mentioned set of complexes which appear in
the explicit form of the minimal resolution given by $\goset$.

Our result also allows to obtain (minimal) resolutions of every
coherent sheaf on $\P$ (although they are no more unique in general),
and it is applied to prove a generalization of Horrock's splitting
criterion for vector bundles on $\P^n$.

Much of the contents of this chapter is summarized (without proofs) in
\cite{Ca2}.

                \section{Koszul complex and sheaves of differentials}

Let $\md_{\p}:=\md_{\p/\K}$ be the $\p=\p(\w)$--module of
$\K$--differentials of $\p$ and $\all j\in\Z$ let $\md_{\p}^j:=
\Lambda^{j}(\md_{\p})$ (of course, $\md_{\p}^j=0$ if $j<0$ or $j>n+1$);
denoting by $d:\p\to\md_{\p}$ the universal derivation, $\md_{\p}^j$ is a free
$\p$--module with basis
\[\{dx_I\st I\subseteq\{0,\dots,n\};\card{I}=j\}\]
(where $dx_I:=dx_{i_1}\wedge\cdots\wedge dx_{i_j}$ if $I=\{i_1,\dots,i_j\}$
with $0\le i_1<\cdots<i_j\le n$), which is naturally graded by $\deg(dx_I)=
\sw{\w_I}$ (so that $\md_{\p}^j=\bigoplus_{\card{I}=j}\p dx_I\iso
\bigoplus_{\card{I}=j}\p(-\sw{\w_I})$).
Setting $\ko^{-j}:=\md_{\p}^j$ and considering the morphisms (of degree zero)
\[\begin{split}
\ko^{-j}=\md_{\p}^j & \to \ko^{1-j}=\md_{\p}^{j-1}\\
dx_I &
\mapsto\sum_{i\in I}(-1)^{\card{\{i'\in I\st i'<i\}}}x_idx_{I\minus\{i\}}
\end{split}\]
we obtain a complex $\cp{\ko}=\cp{\ko}_{(\w)}\in
C^b(\gfmo{\p})$.\index{K@$\cp{\ko}$, $\cp{\ko}_{(\w)}$} Since $\cp{\ko}$ is
just the Koszul complex of the regular sequence $(x_0,\dots, x_n)$, it is exact
everywhere, except that $H^0(\cp{\ko})=\p/(x_0,\dots,x_n) \iso\K$ (in degree
zero). We will consider also the sheafification $\cp{\sko}
=\cp{\sko}_{(\w)}:=(\cp{\ko})\gsh{}$:\index{K@$\cp{\sko}$, $\cp{\sko}_{(\w)}$}
since $\gsh{}$ is an exact functor and $\K\gsh{}=0$, $\cp{\sko}\in
C^b(\gcoh{\gP})$ is an exact complex.

        \begin{lemm}\label{koselfdual}
There are natural isomorphisms of complexes
\begin{align*}
& {\cp{\ko}}\dual\iso\cp{\ko}(\sw{\w})[-n-1], &
{\cp{\sko}}\dual\iso\cp{\sko}(\sw{\w})[-n-1].
\end{align*}
        \end{lemm}

        \begin{proof}
Clearly it is enough to consider the case of modules. Let $\cp{L}:={\cp{\ko}}
\dual(-\sw{\w})[n+1]$: $\all j\in\Z$ we have by definition
\[L^j=(\ko^{-n-1-j})\dual(-\sw{\w})=
(\bigoplus_{\card{I}=n+1+j}\p dx_I)\dual(-\sw{\w})=
(\bigoplus_{\card{I}=-j}\p dx_{\compl{I}})\dual(-\sw{\w}).\]
Since $(\p dx_{\compl{I}})\dual(-\sw{\w})\iso\p dx_I$, we see that there is a
natural isomorphism $L^j\iso\ko^j$. It is easy to see that, with these
identifications, the differential of $\cp{L}$ is given by
\[\cdiff{L}^j(dx_I)=(-1)^{n+1}
\sum_{i\in I}(-1)^{\card{\{j\notin I\st j<i\}}}x_idx_{I\minus\{i\}},\]
and then it is straightforward to check that the maps $f^j:\ko^j\to L^j$
defined by
\[f^j(dx_I):=(-1)^{j(n+1)+\sum_{i\in I}i}dx_I\]
determine an isomorphism of complexes $\cp{f}:\cp{\ko}\to\cp{L}$.
        \end{proof}

$\all j\in Z$ let $\smd^j:=\ker\cdiff{\ko}^{-j}$\index{Syz^j@$\smd^j$} and
$\sd^j_{\gP}:=(\smd^j)\gsh{}=\ker\cdiff{\sko}^{-j}=
\im\cdiff{\sko}^{-j-1}$.\index{Omega_P@$\sd^j_{\gP}$} Therefore we have short
exact sequences
\[\tag*{$\esd{j}$}0\to\sd^j_{\gP}\mor{\omd{j}}\sko^{-j}\iso
\bigoplus_{\card{I}=j}\so_{\gP}(-\sw{\w_I})\mor{\dom{j}}\sd^{j-1}_{\gP}\to
0\]\index{S_j@$\esd{j}$}\index{alpha^j@$\omd{j}$}\index{beta^j@$\dom{j}$} such
that $\cdiff{\sko}^{-j}=\omd{j-1}\comp\dom{j}$.

        \begin{rema}
$\sd^j_{\P}:=(\smd^j)\sh{}=(\sd^j_{\gP})_0$ is the sheaf of regular
differential $j$--forms (see \cite{Do} or \cite{BR}). It is reflexive but not
locally free in general.
        \end{rema}

        \begin{lemm}\label{Omegadual}
$\sd^j_{\gP}$ is a locally free sheaf of rank $\binom{n}{j}$. In particular,
$\sd^0_{\gP}\iso\so_{\gP}$ and $\sd^n_{\gP}\iso\so_{\gP}(-\sw{\w})$. Moreover,
$(\sd_{\gP}^j)\dual\iso\sd_{\gP}^{n-j}(\sw{\w})$.
        \end{lemm}

        \begin{proof}
The first statement follows easily by induction on $j$, using the exact
sequences $\esd{j}$ (notice that $\sd^j_{\gP}=0$ if $j<0$ or $j>n$). It
is also clear from the definition that $\sd^0_{\gP}\iso\so_{\gP}$ and
$\sd^n_{\gP}\iso\so_{\gP}(-\sw{\w})$. The last statement follows immediately
from \ref{koselfdual}.
        \end{proof}

        \begin{prop}\label{gPdualsh}
$\ds_{\gP}:=\sd^n_{\gP}\iso\so_{\gP}(-\sw{\w})$ is the dualizing sheaf on
$\gP$. Moreover, $\all\s{F}\in\gcoh{\gP}$ and $\all i>0$ there are natural
isomorphisms $\Ext_{\gP}^i(\s{F},\ds_{\gP})\iso H^{n-i}(\gP,\s{F})\dual$.
        \end{prop}

        \begin{proof}
One could easily adapt the standard proof for $\P^n$ (see
\cite[III, thm. 7.1]{H1}). Alternatively, observe that $\so_{\P}(-\sw{\w})$ is
the dualizing sheaf on $\P$ (by \cite[cor. 6B.8]{BR}) and that
$\so_{\gP}(-\sw{\w}))\iso\gsHom_{\P}(\so_{\gP},\so_{\P}(-\sw{\w}))$: this is a
consequence of the fact that $\so_{\P}(-m-\sw{\w}))\iso
\sHom_{\P}(\so(m),\so(-\sw{\w}))$ $\all m\in\Z$ (we can clearly assume
$\gcd(\w_0,\dots,\w_n)=1$, and then this is true by \cite[lemma 4.1]{D} if $\w$
is normalized, whereas in general, if $\P(\w)\iso\P(\w')$ with $\w'$
normalized, it is easy to check, using \cite[prop. 3C.7]{BR}, that $\exi m'\in
\Z$ such that $\so_{\P(\w)}(m)\iso\so_{\P(\w')}(m')$, $\so_{\P(\w)}(-\sw{\w})
\iso\so_{\P(\w')}(-\sw{\w'})$ and $\so_{\P(\w)}(-m-\sw{\w})\iso
\so_{\P(\w')}(-m'-\sw{\w'})$). Since moreover
$\sExt^i_{\P}(\so(m),\so(-\sw{\w}))=0$ for $i>0$ and $\all m\in\Z$ (as usual,
one can reduce to the case $\w$ normalized, and then this is proved in
\cite[prop. 5.4]{D}), everything follows from \ref{Serredual}.
        \end{proof}

        \begin{lemm}\label{mondimfor}
$\all\w\in\N_+^{n+1}$ and $\all l\in\Z$ the following formula holds:
\[\sum_{I\subseteq\{0,\dots,n\}}(-1)^{\card{I}}\pd_{l-\sw{\w_I}}(\w)=
\delta_{l,0}.\]
        \end{lemm}

        \begin{proof}
The statement being trivial for $l\le0$, we can assume $l>0$. In this case
$\all I\subseteq\{0,\dots,n\}$ by \ref{lbcohom} we have $\pd_{l-\sw{\w_I}}=
h^0(\gP,\so(l-\sw{\w_I}))=\chi(\gP,\so(l-\sw{\w_I}))$. Therefore, as
$\cp{\sko}(l)$ is an exact complex,
\begin{multline*}
\sum_I(-1)^{\card{I}}\pd_{l-\sw{\w_I}}=
\sum_{j=0}^{n+1}(-1)^j\sum_{\card{I}=j}\chi(\gP,\so(l-\sw{\w_I}))\\
=\sum_{j=0}^{n+1}(-1)^j\chi(\gP,\sko^{-j}(l))=\chi(\gP,\cp{\sko}(l))=0.
\end{multline*}
        \end{proof}

        \begin{lemm}\label{omegacohom}
$\all j,l\in\Z$ the cohomology of $\sd_{\gP}^j(l)$ is given by the following
formulas:
\begin{align*}
h^i(\gP,\sd^j(l)) & =\delta_{i,j}\delta_{l,0} & \text{if $0<i<n$;}\\
h^0(\gP,\sd^j(l)) & =\sum_{\card{I}\le j}(-1)^{j-\card{I}}
\pd_{l-\sw{\w_I}}-(-1)^j\delta_{l,0}(1-\delta_{j,0})\\
& =\sum_{\card{I}>j}(-1)^{j+1-\card{I}}
\pd_{l-\sw{\w_I}}+(-1)^j\delta_{l,0}\delta_{j,0};\\
h^n(\gP,\sd^j(l)) & =h^0(\gP,\sd^{n-j}(-l)).
\end{align*}
        \end{lemm}

        \begin{proof}
It can be easily proved using the exact sequences $\esd{j}$ and \ref{lbcohom}.

Alternatively, since $H^i(\gP,\sd^j(l))\iso H^i(\P,\sd^j(l))$ by
\ref{goodcohom}, the result follows from \cite[thm. 2.3.2]{Do}. Notice that
the last equality in the second equation follows from \ref{mondimfor} and the
third equation is a consequence of Serre duality and of the fact that
$(\sd_{\gP}^j)\dual\iso\sd_{\gP}^{n-j}(\sw{\w})$.
        \end{proof}

        \begin{coro}\label{omcohvan}
$\all j\in\Z$ we have $h^0(\gP,\sd^j(l))=0$ if $l<j$ and $h^n(\gP,\sd^j(l))=0$
if $l>j-n$. If moreover $j\ne0,n$, then also $h^0(\gP,\sd^j(j))=
h^n(\gP,\sd^j(j-n))=0$.
        \end{coro}

                \section{The theorem as equivalence of categories}

        \begin{lemm}\label{extvan}
$\all j,j',l,l'\in\Z$ and $\all i>0$ we have the following results:
\begin{enumerate}

\item $\Ext_{\gP}^i(\so(l),\so(l'))=0$ if $l'>l-\sw{\w}$;

\item $\Ext_{\gP}^i(\so(l),\sd^j(j))=0$ if $l<i$;

\item $\Ext_{\gP}^i(\sd^j(j),\so(l))=0$ if $l>n-\sw{\w}-i$;

\item $\Ext_{\gP}^i(\sd^j(l),\sd^{j'}(j'))=0$ if $l<j+i$;

\item $\Ext_{\gP}^i(\sd^j(l),\sd^{j'}(j'))=0$ if $j'<i<n-j$.

\end{enumerate}
        \end{lemm}

        \begin{proof}
$1$, $2$ and $3$ follow from \ref{omegacohom}, since
$\Ext_{\gP}^i(\so(l),\so(l'))\iso H^i(\gP,\so(l'-l))$,
$\Ext_{\gP}^i(\so(l),\sd^j(j))\iso H^i(\gP,\sd^j(j-l))$ and (in view of
\ref{Omegadual}) $\Ext_{\gP}^i(\sd^j(j),\so(l))\iso
H^i(\gP,\sd^{n-j}(\sw{\w}+l-j))$. Now we prove $4$ and $5$ by induction on
$j\ge0$. If $j=0$, $4$ coincides with $2$ and $5$ follows again from
\ref{omegacohom}. If $j>0$, applying $\Hom_{\gP}(-,\sd^{j'}(j'))$ to
$\esd{j}(l)$ yields the exact sequence
\[\Ext_{\gP}^i(\bigoplus_{\card{I}=j}\so_{\gP}(l-\sw{\w_I}),\sd^{j'}(j'))\to
\Ext_{\gP}^i(\sd^j(l),\sd^{j'}(j'))\to
\Ext_{\gP}^{i+1}(\sd^{j-1}(l),\sd^{j'}(j')).\]
We have
$\Ext_{\gP}^i(\bigoplus_{\card{I}=j}\so_{\gP}(l-\sw{\w_I}),\sd^{j'}(j'))=0$ by
what we have just proved ($\sw{\w_I}\ge j$ if $\card{I}=j$) and
$\Ext_{\gP}^{i+1}(\sd^{j-1}(l),\sd^{j'}(j'))=0$ by the inductive hypothesis,
whence $\Ext_{\gP}^i(\sd^j(l),\sd^{j'}(j'))=0$.
        \end{proof}

    \begin{coro}\label{Beilextvan}
If $\s{A},\s{B}\in\goset$ or $\s{A},\s{B}\in\gdset$, then
$\Ext_{\gP}^i(\s{A},\s{B})=0$ for $i>0$.
    \end{coro}

        \begin{lemm}\label{Homsmd=Homsd}
$\all j,j',l,l'\in\Z$ the functor $\gsh{}$ induces an isomorphism
\[\Hom_{\p}(\smd^j(l),\smd^{j'}(l'))\isomor\Hom_{\gP}(\sd^j(l),\sd^{j'}(l')).\]
        \end{lemm}

        \begin{proof}
First we prove the case $j=0$ (and $l=0$, which we can always assume), namely
that $\smd^{j'}(l')_0\iso\Hom_{\p}(\p,\smd^{j'}(l'))\isomor
\Hom_{\gP}(\so,\sd^{j'}(l'))\iso H^0(\gP,\sd^{j'}(l'))$. The statement being
trivial for $j'<0$, we can assume $j'\ge0$ and we proceed by induction on $j'$.
Since $\smd^0\iso\p$, the case $j'=0$ follows from \ref{lbcohom}. If $j'>0$,
from the exact sequence
\[0\to\smd^{j'}(l')\to\ko^{-j'}(l')\iso\bigoplus_{\card{I}=j'}\p(l'-\sw{\w_I})
\to\smd^{j'-1}(l')\]
we obtain a commutative diagram with exact rows
\[\begin{CD}
0 @>>> H_{\p}(\smd^{j'}(l')) @>>> H_{\p}(\ko^{-j'}(l')) @>>>
H_{\p}(\smd^{j'-1}(l'))\\
@. @VVV @VVV @VVV \\
0 @>>> H_{\gP}(\sd^{j'}(l')) @>>> H_{\gP}(\sko^{-j'}(l'))
@>>> H_{\gP}(\sd^{j'-1}(l'))
\end{CD}\]
(where $H_{\p}(-):=\Hom_{\p}(\p,-)$ and $H_{\gP}(-):=\Hom_{\gP}(\so,-)$). Since
the vertical arrows on the right and in the middle are isomorphisms by
induction, the one on the left is also an isomorphism by the five--lemma.

Now we pass to the general case. As before, we can assume $j\le n$ and we
proceed by descending induction on $j$, the case $j=n$ following from what we
have just proved, since $\smd^n(l)\iso\p(l-\sw{\w})$. If $j<n$ the exact
sequence
\[0\to\smd^{j+1}(l)\to\ko^{-j-1}(l)\iso\bigoplus_{\card{I}=j+1}\p(l-\sw{\w_I})
\to\smd^j(l)\]
yields a commutative diagram with exact rows
\[\begin{CD}
0 @>>> H'_{\p}(\smd^j(l)) @>>> H'_{\p}(\ko^{-j-1}(l)) @>>>
H'_{\p}(\smd^{j+1}(l))\\
@. @VVV @VVV @VVV \\
0 @>>> H'_{\gP}(\sd^j(l)) @>>> H'_{\gP}(\sko^{-j-1}(l)) @>>>
H'_{\gP}(\sd^{j+1}(l))
\end{CD}\]
(where $H'_{\p}(-):=\Hom_{\p}(-,\smd^{j'}(l'))$ and
$H'_{\gP}(-):=\Hom_{\gP}(-,\sd^{j'}(l'))$). The result then follows, as before,
from the inductive hypothesis using the five--lemma.
        \end{proof}

Given $S\subseteq\Z$, we will denote by
$\mco{S}=\mco{S}(\p)$\index{M_S@$\mco{S}$} the full subcategory of $\gfmo{\p}$
whose objects are of the form $\bigoplus_{l\in S} \p(l)^{a_l}$ ($a_l\in\N$) and
by $\mcom{S}$\index{M_S @$\mcom{S}$} the full subcategory of $\gfmo{\p}$ whose
objects are of the form $\bigoplus_{0\le j\le n}\smd^j(j)^{b_j}\oplus X$ with
$b_j\in\N$ and $X\in\mco{S}$ (notice that $\mcom{S}=
\mcom{S\cup\{n-\sw{w},0\}}$). Similarly, $\co{S}$\index{O_S@$\co{S}$}
(respectively $\com{S}$)\index{O_S @$\com{S}$} will denote the full subcategory
of $\gcoh{\gP}$ whose objects are of the form $\bigoplus_{l\in S}\so(l)^{a_l}$
(respectively $\bigoplus_{0\le j\le n}\sd^j(j)^{b_j}\oplus X$ with
$X\in\co{S}$); we will write $\co{}$\index{O@$\co{}$} and $\com{}$\index{O
@$\com{}$} instead of $\co{\Z}$ and $\com{\Z}$. By \ref{Homsmd=Homsd} it is
clear that $\gsh{\p}$ induces an equivalence of categories between $\mco{S}$
(respectively $\mcom{S}$) and $\co{S}$ (respectively $\com{S}$). If $a,b\in\Z$,
we set $]a,b[:=\{m\in\Z\st a<m<b\}$ and $]a,b]:=\{m\in\Z\st a<m\le b\}$ (and
similarly for $[a,b[$ and $[a,b]$).

        \begin{theo}\label{wbeilthm}
The functors $\gsh{\p}:\mco{\into}\to\gcoh{\gP}$ and $\gsh{\p}:\mcom{\intd}\to
\gcoh{\gP}$ extend naturally to exact functors between triangulated categories
\begin{align*}
& \gfo:K^b(\mco{\into})\to D^b(\gcoh{\gP}), &
\gfd:K^b(\mcom{\intd})\to D^b(\gcoh{\gP}).
\end{align*}
$\gfo$ and $\gfd$ are equivalences of categories.
        \end{theo}

        \begin{proof}
Clearly $\gsh{\p}$ induces exact and fully faithful functors from
$K^b(\mco{\into})$ (respectively $K^b(\mcom{\intd})$) to $K^b(\gcoh{\gP})$,
which, composed with the natural exact functor $K^b(\gcoh{\gP})\to
D^b(\gcoh{\gP})$, define the exact functors $\gfo$ and $\gfd$. In order to
prove that they are fully faithful it is therefore enough to show that given
$\cp{\s{F}},\cp{\s{G}}\in K^b(\co{\into})$ (respectively $K^b(\com{\intd})$),
the natural map $\alpha:\Hom_{K^b(\gcoh{\gP})}(\cp{\s{F}},\cp{\s{G}})\to
\Hom_{D^b(\gcoh{\gP})}(\cp{\s{F}},\cp{\s{G}})$ is an isomorphism. Consider the
natural commutative diagram
\[\begin{CD}
\Hom_{K^b(\gcoh{\gP})}(\cp{\s{F}},\cp{\s{G}}) @>{\alpha}>>
\Hom_{D^b(\gcoh{\gP})}(\cp{\s{F}},\cp{\s{G}}) \\
@VV{\beta}V @VV{\beta'}V \\
\Hom_{K^b(\gqco{\gP})}(\cp{\s{F}},\cp{\s{G}}) @>{\alpha'}>>
\Hom_{D^b(\gqco{\gP})}(\cp{\s{F}},\cp{\s{G}})
\end{CD}\]
where $\beta$ is obviously an isomorphism. By \ref{gcgqgc} the abelian
categories $\gqco{\gP}$ and $\gcoh{\gP}$ satisfy the hypotheses of
\ref{HomDA=HomDB}, whence $\beta'$ is an isomorphism. Moreover, it follows from
\ref{Beilextvan} that $\Ext_{\gP}^i(\s{F}^l,\s{G}^m)=0$ for $i>0$ and $\all
l,m\in\Z$, so that $\alpha'$ is an isomorphism by \ref{HomKA=HomDA}. Thus
$\alpha$ is an isomorphism, too.

It remains to prove that $\gfo$ and $\gfd$ are essentially surjective. As
$D^b(\gcoh{\gP})$ is generated (as a triangulated category) by the complexes
concentrated in position $0$ and since every object of $\gcoh{\gP}$ has a
finite resolution with sheaves of the form $\so(j)$ (by Hilbert's syzygy
theorem), it is enough to show that the $\so(j)$ ($j\in\Z$) are in the
essential images $\im\gfo$ and $\im\gfd$. First we prove by descending
induction on $j$ that $\so(j)\in\im\gfd$ for $-\sw{\w}<j\le0$ (this is true by
hypothesis for $n-\sw{\w}\le j\le0$). Now, if $-\sw{\w}<j<n-\sw{\w}$,
$\so(j)[n+1]$ is quasi--isomorphic to the complex
\[0\to\sko^{-n}(\sw{\w}+j)\to\cdots\to\sko^{-\sw{\w}-j-1}(\sw{\w}+j)\to
\sd^{\sw{\w}+j}(\sw{\w}+j)\to0,\] every term of which is in $\im\gfd$ by
induction. Therefore $\im\gfd\supseteq \im\gfo$ and we can consider only the
case of $\gfo$. Again, since $\so(j)$ is quasi--isomorphic to
$\sko^{<0}(j)[-1]$, it follows by induction that $\so(j)\in \im\gfo$ for $j>0$
(in a similar way one proceeds for $j\le-\sw{\w}$).
        \end{proof}

        \begin{rema}
If $\w=(1,\dots,1)$, this result coincides with Beilinson's theorem (taking
into account that $\gcoh{\gP^n}$ is equivalent to $\coh{\P^n}$ by
\ref{shequivcrit} and \ref{stdmori}), except that, in order to identify $\gfd$
with $\fde$, it is also necessary to consider the equivalence between
$\mcom{\intd}=\mcom{\emptyset}$ and $\mco{[-n,0]}(\Lambda)$ (defined on objects
by $\smd^j(j)\mapsto\Lambda(-j)$). That this is an equivalence follows from the
fact that $\Hom_{\P^n}(\sd^j(j),\sd^{j'}(j'))\iso\Lambda^{j-j'}(\p_1\dual)$.
Since this is not true in the weighted case (although it still holds in a
weaker form, see \ref{omegahom}), and also because of the presence of $\p(l)$
for $n-\sw{\w}<l<0$, we don't think it is possible to give in general an
alternative simple description of $\mcom{\intd}$.
        \end{rema}

    \begin{rema}\label{Oexc}
By \ref{Beilextvan}, together with the fact that $\Hom_{\gP}(\so(j),\so(l))=0$
for $j>l$ and $\Hom_{\gP}(\so(j),\so(j))\iso\K$, the sequence
$(\so(1-\sw{\w}),\dots,\so)$ (hence also $(\so(k),\dots,\so(k+\sw{\w}-1))$
$\all k\in\Z$) of $D^b(\gcoh{\gP})$ is strong exceptional (see \ref{exc}), and
also full, since $\goset$ generates $D^b(\gcoh{\gP})$, as we proved before. On
the other hand, while $(\sd^n(n),\dots,\sd^0(0)=\so)$ is a full and strong
exceptional sequence if $\w=(1,\dots,1)$, the objects of $\gdset$ usually do
not form an exceptional sequence, because in general
$\dim_{\K}\Hom_{\gP}(\sd^j(j),\sd^j(j))>1$ (so that $\sd^j(j)$ is not
exceptional), as it is easy to see. In \cite{AKO} also the full (but not
strong) exceptional sequence which is left dual (see \ref{dualmut}) of
$(\so,\dots,\so(\sw{\w}-1))$ plays an important role; later we will study it
(see \ref{Odual} and \ref{notstrong}).
    \end{rema}

    \begin{rema}\label{Dmod}
Both in \cite{B} and in \cite{AKO} a different, although closely related,
result is proved, namely (translating to our setting) that there is an
equivalence of categories between $D^b(\gcoh{\gP})$ and $D^b(\fmo{B})$, where
$B$ is the (noncommutative) endomorphism algebra of $\bigoplus_{0\le
j<\sw{\w}}\so(j)$.
    \end{rema}

                \section{Morphisms between sheaves of differentials}

Given $\s{F},\s{G}\in\gcoh{\gP}$ and $S\subseteq\Z$, let
$\ohom{S}(\s{F},\s{G})$\index{Hom^O_S@$\ohom{S}$} be the subset of
$\Hom_{\gP}(\s{F},\s{G})$ consisting of those morphisms which factor through an
object of $\co{S}$. Since the composition of two morphisms of $\gcoh{\gP}$, one
of which is in $\ohom{S}$, is again in $\ohom{S}$, we can consider the quotient
category $\coq{S}= \gcoh{\gP}/\co{S}$\index{Q_S@$\coq{S}$} with the same
objects as $\gcoh{\gP}$ and with
$\Hom_{\coq{S}}(\s{F},\s{G}):=\Hom_{\gP}(\s{F},\s{G})/\ohom{S}(\s{F},\s{G})$
(as usual, $\coq{\Z}$ will be simply denoted by $\coq{}$).\index{Q@$\coq{}$}

        \begin{rema}\label{coqadd}
It is easy to see that $\all\s{F},\s{G}\in\gcoh{\gP}$ the set
$\ohom{S}(\s{F},\s{G})$ is a subgroup of $\Hom_{\gP}(\s{F},\s{G})$, and that
therefore $\coq{S}$ is an additive category.
        \end{rema}

        \begin{defi}\label{mindef}
A morphism in $\gcoh{\gP}$
\[f:\bigoplus_{0\le j<n}\bigoplus_{l\in\Z}\sd^j(l)^{a^j_l}\to
\bigoplus_{0\le j<n}\bigoplus_{l\in\Z}\sd^j(l)^{b^j_l}\] is called {\em minimal
with respect to}\index{minimal!with respect to} $\sd^j(l)$ ($0\le j<n$ and
$l\in\Z$) if the component of $f$ from $\sd^j(l)^{a^j_l}$ to $\sd^j(l)^{b^j_l}$
is $0$ in case $j=0$ and if it is in
$\ohom{}(\sd^j(j)^{a^j_l},\sd^j(j)^{b^j_l})$ (i.e. its image in
$\Hom_{\coq{}}(\sd^j(j)^{a^j_l},\sd^j(j)^{b^j_l})$ is $0$) in case $0<j<n$. $f$
is {\em minimal}\index{minimal!morphism} if it is minimal with respect to each
$\sd^j(l)$.

A complex $\cp{\s{L}}\in C(\gcoh{\gP})$ such that each $\s{L}^i$ is a direct
sum of sheaves of the form $\sd^j(l)$ is {\em minimal}\index{minimal!complex}
if $\cdiff{\s{L}}^i$ is minimal $\all i\in\Z$.
        \end{defi}

        \begin{rema}
$\ohom{}(\sd^j(j),\sd^{j'}(j'))=\ohom{\intd}(\sd^j(j),\sd^{j'}(j'))$ if $0<j,j'
<n$: this follows from the fact that (by Serre duality and \ref{omcohvan})
\[\Hom_{\gP}(\sd^j(j),\so(l))\iso H^n(\gP,\sd^j(j-\sw{\w}-l))\dual=0\]
if $l\le n-\sw{\w}$ and $\Hom_{\gP}(\so(l),\sd^{j'}(j'))\iso
H^0(\gP,\sd^{j'}(j'-l))=0$ if $l\ge0$.
        \end{rema}

        \begin{lemm}\label{indec}
Let $\so_{\gP}(l)\mor{f}\sd_{\gP}^j\mor{g}\so_{\gP}(l)$ be morphisms with
$l\in\Z$ and $0<j<n$. Then $g\comp f=0$.
        \end{lemm}

        \begin{proof}
Applying the functor $\Hom_{\gP}(\so(l),-)$ to $\esd{j+1}$ we obtain the exact
sequence
\[\Hom_{\gP}(\so(l),\sko^{-j-1})\to\Hom_{\gP}(\so(l),\sd^j)\to
\Ext^1_{\gP}(\so(l),\sd^{j+1}).\]
Since $\Ext^1_{\gP}(\so(l),\sd^{j+1})\iso H^1(\gP,\sd^{j+1}(-l))=0$ by
\ref{omegacohom}, $f$ factors as $f=\dom{j+1}\comp f'$ for some $f':\so(l)
\to\sko^{-j-1}$. Similarly, applying $\Hom_{\gP}(-,\so(l))$ to $\esd{j}$ one
can prove that $g=g'\comp\omd{j}$ for some $g':\sko^{-j}\to\so(l)$. Hence
$g\comp f=g'\comp\omd{j}\comp\dom{j+1}\comp f'=g'\comp\cdiff{\sko}^{-j-1}\comp
f'$. Now $g'$, $\cdiff{\sko}^{-j-1}$ and $f'$ are morphism in $\co{}$ and
$\cdiff{\sko}^{-j-1}$ is minimal, so that their composition $g\comp f$ is also
minimal (this follows immediately from the fact that
$\Hom_{\gP}(\so(l),\so(l'))=0$ if $l>l'$), i.e. $g\comp f=0$.
        \end{proof}

Denoting by $\pi^j:\sko^{-j}(j)\iso\bigoplus_{\card{I}=j}\so_{\gP}(j-\sw{\w_I})
\to\bigoplus_{\card{I}=j,\sw{\w_I}>j}\so_{\gP}(j-\sw{\w_I})$ the natural
projection, and observing that $n-\sw{\w}<j-\sw{\w_I}\le0$ if $\card{I}=j$, we
will consider the map
\[\pi^j\comp\omd{j}(j):\sd^j(j)\to
\bigoplus_{\card{I}=j,\sw{\w_I}>j}\so_{\gP}(j-\sw{\w_I}).\]

        \begin{lemm}\label{omhomfact}
$\phi\in\Hom_{\gP}(\sd^j(j),\sd^{j'}(j'))$ is in
$\ohom{\intd}(\sd^j(j),\sd^{j'}(j'))$ if and only if it factors through
$\pi^j\comp\omd{j}(j)$. Moreover, if $j'>0$, this is the case if and only if
$\phi$ factors through $\omd{j}(j)$.
        \end{lemm}

        \begin{proof}
Assume that $\phi\in\ohom{\intd}(\sd^j(j),\sd^{j'}(j'))$ (the other implication
being obvious), say $\phi=\phi'\comp\psi$ with $\sd_{\gP}^j(j)\mor{\psi}\s{F}:=
\bigoplus_{n-\sw{\w}<l<0}\so_{\gP}(l)^{a_l}\mor{\phi'}\sd_{\gP}^{j'}(j')$.
Applying $\Hom_{\gP}(-,\s{F})$ to $\esd{j}(j)$, we obtain the exact sequence
\[\Hom_{\gP}(\bigoplus_{\card{I}=j}\so(j-\sw{\w_I}),\s{F})\to
\Hom_{\gP}(\sd^j(j),\s{F})\to\Ext_{\gP}^1(\sd^{j-1}(j),\s{F}).\]
Since, for $l>n-\sw{\w}$, $\Ext_{\gP}^1(\sd^{j-1}(j),\so(l))=0$ by
\ref{extvan}, we see that $\psi=\rho\comp\omd{j}(j)$ for some $\rho\in
\Hom_{\gP}(\bigoplus_{\card{I}=j}\so(j-\sw{\w_I}),\s{F})$. On the other hand,
since $\Hom_{\gP}(\so,\so(l))=0$ for $l<0$, it is clear that $\rho$ factors
through $\pi^j$, say $\rho=\rho'\comp\pi^j$. Therefore $\phi=\phi'\comp\rho'
\comp\pi^j\comp\omd{j}(j)$. The last statement follows immediately from the
fact that, if $j'>0$, then $\Hom_{\gP}(\so,\sd^{j'}(j'))=0$ (by
\ref{omcohvan}).
        \end{proof}

Unlike the case of $\gP^n$, it is not true in general that
$\Hom_{\gP}(\sd^j(j),\sd^{j'}(j'))$ (for $0\le j,j'\le n$) only depends on
$j-j'$ (for instance, it can happen that $\Hom_{\gP}(\sd^j(j),\sd^{j'}(j'))\ne
0$ for some $j<j'$, whereas $\Hom_{\gP}(\so,\sd^{j'-j}(j'-j))=0$ by
\ref{omcohvan}). On the other hand, we have the following result.

        \begin{prop}\label{omegahom}
For $0\le j,j'\le n$ there are natural isomorphisms
\[\Hom_{\coq{\intd}}(\sd^j(j),\sd^{j'}(j'))\iso\Lambda^{j-j'}(\p_1\dual).\]
        \end{prop}

        \begin{proof}
We will denote by $\hat{\sko}^{-j}$ the subsheaf
$\bigoplus_{\card{I}=j=\sw{\w_I}}\so(-\sw{\w_I})$ of $\sko^{-j}$; observe that
it is naturally isomorphic to $\so(-j)\otimes_{\K}\Lambda^j(\p_1)$.

First we consider the case $j'=0$: applying $\Hom_{\gP}(-,\so)$ to $\esd{j}(j)$
yields the exact sequence
\[0\to\Hom_{\gP}(\sd^{j-1}(j),\so)\mor{\dom{j}(j)^*}
\Hom_{\gP}(\sko^{-j}(j),\so)\mor{\omd{j}(j)^*}
\Hom_{\gP}(\sd^j(j),\so)\to0\]
(because $\Ext_{\gP}^1(\sd^{j-1}(j),\so)=0$ by \ref{extvan}). It follows from
\ref{omhomfact} that
\[\Hom_{\coq{\intd}}(\sd^j(j),\so)\iso\Hom_{\gP}(\sko^{-j}(j),\so)/
(\im\dom{j}(j)^*+\im(\pi^j)^*).\]
Since $\all\phi\in\Hom_{\gP}(\sd^{j-1}(j),\so)$ we have
$\dom{j}(j)^*(\phi)\rest{\hat{\sko}^{-j}(j)}=0$ (by \ref{indec} if $j>1$ and
because $\phi=0$ if $j\le1$), we obtain
\[\Hom_{\coq{\intd}}(\sd^j(j),\so)\iso\Hom_{\gP}(\hat{\sko}^{-j}(j),\so)\iso
\Hom_{\gP}(\so\otimes_{\K}\Lambda^j(\p_1),\so)\iso\Lambda^j(\p_1\dual).\]
By duality this also proves the case $j=n$: indeed, we have more generally
\[\Hom_{\coq{\intd}}(\s{F},\s{G})\iso
\Hom_{\coq{\intd}}(\s{G}\dual(n-\sw{\w}),\s{F}\dual(n-\sw{\w}))\]
$\all\s{F},\s{G}\in\gcoh{\gP}$ (because the same is true for $\Hom_{\gP}$ and
$\ohom{\intd}$). Taking into account that $\sd^j(j)\dual(n-\sw{\w})\iso
\sd^{n-j}(n-j)$, this implies
\[\Hom_{\coq{\intd}}(\sd^j(j),\sd^{j'}(j'))\iso
\Hom_{\coq{\intd}}(\sd^{n-j'}(n-j'),\sd^{n-j}(n-j)),\]
whence $\Hom_{\coq{\intd}}(\sd^n(n),\sd^{j'}(j'))\iso\Lambda^{n-j'}(\p_1\dual)$
by the case already proved.

Therefore we can assume $j'>0$ and $j<n$. As before,
$\Hom_{\gP}(-,\sd^{j'}(j'))$ applied to $\esd{j}(j)$ yields the exact sequence
(observing that $\Ext_{\gP}^1(\sko^{-j}(j),\sd^{j'}(j'))=0$ by \ref{extvan})
\[\Hom_{\gP}(\sko^{-j}(j),\sd^{j'}(j'))\mor{\omd{j}(j)^*}
\Hom_{\gP}(\sd^j(j),\sd^{j'}(j'))\to\Ext_{\gP}^1(\sd^{j-1}(j),\sd^{j'}(j'))\to
0.\]
From \ref{omhomfact} it follows that $\Hom_{\coq{\intd}}(\sd^j(j),\sd^{j'}(j'))
\iso\Ext_{\gP}^1(\sd^{j-1}(j),\sd^{j'}(j'))$, and we claim that
\[\Hom_{\coq{\intd}}(\sd^j(j),\sd^{j'}(j'))\iso
\Ext_{\gP}^{j'}(\sd^{j-j'}(j),\sd^{j'}(j')).\] In fact, if $0<i<j'$, applying
$\Hom_{\gP}(-,\sd^{j'}(j'))$ to $\esd{j-i}(j)$ we obtain the exact sequence
(taking into account that $\Ext_{\gP}^i(\sko^{i-j}(j),\sd^{j'}(j'))=0$ by
\ref{omegacohom})
\[0\to\Ext_{\gP}^i(\sd^{j-i}(j),\sd^{j'}(j'))\to
\Ext_{\gP}^{i+1}(\sd^{j-i-1}(j),\sd^{j'}(j'))\to
\Ext_{\gP}^{i+1}(\sko^{i-j}(j),\sd^{j'}(j')).\]
Since also $\Ext_{\gP}^{i+1}(\sko^{i-j}(j),\sd^{j'}(j'))=0$ by \ref{extvan},
the claim follows by induction on $i$.
Applying $\Hom_{\gP}(-,\sd^{j'}(j'))$ to $\esd{j-j'-i}(j)$ (for $i=0,1$) gives
the exact sequence
\begin{multline*}
\Ext_{\gP}^{j'}(\sd^{j-j'-1-i}(j),\sd^{j'}(j'))\mor{\dom{j-j'-i}(j)^*}
\Ext_{\gP}^{j'}(\sko^{j'-j+i}(j),\sd^{j'}(j'))\\
\mor{\omd{j-j'-i}(j)^*}\Ext_{\gP}^{j'}(\sd^{j-j'-i}(j),\sd^{j'}(j'))\to
\Ext_{\gP}^{j'+1}(\sd^{j-j'-1-i}(j),\sd^{j'}(j'))=0
\end{multline*}
(the last term is $0$ by part 5 of \ref{extvan}, as $j<n$). Now
$\dom{j-j'}(j)^*\comp\omd{j-j'-1}(j)^*=(\cdiff{\sko}^{j'-j})^*=0$ (this follows
from the fact that $\cdiff{\sko}^{j'-j}$ is a minimal morphism and
$\Ext_{\gP}^{j'}(\so(l),\sd^{j'}(j'))\iso H^{j'}(\gP,\sd^{j'}(j'-l))=0$ for $j'
\ne l$ if $0<j'<n$, as we can assume here), whence also $\dom{j-j'}(j)^*=0$,
because $\omd{j-j'-1}(j)^*$ is surjective. Thus $\omd{j-j'}(j)^*$ is an
isomorphism and
\[\Hom_{\coq{\intd}}(\sd^j(j),\sd^{j'}(j'))\iso
\Ext_{\gP}^{j'}(\sko^{j'-j}(j),\sd^{j'}(j')).\]
On the other hand, by \ref{omegacohom} we have
\begin{multline*}
\Ext_{\gP}^{j'}(\sko^{j'-j}(j),\sd^{j'}(j'))\iso
H^{j'}(\gP,\sd^{j'}(j')\otimes\sko^{j'-j}(j)\dual)\\
\iso H^{j'}(\gP,\sd^{j'}(j')\otimes\hat{\sko}^{j'-j}(j)\dual)\iso
H^{j'}(\gP,\sd^{j'}\otimes_{\K}\Lambda^{j-j'}(\p_1)\dual)\iso
\Lambda^{j-j'}(\p_1\dual).
\end{multline*}
        \end{proof}

        \begin{rema}
One can check that, under the isomorphisms given by \ref{omegahom}, composition
of morphisms in $\coq{\intd}$ corresponds to multiplication in
$\Lambda(\p_1\dual)$.
        \end{rema}

        \begin{coro}\label{compminmin}
If $X\mor{f}Y\mor{g}Z$ are morphisms in $\com{}$ with $f$ or $g$ minimal, then
$g\comp f$ is minimal, too.
        \end{coro}

        \begin{proof}
It follows from \ref{indec} and the fact that
$\Hom_{\coq{\intd}}(\sd^j(j),\sd^{j'}(j'))=0$ if $j<j'$ and
$\Hom_{\gP}(\so(l),\so(l'))=0$ if $l>l'$.
        \end{proof}

                \section{Uniqueness of the minimal resolution}

        \begin{lemm}\label{minnil}
If $f:\s{L}\to \s{L}$ is a minimal morphism in $\com{}$, then $f$ is nilpotent.
        \end{lemm}

        \begin{proof}
Let $\s{L}=\bigoplus_{l\in\Z}\so(l)^{a_l}\bigoplus_{0<j<n}\sd^j(j)^{b_j}$.
First notice that $g:=f^{n-1}\in\ohom{}(\s{L},\s{L})$, i.e. $g=0$ in
$\Hom_{\coq{}}(\s{L},\s{L})$: this follows easily from the fact that for $0<j,
j'<n$ the component of $f$ from $\sd^j(j)^{b_j}$ to $\sd^{j'}(j')^{b_{j'}}$ is
zero in $\Hom_{\coq{}}$ if $j\le j'$ (by \ref{omegahom} if $j<j'$ and because
$f$ is minimal if $j=j'$). Therefore we can write $g=g''\comp g'$ with
$\s{L}\mor{g'}\s{M}\mor{g''}\s{L}$ and $\s{M}\in\co{}$, and then the morphism
$h:=(g'\comp g'')^2=g'\comp g\comp g'':\s{M}\to\s{M}$ is minimal by
\ref{compminmin}. Now it is clear that $h$ is nilpotent (because $\all l,l'\in
\Z$ the component of $f$ from $\so(l)^{a_l}$ to $\so(l')^{a_{l'}}$ is zero if
$l\ge l'$), say $h^m=0$. Then $g^{2m+1}=g''\comp h^m\comp g'=0$, whence also
$f$ is nilpotent.
        \end{proof}

        \begin{coro}
$\sd_{\gP}^j(l)$ is indecomposable for $0\le j\le n$ and $\all l\in\Z$.
        \end{coro}

        \begin{proof}
Since the statement is independent of $l$, we can take $l=j$. By \ref{omegahom}
\[\Hom_{\coq{\intd}}(\sd^j(j),\sd^j(j))\iso\K\]
and it is clearly generated by $\id_{\sd^j(j)}$. Now assume that $\sd^j(j)\iso
\s{F}\oplus\s{G}$ with $\s{F}\ne0\ne\s{G}$, and let $f:\sd^j(j)\to\sd^j(j)$ be
the morphism corresponding to $\id_{\s{F}}$. Then $\all\lambda\in\K$ the
morphism $f-\lambda\id_{\sd^j(j)}$, being not nilpotent, is not minimal, i.e.
its image in $\Hom_{\coq{\intd}}(\sd^j(j),\sd^j(j))$ is not zero. We have thus
obtained a contradiction.
        \end{proof}

        \begin{prop}\label{minuni}
Every bounded complex of $\com{}$ is isomorphic in $K^b(\com{})$ to a minimal
complex, which is unique up to isomorphism in $C^b(\com{})$. Moreover, a
morphism between two minimal complexes is an isomorphism in $C^b(\com{})$ if
and only if it is an isomorphism in $K^b(\com{})$.
        \end{prop}

        \begin{proof}
We will show that every $\cp{\s{L}}\in K^b(\com{})$ is isomorphic to a minimal
complex by induction on $\rk{\cp{\s{L}}}:=\sum_{i\in\Z}\rk{\s{L}^i}$, the case
$\rk{\cp{\s{L}}}=0$ being trivial. We can assume $\cp{\s{L}}$ is not minimal,
and then $\cp{\s{L}}$ must be of the form (for some $i\in\Z$)
\[\cdots\to\s{L}^{i-1}
\mor{\begin{pmatrix}
a^{i-1} \\
b^{i-1}
\end{pmatrix}}
\s{L}^i=\begin{matrix}
\tilde{\s{L}}^i \\
\oplus \\
\s{A}
\end{matrix}
\mor{\begin{pmatrix}
a^i & c^i \\
b^i & s
\end{pmatrix}}
\s{L}^{i+1}=\begin{matrix}
\tilde{\s{L}}^{i+1} \\
\oplus \\
\s{A}
\end{matrix}
\mor{\begin{pmatrix}
a^{i+1} & c^{i+1}
\end{pmatrix}}
\s{L}^{i+2}\to\cdots\]
with $\s{A}=\sd^j(j)$ for some $0<j<n$ or $\s{A}=\so(l)$ for some $l\in\Z$ and
with $s$ an isomorphism (this is clear if $\s{A}=\so(l)$, whereas if $\s{A}=
\sd^j(j)$, just notice that by \ref{omegahom} $s$ must be of the form
$\lambda(\id+f)$ with $\lambda\in\K^*$ and $f$ nilpotent by
\ref{minnil}). The condition $\cdiff{\s{L}}^i\comp\cdiff{\s{L}}^{i-1}=0$ is
equivalent to $b^{i-1}=-s^{-1}\comp b^i\comp a^{i-1}$ and $(a^i-c^i\comp s^{-1}
\comp b^i)\comp a^{i-1}=0$ and the condition $\cdiff{\s{L}}^{i+1}\comp
\cdiff{\s{L}}^i=0$ is equivalent to $c^{i+1}=-a^{i+1}\comp c^i\comp s^{-1}$ and
$a^{i+1}\comp(a^i-c^i\comp s^{-1}\comp b^i)=0$. It follows that
\[\cdots\to\s{L}^{i-1}\mor{a^{i-1}}\tilde{\s{L}}^i\mor{a^i-c^i\comp s^{-1}\comp
b^i}\tilde{\s{L}}^{i+1}\mor{a^{i+1}}\s{L}^{i+2}\to\cdots\]
is a complex, denoted by $\cp{\tilde{\s{L}}}$. It is easy to see that the maps
\[\begin{CD}
\cp{\tilde{\s{L}}}= @. \cdots @>>> \tilde{\s{L}}^{i-1}
@>\cdiff{\tilde{\s{L}}}^{i-1}>> \tilde{\s{L}}^i @>\cdiff{\tilde{\s{L}}}^i>>
\tilde{\s{L}}^{i+1} @>\cdiff{\tilde{\s{L}}}^{i+1}>> \tilde{\s{L}}^{i+2} @>>>
\cdots \\
@VV\cp{f}V @. @VV\id V
@VV\left(\begin{smallmatrix}
\id \\
-s^{-1}\comp b^i
\end{smallmatrix}\right)V
@VV\left(\begin{smallmatrix}
\scriptstyle \id \\
\scriptstyle 0
\end{smallmatrix}\right)V
@VV\id V @. \\
\cp{\s{L}}= @. \cdots @>>> \s{L}^{i-1} @>\cdiff{\s{L}}^{i-1}>>
\tilde{\s{L}}^i\oplus\s{A}
@>\cdiff{\s{L}}^i>> \tilde{\s{L}}^{i+1}\oplus\s{A}
@>\cdiff{\tilde{\s{L}}}^{i+1}>>\tilde{\s{L}}^{i+2} @>>> \cdots \\
@VV\cp{g}V @. @VV\id V
@VV{\left(\begin{smallmatrix}
\id & 0
\end{smallmatrix}\right)}V
@VV{\left(\begin{smallmatrix}
\id & -c^i\comp s^{-1}
\end{smallmatrix}\right)}V
@VV\id V @. \\
\cp{\tilde{\s{L}}}= @. \cdots @>>> \tilde{\s{L}}^{i-1}
@>\cdiff{\tilde{\s{L}}}^{i-1}>> \tilde{\s{L}}^i @>\cdiff{\tilde{\s{L}}}^i>>
\tilde{\s{L}}^{i+1} @>\cdiff{\tilde{\s{L}}}^{i+1}>> \tilde{\s{L}}^{i+2} @>>>
\cdots
\end{CD}\]
define morphisms of complexes, which are inverse isomorphisms in $K^b(\com{})$.
Indeed, $\cp{g}\comp\cp{f}=\id_{\cp{\tilde{\s{L}}}}$ and $\all j\in\Z$
$\id_{\s{L}^j}-f^j\comp g^j=h^{j+1}\comp\cdiff{\s{L}}^j+\cdiff{\s{L}}^{j-1}
\comp h^j$, where $h^j:\s{L}^j\to \s{L}^{j-1}$ is defined by
\begin{align*}
& h^{i+1}:=\begin{pmatrix}
0 & 0\\
0 & s^{-1}
\end{pmatrix}, & h^j=0\text{ if }j\ne i+1.
\end{align*}
Since $\rk{\cp{\tilde{\s{L}}}}<\rk{\cp{\s{L}}}$, by induction
$\cp{\tilde{\s{L}}}$ (and then also $\cp{\s{L}}$) is isomorphic to a minimal
complex.

It remains to prove the last statement (which clearly implies the uniqueness of
the minimal complex in $C^b(\com{})$), namely that if $\cp{f}:\cp{\s{L}}\to
\cp{\s{M}}$ is a morphism in $C^b(\com{})$ inducing an isomorphism in
$K^b(\com{})$ with $\cp{\s{L}}$ and $\cp{\s{M}}$ minimal, then $\cp{f}$ is an
isomorphism in $C^b(\com{})$. Let $\cp{g}:\cp{\s{M}}\to\cp{\s{L}}$ in
$C^b(\com{})$ be such that $\cp{f}$ and $\cp{g}$ are inverse isomorphisms in
$K^b(\com{})$. $\cp{g}\comp\cp{f}$ being homotopic to $\id_{\cp{\s{L}}}$ means
that $\all i\in\Z$ there exist maps $h^i:\s{L}^i\to\s{L}^{i-1}$ such that,
setting $a^i:=\cdiff{\s{L}}^{i-1}\comp h^i+h^{i+1}\comp\cdiff{\s{L}}^i$, we
have $g^i\comp f^i=\id_{\s{L}^i}-a^i$. By \ref{compminmin} each $a^i$ is
minimal, and then, by \ref{minnil}, $\exi m\in\N$ such that $(a^i)^m=0$. Thus,
setting $s^i:=\sum_{j=0}^{m-1}(a^i)^j$, we obtain
\[\id_{\s{L}^i}=(\id_{\s{L}^i})^m-(a^i)^m=(\id_{\s{L}^i}-a^i)\comp s^i=s^i\comp
(\id_{\s{L}^i}-a^i)=(g^i\comp f^i)\comp s^i=s^i\comp(g^i\comp f^i),\]
so that $g^i\comp f^i$ is an isomorphism. In the same way, the fact that
$\cp{f}\comp\cp{g}$ is homotopic to $\id_{\cp{\s{M}}}$ implies that also $f^i
\comp g^i$ is an isomorphism, whence $f^i$ and $g^i$ are isomorphisms $\all i
\in\Z$.
        \end{proof}

    \begin{rema}\label{minunit}
With obvious modifications, the same result holds, of course, also for a
complex a twist of which is in $K^b(\com{})$.
    \end{rema}

        \begin{rema}\label{minuniwrt}
It is clear from the proof of \ref{minuni} that if $\cp{\s{L}},\cp{\s{M}}\in
K^b(\com{})$ are isomorphic, $\cp{\s{M}}$ is minimal and $\cp{\s{L}}$ is
minimal with respect to some $\s{A}$ (with $\s{A}=\sd^j(j)$ for some $0<j<n$ or
$\s{A}=\so(l)$ for some $l\in\Z$), then the coefficients of $\s{A}$ in
$\s{L}^i$ and $\s{M}^i$ are the same $\all i\in\Z$.
        \end{rema}

        \begin{coro}\label{wbeiluni}
Given $\cp{\s{F}}\in D^b(\gcoh{\gP})$ there exist unique (up to isomorphisms)
minimal complexes $\cp{\gbro}\in C^b(\co{\into})$ and $\cp{\gbrd}\in
C^b(\com{\intd})$ such that $\cp{\s{F}}\iso\cp{\gbro}\iso\cp{\gbrd}$ in
$D^b(\gcoh{\gP})$.
        \end{coro}

        \begin{proof}
It follows immediately from \ref{wbeilthm} and \ref{minuni}.
        \end{proof}

                \section{Explicit form of the minimal resolution}

        \begin{lemm}\label{forcrit}
Let $\s{E}\ind{1},\dots,\s{E}\ind{r}$ be distinct vector bundles on $\gP$ of
the form $\so(l)$ ($l\in\Z$) or $\sd^j(j)$ ($0<j<n$) and let $\cp{\s{U}}$ be a
minimal complex in $C^b(\com{})$ with $\s{U}^i=\bigoplus_{k=1}^r\s{E}\ind{k}
^{a^i_k}$. Assume that there exist $h\in\{1,\dots,r\}$ and $\cp{\s{L}}\in
C^b(\gcoh{\gP})$ such that $h^i(\gP,\cp{\s{L}}\otimes\s{E}\ind{k})=\delta_{h,k}
\delta_{i,0}$. Then $\all i\in\Z$ there is a natural isomorphism
$\s{E}\ind{h}^{a^i_h}\iso\s{E}\ind{h}\otimes_{\K}
H^i(\gP,\cp{\s{L}}\otimes\cp{\s{U}})$.

Moreover, if $\cp{\phi}:\cp{\s{U}}\to\cp{\s{V}}$ is a morphism in $C^b(\com{})$
with $\cp{\s{V}}$ also minimal and $\s{V}^i=\bigoplus_{k=1}^r
\s{E}\ind{k}^{b^i_k}$, let $\phi^i\ind{h}$ be the component of $\phi^i$ from
$\s{E}\ind{h}^{a^i_h}$ to $\s{E}\ind{h}^{b^i_h}$. With the above
identifications, if $\s{E}\ind{h}=\so(l)$ for some $l\in\Z$, then
$\phi^i\ind{h}$ is the natural map $\id_{\so(l)}\otimes_{\K}
H^i(\gP,\cp{\s{L}}\otimes\cp{\phi})$, whereas if $\s{E}\ind{h}=\sd^j(j)$ with
$0<j<n$, then the image of $\phi^i\ind{h}$ in $\Hom_{\coq{}}$ coincides with
that of $\id_{\sd^j(j)}\otimes_{\K}H^i(\gP,\cp{\s{L}}\otimes\cp{\phi})$.
        \end{lemm}

        \begin{proof}
$\all m\in\Z$ by \ref{homotdist} there is a morphism of distinguished triangles
in $K^b(\gcoh{\gP})$
\begin{equation}\label{cd}
\begin{CD}
\s{U}^{<m}[-1] @>{\cp{\tbo{\s{U}^{\le m}}{m}}}>> \s{U}^m[-m]
@>{\cp{i}}>> \s{U}^{\le m} @>{\cp{p}}>> \s{U}^{<m} \\
@VV{\phi^{<m}[-1]}V @VV{\phi^m[-m]}V @VV{\phi^{\le m}}V @VV{\phi^{<m}}V \\
\s{V}^{<m}[-1] @>{\cp{\tbo{\s{V}^{\le m}}{m}}}>> \s{V}^m[-m] @>{\cp{j}}>>
\s{V}^{\le m} @>{\cp{q}}>> \s{V}^{<m}
\end{CD}
\end{equation}
(where $\cp{i}$ and $\cp{j}$ are the natural inclusions and $\cp{p}$ and
$\cp{q}$ the natural projections).

Let's denote for simplicity by $F$ the functor
$\gsec(\gP,\cp{\s{L}}\otimes-)_0$ (so that $R^iF=H^i(\gP,\cp{\s{L}}\otimes-)$).
We are going to show by induction on $m$ that
\begin{align*}
& R^iF(\s{U}^{<m})\iso\begin{cases}
R^0F(\s{U}^i) & \text{if $i<m$} \\
0 & \text{if $i\ge m$}
\end{cases}, &
R^iF(\cp{\tbo{\s{U}^{\le m}}{m}})=0
\end{align*}
$\all i,m\in\Z$ (and similarly for $\cp{\s{V}}$, of course), and that
$R^iF(\phi^{<m})$ can be identified with $R^0F(\phi^i)$ if $i<m$.

Since the above equations are trivially true for $m<<0$, we have only to prove
the inductive step from $m$ to $m+1$. Applying $F$ to \eqref{cd} and taking
into account that, by induction, $R^iF(\cp{\tbo{\s{U}^{\le m}}{m}})
=R^iF(\cp{\tbo{\s{V}^{\le m}}{m}})=0$, we obtain $\all i\in\Z$ a commutative
diagram with exact rows
\[\begin{CD}
0 @>>> R^{i-m}F(\s{U}^m) @>R^iF(\cp{i})>> R^iF(\s{U}^{\le m}) @>R^iF(\cp{p})>>
R^iF(\s{U}^{<m}) @>>> 0 @. \\
@. @VVR^{i-m}F(\phi^m)V @VVR^iF(\phi^{\le m})V @VVR^iF(\phi^{<m})V @. \\
0 @>>> R^{i-m}F(\s{V}^m) @>R^iF(\cp{j})>> R^iF(\s{V}^{\le m}) @>R^iF(\cp{q})>>
R^iF(\s{V}^{<m}) @>>> 0 @. .
\end{CD}\]
For $i\ne m$ we have $R^{i-m}F(\s{U}^m)=R^{i-m}F(\s{V}^m)=0$ by hypothesis, so
that $R^iF(\cp{p})$ and $R^iF(\cp{q})$ are isomorphisms (hence $R^iF(\phi^{\le
m})$ can be identified with $R^iF(\phi^{<m})$). On the other hand, by the
inductive hypothesis $R^mF(\s{U}^{<m})=R^mF(\s{V}^{<m})=0$, whence
$R^mF(\cp{i})$ and $R^mF(\cp{j})$ are isomorphisms. It follows that (as wanted)
\[R^iF(\s{U}^{\le m})\iso\begin{cases}
R^0F(\s{U}^i) & \text{if $i\le m$} \\
0 & \text{if $i>m$}
\end{cases}\]
(and similarly for $\cp{\s{V}}$) and that $R^iF(\phi^{\le m})$ can be
identified with $R^0F(\phi^i)$ if $i\le m$. It remains to prove that
$R^iF(\cp{\tbo{\s{U}^{\le m+1}}{m+1}})=0$. Now, for $i \ne m+1$ this is true
because the target of the map is $R^{i-m-1}F(\s{U}^{m+1}) =0$ by hypothesis,
whereas for $i=m+1$ observe that $\cdiff{\s{U}}^m= \cp{\tbo{\s{U}^{\le
m+1}}{m+1}}[m+1]\comp\cp{i}[m]:\s{U}^m\to\s{U}^{m+1}$, so that
\[R^0F(\cdiff{\s{U}}^m)=R^{m+1}F(\cp{\tbo{\s{U}^{\le m+1}}{m+1}})\comp
R^mF(\cp{i}).\]
As $R^0F(\cdiff{\s{U}}^m)=0$ because $\cdiff{\s{U}}^m$ is
minimal (remember that $R^0F(\s{E}\ind{h})\iso\K$ and that each component of
$\cdiff{\s{U}}^m$ from $\s{E}\ind{h}$ to $\s{E}\ind{h}$ is nilpotent by
\ref{minnil}) and as $R^mF(\cp{i})$ is an isomorphism, it follows that also
$R^{m+1}F(\cp{\tbo{\s{U}^{\le m+1}}{m+1}})=0$.

$\all i\in\Z$, taking $m>i$ such that $\s{U}^{<m}=\cp{\s{U}}$, we have
therefore natural isomorphisms
\[H^i(\gP,\cp{\s{L}}\otimes\cp{\s{U}})=R^iF(\cp{\s{U}})\iso R^iF(\s{U}^{<m})
\iso R^0F(\s{U}^i)=H^0(\gP,\cp{\s{L}}\otimes\s{U}^i)\] (and similarly for
$\cp{\s{V}}$) such that $H^i(\gP,\cp{\s{L}}\otimes\cp{\phi})$ can be identified
with $H^0(\gP,\cp{\s{L}}\otimes\phi^i)$. Then everything follows, as the
hypothesis on $\cp{\s{L}}$ implies that $\s{E}\ind{h}^{a^i_h}\iso
\s{E}\ind{h}\otimes_{\K}H^0(\gP,\cp{\s{L}}\otimes\s{U}^i)$ and that in this way
$\phi^i\ind{h}$ if $\s{E}\ind{h}=\so(l)$ (respectively the image of
$\phi^i\ind{h}$ in $\Hom_{\coq{}}$ if $\s{E}\ind{h}=\sd^j(j)$) can be
identified with $\id_{\s{E}\ind{h}}\otimes_{\K}
H^0(\gP,\cp{\s{L}}\otimes\phi^i)$.
        \end{proof}

    \begin{defi}
For $-\sw{\w}<l\le0$ let $\cp{\bo{l}}$\index{M_(l)@$\cp{\bo{l}}$} be the
subcomplex of $\cp{\sko}(-l)$ defined by
\[\bo{l}^j:=\bigoplus_{\card{I}=-j,\sw{\w_I}\le -l}\so(-l-\sw{\w_I})\subseteq
\bigoplus_{\card{I}=-j}\so(-l-\sw{\w_I})=\sko^j(-l)\] $\all j\in\Z$ (it is a
subcomplex because $\diff{\cp{\sko}(-l)}^j(\bo{l}^j) \subseteq\bo{l}^{j+1}$).
Similarly, for $n-\sw{\w}<l<0$ we denote by
$\cp{\bd{l}}$\index{N_(l)@$\cp{\bd{l}}$} the subcomplex of $\cp{\sko}(-l)$
defined by
\[\bd{l}^j:=\bigoplus_{\card{I}=-j,\sw{\w_I}<-l-j}\so(-l-\sw{\w_I})\subseteq
\bigoplus_{\card{I}=-j}\so(-l-\sw{\w_I})=\sko^j(-l)\] $\all j\in\Z$. Let
moreover $\cp{\abo{l}}$\index{M_(l) @$\cp{\abo{l}}$} and
$\cp{\abd{l}}$\index{N_(l) @$\cp{\abd{l}}$} be the complexes defined by the
exact sequences
\begin{gather*}
0\to\cp{\bo{l}}\to\cp{\sko}(-l)\to\cp{\abo{l}}[1]\to0,\\
0\to\cp{\bd{l}}\to\cp{\sko}(-l)\to\cp{\abd{l}}[1]\to0.
\end{gather*}
    \end{defi}

    \begin{rema}
By definition, we have
\begin{align*}
& \abo{l}^j\iso\bigoplus_{\card{I}=1-j,\sw{\w_I}>-l}\so(-l-\sw{\w_I}), &
\abd{l}^j\iso\bigoplus_{\card{I}=1-j,\sw{\w_I}>-l-j}\so(-l-\sw{\w_I}).
\end{align*}
    \end{rema}

        \begin{rema}\label{bobdvan}
$\bo{l}^j=\abo{l}^j=\bd{l}^j=\abd{l}^j=0$ if $j<-n$ or $j>0$. It is also
immediate to see that each $\bo{l}^j$ (respectively $\abo{l}^j$) contains
only terms of the form $\so(i)$ with $0\le i\le-l$ (respectively
$-l-\sw{\w}\le i<0$). Similarly, each $\bd{l}^j$ (respectively
$\abd{l}^j$) contains only terms of the form $\so(i)$ with $j<i\le-l$
(respectively $-l-\sw{\w}\le i<j$).
        \end{rema}

    \begin{rema}\label{altiso}
Since $\cp{\sko}(-l)\iso0$ in $D^b(\gcoh{\gP})$, the distinguished triangles
$\tri{\cp{\abo{l}}}{\cp{\bo{l}}}{\cp{\sko}(-l)}{}{}{}$ and
$\tri{\cp{\abd{l}}}{\cp{\bd{l}}}{\cp{\sko}(-l)}{}{}{}$ imply that $\cp{\abo{l}}
\iso\cp{\bo{l}}$ and $\cp{\abd{l}}\iso\cp{\bd{l}}$ in $D^b(\gcoh{\gP})$.
    \end{rema}

        \begin{lemm}\label{coeffsymm}
There are the following isomorphisms in $D^b(\gcoh{\gP})$:
\begin{enumerate}

\item $(\cp{\bo{l}})\dual\iso\cp{\bo{1-\sw{\w}-l}}(1)[-n]$ for $-\sw{\w}<l\le
0$;

\item $(\cp{\bd{l}})\dual\iso\cp{\bd{n-\sw{\w}-l}}(n)[-n]$ for $n-\sw{\w}<l<0$.
\end{enumerate}
    \end{lemm}

    \begin{proof}
Since the proofs are very similar, we prove only $1$. Dualizing the surjection
$\cp{\sko}(-l)\epi\cp{\abo{l}}[1]$ gives an injective map
$(\cp{\abo{l}})\dual[-1]\mono(\cp{\sko}(-l))\dual$. Since $(\cp{\sko}(-l))
\dual\iso\cp{\sko}(l+\sw{\w})[-n-1]$ by \ref{koselfdual}, we see that
$(\cp{\abo{l}})\dual(-1)[n]$ can be identified with a subcomplex of
$\cp{\sko}(-(1-\sw{\w}-l))$. Therefore, the natural isomorphisms ($\all j\in
\Z$)
\begin{multline*}
(\cp{\abo{l}})\dual(-1)[n]^j\iso\abo{l}^{-n-j}(1)\dual\iso
\bigoplus_{\card{I}=n+j+1,\sw{\w_I}>-l}\so(1-l-\sw{\w_I})\dual\\
\iso\bigoplus_{\card{I}=n+j+1,\sw{\w_I}>-l}\so(l+\sw{\w_I}-1)\iso
\bigoplus_{\card{\compl{I}}=-j,\sw{\w_{\compl{I}}}<\sw{\w}+l}
\so(l+\sw{\w}-1-\sw{\w_{\compl{I}}})\\
\iso\bigoplus_{\card{J}=-j,\sw{\w_J}\le\sw{\w}+l-1}\so(\sw{\w}+l-1-\sw{\w_J})
\iso\bo{1-\sw{\w}-l}^j
\end{multline*}
imply that $(\cp{\abo{l}})\dual(-1)[n]\iso\cp{\bo{1-\sw{\w}-l}}$ in
$C^b(\gcoh{\gP})$. Therefore, in view of \ref{altiso}, there is an isomorphism
$(\cp{\bo{l}})\dual\iso(\cp{\abo{l}})\dual\iso\cp{\bo{1-\sw{\w}-l}}(1)[-n]$ in
$D^b(\gcoh{\gP})$.
    \end{proof}

    \begin{lemm}\label{vancrit}
Let $\cp{\s{L}},\cp{\s{F}}\in D^b(\gcoh{\gP})$ and $a\in\Z$. If
$H^i(\gP,\s{L}^m\lotimes\cp{\s{F}})=0$ for $i>a-m$ (respectively
$i<a-m$)\footnote{$\lotimes$ denotes the left derived functor of $\otimes$.}
and $\all m\in\Z$, then $H^i(\gP,\cp{\s{L}}\lotimes\cp{\s{F}})=0$ for $i>a$
(respectively $i<a$).
    \end{lemm}

    \begin{proof}
We consider only the case with $>$ (the other one is similar). $\cp{\s{L}}$
being bounded, it is enough to prove that
$H^i(\gP,\s{L}^{<m}\lotimes\cp{\s{F}})=0$ for $i>a$ and $\all m\in\Z$. Since
this is obviously true for $m<<0$, we can proceed by induction on $m$. Applying
the functor $H^0(\gP,-\lotimes\cp{\s{F}})$ to the distinguished triangle
$\s{L}^{<m}[-1]\to\s{L}^m[-m]\to\s{L}^{\le m}\to\s{L}^{<m}$ we obtain
$\all i\in\Z$ an exact sequence
\[H^{i-m}(\gP,\s{L}^m\lotimes\cp{\s{F}})\to
H^i(\gP,\s{L}^{\le m}\lotimes\cp{\s{F}})\to
H^i(\gP,\s{L}^{<m}\lotimes\cp{\s{F}}).\]
Now, if $i>a$, $H^i(\gP,\s{L}^{<m}\lotimes\cp{\s{F}})=0$ by induction and
$H^{i-m}(\gP,\s{L}^m\lotimes\cp{\s{F}})=0$ by hypothesis ($i-m>a-m$), whence
also $H^i(\gP,\s{L}^{\le m}\lotimes\cp{\s{F}})=0$.
    \end{proof}

    \begin{theo}\label{wbeilfor}
Given $\cp{\s{F}}\in D^b(\gcoh{\gP})$ let $\cp{\gbro}=\cp{\gbro}(\cp{\s{F}})\in
C^b(\co{\into})$\index{X@$\cp{\gbro}$} and
$\cp{\gbrd}=\cp{\gbrd}(\cp{\s{F}})\in C^b(\com{\intd})$\index{Y@$\cp{\gbrd}$}
be the minimal complexes (which exist unique up to isomorphisms by
\ref{wbeiluni}) such that $\cp{\s{F}}\iso\cp{\gbro}\iso\cp{\gbrd}$ in
$D^b(\gcoh{\gP})$. Then $\all i\in\Z$ we have
\begin{align*}
\gbro^i & =\bigoplus_{-\sw{\w}<j\le0}\so(j)\otimes_{\K}
H^i(\gP,\cp{\s{F}}\otimes\cp{\bo{j}}),\\
\gbrd^i & =\bigoplus_{n-\sw{\w}<j<0}\so(j)\otimes_{\K}
H^i(\gP,\cp{\s{F}}\otimes\cp{\bd{j}})\bigoplus_{0\le j\le
n}\sd^j(j)\otimes_{\K}H^{i+j}(\gP,\cp{\s{F}}(-j)).
\end{align*}
Moreover, given $\cp{\phi}:\cp{\s{F}}\to\cp{\s{G}}$ in $D^b(\gcoh{\gP})$, let
$\cp{\gbro}(\cp{\phi}):\cp{\gbro}(\cp{\s{F}})\to \cp{\gbro}(\cp{\s{G}})$ and
$\cp{\gbrd}(\cp{\phi}):\cp{\gbrd}(\cp{\s{F}})\to \cp{\gbrd}(\cp{\s{G}})$ be the
induced morphisms in $K^b(\gcoh{\gP})$. Then $\all i\in\Z$ and for
$-\sw{\w}<j\le0$ the component of $\gbro^i(\cp{\phi})$ from $\so(j)\otimes_{\K}
H^i(\gP,\cp{\s{F}}\otimes\cp{\bo{j}})$ to $\so(j)\otimes_{\K}
H^i(\gP,\cp{\s{G}}\otimes\cp{\bo{j}})$ (respectively for $n-\sw{\w}<j<0$ the
component of $\gbrd^i(\cp{\phi})$ from $\so(j)\otimes_{\K}
H^i(\gP,\cp{\s{F}}\otimes\cp{\bd{j}})$ to $\so(j)\otimes_{\K}
H^i(\gP,\cp{\s{G}}\otimes\cp{\bd{j}})$) is given by the natural map
$\id_{\so(j)}\otimes_{\K}H^i(\gP,\cp{\phi}\otimes\cp{\bo{j}})$ (respectively
$\id_{\so(j)}\otimes_{\K}H^i(\gP,\cp{\phi}\otimes\cp{\bd{j}})$), whereas for
$0\le j\le n$ the image of the component of $\gbrd^i(\cp{\phi})$ from
$\sd^j(j)\otimes_{\K}H^{i+j}(\gP,\cp{\s{F}}(-j))$ to
$\sd^j(j)\otimes_{\K}H^{i+j}(\gP,\cp{\s{G}}(-j))$ in $\Hom_{\coq{\intd}}$
coincides with the image of the natural map $\id_{\sd^j(j)}\otimes_{\K}
H^i(\gP,\cp{\phi}(-j))$.
    \end{theo}

    \begin{proof}
Let $\s{E}\ind{j}:=\so(j)$ if $n-\sw{\w}<j<0$ and $\s{E}\ind{j}:=\sd^j(j)$ if
$0\le j\le n$. Define moreover $\cp{\bd{l}}:=\so(-l)[l]$ for $0\le l\le n$.
Then, by \ref{forcrit}, it is enough to prove that the following formulas hold
$\all i\in\Z$:
\begin{align*}
& h^i(\gP,\cp{\bo{l}}(j))=\delta_{j,l}\delta_{i,0} & & -\sw{\w}<j,l\le0 \\
& h^i(\gP,\cp{\bd{l}}\otimes\s{E}\ind{j})=\delta_{j,l}\delta_{i,0} & &
n-\sw{\w}<j,l\le n.
\end{align*}
Notice first that the case of $\cp{\bd{l}}$ for $0\le l\le n$ follows
immediately from \ref{omegacohom}.

Since $h^i(\gP,\bo{l}^m(j))=0$ for $i>-m$ and $\all m\in\Z$ (by \ref{lbcohom}
and \ref{bobdvan}), we have $h^i(\gP,\cp{\bo{l}}(j))=0$ if $i>0$ by
\ref{vancrit}. Moreover, also $h^i(\gP,\abo{l}^m(j))=0$ for $i<-m$ and $\all
m\in\Z$, and so $h^i(\gP,\cp{\bo{l}}(j))=h^i(\gP,\cp{\abo{l}}(j))=0$ if $i<0$,
again by \ref{vancrit}. In a similar way (using \ref{omegacohom} instead of
\ref{lbcohom} in the case $0\le j\le n$) one can prove that
$h^i(\gP,\cp{\bd{l}}\otimes\s{E}\ind{j})=0$ if $i\ne0$ for $n-\sw{\w}<l<0$ and
$n-\sw{\w}<j\le n$.

Therefore it is enough prove that $\chi(\gP,\cp{\bo{l}}(j))=\delta_{j,l}$ for
$-\sw{\w}<j,l\le0$ and $\chi(\gP,\cp{\bd{l}}\otimes\s{E}\ind{j})=\delta_{j,l}$
for $n-\sw{\w}<l<0$ and $n-\sw{\w}<j\le n$. Now, if $-\sw{\w}<j,l\le0$, by
\ref{lbcohom} and \ref{bobdvan} we have $\chi(\gP,\bo{l}^m(j))=
h^0(\gP,\bo{l}^m(j))=h^0(\gP,\sko^m(j-l))$ $\all m\in\Z$, and so
\begin{multline*}
\chi(\gP,\cp{\bo{l}}(j))=\sum_m(-1)^m\chi(\gP,\bo{l}^m(j))=
\sum_m(-1)^mh^0(\gP,\sko^m(j-l))\\
=\sum_m(-1)^m\sum_{\card{I}=-m}\pd_{j-l-\sw{\w_I}}=
\sum_I(-1)^{\card{I}}\pd_{j-l-\sw{\w_I}}=\delta_{j,l},
\end{multline*}
the last equality following from \ref{mondimfor}. Similarly one can prove that
$\chi(\gP,\cp{\bd{l}}(j))=\delta_{j,l}$ if $n-\sw{\w}<l,j<0$. It remains to
show that $\chi(\gP,\cp{\bd{l}}\otimes\sd^j(j))=0$ if $n-\sw{\w}<l< 0$ and
$0\le j\le n$. Taking into account that $\sd^j(j)\iso \sko^{\ge-j}(j)[-j]$ in
$D^b(\coh{\gP})$, we have
\begin{multline*}
\chi(\gP,\cp{\bd{l}}\otimes\sd^j(j))=
\chi(\gP,\cp{\bd{l}}\otimes\sko^{\ge-j}(j)[-j])\\
=\sum_{\card{I}\le j}(-1)^{\card{I}-j}\chi(\gP,\cp{\bd{l}}(j-\sw{\w_I}))=
\sum_{\card{I}\le j}\sum_m(-1)^{\card{I}-j+m}\chi(\gP,\bd{l}^m(j-\sw{\w_I}))\\
=\sum_{\card{I}\le j}\sum_{\sw{\w_J}<\card{J}-l}(-1)^{\card{I}-j+\card{J}}
\pd_{j-l-\sw{\w_I}-\sw{\w_J}}
\end{multline*}
(using the fact that $\chi(\gP,\bd{l}^m(j-\sw{\w_I}))=
h^0(\gP,\bd{l}^m(j-\sw{\w_I}))$, since $j-l-\sw{\w_I}-\sw{\w_J}>-\sw{\w}$ if
$\card{I}\le j$ and $\sw{\w_J}<\card{J}-l$). Hence we can write
$\chi(\gP,\cp{\bd{l}}\otimes\sd^j(j))=x-y$ with
\begin{align*}
x & :=\sum_{\card{I}\le j}\sum_J(-1)^{\card{I}-j+\card{J}}
\pd_{j-l-\sw{\w_I}-\sw{\w_J}},\\
y & :=\sum_{\card{I}\le j}\sum_{\sw{\w_J}\ge\card{J}-l}
(-1)^{\card{I}-j+\card{J}}\pd_{j-l-\sw{\w_I}-\sw{\w_J}}.
\end{align*}
Using \ref{mondimfor} we obtain
\begin{multline*}
x=\sum_{\card{I}\le j}(-1)^{\card{I}-j}\sum_J(-1)^{\card{J}}
\pd_{j-l-\sw{\w_I}-\sw{\w_J}}\\
=\sum_{\card{I}\le j}(-1)^{\card{I}-j}\delta_{\sw{\w_I},j-l}=
(-1)^j\sum_{\card{I}\le j,\sw{\w_I}=j-l}(-1)^{\card{I}}.
\end{multline*}
On the other hand, observing that $j-l-\sw{\w_I}-\sw{\w_J}<0$ if $\card{I}>j$
and $\sw{\w_J}\ge\card{J}-l$, again by \ref{mondimfor} we have
\begin{multline*}
y=\sum_{\sw{\w_J}\ge\card{J}-l}(-1)^{j+\card{J}}\sum_I(-1)^{\card{I}}
\pd_{j-l-\sw{\w_I}-\sw{\w_J}}=
\sum_{\sw{\w_J}\ge\card{J}-l}(-1)^{j+\card{J}}\delta_{\sw{\w_J},j-l}\\
=(-1)^j\sum_{\sw{\w_J}\ge\card{J}-l,\sw{\w_J}=j-l}(-1)^{\card{J}}=
(-1)^j\sum_{\card{J}\le j,\sw{\w_J}=j-l}(-1)^{\card{J}}.
\end{multline*}
Thus $x=y$ and $\chi(\gP,\cp{\bd{l}}\otimes\sd^j(j))=0$.
    \end{proof}

    \begin{rema}
Of course, in the particular case $\w=(1,\dots,1)$, the above resolutions
coincide with the usual Beilinson's resolutions on $\P^n$, taking into account
that $\gcoh{\gP^n}$ is equivalent to $\coh{\P^n}$ by \ref{shequivcrit} and
\ref{stdmori}. Indeed, in this case $n-\sw{\w}=-1$, and so
$\gbrd^i(\cp{\s{F}})=\bigoplus_{0\le j\le n}
\sd^j(j)\otimes_{\K}H^{i+j}(\gP^n,\cp{\s{F}}(-j))$. On the other hand,
$\cp{\bo{j}}=\sko^{\ge j}(-j)$, which is isomorphic in $D^b(\gcoh{\gP^n})$ to
$\sd^{-j}(-j)[-j]$ for $-n=1-\sw{\w}\le j\le0$, whence $\gbro^i(\cp{\s{F}})=
\bigoplus_{0\le j\le n}\so(-j)\otimes_{\K}
H^{i+j}(\gP^n,\cp{\s{F}}\otimes\sd^j(j))$.
    \end{rema}

The complexes $\cp{\bo{j}}$ and $\cp{\bd{j}}$ (including $\so(-j)[j]$ for $0\le
j\le n$) which appear in the resolutions of \ref{wbeilfor} are uniquely
determined as objects of $D^b(\gcoh{\gP})$. More precisely, we have the
following result.

    \begin{prop}
If $\cp{\s{P}}\in D^b(\gcoh{\gP})$ is such that for some $-\sw{\w}<a\le0$ and
$\all\cp{\s{F}}\in D^b(\gcoh{\gP})$ the coefficient of $\so(a)$ in
$\gbro^i(\cp{\s{F}})$ is given by $h^i(\gP,\cp{\s{F}}\otimes\cp{\s{P}})$, then
$\cp{\s{P}}\iso\cp{\bo{a}}$ in $D^b(\gcoh{\gP})$. Similarly, if $\cp{\s{Q}}\in
D^b(\gcoh{\gP})$ is such that for some $n-\sw{\w}<b<0$ (respectively $0\le b\le
n$) and $\all\cp{\s{F}}\in D^b(\gcoh{\gP})$ the coefficient of $\so(b)$
(respectively $\sd^b(b)$) in $\gbrd^i(\cp{\s{F}})$ is given by
$h^i(\gP,\cp{\s{F}}\otimes\cp{\s{Q}})$, then $\cp{\s{Q}}\iso\cp{\bd{a}}$
(respectively $\cp{\s{Q}}\iso\so(-b)[b]$) in $D^b(\gcoh{\gP})$.
    \end{prop}

    \begin{proof}
Possibly after substituting $\cp{\s{P}}$ with
$\cp{\gbro}(\cp{\s{P}}(1-\sw{\w}))(\sw{\w}-1)$, we can assume that $\cp{\s{P}}$
is minimal in $C^b(\co{[0,\sw{\w}[})$. Then, by \ref{wbeilfor} and by the
hypothesis on $\cp{\s{P}}$, we have $\all i\in\Z$:
\begin{multline*}
\s{P}^i(1-\sw{\w})=\bigoplus_{-\sw{\w}<j\le0}\so(j)\otimes_{\K}
H^i(\gP,\cp{\s{P}}(1-\sw{\w})\otimes\cp{\bo{j}})\\
=\bigoplus_{-\sw{\w}<j\le0}\so(j)\otimes_{\K}
H^i(\gP,\cp{\bo{a}}(1-\sw{\w})\otimes\cp{\bo{j}}).
\end{multline*}
As $\cp{\bo{a}}$ is also minimal in $C^b(\co{[0,\sw{\w}[})$ (so that
$\cp{\bo{a}}=\cp{\gbro}(\cp{\bo{a}}(1-\sw{\w}))(\sw{\w}-1)$), the last term is
equal to $\bo{a}^i(1-\sw{\w})$, whence $\s{P}^i=\bo{a}^i$ $\all i\in\Z$.
Similarly, $\cp{\tilde{\s{P}}}:=\cp{\gbro}(\cp{\s{P}}(1))(-1)$ (which is
minimal in $C^b(\co{[-\sw{\w},0[})$) is such that $\tilde{\s{P}}^i=\abo{a}^i$
$\all i\in\Z$.

Since $\Ext_{\gP}^i(\tilde{\s{P}}^l,\s{P}^m)=0$ for $i>0$ and $\all l,m\in\Z$,
it follows from \ref{HomKA=HomDA} that
$\Hom_{K^b(\gcoh{\gP})}(\cp{\tilde{\s{P}}},\cp{\s{P}})\isomor
\Hom_{D^b(\gcoh{\gP})}(\cp{\tilde{\s{P}}},\cp{\s{P}})$. Choosing $\alpha:
\cp{\tilde{\s{P}}}\to\cp{\s{P}}$ whose image in $D^b(\gcoh{\gP})$ is an
isomorphism, it is clear that $\cp{\s{L}}:=\mc{\alpha}$ satisfies the following
properties: $\cp{\s{L}}\iso0$ in $D^b(\gcoh{\gP})$ (i.e., $\cp{\s{L}}$ is an
exact complex), $\cp{\s{L}}$ is minimal and $\s{L}^i=\sko^i(-a)$ $\all i\in\Z$.
It is easy to see that this implies that $\cp{\s{L}}\iso\cp{\sko}(-a)$ in
$C^b(\gcoh{\gP})$ ($\s{L}^{\ge-n}$ and $\sko^{\ge-n}(-a)$ are isomorphic as
complexes because they are both $\cp{\gbro}$ resolutions of
$\so(-a-\sw{\w})[n+1]$, and it is immediate to check that an isomorphism
$\s{L}^{\ge-n}\to\sko^{\ge-n}(-a)$ extends to an isomorphism
$\cp{\s{L}}\to\cp{\sko}(-a)$). Therefore $\cp{\s{P}}$, which is a subcomplex of
$\cp{\s{L}}\iso\cp{\sko}(-a)$ by definition of mapping cone, can be identified
with $\cp{\bo{a}}$.

In a similar way (always using the $\cp{\gbro}$ resolution) one can prove that
$\cp{\s{Q}}\iso\cp{\bd{b}}$ or $\cp{\s{Q}}\iso\so(-b)[b]$.
    \end{proof}

    \begin{prop}\label{Odual}
The left dual of $(\so,\dots\so(\sw{\w}-1))$ (which is a full and strong
exceptional sequence of $D^b(\gcoh{\gP})$ by \ref{Oexc}) is the (full
exceptional by \ref{dualmut}) sequence
$(\cp{\bo{1-\sw{\w}}}[1-\sw{\w}],\dots,\cp{\bo{l}}[l],\dots\cp{\bo{0}})$.
    \end{prop}

    \begin{proof}
We will prove more generally that $L^{(j)}\so(i)\iso\cp{\bo{-j}}(i-j)[-j]$ for
$0\le j\le i<\sw{\w}$. We proceed by induction on $j$, the case $j=0$ being
clear, since $\cp{\bo{0}}\iso\so$. For $j>0$, by definition, there is a
distinguished triangle in $D^b(\gcoh{\gP})$
\[L^{(j)}\so(i)\to\cp{\s{H}}\mor{u}L^{(j-1)}\so(i)\to(L^{(j)}\so(i))[1],\]
where $\cp{\s{H}}:=\bigoplus_{k\in\Z}
(\Hom_{D^b(\gcoh{\gP})}(\so(i-j)[k],L^{(j-1)}\so(i))\otimes_{\K}\so(i-j)[k])$
and $u$ is the natural morphism. We can assume
$L^{(j-1)}\so(i)=\cp{\bo{1-j}}(i+1-j)[1-j]$ by induction, whence we have
\begin{multline*}
\s{H}^p=
\Hom_{D^b(\gcoh{\gP})}(\so(i-j)[-p],L^{(j-1)}\so(i))\otimes_{\K}\so(i-j)\\
\iso H^p(\gP,(L^{(j-1)}\so(i))(j-i))\otimes_{\K}\so(i-j)\iso
H^{p+1-j}(\gP,\cp{\bo{1-j}}(1))\otimes_{\K}\so(i-j)
\end{multline*}
$\all p\in\Z$. Now, $h^{p+1-j}(\gP,\cp{\bo{1-j}}(1))$ is the coefficient of
$\so(1-j)$ in $\gbro^{p+1-j}(\so(1))$ by \ref{wbeilfor}, and, as
$\cp{\gbro}(\so(1))$ is clearly given by $\sko^{<0}(1)[-1]$, we obtain
\[h^{p+1-j}(\gP,\cp{\bo{1-j}}(1))=\card{\{I\st\card{I}=j-p,\sw{\w_I}=j\}}.\]
Using \ref{HomKA=HomDA} it is immediate to see that $u$ is actually a morphism
of $K^b(\gcoh{\gP})$; then we can assume that $L^{(j)}\so(i)=\MC{u}[-1]$,
so that we have
\begin{multline*}
(L^{(j)}\so(i))^p=\s{H}^p\oplus(L^{(j-1)}\so(i))^{p-1}=
\s{H}^p\oplus\bo{1-j}^{p-j}(i+1-j)\\
\iso\so(i-j)^{\card{\{I\st\card{I}=j-p,\sw{\w_I}=j\}}}
\bigoplus_{\card{I}=j-p,\sw{\w_I}\le j-1}\so(i-\sw{\w_I})\iso
\bigoplus_{\card{I}=j-p,\sw{\w_I}\le j}\so(i-\sw{\w_I})\\
=\bo{-j}^{p-j}(i-j)=\cp{\bo{-j}}(i-j)[-j]^p.
\end{multline*}
It is also easy to check that the differential of $L^{(j)}\so(i)$ can
be identified with that of $\cp{\bo{-j}}(i-j)[-j]$.
    \end{proof}

    \begin{rema}\label{notstrong}
One can also prove directly that
$(\cp{\bo{1-\sw{\w}}},\dots,\cp{\bo{0}})$ is a (full) exceptional
sequence of $D^b(\gcoh{\gP})$ (see \ref{shiftexc}). Indeed, for
$-\sw{\w}<j,l\le0$ and $\all k\in\Z$ we have (taking into account
\ref{coeffsymm})
\begin{multline*}
\Hom_{D^b(\gcoh{\gP})}(\cp{\bo{l}},\cp{\bo{j}}[k])\iso
H^k(\gP,(\cp{\bo{l}})\dual\otimes\cp{\bo{j}})\\
\iso H^k(\gP,\cp{\bo{1-\sw{\w}-l}}(1)[-n]\otimes\cp{\bo{j}})\iso
H^{k-n}(\gP,\cp{\abo{1-\sw{\w}-l}}(1)\otimes\cp{\bo{j}}).
\end{multline*}
Now, by \ref{wbeilfor} the dimension of this $\K$--vector space is the
coefficient of $\so(j)$ in $\gbro^{k-n}(\cp{\abo{1-\sw{\w}-l}}(1))$. But
$\cp{\abo{1-\sw{\w}-l}}(1)$ is minimal in $C^b(\co{]-\sw{\w},0]})$ (by
\ref{bobdvan}), so that it coincides with
$\cp{\gbro}(\cp{\abo{1-\sw{\w}-l}}(1))$. Therefore
\[\gbro^{k-n}(\cp{\abo{1-\sw{\w}-l}}(1))\iso\abo{1-\sw{\w}-l}^{k-n}(1)\iso
\bigoplus_{\card{I}=n-k+1,\sw{\w_I}\ge\sw{\w}+l}\so(\sw{\w}+l-\sw{\w_I}),\]
and we conclude that
\begin{multline*}
\dim_{\K}\Hom_{D^b(\gcoh{\gP})}(\cp{\bo{l}},\cp{\bo{j}}[k])=
\card{\{I\st\card{I}=n-k+1,\sw{\w_I}=\sw{\w}+l-j\}}\\
=\card{\{I\st\card{\compl{I}}=k,\sw{\w_{\compl{I}}}=j-l\}}=
\card{\{J\st\card{J}=k,\sw{\w_J}=j-l\}}.
\end{multline*}
As this number is $0$ if $l>j$ and $\delta_{k,0}$ if $l=j$, we see
that $(\cp{\bo{1-\sw{\w}}},\dots,\cp{\bo{0}})$ is really an exceptional
sequence (it is also full, since it clearly generates the same
triangulated subcategory as $\{\so,\dots,\so(\sw{\w}-1)\}$). Moreover,
the above formula shows that, for general $\w$, there do not exist
integers $k_i$ (for $-\sw{\w}<i\le0$) such that
$(\cp{\bo{1-\sw{\w}}}[k_{1-\sw{\w}}],\dots,\cp{\bo{0}}[k_0])$ is also
strong. Indeed, a necessary condition for the existence of such $k_i$
is clearly that $\sw{\w_J}=\sw{\w_{J'}}$ implies $\card{J}=\card{J'}$,
and it is not difficult to prove that a necessary and sufficient
condition is $\w_0=\cdots=\w_n$ (in this case one can take
$k_i=q_i+a_{r_i}$, where $i=q_i\w_0-r_i$ with $0\le r_i<\w_0$ and
$a_0,\dots,a_{\w_0-1}\in\Z$ are arbitrary).
    \end{rema}

    \begin{rema}
Under the equivalence of categories mentioned in \ref{Dmod}, the
complexes $\cp{\bo{l}}$ correspond to the $B$--modules denoted by
$Q_{-l}$ in \cite{AKO}.
    \end{rema}

         \section{Some computations of the minimal resolution}

The explicit computation of $H^i(\gP,\cp{\s{F}}\otimes\cp{\bo{j}})$ or
$H^i(\gP,\cp{\s{F}}\otimes\cp{\bd{j}})$ is quite difficult in general. In the
latter case, however, one can always reduce to consider a similar problem on
the weighted projective space given by the weights $\w_i>1$, as we are going to
see.

        \begin{lemm}\label{Omegadecomp}
Let $\w=(\w_0,\dots,\w_n)$ be such that $\w_{m+1}=\cdots=\w_n=1$ for some $m\le
n$ and let $\w':=(\w_0,\dots,\w_m)$. Then, denoting by $\iota:\gP(\w')\to
\gP(\w)$ the closed immersion of graded schemes induced by the natural
epimorphism of graded rings $\p(\w)\to\p(\w')$, $\all j\in\Z$ there is an
isomorphism of graded sheaves on $\gP(\w')$
\[\iota^*\sd_{\gP(\w)}^j(j)\iso
\bigoplus_{0\le k\le j}\sd_{\gP(\w')}^k(k)^{\binom{n-m}{j-k}}.\]
        \end{lemm}

        \begin{proof}
Applying $\iota^*$ to the exact sequence
\[\sko_{(\w)}^{-j-2}\iso\bigoplus_{\card{I}=j+2}\so_{\gP(\w)}(-\sw{\w_I})
\mor{\diff{\cp{\sko}_{(\w)}}^{-j-2}}\sko_{(\w)}^{-j-1}\iso
\bigoplus_{\card{I}=j+1}\so_{\gP(\w)}(-\sw{\w_I})\to\sd_{\gP(\w)}^j\to0\]
and twisting by $j$ gives the exact sequence
\[\iota^*\sko_{(\w)}^{-j-2}(j)
\mor{\varphi:=\iota^*\diff{\cp{\sko}_{(\w)}(j)}^{-j-2}}
\iota^*\sko_{(\w)}^{-j-1}(j)\to\iota^*\sd_{\gP(\w)}^j(j)\to0.\]
Since $\iota^*\so_{\gP(\w)}(l)\iso\so_{\gP(\w')}(l)$ $\all l\in\Z$, we have
$\iota^*\sko_{(\w)}^{-j-h}(j)\iso\bigoplus_{J\subseteq\{m+1,\dots,n\}}
\s{A}^h(J)$, where
\[\s{A}^h(J):=\bigoplus_{\card{I}=j+h-\card{J}}
\so_{\gP(\w')}(j-\card{J}-\sw{\w'_I})\iso
\sko_{(\w')}^{\card{J}-j-h}(j-\card{J}).\]
As $\iota^*(x_l)=0$ for $l>m$, it is clear that the component of $\varphi$ from
$\s{A}^2(J)$ to $\s{A}^1(J')$ is $0$ if $J\ne J'$ and can be identified with
$\diff{\cp{\sko}_{(\w')}(j-\card{J})}^{\card{J}-j-2}$ if $J=J'$.
Therefore
\begin{multline*}
\iota^*\sd_{\gP(\w)}^j(j)\iso\bigoplus_{J\subseteq\{m+1,\dots,n\}}
\cok\diff{\cp{\sko}_{(\w')}(j-\card{J})}^{\card{J}-j-2}\\
\iso\bigoplus_{J\subseteq\{m+1,\dots,n\}}
\sd_{\gP(\w')}^{j-\card{J}}(j-\card{J})
\iso\bigoplus_{0\le k\le j}\sd_{\gP(\w')}^k(k)^{\binom{n-m}{j-k}}.
\end{multline*}
        \end{proof}

        \begin{coro}\label{weightred}
With the same notation used in \ref{Omegadecomp}, we have
\[h^i(\gP(\w),\cp{\s{F}}\otimes\cp{\bd{j}})=
h^i(\gP(\w'),L\iota^*(\cp{\s{F}})\otimes\cp{\bd{j}})\]
$\all i\in\Z$, $\all\cp{\s{F}}\in D^b(\gcoh{\gP(\w)})$ and for
$n-\sw{\w}=m-\sw{\w'}<j<0$.
        \end{coro}

        \begin{proof}
Since $\iota^*\cp{\gbrd}(\cp{\s{F}})\iso L\iota^*(\cp{\gbrd}(\cp{\s{F}}))\iso
L\iota^*(\cp{\s{F}})$ in $D^b(\gcoh{\gP(\w')})$, by \ref{Omegadecomp} the
resolution $\iota^*\cp{\gbrd}(\cp{\s{F}})$ of $L\iota^*(\cp{\s{F}})$ only
involves the sheaves of the form $\sd_{\gP(\w')}^j(j)$ for $0\le j\le m$ and
$\so_{\gP(\w')}(l)$ for $m-\sw{\w'}<l<0$. As $\iota^*\cp{\gbrd}(\cp{\s{F}})$ is
clearly minimal with respect to each $\so_{\gP(\w')}(l)$, the result follows
from \ref{wbeilfor} and \ref{minuniwrt}.
        \end{proof}

The following result (which will be used in chapter $3$) shows that,
up to twist, the Beilinson resolutions commute with taking the dual.

        \begin{lemm}\label{symbeilres}
$\all\cp{\s{F}}\in D^b(\gcoh{\gP})$ there are isomorphisms in $C^b(\gcoh{\gP})$
\begin{align*}
\cp{\gbro}(\cp{\s{F}})\dual(1-\sw{\w}) & \iso
\cp{\gbro}(R\gsHom_{\gP}(\cp{\s{F}},\so(1-\sw{\w}))),\\
\cp{\gbrd}(\cp{\s{F}})\dual(n-\sw{\w}) & \iso
\cp{\gbrd}(R\gsHom_{\gP}(\cp{\s{F}},\so(n-\sw{\w}))).
\end{align*}
        \end{lemm}

        \begin{proof}
It is clear that $\cp{\gbro}(\cp{\s{F}})\dual(1-\sw{\w})\iso
R\gsHom_{\gP}(\cp{\s{F}},\so(1-\sw{\w})):=\cp{\s{G}}$ in $D^b(\gcoh{\gP})$.
Since $\so(j)\dual(1-\sw{\w})\iso\so(1-\sw{\w}-j)$,
$\cp{\gbro}(\cp{\s{F}})\dual(1-\sw{\w})$ only involves sheaves of the form
$\so(j)$ for $-\sw{\w}<j\le0$. As it is also clearly minimal, it must coincide
with $\cp{\gbro}(\cp{\s{G}})$. Similarly, the statement about $\cp{\gbrd}$
follows from the fact that $\so(j)\dual(n-\sw{\w})\iso\so(n-\sw{\w}-j)$ and
$\sd^j(j)\dual(n-\sw{\w})\iso\sd^{n-j}(n-j)$.
        \end{proof}

In the case of $\P^n$ it is known (see \cite{AO} and \cite{EFS}) that in the
Beilinson resolution with the $\sd^j(j)$ every component of the differential
from $\sd^j(j)$ to $\sd^{j-1}(j-1)$ is given by the ``natural'' map. We are
going to see that, in a weaker form, the same fact is true in the weighted
case, as well. Fix $0<j\le n$: by \ref{omegahom} there is a natural isomorphism
\[\Hom_{\coq{\intd}}(\sd^j(j),\sd^{j-1}(j-1)\otimes_{\K}\p_1)\iso
\p_1\dual\otimes_{\K}\p_1\iso\Hom_{\K}(\p_1,\p_1).\]
Let's denote by $\upsilon$ the morphism corresponding to $\id_{\p_1}$.

    \begin{prop}\label{beildiff}
$\all\cp{\s{F}}\in D^b(\gcoh{\gP})$, $\all i\in\Z$ and for  $0<j\le n$ the
image in $\Hom_{\coq{\intd}}$ of the component
\[f^i_j(\cp{\s{F}}):\sd^j(j)\otimes_{\K}H^{i+j}(\gP,\cp{\s{F}}(-j))\to
\sd^{j-1}(j-1)\otimes_{\K}H^{i+j}(\gP,\cp{\s{F}}(1-j))\]
of $\diff{\cp{\gbrd}(\cp{\s{F}})}^i$ can be identified with the natural map
\[\xymatrix{\sd^j(j)\otimes_{\K}H^{i+j}(\gP,\cp{\s{F}}(-j))
\ar[r]^{f^i_j(\cp{\s{F}})} \ar[d]_{\upsilon\otimes\id} &
\sd^{j-1}(j-1)\otimes_{\K}H^{i+j}(\gP,\cp{\s{F}}(1-j)) \\
\sd^{j-1}(j-1)\otimes_{\K}\p_1\otimes_{\K}H^{i+j}(\gP,\cp{\s{F}}(-j))
\ar[ur]_{\id\otimes m}}\]
where $m:\p_1\otimes_{\K}H^{i+j}(\gP,\cp{\s{F}}(-j))\to
H^{i+j}(\gP,\cp{\s{F}}(1-j))$ is the multiplication map.\footnote{Remember
that, by \ref{lbcohom}, $\p_1$ is naturally isomorphic to $H^0(\gP,\so(1))$.}
    \end{prop}

    \begin{proof}
Given $\varphi\in H^{i+j}(\gP,\cp{\s{F}}(-j))=
R^{i+j}\gsec(\gP,\cp{\s{F}}(-j))_0$, by \ref{Yoneda} (since $\gsec(\gP,-(-j))_0
\iso\Hom_{\gP}(\so(j),-)$) we can consider $\varphi$ as a morphism from
$\so(j)[-i-j]$ to $\cp{\s{F}}$. It follows from \ref{wbeilfor} that there is a
commutative diagram in $\Hom_{\coq{\intd}}$
\[\begin{CD}
\sd^j(j)\otimes_{\K}H^0(\gP,\so)\iso\sd^j(j) @>f^i_j(\so(j)[-i-j])>>
\sd^{j-1}(j-1)\otimes_{\K}H^0(\gP,\so(1)) \\
@VV{\id\otimes\varphi}V @VV{\id\otimes\varphi}V \\
\sd^j(j)\otimes_{\K}H^{i+j}(\gP,\cp{\s{F}}(-j)) @>f^i_j(\cp{\s{F}})>>
\sd^{j-1}(j-1)\otimes_{\K}H^{i+j}(\gP,\cp{\s{F}}(1-j)).
\end{CD}\]
It is then clear that it is enough to prove that $f^i_j(\so(j)[-i-j])=
f^{-j}_j(\so(j))=\upsilon$.

Let's assume for simplicity that $\w_i>1$ if and only if $0\le i \le m$ for
some $0\le m\le n$, so that $H^0(\gP,\so(1))\iso\p_1=\langle x_{m+1},\dots,x_n
\rangle$. It is immediate to check that the maps
\[\begin{split}
g^l:\ko^{l-1}(j)=\bigoplus_{\card{I}=1-l}\p(j)dx_I & \to
\ko^l(j-1)\otimes_{\K}\p_1=(\bigoplus_{\card{I}=-l}\p(j-1)dx_I)\otimes_{\K}\p_1
\\
dx_I & \mapsto
\sum_{i\in I,i>m}(-1)^{\card{\{k\in I\st k<i\}}}dx_{I\minus\{i\}}\otimes x_i
\end{split}\]
define a morphism of complexes $\cp{g}:\cp{\ko}(j)[-1]\to\cp{\ko}(j-1)$. By
definition of $\sd^j$ and $\sd^{j-1}$ it is clear that restricting
$(\cp{g})\gsh{}$ yields a morphism of complexes (concentrated in positions
between $-j$ and $0$) $\cp{\psi}:\cp{\s{U}}\to\cp{\s{V}}\otimes_{\K}\p_1$,
where
\begin{align*}
\cp{\s{U}} & =0\to\sd^j(j)\to\sko^{-j}(j)\mor{-\diff{\cp{\sko}(j)}^{-j}}\cdots
\mor{-\diff{\cp{\sko}(j)}^{-2}}\sko^{-1}(j)\to0,\\
\cp{\s{V}} & =0\to\sd^{j-1}(j-1)\to\sko^{1-j}(j-1)
\mor{\diff{\cp{\sko}(j-1)}^{1-j}}\cdots\mor{\diff{\cp{\sko}(j-1)}^{-1}}
\sko^0(j-1)\to0.
\end{align*}
As $\cp{\s{U}}\iso\so(j)$ and $\cp{\s{V}}\iso0$ in $D^b(\gcoh{\gP})$, by
\ref{trcatprop} $\MC{\cp{\psi}[-1]}\iso\cp{\s{U}}\iso\so(j)$. Denoting by
$\cp{\s{W}}$ the minimal complex isomorphic to $\MC{\cp{\psi}[-1]}$ in
$K^b(\com{})$, it is then easy to see that $\s{W}^i=0$ for $i<-j$ or $i>0$,
$\s{W}^{-j}=\sd^j(j)$, $\s{W}^{1-j}=\sd^{j-1}(j-1)\otimes_{\K}\p_1\oplus\s{L}$
with $\s{L}\in\co{[n-\sw{\w},0]}$, that the component of $\cdiff{\s{W}}^{-j}$
from $\sd^j(j)$ to $\sd^{j-1}(j-1)\otimes_{\K}\p_1$ is given by $\psi^{-j}$ and
that $\s{W}^{\ge2-j}\in K^b(\co{]n-\sw{\w},j-2]})$ (notice that $\psi^0$ is
surjective, so that no term of the form $\so(j-1)$ survives). Applying
$\cp{\gbrd}$ to the distinguished triangle $\s{W}^{<2-j}[-1]\to\s{W}^{\ge2-j}
\to\cp{\s{W}}\to\s{W}^{<2-j}$ yields a distinguished triangle in
$K^b(\com{\intd})$
\[\cp{\gbrd}(\s{W}^{<2-j})[-1]\to\cp{\gbrd}(\s{W}^{\ge2-j})\to
\cp{\gbrd}(\cp{\s{W}})\to\cp{\gbrd}(\s{W}^{<2-j}).\] For $j-1\le l\le n$ and
$\all i\in\Z$, as $H^i(\gP,\s{W}^m(-l))=0$ for $m\ge 2-j$, by \ref{vancrit} we
have $H^i(\gP,\s{W}^{\ge2-j}(-l))=0$, whence $\cp{\gbrd}(\s{W}^{\ge2-j})$
contains no term of the form $\sd^l(l)$ for $j-1 \le l\le n$. Since moreover
$\s{W}^{<2-j}$ can be identified with $\cp{\gbrd}(\s{W}^{<2-j})$ and
$\cp{\gbrd}(\so(j))$ with $\cp{\gbrd}(\cp{\s{W}})$, it follows that
$f^{-j}_j(\so(j))$ is given by $\psi^{-j}$. It is not difficult to check that
its image in $\Hom_{\coq{\intd}}$ is $\upsilon$.
    \end{proof}

                \section{The theorem on $\P(\w)$}

From \ref{wbeilfor} we can obtain also a weaker version of the theorem for
$\coh{\P}$.

        \begin{prop}
Each of the two sets of coherent sheaves on $\P$
\begin{align*}
& \{\so(j)\st-\sw{\w}<j\le0\}, &
\{\so(j)\st n-\sw{\w}<j<0\}\cup\{\sd^j(j)\st0\le j\le n\}
\end{align*}
generates $D^b(\coh{\P})$ as a triangulated category.
        \end{prop}

        \begin{proof}
Since $\usmf{\gP}:\gcoh{\gP}\to\coh{\P}$ is exact and essentially surjective by
\ref{esssurj}, $\all\s{G}\in\coh{\P}$ there exists $\s{F}\in\gcoh{\gP}$ such
that $\s{F}_0\iso \s{G}$, and then $\cp{\gbro}(\s{F})_0$ and
$\cp{\gbrd}(\s{F})_0$ are two resolutions of $\s{G}$. The result follows, as
$\sd_{\gP}^j(l)_0\iso \sd_{\P}^j(l)$ $\all j,l\in\Z$.
        \end{proof}

        \begin{rema}
It is clear that, as in \ref{wbeilthm}, one can define functors
\begin{align*}
& \fo:K^b(\mco{\into})\to D^b(\coh{\P}), &
\fd:K^b(\mcom{\intd})\to D^b(\coh{\P})
\end{align*}
and the above result implies that they are essentially surjective, but, of
course, they are not fully faithful if $\w\ne(1,\dots,1)$. Indeed, even
assuming $\w$ normalized, one still has
\[\Hom_{\p}(\smd^j(l),\smd^{j'}(l'))\isomor\Hom_{\P}(\sd^j(l),\sd^{j'}(l'))\]
(this can be proved as in \ref{Homsmd=Homsd}, since the natural map $\p_l\to
\Hom_{\P}(\so(l'),\so(l+l'))$ is an isomorphism $\all l,l'\in\Z$ by \cite[lemma
4.1]{D}), but (as we already observed in \cite{Ca1}) the vanishing result
analogous to \ref{Beilextvan} is not true.
        \end{rema}

        \begin{rema}
In the notation of the above proof, $\cp{\gbro}(\s{F})_0$ and
$\cp{\gbrd}(\s{F})_0$ are resolutions of $\s{G}$, and they are minimal (at
least assuming $\w$ normalized), but they are not unique in general. To see
this, just observe that if $\w\ne(1,\dots,1)$, there exists $0\ne\s{F}\in
\gcoh{\gP}$ such that $\s{F}_0=0$, and clearly $\cp{\gbro}(\s{F})_0\ne0\ne
\cp{\gbrd}(\s{F})_0$.

On the other hand, there are at least two functorial choices for $\s{F}$,
namely $\s{G}\otimes_{\so_{\P}}\so_{\gP}$ and $\gsHom_{\P}(\so_{\gP},\s{G})$,
which don't coincide in general (anyway, they are at least both $0$ if
$\s{G}=0$!).
        \end{rema}

        \begin{rema}
$\usmf{\gP}:\gcoh{\gP}\to\coh{\P}$, being exact and essentially surjective,
clearly induces a surjective morphism of Grothendieck groups $K_0(\gP)\epi
K_0(\P)$ (where, of course, $K_0(\gP)$ is defined from $\gcoh{\gP}$).
\ref{wbeilthm} easily implies that $K_0(\gP)\iso\Z^{\sw{\w}}$ (a basis is given
by the elements of $\goset$ or of $\gdset$), and it is proved in \cite{A} that
$K_0(\P)\iso\Z^{n+1}$.
        \end{rema}

We conclude by remarking that the two Beilinson--type resolutions described in
\cite{Ca1} are both (non minimal) of the above forms. This is clear for the
first one, whereas for the second one, the point is that the sheaves of
logarithmic differentials introduced there actually decompose as direct sums of
sheaves $\so_{\P}(j)$. Indeed, in the notation of \cite{Ca1}, it is easy to see
that if $\emptyset\ne I\subseteq\{0,\dots,n\}$, there is an isomorphism of
$\p$--modules
\[\bar{\sd}^j_{\p}(\log x^I)\iso
\bigoplus_{J\cap I=\emptyset}\p(-\sw{\w_J})^{\binom{\card{I}-1}{j-\card{J}}}.\]
It follows that if $0\ne\chi\in\mu^*_{\w}$, there is an isomorphism of sheaves
on $\P$
\[\bar{\sd}^j_{\P}(\log x^I)(j-|\chi|)\iso\bigoplus_{J\cap I(\chi)=\emptyset}
\so_{\P}(j-|\chi|-\sw{\w_J})^{\binom{\card{I(\chi)}-1}{j-\card{J}}}.\] Notice
that $\binom{\card{I(\chi)}-1}{j-\card{J}}>0$ if and only if $j+1-
\card{I(\chi)}\le\card{J}\le j$, and that in this case $n-\sw{\w}<j-|\chi|-
\sw{\w_J}<0$ (since $J\cap I(\chi)=\emptyset$). Therefore the sheaves not of
the form $\sd_{\P}^j(j)$ (for $0\le j\le n$) which appear in the second
resolution of \cite[thm. 4.1]{Ca1}, all decompose as direct sums of sheaves of
the form $\so_{\P}(j)$ for $n-\sw{\w}<j<0$.

                \section{Application: a splitting criterion}

Beilinson's theorem can be used to give a simple proof of Horrock's splitting
criterion for vector bundles (or, more generally, torsion--free sheaves) on
$\P^n$ (see \cite{H}). Here we will show that also this result can be easily
extended to the weighted case. Of course, $\s{F}\in\gsmo{\gP}$ will be called
{\em torsion--free} if $\all\I{p}\in\gP$, given $\sigma\in\so_{\gP,\I{p}}$ and
$\tau\in\s{F}_{\I{p}}$, $\sigma\tau=0$ implies $\sigma=0$ or $\tau=0$.

In the following we will assume that the dimension $n$ of $\gP=\gP(\w)$ is
$>0$.

        \begin{lemm}\label{torcohom}
Let $0\ne\s{F}\in\gcoh{\gP}$ be a torsion--free sheaf. Then $\exi k\in\Z$ such
that $h^0(\gP,\s{F}(k))\ne0$ and $h^0(\gP,\s{F}(j))=0$ for $j<k$.
        \end{lemm}

        \begin{proof}
Since $\exi m\in\Z$ such that $h^0(\gP,\s{F}(m))\ne0$ (because
$\gsec(\s{F})\gsh{}\iso\s{F}\ne0$ by \ref{gmogqco}, and so $\gsec(\s{F})\ne0$),
it is enough to prove that $h^0(\gP,\s{F}(j))=0$ for $j<<0$. Notice that if $0
\ne\sigma\in H^0(\gP,\s{F}(j))$, then $\all i\in\N$ multiplication by $\sigma$
induces an injective map $H^0(\gP,\so(i))\mono H^0(\gP,\s{F}(i+j))$ (because
$\s{F}$ is torsion--free), whence $h^0(\gP,\s{F}(i+j))\ge\pd_{i}$. Setting $d:=
\gcd(\w_0,\dots,\w_n)$, let $c\in\N$ be such that $h^0(\gP,\s{F}(j))\le c$ for
$0\le j<d$. It is clear that there exists $l\in\N$ such that $\pd_{dl'}>c$ if
$l'\ge l$, and then $h^0(\gP,\s{F}(j))=0$ if $j\le-dl$.
        \end{proof}

        \begin{prop}\label{splitcrit}
Let $\s{F}\in\gcoh{\gP}$ be a torsion--free sheaf. Then $\s{F}$ decomposes as a
direct sum of line bundles $\so(l)$ if and only if $H^i(\gP,\s{F}(j))=0$ for
$0<i<n$ and $\all j\in\Z$.
        \end{prop}

        \begin{proof}
The other implication following immediately from \ref{lbcohom}, we can assume
that $0\ne\s{F}\in\gcoh{\gP}$ is a torsion--free sheaf. By \ref{torcohom} $\exi
k\in\Z$ such that $h^0(\gP,\s{F}(k))\ne0$ and $h^0(\gP,\s{F}(j))=0$ for $j<k$,
and we claim that $\gbrd^i(\s{F}(k))=0$ for $i<0$. Indeed, by the choice of $k$
it is clear that if $i<0$ and $0\le j\le n$ then $H^{i+j}(\gP,\s{F}(k-j))=0$.
Moreover, for $n-\sw{\w}<j<0$ and $\all m\in\Z$ we have
$H^i(\gP,\s{F}(k)\otimes\abd{j}^m)=0$ if $i<-m$, whence by \ref{vancrit}
$H^i(\gP,\s{F}\otimes\cp{\bd{j}})=H^i(\gP,\s{F}\otimes\cp{\abd{j}})=0$ for
$i<0$, and the claim follows from \ref{wbeilfor}. Therefore there is an
injective morphism of graded sheaves on $\gP$
\[\phi=\begin{pmatrix}
\phi^0 \\
\phi'
\end{pmatrix}:\s{F}(k)\mono
\gbrd^0(\s{F}(k))=\so\otimes_{\K}H^0(\gP,\s{F}(k))\oplus\s{E}.\]
As $H^0(\gP,\s{E})=0$, we see that $H^0(\gP,\phi^0)$ is an isomorphism, and
then it is clear that there exists $\psi:\so\otimes_{\K}H^0(\gP,\s{F}(k))\to
\s{F}(k)$ such that $\phi^0\comp\psi=\id$, which implies that $\s{F}\iso\so(-k)
\otimes_{\K}H^0(\gP,\s{F}(k))\oplus\s{F}'$. If $\s{F}'\ne0$, we can repeat the
procedure with $\s{F}'$ in place of $\s{F}$, and so on: we conclude in a finite
number of steps because $\s{F}$ is coherent.
        \end{proof}

The above result can be formulated also for ordinary sheaves on $\P$. Observe
that if $\s{F}\in\gcoh{\gP}$ is torsion--free, then $\s{F}_0\in\coh{\P}$ is
also torsion--free, but the converse is not true in general. On the other hand,
we have the following result.

        \begin{lemm}\label{torgtor}
If $\s{G}\in\coh{\P}$ is torsion--free, then $\s{F}:=
\gsHom_{\P}(\so_{\gP},\s{G})\in\gcoh{\gP}$ is torsion--free, too.
        \end{lemm}

        \begin{proof}
Given $\I{p}\in\gP$, let $R:=\so_{\gP,\I{p}}$ and $M:=\s{G}_{\I{p}}\in
\fmo{R_0}$. We have to prove that if $M$ is a torsion--free $R_0$--module,
then, given $f\in\Hom_{R_0}(R_i,M)$ and $r\in R_j$ (for some $i,j\in\Z$),
$f(rR_{i-j})=0$ implies $f=0$ or $r=0$. Assume on the contrary $f\ne0$ (say
$f(s)\ne0$ for some $s\in R_i$) and $r\ne0$. As $R$ is a {\qstd} domain, $\exi
t\in R_d$ invertible for some $d>0$, and $r^d=r't^j$ for some $0\ne r'\in R_0$.
Since $f(rR_{i-j})=0$, we have also $f(r^dR_{i-dj})=f(r't^jR_{i-dj})=0$, whence
$f(r'R_i)=0$ (because $t^j:R_{i-dj}\isomor R_i$). In particular, $0=f(r's)=
r'f(s)$, which implies $r'=0$ (whence $r=0$) or $f(s)=0$, a contradiction.
        \end{proof}

        \begin{coro}
Let $\s{G}\in\coh{\P}$ be a torsion--free sheaf. Then $\s{G}$ decomposes as a
direct sum of sheaves of the form $\so(l)$ ($l\in\Z$) if and only if
$H^i(\P,\sHom_{\P}(\so(j),\s{G}))=0$ for $0<i<n$ and $\all j\in\Z$.
        \end{coro}

        \begin{proof}
Assume first that $\s{G}\iso\bigoplus_m\so_{\P}(m)^{a_m}$ and notice that $\all
j,l\in\Z$ $\exi k\in\Z$ such that $\sHom_{\P}(\so(j),\so(l))\iso\so(k)$ (see
the proof of \ref{gPdualsh}). Then the required vanishing result follows from
the fact that $H^i(\P,\so(k))=0$ for $0<i<n$.

Conversely, let $\s{F}:=\gsHom_{\P}(\so_{\gP},\s{G})\in\gcoh{\gP}$: then by
\ref{goodcohom} we have
\[H^i(\gP,\s{F}(j))\iso H^i(\P,\s{F}_j)\iso
H^i(\P,\sHom_{\P}(\so(-j),\s{G}))=0\]
for $0<i<n$ and $\all j\in\Z$. As $\s{F}$ is torsion--free by \ref{torgtor}, it
follows from \ref{splitcrit} that $\s{F}\iso\bigoplus_m\so_{\gP}(m)^{a_m}$,
whence $\s{G}\iso\s{F}_0\iso\bigoplus_m\so_{\P}(m)^{a_m}$.
        \end{proof}

    \begin{rema}
If $\P=\P^n$, this result is a particular case of the equivalence (described
in \cite{W}) between the stable category of vector bundles on $\P^n$ and a
suitable subcategory of $D^b(\gmo{\p})$. We believe that also this more general
result can be extended to the (graded) weighted case.
    \end{rema}

                        \chapter{The theorem on weighted canonical projections}

In this chapter we prove the main result of this paper, namely an extension to
the (graded) weighted case of the theorem by Catanese and Schreyer on good
birational canonical projections of surfaces of general type proved in
\cite{CS} (where the case, which we will not consider, of good
projections of degree $2$ is also treated).

We will work over a fixed algebraically closed field $\K$ of characteristic
$\ne2$.\footnote{This assumption is needed both in the proof of \ref{symmres}
and for some standard results about pluricanonical maps of surfaces of general
type, which are well known in characteristic $0$ (see e.g. \cite{BPV}) and
whose proofs can be found in \cite{Ek} in the case of positive characteristic.}
Using standard notation, if $\S$\index{S_0@$\S$} is a minimal surface of
general type, $K_{\S}$\index{K_S_0@$K_{\S}$} will denote a canonical divisor on
$\S$ and $\ds_{\S}:=\so_{\S}(K_{\S})\iso\sd^2_{\S}$\index{omega_S_0@$\ds_{\S}$}
the dualizing sheaf. We recall that the basic numerical invariants of $\S$ are:

the geometric genus
$p_g(\S):=h^0(\S,\ds_{\S})=h^2(\S,\so_{\S})$;\index{p_g(S_0)@$p_g(\S)$}

the irregularity
$q(\S):=h^1(\S,\so_{\S})=h^1(\S,\ds_{\S})$;\index{q(S_0)@$q(\S)$}

the self--intersection of the canonical divisor
$K^2_{\S}$.\index{K_S_0^2@$K^2_{\S}$}

$\S$ being of general type, the plurigenera $h^0(\S,\ds_{\S}^m)$ (for $m>1$)
are given by $\chi(\so_{\S})+\frac{m(m-1)}{2}K^2_{\S}$ (where $\chi(\so_{\S}):=
\chi(\S,\so_{\S})=1+p_g(\S)-q(\S)$).\index{chi(S_0,O_S_0)@$\chi(\S,\so_{\S})$}
Thus all the cohomology groups of all the powers of $\ds_{\S}$ are completely
determined by $p_g(\S)$, $q(\S)$ and $K^2_{\S}$: indeed, we have also
$H^0(\S,\ds_{\S}^m)=0$ if $m<0$, $H^1(\S,\ds_{\S}^m)=0$ if $m\ne0,1$ and, by
Serre duality, $H^2(\S,\ds_{\S}^m) \iso H^0(\S,\ds_{\S}^{1-m})\dual$ $\all
m\in\Z$. Notice that, in particular,
$\chi(\S,\ds_{\S}^m)=\chi(\so_{\S})+\frac{m(m-1)}{2}K^2_{\S}$ $\all m\in\Z$.

After extending in a natural way the notion of weighted canonical
projection to the graded setting, we will show that the datum of a
{\gbwcp} $\gph:\gS\to\gY\subset\gP$ (where $\gS$ is the {\std} graded
scheme $\stdsch{\S}{\ds_{\S}}$ and $\gP=\gP(\w)$ is a $3$--dimensional
weighted projective space) is equivalent to the datum of a symmetric
minimal morphism
$\alpha:(\so\oplus\s{E})\dual(-1-\sw{\w})\to\so\oplus\s{E}$ of vector bundles
on $\gP$, satisfying some conditions, which are just the graded
weighted versions of those appearing in the theorem by Catanese and
Schreyer described in the introduction.

To this purpose, we first prove that for every good weighted canonical
projection $\gph:\gS\to\gY\subset\gP$ there exists a short exact
sequence in $\gcoh{\gP}$
\[0\to(\so\oplus\s{E})\dual(-1-\sw{\w})\mor{\alpha=\begin{pmatrix}
\alpha^{(1)} \\
\alpha'
\end{pmatrix}}\so\oplus\s{E}\to
\gph_*\so_{\gS}\to0,\]
where $\s{E}=\s{E}(\gph)$ is a direct sum of vector bundles of the form
$\sd^j(j-2)$ for $j=0,1,2$ and $\so(j-2)$ for $3-\sw{\w}<j<0$, and
$\alpha$ is a minimal (in the sense of \ref{mindef}) but
possibly not symmetric morphism. This is obtained as follows:
setting $\s{F}:=(\gph_*\so_{\gS})(2)$, the Beilinson resolution
$\cp{\gbrd}(\s{F})$ provided by \ref{wbeilfor} turns out to be a
symmetric complex, which however is not yet the one we are looking
for. Indeed, the most delicate point is to prove that
$\cp{\gbrd}(\s{F})$ is quasi--isomorphic to a length $1$ complex,
which (twisted by $\so(-2)$) yields the wanted exact sequence. For
this we use the functoriality of the Beilinson resolution, applied to
the natural morphism $\so(2)\to\s{F}$: here difficulties come from
the fact that in general $\cp{\gbrd}(\so(2))$ is not a subcomplex of
$\cp{\gbrd}(\s{F})$ (as it is if $\w=(1,1,1,1)$ and $\deg(\Y)>2$).
Anyway, eventually we manage (again, under the mild assumption
$\deg(\Y)>2$) to determine the explicit form of the vector bundle
$\s{E}$, in  which the coefficients of the $\sd^j(j-2)$ depend only on
$p_g(\S)$, $q(\S)$, $K^2_{\S}$ and $\w$ (the examples of chapter $4$
will show that instead the complicated expressions we get for the
coefficients of the $\so(j-2)$ really depend on $\gph$).

We also prove that an exact sequence as above is unique up to
isomorphism, and this fact implies that $\alpha$ can be chosen to be
symmetric, provided $\gph$ is birational. Then, as in the case of
projections to $\P^3$, the main result needed to complete the proof of
the theorem is that the rank condition
$\s{I}_{\rk{\s{E}}}(\alpha)=\s{I}_{\rk{\s{E}}}(\alpha')$ is
satisfied if and only if $\cok\alpha$ is a graded sheaf of commutative
$\so_{\gY}$--algebras.

                \section{Weighted canonical projections}

Let $\S$ be a minimal surface of general type. The canonical ring of $\S$ is
the positively graded ring
\[\R=\R(\S,\ds_{\S}):=\bigoplus_{d\ge0}H^{0}(\S,\ds_{\S}^d).\]\index{R@$\R$}
It is well known that $\R$ is a noetherian domain, so that $\R\in\prng$. The
canonical model of $\S$ is $\X:=\proj\R$.\index{X_0@$\X$} $\X$ is a surface
with only rational double points as singularities and there is a natural
birational morphism $\pr:\S\to\X$,\index{pi_0@$\pr$} which contracts
$(-2)$--curves (i.e. curves $C\iso\P^1$ with $C^2=-2$, $CK_{\S}=0$).
$\ds_{\X}\iso\so_{\X}(1)$ is an ample invertible sheaf ($\X$ is Gorenstein),
and there are isomorphisms $\all m\in\Z$:
\begin{align*}
& \pr^*\ds_{\X}^m\iso\ds_{\S}^m, & {\pr}_*\ds_{\S}^m\iso\ds_{\X}^m.
\end{align*}
The singularities of $\X$ being rational means that $R^1{\pr}_*(\so_{\S})=0$,
and this implies that $H^i(\S,\ds_{\S}^m)\iso H^i(\X,\ds_{\X}^m)$
$\all m\in\Z$.

Let $\gS$\index{S@$\gS$} be the {\std} graded scheme $\stdsch{\S}{\ds_{\S}}$
(i.e. $\gS$ has the same topological space as $\S$ and $\so_{\gS}=
\bigoplus_{d\in\Z}\ds_{\S}^d$) and let $\gX:=\gproj\R$.\index{X@$\gX$} Although
$\R$ is usually not generated by $\R_1$ as an $\R_0=\K$--algebra, we have the
following result.

        \begin{lemm}\label{canmodstd}
$\gX$ is a {\std} graded scheme.
        \end{lemm}

        \begin{proof}
Since the linear system $\linsys{mK_{\X}}$ is base point free for $m>>0$ ($m\ge
4$), $\all\I{p}\in\gX=\gproj\R$ there exist $s\in\R_m\minus\I{p}$ and $t\in
\R_{m+1}\minus\I{p}$ for some $m$. Then $t/s\in(\R_{\I{p}})_1$ is an invertible
element, i.e. $\R_{\I{p}}$ is \std.
        \end{proof}

By \ref{stdchar} the natural isomorphism $\pr^*\ds_{\X}\iso\ds_{\S}$ induces an
extension of $\pr$ to a morphism of graded schemes
$\gpr:\gS\to\gX$\index{pi@$\gpr$} such that $\gpr\mrs$ is an isomorphism.
Notice moreover that, since $\ds_{\X}\iso \so_{\X}(1)$, \ref{canmodstd} and
\ref{Serredual} imply that also $\ds_{\gX} \iso\so_{\gX}(1)$; by definition of
$\gS$ it is also clear that $\ds_{\gS}\iso \so_{\gS}(1)$.

Let $\w\in\N_+^{n+1}$ for some $n>0$. By \ref{projmor} the choice of sections
\begin{align*}
& \sigma_i\in H^0(\gS,\so_{\gS}(\w_i))\iso H^0(\S,\ds_{\S}^{\w_i}) &
i=0,\dots,n
\end{align*}
determines a rational map $\gph:\gS\rat\gP(\w)$\index{phi@$\gph$} (and also
$\ph:\S\rat\P(\w)$),\index{phi_0@$\ph$} called a {\em weighted canonical
projection}.\index{weighted canonical projection} $\gph$ always factors through
$\gpr:\gS\to\gX$, i.e. $\gph=\gps\comp\gpr$\index{psi@$\gps$} (and
$\ph=\ps\comp\pr$),\index{psi_0@$\ps$} where the rational map $\gps:\gX
\to\gP(\w)$ is the one induced by the morphism of graded rings\index{rho@$\rh$}
\[\begin{split}
\rh:\p(\w) & \to\R\\
x_i & \mapsto\sigma_i
\end{split}\]
determined by the chosen sections.

        \begin{defi}
A weighted canonical projection $\gph$ (or $\ph$) is called {\em
good}\index{good!weighted canonical projection} if it is a morphism; it is
called {\em birational}\index{birational weighted canonical projection} if it
is birational onto the image.
        \end{defi}

        \begin{rema}
$\ph$ is good if and only if $\gph$ is good if and only if $\gps$ (and $\ps$)
is a morphism. By \ref{prngmor} this happens if and only if $\rh$ is a finite
morphism, and in this case $\gps$ and $\ps$ are also finite.
        \end{rema}

        \begin{rema}\label{agalg}
$\gph$ (respectively $\ph$) is good birational if and only if $\gps$
(respectively $\ps$) is a morphism which is birational onto the image. This
happens if and only if $\rh$ is finite and $\exi p\in\p_+\hom$ such that
$p\notin\ker\rh$ and $\rh_p: \p_p\to\R_p$ (respectively
$\rh_{(p)}:\p_{(p)}\to\R_{(p)}$) is surjective. It is thus clear that if $\gph$
is good birational, then $\ph$ is good birational, too, but the converse is not
true in general (see \ref{gphagphag}).
        \end{rema}

        \begin{rema}\label{psnorm}
If $\ph$ is good birational, then $\ps$ is the normalization morphism (because
it is finite and birational and $\X$ is normal).
        \end{rema}

Keeping the notation introduced so far, from now on we will assume that $\gph:
\gS\to\gP=\gP(\w)$ is a good weighted canonical projection with $\w\in\N_+^4$
(i.e. $n=\dim\P=3$). We will denote by $\gY$\index{Y@$\gY$} (respectively
$\Y$)\index{Y_0@$\Y$} the image of $\gph$ in $\gP$ (respectively of $\ph$ in
$\P$). As $\gY$ is a hypersurface in $\gP$, we have $\gY=\gproj\p/(\f)$ (and
$\Y=\proj\p/(\f)$), where $\f\in\p=\p(\w)$\index{f@$\f$} is a homogeneous
irreducible polynomial which generates $\ker\rh$.\footnote{It is easy to see
that $\deg(\f)=(K^2_{\S}\prod_{i=0}^3\w_i)/\deg(\ps)$.}

Now we are going to prove a result on relative duality for the
morphism $\gps:\gX\to\gP$, which will be needed to obtain a symmetric
resolution of $\gps_*\so_{\gX}\iso\gph_*\so_{\gS}$. For later use, we
will prove it under more general hypotheses. So, let $\gamma:Z\to\gP$
be a finite morphism, where $Z$ is a (necessarily {\qstd} by \ref{YgoodXgood})
graded scheme with $Z_0$ of dimension $2$.\footnote{Actually, with obvious
modifications, everything holds for arbitrary dimensions of $Z$ and $\gP$.
Indeed, the proofs of the following results are essentially a rewriting of
those of \cite[III, lemma 7.3, lemma 7.4, prop. 7.5]{H1}, where it is proved
that every closed subscheme of $\P^n$ has a dualizing sheaf.} $\all i\in\N$ let
$\s{E}^i:=\gsExt^i_{\gP}(\gamma_*\so_Z,\ds_{\gP})$, and observe that
$\s{E}^i\in\gcoh{\gP}$ by \ref{gsExtgcoh} ($\gamma_*\so_Z$ is coherent because
$\gamma$ is a finite morphism).

        \begin{lemm}\label{gsExtvan}
If $H^{3-i}(Z,\so_Z(j))=0$ for $j<<0$, then $\s{E}^i=0$. In particular,
$\s{E}^0=0$.
        \end{lemm}

        \begin{proof}
As $\s{E}^i\iso\gsec(\s{E}^i)\gsh{}$ by \ref{gmogqco}, in order to prove that
$\s{E}^i=0$ it is enough to show that $H^0(\gP,\s{E}^i(d))=0$ for $d>>0$. By
\ref{secgsExt} $\all i\in\Z$ $\exi m\in\N$ such that
\[H^0(\gP,\s{E}^i(d))\iso\Ext^i_{\gP}(\gamma_*\so_Z,\ds_{\gP}(d))\iso
\Ext^i_{\gP}(\gamma_*\so_Z(-d),\ds_{\gP})\]
$\all d>m$. By Serre duality and since $\gamma$ is an affine morphism we have
\[\Ext^i_{\gP}(\gamma_*\so_Z(-d),\ds_{\gP})\iso
H^{3-i}(\gP,\gamma_*\so_Z(-d))\dual\iso
H^{3-i}(Z,\so_Z(-d))\dual,\]
and everything follows.
        \end{proof}

As clearly $\s{E}^1\in\gcoh{\gP}\cap\gsmo{\gamma_*\so_Z}$, by \ref{gSpecA}
there exists unique (up to isomorphism) $\s{D}\in\gcoh{Z}$ such that $\s{E}^1
\iso\gamma_*\s{D}$.

        \begin{lemm}\label{HomExt}
$\all\s{F}\in\gcoh{Z}$ there is a functorial isomorphism
\[\Hom_Z(\s{F},\s{D})\iso\Ext^1_{\gP}(\gamma_*\s{F},\ds_{\gP}).\]
        \end{lemm}

        \begin{proof}
$\all\s{G}\in\gqco{\gP}$ let $\gamma^!\s{G}\in\gqco{Z}$ be the graded sheaf
(which exists unique up to isomorphism by \ref{gSpecA}) such that
$\gamma_*\gamma^!\s{G}\iso\gsHom_{\gP}(\gamma_*\so_Z,\s{G})$. It is easy to see
that $\gamma^!:\gqco{\gP}\to\gqco{Z}$ is right adjoint of
$\gamma_*:\gqco{Z}\to\gqco{\gP}$.

Taking an injective resolution $\cp{\s{I}}$ of $\ds_{\gP}$ in $\gqco{\gP}$, we
have natural isomorphisms
\[\gamma_*H^i(\gamma^!\cp{\s{I}})\iso H^i(\gamma_*\gamma^!\cp{\s{I}})\iso
H^i(\gsHom_{\gP}(\gamma_*\so_Z,\cp{\s{I}}))\iso
\gsExt^i_{\gP}(\gamma_*\so_{Z},\ds_{\gP})=\s{E}^i\]
$\all i\in\N$. In particular, it follows from \ref{gsExtvan} that
$H^0(\gamma^!\cp{\s{I}})=0$. Since moreover $\gamma^!\s{I}^j$ is injective in
$\gqco{Z}$ ($\Hom_{\gqco{Z}}(-,\gamma^!\s{I}^j)\iso
\Hom_{\gqco{\gP}}(\gamma_*-,\s{I}^j)$ is exact because composition of two exact
functors), $\gamma^!\cp{\s{I}}$ splits as $\cp{\s{J}}\oplus\cp{\s{L}}$ with
$\cp{\s{L}}$ exact and $\s{J}^0=0$. Therefore
\[\gamma_*\s{D}\iso\s{E}^1\iso\gamma_*H^1(\gamma^!\cp{\s{I}})\iso
\gamma_*H^1(\cp{\s{J}}),\]
whence $\s{D}\iso H^1(\cp{\s{J}})=\ker\cdiff{\s{J}}^1$ by \ref{gSpecA}. We
have then natural isomorphisms
\begin{multline*}
\Hom_{Z}(\s{F},\s{D})\iso\ker\Hom_{Z}(\s{F},\cdiff{\s{J}}^1)\iso
H^1(\Hom_{Z}(\s{F},\cp{\s{J}}))\iso\\
\iso H^1(\Hom_{Z}(\s{F},\gamma^!\cp{\s{I}}))
\iso H^1(\Hom_{\gP}(\gamma_*\s{F},\cp{\s{I}}))\iso
\Ext^1_{\gP}(\gamma_*\s{F},\ds_{\gP}).
\end{multline*}
        \end{proof}

        \begin{prop}\label{gsExtdual}
Let $\gamma:Z\to\gP$ be a finite morphism of graded schemes with $Z_0$ of
dimension $2$. Then $Z$ has a dualizing sheaf $\ds_Z$ and
$\gamma_*\ds_{Z}\iso\gsExt^1_{\gP}(\gamma_*\so_{Z},\ds_{\gP})$.
        \end{prop}

        \begin{proof}
Using \ref{HomExt}, Serre duality on $\gP$ and the fact that $\gamma$ is a
finite morphism, we see that $\all\s{F}\in\gcoh{Z}$ there are natural
isomorphisms
\[\Hom_{Z}(\s{F},\s{D})\iso\Ext^1_{\gP}(\gamma_*\s{F},\ds_{\gP})\iso
H^2(\gP,\gamma_*\s{F})\dual\iso H^2(Z,\s{F})\dual.\]
Then $\s{D}\iso\ds_Z$ by definition of dualizing sheaf, whence
$\gamma_*\ds_{Z}\iso\gsExt^1_{\gP}(\gamma_*\so_{Z},\ds_{\gP})$.
        \end{proof}

                \section{Construction of the resolution}

Keeping the notation introduced so far, in this section we set
$\s{F}:=(\gps_*\so_{\gX})(2)\iso\gps_*(\so_{\gX}(2))\in\gcoh{\gP}$
(the twist is necessary in order to have a symmetric Beilinson resolution).

        \begin{lemm}\label{Fsym}
$R\gsHom_{\gP}(\s{F},\so(3-\sw{\w}))\iso\s{F}[-1]$ in $D^b(\gcoh{\gP})$. More
explicitly,
\[\gsExt^i_{\gP}(\s{F},\so_{\gP}(3-\sw{\w}))\iso
\begin{cases}
\s{F} & \text{if $i=1$} \\
0 & \text{if $i\ne1$}
\end{cases}.\]
        \end{lemm}

        \begin{proof}
As $H^{3-i}(\gX,\so_{\gX}(j))\iso H^{3-i}(\X,\ds_{\X}^j)=0$ for $j<0$ if
$i\ne1$, by \ref{gPdualsh} and \ref{gsExtvan}
$\gsExt^i_{\gP}(\s{F},\so_{\gP}(3-\sw{\w}))\iso
\gsExt^i_{\gP}(\gps_*\so_{\gX},\ds_{\gP})(1)=0$ if $i\ne1$, and by
\ref{gsExtdual}
\[\gsExt^1_{\gP}(\s{F},\so_{\gP}(3-\sw{\w}))\iso
\gsExt^1_{\gP}(\gps_*\so_{\gX},\ds_{\gP}(1))\iso\gps_*\ds_{\gX}(1)\iso\s{F}.\]
        \end{proof}

Let $\cp{\gbrd}:=\cp{\gbrd}(\s{F})$ and, for $3-\sw{\w}<j<0$ and $\all i\in\Z$
let $\cbrd^i_j:=h^i(\gP,\s{F}\otimes\cp{\bd{j}})$ be the coefficient of
$\so(j)$ in $\gbrd^i$.

        \begin{lemm}\label{Fbdres}
$(\cp{\gbrd})\dual(3-\sw{\w})\iso\cp{\gbrd}[-1]$ in $C^b(\gcoh{\gP})$.
Moreover, for $3-\sw{\w}<j<0$ we have $\cbrd^i_j=0$ if $i>0$ or $i<-1$ and
$\cbrd^{-1}_{3-\sw{\w}-j}=\cbrd^0_j:=\cbrd_j$.\index{y_j@$\cbrd_j$}
        \end{lemm}

        \begin{proof}
By \ref{Fsym} and \ref{symbeilres}  $(\cp{\gbrd})\dual(3-\sw{\w})\iso
\cp{\gbrd}[-1]$ in $C^b(\gcoh{\gP})$. Therefore $\cbrd^i_j=
\cbrd^{-i-1}_{3-\sw{\w}-j}$ for $3-\sw{\w}<j<0$ and $\all i\in\Z$, so that it
remains to prove that $\cbrd^i_j=0$ if $i>0$.  Taking into account that
$H^i(\gP,\s{F}(l))\iso H^i(\X,\ds_{\X}^{l+2})$, it is immediate to check that
$H^i(\gP,\s{F}\otimes\bd{j}^m)=0$ for $3-\sw{\w}<j<0$ and for $i>-m$. The
result then follows from \ref{vancrit}.
        \end{proof}

By \ref{wbeilfor} $\cp{\gbrd}$ has therefore the following form:
\[0\to\sd^2(2)\to\begin{matrix}
\so(3-\sw{\w})^{\chi(\so_{\S})+K^2_{\S}}\\
\oplus\sd^2(2)^{q(\S)}\oplus\sd^1(1)^{p_g(\S)}\\
\bigoplus_{3-\sw{\w}<j<0}\so(j)^{\cbrd_{3-\sw{\w}-j}}\end{matrix}
\to\begin{matrix}
\so^{\chi(\so_{\S})+K^2_{\S}}\\
\oplus\sd^1(1)^{q(\S)}\oplus\sd^2(2)^{p_g(\S)}\\
\bigoplus_{3-\sw{\w}<j<0}\so(j)^{\cbrd_j}\end{matrix}
\to\sd^1(1)\to0.\]

        \begin{lemm}\label{O2res}
Here we don't assume that the dimension $n$ of $\gP(\w)$ is $3$.
The Beilinson resolution $\cp{\gbrd}(\so(2))$ is given by the following complex
\[0\to\sd^2(2)\to\sd^1(1)^{\pd_1}\bigoplus_{n-\sw{\w}<j<0}\so(j)^{z^{-1}_j}
\to\so^{\pd_2}\bigoplus_{n-\sw{\w}<j<0}\so(j)^{z^0_j}\to0,\] where
$z^{-1}_j=\card{\{i<i'\st\w_i,\w_{i'}>1,\w_i+\w_{i'}=2-j\}}$ and
$z^0_j=\card{\{i\st\w_i=2-j\}}$.\index{z^i_j@$z^i_j$}
        \end{lemm}

        \begin{proof}
By \ref{wbeilfor} the only non trivial part is to show that
$h^i(\gP(\w),\cp{\bd{j}}(2))=z^i_j$ for $n-\sw{\w}<j<0$ and $\all i\in\Z$
(where $z^i_j:=0$ if $i\ne0,1$). We can assume that there exists $m\le n$ such
that $\w_i=1$ if and only if $i>m$. Then, setting $\w':=(\w_0,\dots,\w_m)$, by
\ref{weightred} we have
$h^i(\gP(\w),\cp{\bd{j}}(2))=h^i(\gP(\w'),\cp{\bd{j}}(2))$ (since
$L\iota^*\so_{\gP(\w)}(2)\iso\iota^*\so_{\gP(\w)}(2)\iso\so_{\gP(\w')}(2)$).
The complex on $\gP(\w')$ (obtained truncating $\cp{\sko}_{(\w')}(2)[-1]$)
\[0\to\sd^2(2)\to\sko^{-2}(2)\iso\bigoplus_{i<i'}\so(2-\w'_i-\w'_{i'})
\to\sko^{-1}(2)\iso\bigoplus_{i}\so(2-\w'_i)\to0\]
is a minimal resolution of
$\so_{\gP(\w')}(2)$ and is in $C^b(\com{]m-\sw{\w'},0[})$. Therefore, by the
uniqueness of the minimal Beilinson resolution, it must coincide with
$\cp{\gbrd}(\so_{\gP(\w')}(2))$. It follows that
$h^i(\gP(\w'),\cp{\bd{j}}(2))=z^i_j$.
        \end{proof}

        \begin{lemm}\label{coeffinj}
For $3-\sw{\w}<j<0$ the natural map
\[H^{-1}(\gP,\gps\mrs(2)\otimes\cp{\bd{j}}):H^{-1}(\gP,\cp{\bd{j}}(2))\iso
\K^{z^{-1}_j}\to H^{-1}(\gP,\s{F}\otimes\cp{\bd{j}})\iso
\K^{\cbrd_{3-\sw{\w}-j}}\]
is injective.
        \end{lemm}

        \begin{proof}
Let $\cp{\s{A}}:=\MC{\gps\mrs(2)}$: then there is an exact sequence
\[\cdots\to H^{-2}(\gP,\cp{\s{A}}\otimes\cp{\bd{j}})\to
H^{-1}(\gP,\cp{\bd{j}}(2))\to H^{-1}(\gP,\s{F}\otimes\cp{\bd{j}})\to\cdots,\]
so that it is enough to prove that $H^{-2}(\gP,\cp{\s{A}}\otimes\cp{\bd{j}})=0$
for $3-\sw{\w}<j<0$. By \ref{vancrit} this will follow if we show that
$H^i(\gP,\cp{\s{A}}\otimes\abd{j}^m)=0$ for $i<-1-m$ and $\all m\in\Z$. By the
exact sequence
\[H^i(\gP,\s{F}\otimes\abd{j}^m)\to H^i(\gP,\cp{\s{A}}\otimes\abd{j}^m)\to
H^{i+1}(\gP,\abd{j}^m(2))\to H^{i+1}(\gP,\s{F}\otimes\abd{j}^m)\]
this is a consequence of the fact that $H^i(\gP,\s{F}\otimes\abd{j}^m)=0$ and
$H^{i+1}(\gP,\gps\mrs(2)\otimes\abd{j}^m)$ is injective $\all m\in\Z$ and
$\all i<-1-m$ (this last fact is straightforward to check, if one remembers
that $H^i(\gP,\s{F}(l))\iso H^i(\X,\ds_{\X}^{l+2})$, taking into account
\ref{bobdvan}).
        \end{proof}

Let $\cp{\s{Z}}$ be the exact complex
\[0\to\sd^2(2)\to\sd^1(1)^{\pd_1}
\bigoplus_{3-\sw{\w}<j<0}\so(j)^{z^{-1}_j}\to
\so^{\pd_2}\bigoplus_{3-\sw{\w}<j<0}\so(j)^{z^0_j}\to\so(2)\to0\] obtained by
adding $\so(2)$ to its left resolution $\cp{\gbrd}(\so(2))$ (with
$\s{Z}^1=\so(2)$). The morphism $\cp{\gbrd}(\gps\mrs(2)):\cp{\gbrd}(\so(2))\to
\cp{\gbrd}=\cp{\gbrd}(\s{F})$ extends naturally to a morphism
$\cp{\varphi}:\cp{\s{Z}}\to \cp{\gbrd}$ in $K^b(\gcoh{\gP})$ (observe that
\[\Hom_{\gP}(\s{Z}^1,\gbrd^1)=\Hom_{\gP}(\so(2),\sd^1(1))\iso
H^0(\gP,\sd^1(-1))=0\]
by \ref{omcohvan}). We will moreover denote by $\cp{\tilde{\s{Z}}}$ the (exact)
complex $(\cp{\s{Z}})\dual(3-\sw{\w})$ and by $\cp{\tilde{\varphi}}$ the
morphism
\[\cp{\tilde{\varphi}}:=(\cp{\varphi})\dual(3-\sw{\w}):\cp{\gbrd}[-1]
\iso(\cp{\gbrd})\dual(3-\sw{\w})\to\cp{\tilde{\s{Z}}}\]
in $K^b(\gcoh{\gP}$. Notice that $\tilde{\s{Z}}^{-1}=\so(1-\sw{\w})$ and (by
\ref{symbeilres}) $\tilde{\s{Z}}^{\ge0}\iso\cp{\gbrd}(\so(1-\sw{\w}))$, which
is therefore given by
\[0\to\begin{matrix}
\so(3-\sw{\w})^{\pd_2}\\
\bigoplus_{3-\sw{\w}<j<0}\so(j)^{z^{0}_{3-\sw{\w}-j}}\end{matrix}
\to\begin{matrix}
\sd^2(2)^{\pd_1}\\
\bigoplus_{3-\sw{\w}<j<0}\so(j)^{z^{-1}_{3-\sw{\w}-j}}\end{matrix}
\to\sd^1(1)\to0.\]

        \begin{lemm}\label{HomKb=0}
$\Hom_{K^b(\gcoh{\gP})}(\cp{\s{Z}},\cp{\tilde{\s{Z}}}[1])=0$.
        \end{lemm}

        \begin{proof}
Let $\cp{d}:=\cp{\cdiff{\s{Z}}}$, $\cp{\tilde{d}}:=
\cp{\diff{\cp{\tilde{\s{Z}}}[1]}}$, and assume that $\cp{\alpha}:\cp{\s{Z}}\to
\cp{\tilde{\s{Z}}}[1]$
\[\begin{CD}
0 @>>> \s{Z}^{-2}=\sd^2(2) @>d^{-2}>> \s{Z}^{-1} @>d^{-1}>> \s{Z}^{0} @>d^0>>
\so(2)=\s{Z}^1 @>>> 0\\
@. @VV{\alpha^{-2}}V @VV{\alpha^{-1}}V @VV{\alpha^0}V @VV{\alpha^1}V @. \\
0 @>>> \tilde{\s{Z}}^{-1}=\so(1-\sw{\w}) @>\tilde{d}^{-2}>> \tilde{\s{Z}}^0
@>\tilde{d}^{-1}>> \tilde{\s{Z}}^1 @>\tilde{d}^0>> \sd^1(1)=\tilde{\s{Z}}^2
@>>> 0
\end{CD}\]
is a morphism of complexes. Since $\Hom_{\gP}(\so(2),\sd^1(1))=0$, we have
$\alpha^1=0$, and, by duality, also $\alpha^{-2}=0$. Therefore $\alpha^{-1}
\comp d^{-2}=0$, and so $\alpha^{-1}$ factors through $\im d^{-1}=\ker d^0$. As
$\Ext^1_{\gP}(\so(2),\tilde{\s{Z}}^0)\iso H^1(\gP,\tilde{\s{Z}}^0(-2))=0$ by
\ref{lbcohom}, it follows that $\alpha^{-1}$ actually factors through $d^{-1}$,
say $\alpha^{-1}=\beta\comp d^{-1}$ for some $\beta:\s{Z}^0\to\tilde{\s{Z}}^0$.
Setting $\gamma:=\alpha^0-\tilde{d}^{-1} \comp\beta$, we have $\gamma\comp
d^{-1}=0$, which implies that $\gamma=\beta' \comp d^0$ for some
$\beta':\so(2)\to\tilde{\s{Z}}^1$. As $\Hom_{\gP}(\so(2),\tilde{\s{Z}}^1)\iso
H^0(\gP,\tilde{\s{Z}}^1(-2))=0$ by \ref{omcohvan}, we see that
$\alpha^0=\tilde{d}^{-1}\comp\beta$, and this proves that $\cp{\alpha}$ is
homotopic to zero.
        \end{proof}

Consider the distinguished triangle of $K^b(\gcoh{\gP})$
\[\cp{\gbrd}[-1]\mor{\cp{\tilde{\varphi}}}\cp{\tilde{\s{Z}}}
\mor{\inc{\tilde{\varphi}}}\mc{\tilde{\varphi}}\mor{\pro{\tilde{\varphi}}}
\cp{\gbrd}\mor{\cp{\tilde{\varphi}}[1]}\cp{\tilde{\s{Z}}}[1].\]
As $\cp{\tilde{\s{Z}}}\iso0$ in $D^b(\gcoh{\gP})$, by \ref{trcatprop}
$[\pro{\tilde{\varphi}}]$ is an isomorphism in $D^b(\gcoh{\gP})$. Moreover,
$\cp{\tilde{\varphi}}[1]\comp\cp{\varphi}=0$ by \ref{HomKb=0}, so that, again
by \ref{trcatprop}, $\cp{\varphi}=\pro{\tilde{\varphi}}\comp\cp{\eta}$ for some
$\cp{\eta}:\cp{\s{Z}}\to\mc{\tilde{\varphi}}$ in $K^b(\gcoh{\gP})$. Again
because $\cp{\s{Z}}\iso0$ in $D^b(\gcoh{\gP})$, we see that $[\inc{\eta}]:
\mc{\tilde{\varphi}}\to\mc{\eta}$ is an isomorphism, i.e. there are
isomorphisms in $D^b(\gcoh{\gP})$
\[\mc{\eta}\iso\mc{\tilde{\varphi}}\iso\cp{\gbrd}\iso\s{F}.\]

        \begin{prop}\label{Fres}
In the notation introduced so far, let $\cp{\s{W}}\in K^b(\com{})$ be the
minimal complex isomorphic to $\mc{\eta}$ ($\cp{\s{W}}$ exists unique up to
isomorphism in $C^b(\com{})$ by \ref{minuni}). Then $\cp{\s{W}}\iso\s{F}$ in
$D^b(\gcoh{\gP})$ and, assuming that $\deg(\f)>2$, we have $\s{W}^i=0$ if
$i\ne-1,0$, $\s{W}^{-1}\iso(\s{W}^0)\dual(3-\sw{\w})$ and
\[\s{W}^0=\so(2)\oplus\so^{\chi(\so_{\S})+K_{\S}^2-\pd_2}\oplus\sd^1(1)^{q(\S)}
\oplus\sd^2(2)^{p_g(\S)-\pd_1}\bigoplus_{3-\sw{\w}<j<0}\so(j)^{\cre_j(\gph)},\]
where $\cre_j(\gph):=\cbrd_j-z^0_j-z^{-1}_{3-\sw{\w}-j}+
k_j+k_{3-\sw{\w}-j}$\index{c_j(phi)@$\cre_j(\gph)$} and (for $3-\sw{\w}<l<0$)
$k_l:=\dim_{\K} \ker H^0(\gP,\gps\mrs(2)\otimes\cp{\bd{l}})$, except that $-1$
must be added to the above formula for $\cre_j(\gph)$ if $\sw{\w}=2\w_i-1$ (for
some $0\le i\le 3$), $j=2-\w_i$, $k_j=1$ and
$\cp{\tilde{\varphi}}[1]\comp\cp{\varphi}\ne0\in C^b(\com{})$.\footnote{With
these assumptions, $\cp{\tilde{\varphi}}[1]\comp \cp{\varphi}\in C^b(\com{})$
only depends on the homotopy equivalence class of $\cp{\varphi}$. We don't know
if it can actually happen that it is not $0$.}
        \end{prop}

        \begin{proof}
By definition of mapping cone the complex $\cp{\s{V}}:=\mc{\eta}$ is given by
\[0\to\s{Z}^{-2}\to\gbrd^{-2}\oplus\s{Z}^{-1}\to\gbrd^{-1}\oplus
\tilde{\s{Z}}^{-1}\oplus\s{Z}^0\to\gbrd^0\oplus\tilde{\s{Z}}^0\oplus\s{Z}^1\to
\gbrd^1\oplus\tilde{\s{Z}}^1\to\tilde{\s{Z}}^2\to0\]
(with $\s{V}^{-3}=\s{Z}^{-2}$ and $\s{V}^2=\tilde{\s{Z}}^2$). Observe that
$\cp{\s{V}}\in C^b(\com{\{1-\sw{\w},2\}\cup]3-\sw{\w},0[})$ and that
$\cp{\s{V}}$ is already minimal with respect to $\so(2)$ and $\so(1-\sw{\w})
\iso\so(2)\dual(3-\sw{\w})$ (because they appear, respectively, only in
$\s{V}^0$ and $\s{V}^{-1}$, both with coefficient $1$). As for $\so$, we have
only $\so\otimes_{\K}H^0(\gP,\so(2))\subset\s{Z}^0\subset\s{V}^{-1}$ and
$\so\otimes_{\K}H^0(\gP,\s{F})\subset\gbrd^0\subset\s{V}^0$, and the map
between them, being the restriction of $\eta^0$ (which, by definition of
$\eta$, is the same as the restriction of $\varphi^0$), is just
$\id_{\so}\otimes H^0(\gP,\gps\mrs(2))$. As $\deg(\f)>2$, this last map is
injective, and so $\so$ only survives in $\s{W}^0$ with coefficient
$h^0(\gP,\s{F})-h^0(\gP,\so(2))=\chi(\so_{\S})+K_{\S}^2-\pd_2$. Now we consider
$\sd^1(1)$: we have
\begin{gather*}
\sd^1(1)\otimes_{\K}H^0(\gP,\so(1))\subset\s{Z}^{-1}\subset\s{V}^{-2},\\
\sd^1(1)\otimes_{\K}H^0(\gP,\s{F}(-1))\subset\gbrd^{-1}\subset\s{V}^{-1},\\
\sd^1(1)\otimes_{\K}H^1(\gP,\s{F}(-1))\subset\gbrd^0\subset\s{V}^0,\\
\sd^1(1)\otimes_{\K}H^2(\gP,\s{F}(-1))\iso
\sd^1(1)\otimes_{\K}H^0(\gP,\gps_*\so_{\gX})\dual\subset\gbrd^1\subset\s{V}^1,
\\
\sd^1(1)\otimes_{\K}H^3(\gP,\so(-\sw{\w}))\iso
\sd^1(1)\otimes_{\K}H^0(\gP,\so)\dual\subset\tilde{\s{Z}}^2=\s{V}^2.
\end{gather*}
The only non trivial maps between them in $\Hom_{\coq{}}$ are the restriction
of $\eta^{-1}$ (which is the same as the restriction of $\varphi^{-1}$) between
the first two terms and the restriction of $\tilde{\varphi}^2=
(\varphi^{-2})\dual(3-\sw{\w})$ between the last two terms. The former is
$\id_{\sd^1(1)}\otimes H^0(\gP,\gps\mrs(1))$ (which is injective) and the
latter is $\id_{\sd^1(1)}\otimes H^0(\gP,\gps\mrs)\dual$ (which is an
isomorphism). It follows that $\sd^1(1)$ only survives in $\s{W}^0$ with
coefficient $h^1(\gP,\s{F}(-1))=q(\S)$ and in $\s{W}^{-1}$ with coefficient
$h^0(\gP,\s{F}(-1))-h^0(\gP,\so(1))=p_g(\S)-\pd_1$ (notice that
$(\sd^1(1)^{p_g(\S)-\pd_1})\dual(3-\sw{\w})\iso\sd^2(2)^{p_g(\S)-\pd_1}$, the
expected term in $\s{W}^0$).
In a completely similar way (or by symmetry, see \ref{unisym}) one can proceed
in the case of $\sd^2(2)$ and $\so(3-\sw{\w})$.

It remains to consider the case of $\so(j)$ for $3-\sw{\w}<j<0$. $\all i\in\Z$
we have $\so(j)\otimes_{\K}V^i_j\subset\s{V}^i$, where
\begin{gather*}
V^{-2}_j=H^{-1}(\gP,\cp{\bd{j}}(2)),\\
V^{-1}_j=H^{-1}(\gP,\s{F}\otimes\cp{\bd{j}})\oplus H^0(\gP,\cp{\bd{j}}(2)),\\
V^0_j=H^0(\gP,\s{F}\otimes\cp{\bd{j}})\oplus
H^0(\gP,\cp{\bd{3-\sw{\w}-j}}(2))\dual,\\
V^1_j=H^{-1}(\gP,\cp{\bd{3-\sw{\w}-j}}(2))\dual
\end{gather*}
and $V^i_j=0$ if $i<-2$ or $i>1$. Setting
\[h^i_j:=H^i(\gP,\gps\mrs(2)\otimes\cp{\bd{j}}):
H^i(\gP,\cp{\bd{j}}(2))\to H^i(\gP,\s{F}\otimes\cp{\bd{j}})\]
and remembering that $H^i(\gP,\s{F}\otimes\cp{\bd{j}})\iso
H^{-1-i}(\gP,\s{F}\otimes\cp{\bd{3-\sw{\w}-j}})\dual$, it is clear that the
maps from $\so(j)\otimes_{\K}V^i_j$ to $\so(j)\otimes_{\K}V^{i+1}_j$ are given
by $\id_{\so(j)}\otimes d^i_j$, where
\begin{align*}
& d^{-2}_j=\begin{pmatrix}
h^{-1}_j \\
0
\end{pmatrix}, &
& d^{-1}_j=\begin{pmatrix}
0 & h^0_j \\
(h^0_{3-\sw{\w}-j})\dual & g
\end{pmatrix}, &
d^0_j=\begin{pmatrix}
(h^{-1}_{3-\sw{\w}-j})\dual & 0
\end{pmatrix}
\end{align*}
for some $g:H^0(\gP,\cp{\bd{j}}(2))\to H^0(\gP,\cp{\bd{3-\sw{\w}-j}}(2))\dual$
(observe that $d^{i+1}_j\comp d^i_j=0$ because $\cp{\s{V}}$ is a complex and
every morphism $\so(j)\to\s{E}\to\so(j)$ is $0$ by \ref{indec} if $\s{E}=
\sd^{j'}(j')$ or $\s{E}=\so(l)$ for some $l\ne j$). Let's denote by $\cre^i_j$
the coefficient of $\so(j)$ in $\s{W}^i$. As $h^{-1}_j$ is injective by
\ref{coeffinj} (whence also $(h^{-1}_{3-\sw{\w}-j})\dual$ is surjective), we
see that $\cre^i_j=0$ if $i<-1$ or $i>0$. Then by symmetry (see \ref{unisym})
it is enough to show that $\cre^0_j=\cre_j(\gph)$. Let's assume first that
$g=0$: in this case we have
\begin{multline*}
\cre^0_j=\dim_{\K}\ker d^0_j/\im d^{-1}_j=\dim V^0_j-\dim V^1_j-
\dim\im h^0_j-\dim\im h^0_{3-\sw{\w}-j}=\\
=\cbrd_j+z^0_{3-\sw{\w}-j}-z^{-1}_{3-\sw{\w}-j}-(z^0_j-k_j)-
(z^0_{3-\sw{\w}-j}-k_{3-\sw{\w}-j})=\cre_j(\gph).
\end{multline*}
Suppose now that $g\ne0$: we must have, in particular, $h^0(\gP,\cp{\bd{j}}(2))
=z^0_j\ne0$ and $h^0(\gP,\cp{\bd{3-\sw{\w}-j}}(2))=z^0_{3-\sw{\w}-j}\ne0$.
Assuming for simplicity $\w_0\le \w_1\le \w_2\le \w_3$, by \ref{O2res} this
happens if and only if $\sw{\w}=2\w_3-1$ (i.e. $\w_3=\w_0+\w_1+\w_2+1$) and
$j=2-\w_3$ (notice that then $j=3-\sw{\w}-j$), in which case $z^0_{2-\w_3}=1$
(whence $k_{2-\w_3}$ can only be $0$ or $1$). If $k_{2-\w_3}=0$ (so that
$h^0_{2-\w_3}$ is injective and $(h^0_{2-\w_3})\dual$ is surjective), it is
clear that $\im d^{-1}_j$ (and of course also $\ker d^0_j$) doesn't depend on
$g$; therefore also in this case $\cre^0_j=\cre_j(\gph)$. On the other hand, if
$k_{2-\w_3}=1$ (i.e. $h^0_{2-\w_3}=0$), $\dim_{\K}\im d^{-1}_j$ increases by
$1$ (with respect to the case $g=0$).

Thus, to conclude we have only to show that (under the hypothesis $\sw{\w}=
2\w_3-1$ and $h^0_{2-\w_3}=0$) $g=0$ if and only if $\cp{\tilde{\varphi}}[1]
\comp\cp{\varphi}=0\in C^b(\com{})$. Now, again by \ref{O2res} it is clear that
the restriction $\xi$ of $\cdiff{\s{V}}^{-1}$ from $\s{Z}^0$ to
$\tilde{\s{Z}}^0$ is just $\id_{\so(2-\w_3)}\otimes g$, whence $g=0$ if and
only if $\xi=0$. Denoting by $\cp{d}$ the differential of $\cp{\s{Z}}$, the
fact that $\Hom_{\gP}(\s{Z}^1,\tilde{\s{Z}}^0)=0$ implies that $\xi=0$ if and
only if $\xi\comp d^{-1}=0$. Moreover, by definition of mapping cone $\xi$ is
also the component of $\eta^0$ from $\s{Z}^0$ to $\tilde{\s{Z}}^0$, and the
condition that $\eta$ is a morphism of complexes easily implies (since
$\Hom_{\gP}(\s{Z}^{-1},\tilde{\s{Z}}^{-1})=0$) that $\xi\comp d^{-1}=
\tilde{\varphi}^0\comp\varphi^{-1}$. Then it is enough to observe that if this
last term is $0$, also $\tilde{\varphi}^1\comp\varphi^0=(\tilde{\varphi}^0\comp
\varphi^{-1})\dual(3-\sw{\w})=0$, so that $\cp{\tilde{\varphi}}[1]\comp
\cp{\varphi}=0$ (in any case $\tilde{\varphi}^{i+1}\comp\varphi^i=0$ if $i\ne
-1,0$).
        \end{proof}

        \begin{rema}\label{unisym}
A minimal complex $\cp{\s{L}}\in C^b(\com{\{1-\sw{\w},2\}\cup]3-\sw{\w},0[})$
with $\s{L}^i=0$ for $i\ne-1,0$ and quasi--isomorphic to $\s{F}$ is unique up
to isomorphism in $C^b(\gcoh{\gP})$ (hence it must coincide with the complex
$\cp{\s{W}}$ defined in \ref{Fres}). Indeed, if $\cp{\s{M}}$ is another complex
with the same properties, since $\Ext^1_{\gP}(\s{L}^0,\s{M}^{-1})=0$ (as it is
easy to check), by \ref{AfullDA} and \ref{HomKA=HomDA} there are natural
isomorphisms
\[\Hom_{\gP}(\s{F},\s{F})\iso\Hom_{D^b(\gcoh{\gP})}(\cp{\s{L}},\cp{\s{M}})\iso
\Hom_{K^b(\gcoh{\gP})}(\cp{\s{L}},\cp{\s{M}}).\]
It follows that $\cp{\s{L}}\iso\cp{\s{M}}$ in $K^b(\gcoh{\gP})$, whence also in
$C^b(\gcoh{\gP})$ by \ref{minuni}.

Moreover, it is clear a priori that $\s{L}$ must satisfy $\s{L}^{-1}\iso
\s{L}^0(3-\sw{\w})$: applying $\gsHom_{\gP}(-,\so(3-\sw{\w}))$ to the exact
sequence $0\to\s{L}^{-1}\to\s{L}^0\to\s{F}\to0$, by \ref{Fsym} we obtain the
exact sequence
\[0\to(\s{L}^0)\dual(3-\sw{\w})\to(\s{L}^{-1})\dual(3-\sw{\w})\to
\gsExt^1_{\gP}(\s{F},\so(3-\sw{\w}))\iso\s{F}\to0.\]
Then the result follows from the uniqueness proved above.
        \end{rema}

                \section{Symmetric resolution}

    \begin{defi}\label{Edef}
If $\gph:\gS\to\gP$ is a good weighted canonical projection, let $\s{E}(\gph)$
be the vector bundle on $\gP$ defined by
\[\s{E}(\gph):=\so(-2)^{\chi(\so_{\S})+K_{\S}^2-\pd_2}\oplus
\sd^1(-1)^{q(\S)}\oplus(\sd^2)^{p_g(\S)-\pd_1}
\bigoplus_{3-\sw{\w}<j<0}\so(j-2)^{\cre_j(\gph)},\]\index{E(phi)@$\s{E}(\gph)$}
with $\cre_j(\gph)$ defined as in \ref{Fres}.
    \end{defi}

If $\gph:\gS\to\gP$ is a good weighted canonical projection and $\deg(\f)>2$,
it follows from \ref{Fres} that there is an exact sequence on $\gP$
\begin{equation}\label{symres}
0\to(\so\oplus\s{E}(\gph))\dual(-1-\sw{\w})\mor{\alpha}\so\oplus\s{E}(\gph)\to
\gph_*\so_{\gS}\iso\gps_*\so_{\gX}\to0,
\end{equation}
where $\alpha$ is a minimal morphism.

From now on we will assume that $\gph$ is good birational (this implies that
$\deg(\f)>2$, as it is easy to see). We are going to prove that, with this
hypothesis, also $\alpha$ in \eqref{symres} can be taken to be symmetric (i.e.,
$\alpha=\alpha\dual(-1-\sw{\w})$).

        \begin{lemm}\label{gphagphag}
$\gph$ is good birational if and only if $\gcd(\w_0,\w_1,\w_2,\w_3)=1$ and
$\ph$ is good birational.
        \end{lemm}

        \begin{proof}
First observe that, since $\gS$ is \std, given $s\in\gS$,
$\gph\mrs_s:\so_{\gY,\gph(s)}\to\so_{\gS,s}$ is an isomorphism if and only if
$(\ph\mrs)_s=(\gph\mrs_s)_0$ is an isomorphism and $\gph(s)$ is \std. Therefore
it is enough to prove that the (open) set $V$ of {\std} points of $\gY$ is
empty if and only if $d:=\gcd(\w_0,\w_1,\w_2,\w_3)\ne1$. By \ref{YgoodXgood}
$V=U\cap\gY$, where $U$ denotes the set of {\std} points of $\gP$. Moreover, by
\ref{stdmori}, $U=\emptyset$ if $d\ne1$, whereas $\gP\minus
U\subseteq\V{(x_0x_1x_2x_3)}$ if $d=1$. To conclude, just notice that, if $\ph$
is birational, then $\gY\nsubseteq\V{(x_0x_1x_2x_3)}$ (otherwise $\Y$ would be
isomorphic to a $2$--dimensional weighted projective space, which is a rational
surface).
        \end{proof}

        \begin{lemm}\label{det=Y}
$\det(\alpha)=\lambda\f$ for some $\lambda\in\K^*$.
        \end{lemm}

    \begin{proof}
Similar to that of \cite[lemma 2.11]{C1}.
    \end{proof}

        \begin{lemm}\label{birend}
$\Hom_{\gP}(\gps_*\so_{\gX},\gps_*\so_{\gX})\iso\K$.
        \end{lemm}

        \begin{proof}
Since $\gps_*\so_{\gX}\iso\R\gsh{\p}$,
$\Hom_{\gP}(\R\gsh{\p},\R\gsh{\p})\iso\Hom_{\p}(\R,\gsec(\R\gsh{\p}))$ (by
\ref{gmogqco}) and $\gsec(\R\gsh{\p})\iso\R$ ($\all m\in\Z$ the natural map
\[\R_m\to\gsec(\R\gsh{\p})_m=H^0(\gP,\R\gsh{\p}(m))\iso H^0(\gX,\so_{\gX}(m))
\iso H^0(\X,\ds_{\X}^m)\]
is an isomorphism), it is enough to prove that $\Hom_{\p}(\R,\R)\iso\K$.

Given $\varphi\in\Hom_{\p}(\R,\R)$, let $\lambda:=\varphi(1)\in\R_0\iso\K$. As
$\gps$ is birational, by \ref{agalg} $\exi p\in\p_+\hom$ such that $p\notin
\ker\rh$ and $\rh_p:\p_p\to\R_p$ is surjective. This means that $\all r\in\R$
we can find $m\in\N$ and $p'\in\p$ such that $p^mr=\rh(p')$. Therefore
\[p^m\varphi(r)=\varphi(p^mr)=\varphi(\rh(p'))=\lambda\rh(p')=p^m\lambda r.\]
As $\R$ is a domain and $\rh(p)\ne0$, this implies that $\varphi(r)=\lambda r$,
whence $\varphi=\lambda\id_{\R}$.
        \end{proof}

        \begin{prop}\label{symmres}
If $\gph$ is a \gbwcp, then there exists an exact sequence on $\gP$
\[0\to(\so\oplus\s{E}(\gph))\dual(-1-\sw{\w})\mor{\alpha}\so\oplus\s{E}(\gph)
\to\gph_*\so_{\gS}\iso\gps_*\so_{\gX}\to0,\]
where $\alpha$ is a minimal morphism and $\alpha=\alpha\dual(-1-\sw{\w})$.
        \end{prop}

        \begin{proof}
Denoting $\so\oplus\s{E}(\gph)$ by $\s{A}$, $\dual(-1-\sw{\w})$ by
$\tdual$\index{$\tdual$} and $\gsExt^1_{\gP}(-,\so(-1-\sw{\w}))$ by $E$, we
already know that there is an exact sequence
\[0\to\s{A}\tdual\mor{\beta}\s{A}\mor{\pi}\gph_*\so_{\gS}\to0\]
with $\beta$ minimal. Applying $\gsHom_{\gP}(-,\so(-1-\sw{\w}))$ we obtain the
exact sequence
\[0\to\s{A}\tdual\mor{\beta\tdual}\s{A}\mor{\pi'}E(\gph_*\so_{\gS})\to0\]
and applying it once more we have also the exact sequence
\[0\to\s{A}\tdual\mor{\beta}\s{A}\mor{\pi''}E(E(\gph_*\so_{\gS}))\to0\]
(since $\dtdual$ is naturally isomorphic to $\id$).
Notice that there exists a unique isomorphism $\iota:\gph_*\so_{\gS}\to
E(E(\gph_*\so_{\gS}))$ such that $\pi''=\iota\comp\pi$. Choosing an isomorphism
$\epsilon:\gph_*\so_{\gS}\isomor E(\gph_*\so_{\gS})$, we know from \ref{unisym}
that there exists (unique up to homotopy) an isomorphism of complexes
\[\begin{CD}
0 @>>> \s{A}\tdual @>{\beta}>> \s{A} @>{\pi}>> \gph_*\so_{\gS} @>>> 0 \\
@. @VV{\eta}V @VV{\xi}V @VV{\epsilon}V @. \\
0 @>>> \s{A}\tdual @>{\beta\tdual}>> \s{A} @>{\pi'}>> E(\gph_*\so_{\gS}) @>>> 0
.
\end{CD}\]
Applying $\gsHom_{\gP}(-,\so(-1-\sw{\w}))$ to the above diagram yields the
commutative diagram
\[\begin{CD}
0 @>>> \s{A}\tdual @>{\beta}>> \s{A} @>{\pi}>> \gph_*\so_{\gS} @>>> 0 \\
@. @VV{\xi\tdual}V @VV{\eta\tdual}V @VV{E(\epsilon)\comp\iota}V @. \\
0 @>>> \s{A}\tdual @>{\beta\tdual}>> \s{A} @>{\pi'}>> E(\gph_*\so_{\gS}) @>>> 0
,
\end{CD}\]
and applying it once more we obtain the commutative diagram
\[\begin{CD}
0 @>>> \s{A}\tdual @>{\beta}>> \s{A} @>{\pi}>> \gph_*\so_{\gS} @>>> 0 \\
@. @VV{\eta}V @VV{\xi}V @VV{E(E(\epsilon)\comp\iota)\comp\iota}V @. \\
0 @>>> \s{A}\tdual @>{\beta\tdual}>> \s{A} @>{\pi'}>> E(\gph_*\so_{\gS}) @>>> 0
.
\end{CD}\]
It follows that $E(E(\epsilon)\comp\iota)\comp\iota=\epsilon$. On the other
hand, by \ref{birend} $E(\epsilon)\comp\iota=\lambda\epsilon$ for some $\lambda
\in\K$. Therefore, as $E$ is $\K$--linear,
\[\epsilon=E(E(\epsilon)\comp\iota)\comp\iota=E(\lambda\epsilon)\comp\iota=
\lambda E(\epsilon)\comp\iota=\lambda^2\epsilon,\]
whence $\lambda=\pm1$. Moreover, since $(\eta,\xi,\epsilon)$ and $(\xi\tdual,
\eta\tdual,\lambda\epsilon)$ are morphisms between the same complexes, so is
also
\[(\eta+\lambda\xi\tdual,\xi+\lambda\eta\tdual,
\epsilon+\lambda^2\epsilon=2\epsilon).\]
As $2\epsilon$ is an isomorphism (here we are using the fact that the
characteristic of $\K$ is not $2$), this morphism is actually an isomorphism by
\ref{unisym}. This implies that $\alpha:=(\xi+\lambda\eta\tdual)\comp\beta$ is
injective and satisfies $\cok\alpha\iso\cok\beta$, so that there is an exact
sequence
\[0\to\s{A}\tdual\mor{\alpha}\s{A}\to\gph_*\so_{\gS}\to0.\]
By \ref{compminmin} it is also clear that $\alpha$ is minimal. Moreover (using
the above commutative diagrams) we have
\[\alpha\tdual=[(\xi+\lambda\eta\tdual)\comp\beta]\tdual=\beta\tdual\comp
(\xi+\lambda\eta\tdual)\tdual=\beta\tdual\comp\xi\tdual+
\lambda\beta\tdual\comp\eta=\eta\tdual\comp\beta+\lambda\xi\comp\beta=
\lambda\alpha.\]
By \ref{det=Y} it cannot happen that $\lambda=-1$ (otherwise $\det(\alpha)$
would be a square), i.e. $\alpha\tdual=\alpha$.
        \end{proof}

The following result will be used only in chapter $4$.

    \begin{coro}\label{resdet}
Let $\gph:\gS\to\gP$ be a {\gbwcp} and assume that $\gamma:\s{A}:=
\so\oplus\s{E}(\gph)\to\gph_*\so_{\gS}$ is a morphism such that
$H^1(\gP,\gamma(1))$ and $H^0(\gP,\gamma(i))$ $\all i\in\Z$ are surjective
(this implies that $p_g(\S)=\pd_1$). Then
\begin{enumerate}

\item $\gamma$ is surjective and $\ker\gamma\iso\s{A}\dual(-1-\sw{\w})$;

\item if $\beta:\s{A}\dual(-1-\sw{\w})\to\s{A}$ is a morphism such that $\gamma
\comp\beta=0$ and $H^0(\gP,\beta(i))$ is injective $\all i\in\Z$, then the
sequence $0\to\s{A}\dual(-1-\sw{\w})\mor{\beta}\s{A}
\mor{\gamma}\gph_*\so_{\gS}\to0$ is exact.

\item if $0\to\s{A}\dual(-1-\sw{\w})\mor{\beta}\s{A}\mor{\gamma}\gph_*\so_{\gS}
\to0$ is an exact sequence, then there exists an isomorphism $\delta:
\s{A}\dual(-1-\sw{\w})\to\s{A}\dual(-1-\sw{\w})$ such that $\alpha:=\beta\comp
\delta$ is symmetric ($\alpha=\alpha\dual(-1-\sw{\w})$) and the sequence
\[0\to\s{A}\dual(-1-\sw{\w})\mor{\alpha}\s{A}\mor{\gamma}\gph_*\so_{\gS}\to0\]
is also exact.
\end{enumerate}
    \end{coro}

    \begin{proof}
$\gamma$ is surjective because $\gamma=\gsec(\gamma)\gsh{}$ by \ref{gmogqco}
and $\gsec(\gamma)$ is surjective by hypothesis. In a similar way we can prove
$2$ (the hypothesis implies, by dimension reasons, that
\[0\to H^0(\gP,\s{A}\dual(i-1-\sw{\w}))\mor{H^0(\gP,\beta(i))}
H^0(\gP,\s{A}(i))\mor{H^0(\gP,\gamma(i))}H^0(\gP,\gph_*\so_{\gS}(i))\to0\]
is exact $\all i\in\Z$). Setting $\s{J}:=\ker\gamma$, we claim that
$\Ext^1_{\gP}(\s{A},\s{J})=0$. Indeed, $\all j\in\Z$ there is an exact sequence
\[\Hom_{\gP}(\so(j),\s{A})\mor{\Hom_{\gP}(\so(j),\gamma)}
\Hom_{\gP}(\so(j),\gph_*\so_{\gS})\to
\Ext^1_{\gP}(\so(j),\s{J})\to\Ext^1_{\gP}(\so(j),\s{A}).\] Since
$\Hom_{\gP}(\so(j),\gamma)=H^0(\gP,\gamma(-j))$ is surjective and
$\Ext^1_{\gP}(\so(j),\s{A})=0$ for $j\ne-1$, we see that
$\Ext^1_{\gP}(\so(j),\s{J})=0$ for $j\ne-1$, whence it remains to show that
$\Ext^1_{\gP}(\sd^1(-1),\s{J})=0$ (recall that $\s{A}$ does not contain
$\sd^2$, as $p_g(\S)=\pd_1$). Applying $\Hom_{\gP}(-,\s{J})$ to $\esd{1}(-1)$
we obtain the exact sequence
\[\Ext^1_{\gP}(\sko^{-1}(-1),\s{J})\to\Ext^1_{\gP}(\sd^1(-1),\s{J})
\to\Ext^2_{\gP}(\so(-1),\s{J})\iso H^2(\gP,\s{J}(1)).\]
By what we have already seen $\Ext^1_{\gP}(\sko^{-1}(-1),\s{J})=0$ and the
exact sequence
\[H^1(\gP,\s{A}(1))\mor{H^1(\gP,\gamma(1))}
H^1(\gP,\gph_*\so_{\gS}(1))\to H^2(\gP,\s{J}(1))\to
H^2(\gP,\s{A}(1))\]
shows that $H^2(\gP,\s{J}(1))=0$, since $H^2(\gP,\s{A}(1))=0$ and
$H^1(\gP,\gamma(1))$ is surjective by hypothesis. Thus
$\Ext^1_{\gP}(\sd^1(-1),\s{J})=0$ and then also $\Ext^1_{\gP}(\s{A},\s{J})=0$.

By \ref{symmres} we know that there is an exact sequence of the form
\[0\to\s{A}\dual(-1-\sw{\w})\mor{\tilde{\alpha}}\s{A}\mor{\tilde{\gamma}}
\gph_*\so_{\gS}\to0\]
with $\tilde{\alpha}$ symmetric (and minimal). Denoting by $\cp{\s{U}}$ the
complex $0\to\s{A}\dual(-1-\sw{\w})\mor{\tilde{\alpha}}\s{A}\to0$ and by
$\cp{\s{V}}$ the complex $0\to\s{J}\to\s{A}\to0$, the same argument used in
\ref{unisym} shows that $\cp{\s{U}}$ and $\cp{\s{V}}$ are isomorphic in
$K^b(\gcoh{\gP})$. Then, in particular, there are morphisms of complexes
$\cp{g}:\cp{\s{U}}\to\cp{\s{V}}$ and $\cp{h}:\cp{\s{V}}\to\cp{\s{U}}$ such that
$\cp{h}\comp\cp{g}$ is homotopic to $\id_{\cp{\s{U}}}$. By \ref{minuni} this
implies that $\cp{h}\comp\cp{g}$ is an isomorphism in $C^b(\gcoh{\gP})$, whence
each $h^i$ is surjective. Being a surjective endomorphism of a vector bundle,
$h^0:\s{A}\to\s{A}$ is actually an isomorphism, and then it is clear that also
$h^{-1}:\s{J}\to\s{A}\dual(-1-\sw{\w})$ is an isomorphism. This proves $1$.

As for $3$, if $0\to\s{A}\dual(-1-\sw{\w})\mor{\beta}\s{A}\mor{\gamma}
\gph_*\so_{\gS}\to0$ is an exact sequence, then by \ref{unisym} there is an
isomorphism of complexes
\[\begin{CD}
0 @>>> \s{A}\dual(-1-\sw{\w}) @>{\tilde{\alpha}}>> \s{A} @>{\tilde{\gamma}}>>
\gph_*\so_{\gS} @>>> 0 \\
@. @VV{b}V @VV{a}V @VV{\id}V @. \\
0 @>>> \s{A}\dual(-1-\sw{\w}) @>{\beta}>> \s{A} @>{\gamma}>> \gph_*\so_{\gS}
@>>> 0.
\end{CD}\]
It follows that the sequence $0\to\s{A}\dual(-1-\sw{\w})
\mor{\alpha:=a\comp\tilde{\alpha}\comp a\dual(-1-\sw{\w})}\s{A}\mor{\gamma}
\gph_*\so_{\gS}\to0$ is exact and $\alpha\dual(-1-\sw{\w})=\alpha$ because
$\tilde{\alpha}\dual(-1-\sw{\w})=\tilde{\alpha}$. Moreover, $\alpha=\beta\comp
b\comp a\dual(-1-\sw{\w})$ and $b\comp a\dual(-1-\sw{\w})$ is an isomorphism.
    \end{proof}

                \section{Rank condition}

We are going to prove a version of the result ``rank condition = ring
condition'' for a morphism of vector bundles on a graded scheme. As in
\cite{CS} (see also \cite{JS}), we prove it under weaker hypotheses than we
will need: in particular, the morphism is not assumed to be symmetric.

Let $Z$ be a graded scheme. If $\s{E}$ is a vector bundle of rank $r$ on $Z$,
the line bundle $\Lambda^r(\s{E})$ will be denoted by $\det(\s{E})$. As usual,
for $0\le k\le r$ the natural bilinear map $\Lambda^k(\s{E})\times
\Lambda^{r-k}(\s{E})\to\det(\s{E})$ is a perfect pairing.

Given a morphism of vector bundles $\alpha:\s{G}\to\s{F}$ on $Z$, $\all k\in\N$
we will denote by $\s{I}_k(\alpha)$\index{I_k(alpha)@$\s{I}_k(\alpha)$} the
graded ideal sheaf of $k\times k$ minors of $\alpha$. More precisely,
$\s{I}_k(\alpha)$ is defined as the image of the natural map
$\Lambda^k(\s{G})\otimes\Lambda^k(\s{F})\dual\to\so_Z$ induced by
$\Lambda^k(\alpha)$. $\s{I}_k(\alpha)$ is also called the $(\rk{\s{F}}-k)^{\rm
th}$ Fitting ideal of $\cok\alpha$.

In the following we will assume that $\rk{\s{G}}=\rk{\s{F}}:=r+1$. The morphism
\[\det(\alpha):=\Lambda^{r+1}(\alpha):\det(\s{G})\to\det(\s{F})\]
will be identified with the corresponding map $\s{L}:=\det(\s{G})\otimes
\det(\s{F})\dual\to\so_Z$ (whose image is $\s{I}_{r+1}(\alpha)$) and also with
the corresponding element of $H^0(Z,\s{L}\dual)$. The morphism
\[\Lambda^r(\alpha):\Lambda^r(\s{G})\iso\s{G}\dual\otimes\det(\s{G})\to
\Lambda^r(\s{F})\iso\s{F}\dual\otimes\det(\s{F})\]
corresponds (as
$\Hom_Z(\s{G}\dual\otimes\det(\s{G}),\s{F}\dual\otimes\det(\s{F}))\iso
\Hom_Z(\s{F}\otimes\s{L},\s{G})$) to a map $\beta:\s{F}\otimes\s{L}\to\s{G}$,
which satisfies the identities
\begin{align*}
& \alpha\comp\beta=\id_{\s{F}}\otimes\det(\alpha):\s{F}\otimes\s{L}\to\s{F}, &
\beta\comp(\alpha\otimes\s{L})=\id_{\s{G}}\otimes\det(\alpha):
\s{G}\otimes\s{L}\to\s{G}.
\end{align*}

        \begin{lemm}\label{Qref}
In the notation introduced so far, let $Y\subseteq Z$ be the closed graded
subscheme defined by $\s{I}_{r+1}(\alpha)$ (i.e., $\so_Y:=
\so_Z/\s{I}_{r+1}(\alpha)$). Let moreover $\s{Q}:=\cok\alpha$, which is in a
natural way a graded $\so_Y$--module. If $\det(\alpha)$ is a {\nzdiv} (meaning
that $\det(\alpha)_z$ is a {\nzdiv} $\all z\in Z$), then $\s{C}:=
\gsHom_Y(\s{Q},\so_Y)$ satisfies $\s{Q}\iso\gsHom_Y(\s{C},\so_Y)$ (as graded
$\so_Y$--modules).
        \end{lemm}

    \begin{proof}
Applying the functor $-_Y:=-\otimes_{\so_Z}\so_Y$ to the exact sequence $\s{G}
\mor{\alpha}\s{F}\to\s{Q}\to0$, and taking into account that
$\s{I}_{r+1}(\alpha)\s{Q}=0$ ($\alpha\comp\beta=\id_{\s{F}}\otimes\det(\alpha)$
implies that $\s{I}_{r+1}(\alpha)\s{F}\subseteq\im\alpha$), we obtain the
exact sequence in $\gcoh{Y}$
\[\s{G}_Y\mor{\alpha_Y}\s{F}_Y\to\s{Q}_Y\iso\s{Q}\to0,\]
and also, applying $\gsHom_Y(-,\so_Y)$, the exact sequence
\[0\to\s{C}\to\s{F}_Y\dual\mor{\alpha_Y\dual}\s{G}_Y\dual.\]
We claim moreover that there is a long exact sequence in $\gcoh{Y}$
\[\cdots\to\s{F}_Y\otimes\s{L}_Y^{n+1}\mor{\beta_Y\otimes\s{L}_Y^n}
\s{G}_Y\otimes\s{L}_Y^n\mor{\alpha_Y\otimes\s{L}_Y^n}\s{F}_Y\otimes\s{L}_Y^n
\mor{\beta_Y\otimes\s{L}_Y^{n-1}}\s{G}_Y\otimes\s{L}_Y^{n-1}\to\cdots\] Indeed,
$(\alpha_Y\otimes\s{L}_Y^n)\comp(\beta_Y\otimes\s{L}_Y^n)=
(\id_{\s{F}\otimes\s{L}^n}\otimes\det(\alpha))_Y=0$ because $\det(\alpha)_Y=0$,
and if (for some $y\in Y$) $\sigma\in(\s{G}_Y\otimes\s{L}_Y^n)_y$ is such that
$(\alpha_Y\otimes\s{L}_Y^n)_y(\sigma)=0$, then (choosing
$\tilde{\sigma}\in(\s{G}\otimes\s{L}^n)_y$ lift of $\sigma$) $\exi\tau\in
(\s{F}\otimes\s{L}^{n+1})_y$ such that$(\alpha\otimes\s{L}^n)_y(\tilde{\sigma})
=(\id_{\s{F}\otimes\s{L}^n}\otimes\det(\alpha))_y(\tau)$, which implies (since
$\alpha\otimes\s{L}^n$ is injective, because $(\alpha\otimes\s{L}^n)\comp
(\beta\otimes\s{L}^n)=\id_{\s{F}\otimes\s{L}^n}\otimes\det(\alpha)$ is
injective by hypothesis) that $\tilde{\sigma}=(\beta\otimes\s{L}^n)_y(\tau)$,
whence $\sigma\in\im(\beta_Y\otimes\s{L}_Y^n)_y$. So the sequence is exact at
$\s{G}_Y\otimes\s{L}_Y^n$, and in a completely similar way one can check that
it is exact at $\s{F}_Y\otimes\s{L}_Y^n$. Besides, observing that
$\s{I}_{r+1}(\alpha)=\s{I}_{r+1}(\alpha\dual)$, from $\alpha\dual$ we obtain in
the same way the long exact sequence
\[\cdots\to\s{G}_Y\dual\otimes\s{L}_Y^{1-n}
\mor{\beta_Y\dual\otimes\s{L}_Y^{1-n}}\s{F}_Y\dual\otimes\s{L}_Y^{-n}
\mor{\alpha_Y\dual\otimes\s{L}_Y^{-n}}\s{G}_Y\dual\otimes\s{L}_Y^{-n}
\mor{\beta_Y\dual\otimes\s{L}_Y^{-n}}\s{F}_Y\dual\otimes\s{L}_Y^{-n-1}\to\cdots
\]
In particular, it follows that $\s{C}\iso\ker\alpha_Y\dual\iso
\im(\beta_Y\dual\otimes\s{L}_Y)$, whence there is an exact sequence
\[\s{F}_Y\dual\otimes\s{L}_Y\mor{\alpha_Y\dual\otimes\s{L}_Y}
\s{G}_Y\dual\otimes\s{L}_Y\to\s{C}\to0,\]
and applying $\gsHom_Y(-,\so_Y)$ to it we obtain the exact sequence
\[0\to\gsHom_Y(\s{C},\so_Y)\to\s{G}_Y\otimes\s{L}_Y^{-1}
\mor{\alpha_Y\otimes\s{L}_Y^{-1}}\s{F}_Y\otimes\s{L}_Y^{-1}.\]
Therefore $\gsHom_Y(\s{C},\so_Y)\iso\ker(\alpha_Y\otimes\s{L}_Y^{-1})\iso
\im(\beta_Y\otimes\s{L}_Y^{-1})\iso\cok\alpha_Y\iso\s{Q}$.
    \end{proof}

        \begin{prop}\label{RC}
In the notation of \ref{Qref}, assume moreover that $\s{F}=\so_Z\oplus\s{E}$
(for some vector bundle $\s{E}$ of rank $r$) and that
\[\alpha=\begin{pmatrix}
\alpha^{(1)} \\
\alpha'
\end{pmatrix}:\s{G}\to\s{F}=\so_Z\oplus\s{E}.\]
If $\det(\alpha)$ is a {\nzdiv} and $\depth(\s{I}_r(\alpha'))\ge2$ (i.e.,
$\depth(\s{I}_r(\alpha')_z)\ge2$ $\all z\in Z$), then the following conditions
are equivalent:
\begin{enumerate}

\item $\s{Q}=\cok\alpha$ carries the structure of a graded sheaf of commutative
$\so_Y$--algebras, with unit given by the image of $1\in\gsec(Z,\so_Z)\subset
\gsec(Z,\s{F})$ in $\s{Q}$;

\item $\s{I}_r(\alpha')=\s{I}_r(\alpha)$.
\end{enumerate}
        \end{prop}

    \begin{proof}
As we can always reduce to make local computations, we will adopt the following
notation at a point $z\in Z$: $R:=\so_{Z,z}$, $F:=\s{F}_z$ and $G:=\s{G}_z$
(since $\s{F}$ and $\s{G}$ are locally free, we have $F\iso
\bigoplus_{1\le i\le r+1}R(m_i)$ and $G\iso\bigoplus_{1\le i\le r+1}R(n_i)$ for
some $m_i,n_i\in\Z$ with $m_1=0$, so that $\s{L}_z\iso R(n-m)$, where $m:=
\sum m_i$ and $n:=\sum n_i$). Then $\alpha_z$ and $\beta_z$ can be identified
with maps $a:G\to F$ and $b:F(n-m)\to G$, which can also be considered as
$(r+1)\times(r+1)$ matrices (with entries $a_{i,j}\in R_{m_i-n_j}$ and $b_{i,j}
\in R_{m-n+n_i-m_j}$). For simplicity $\all l\in\Z$ we will write $a$ and $b$
instead of $a(l)$ and $b(l)$. Thus, setting $d:=\det(a)\in R_{m-n}$, we have
$a\comp b=d\id_F$ and $b\comp a=d\id_G$. $\s{I}_r(\alpha')_z$ and
$\s{I}_r(\alpha)_z$ can be identified with the ideals of $R$
\[I':=(b_{i,1}\st i=1,\dots,r+1)\subseteq I:=(b_{i,j}\st i,j=1,\dots,r+1).\]
As $d\in I'$ is a {\nzdiv} and $\depth(I')\ge2$, there exists $d'\in I'$ such
that $(d,d')$ is a regular sequence. We set moreover $\bar{R}:=R/(d)$ and we
will denote more generally by $\bar{-}$ the functor $-\otimes_R\bar{R}$.
Clearly $Q:=\s{Q}_z=\cok a$ and $C:=\s{C}_z=\gHom_{\bar{R}}(Q,\bar{R})$ satisfy
$Q\iso\bar{Q}$ and $C\iso\bar{C}$. Finally, $\iota:R=R(m_1)\mono F$ and
$p:F\epi Q$ will denote, respectively, the natural inclusion and projection.

First we prove that the natural map $\so_Z\to\s{Q}$ induces an injection $\so_Y
\mono\s{Q}$, which is the same as proving that $\ker(p\comp\iota)=(d)$. Now,
\[p\comp\iota(d)=p(d(1,0,\dots,0))=p(a\comp b(1,0,\dots,0))=0\]
because $p\comp a=0$. On the other hand, if $x\in R$ is such that $p\comp
\iota(x)=0$, then $\exi g=(g_1,\dots,g_{r+1})\in G$ such that $\iota(x)=a(g)$.
Therefore we have
\[dg=b\comp a(g)=b(\iota(x))=x(b_{1,1},\dots,b_{r+1,1}),\]
which shows that $xI'\subseteq(d)$. Then, in particular, $xd'\in(d)$, which
implies that $x\in(d)$, because $(d,d')$ is a regular sequence.

Next we show that the quotient $\s{Q}/\so_Y$ is annihilated by
$\s{I}_r(\alpha')$, i.e that $I'Q/\bar{R}=0$. As $Q/\bar{R}= F/(\im
a+\im\iota)$, we have to prove that $b_{i,1}F\subseteq\im a+\im\iota$ for
$i=1,\dots,r+1$. Let $a_{(i)}$ be the matrix equal to $a$, except that the
first row is substituted by $(0,\dots,1,\dots,0)$ (with $1$ in position $i$):
as $\det(a_{(i)})=\pm b_{i,1}$, $b_{(i)}$ (the matrix of cofactors of
$a_{(i)}$) satisfies $a_{(i)}\comp b_{(i)}=\pm b_{i,1}\id$. Thus $\all f\in F$
we have
\[\pm b_{i,1}f=a_{(i)}\comp b_{(i)}(f)=a\comp b_{(i)}(f)+(a-a_{(i)})\comp
b_{(i)}(f).\]
Since $\im(a-a_{(i)})\subseteq\im\iota$, we see that $b_{i,1}f\in\im a+
\im\iota$.

Now we can prove that $\gsHom_Y(\s{Q}/\so_Y,\so_Y)=0$, i.e. that
$\gHom_{\bar{R}}(Q/\bar{R},\bar{R})=0$. Indeed, an element of
$\gHom_{\bar{R}}(Q/\bar{R},\bar{R})$ can be identified with an element $\varphi
\in\gHom_{\bar{R}}(Q,\bar{R})$ such that $\varphi\rest{\bar{R}}=0$. Then
(remembering that $I'Q\subseteq\bar{R}$) $\all q\in Q$ we have $0=\varphi(d'q)=
d'\varphi(q)$, from which it follows that $\varphi(q)=0\in\bar{R}$ (again,
because $(d,d')$ is a regular sequence).

Thus the natural map $\s{C}=\gsHom_Y(\s{Q},\so_Y)\to\gsHom_Y(\so_Y,\so_Y)\iso
\so_Y$ is injective, so that $\s{C}$ can be identified with its image, which we
claim to be the graded sheaf of ideals $\s{I}_r(\alpha')/\s{I}_{r+1}(\alpha)$
of $\so_Y=\so_Z/\s{I}_{r+1}(\alpha)$. To see this, it is enough to check that
if $c\in C=\gHom_{\bar{R}}(Q,\bar{R})\subseteq\bar{R}$, then $c(1)\in\bar{I'}$,
and that every element of $I'$ is of the form $c(1)$ for some $c\in C$. As $C$
is naturally isomorphic to $\im(\bar{b}\dual:\bar{G}(m-n)\dual\to
\bar{F}\dual)$, every $c\in C$ corresponds to an element of the form
$\bar{b}\dual(\varphi)=\varphi\comp\bar{b}$ for some $\varphi\in
\gHom_{\bar{R}}(\bar{G}(m-n),\bar{R})$ ($\varphi\comp\bar{b}:\bar{F}\to\bar{R}$
is null on $\im\bar{a}$, hence it induces $c:Q=\bar{F}/\im\bar{a}\to
\bar{R}$). Then we have by definition
\[c(1)=\varphi\comp\bar{b}(1,0,\dots,0)=
\varphi(\bar{b}_{1,1},\dots,\bar{b}_{r+1,1})=\sum_{i=1}^{r+1}\bar{b}_{i,1}x_i\]
(where $x_i:=\varphi(0,\dots,1,\dots,0)$, with $1$ in position $i$). Since the
$x_i\in\bar{R}$ can be arbitrary, the claim follows.

We are now ready to prove that the condition $\s{I}_r(\alpha')=\s{I}_r(\alpha)$
is satisfied if and only if the natural inclusion $\gsHom_Y(\s{C},\s{C})
\subseteq\gsHom_Y(\s{C},\so_Y)$ is an equality, i.e. that $I'=I$ if and only if
$\gHom_{\bar{R}}(C,C)=\gHom_{\bar{R}}(C,\bar{R})$. Given $q\in Q\iso
\gHom_{\bar{R}}(C,\bar{R})$ and $c\in C$, by \ref{Qref} we have $q(c)=c(q)$.
Since $Q\iso\im(\bar{b}:\bar{F}\to\bar{G}(m-n))$ and $C\iso\im\bar{b}\dual$, we
can find $f\in\bar{F}$ and $\varphi\in\gHom_{\bar{R}}(\bar{G}(m-n),\bar{R})$
such that $q$ and $c$ correspond to $\bar{b}(f)$ and $\bar{b}\dual(\varphi)=
\varphi\comp\bar{b}$. This implies that $q(c)=c(q)=\varphi\comp\bar{b}(f)$.
Therefore the condition $\gHom_{\bar{R}}(C,C)=\gHom_{\bar{R}}(C,\bar{R})$ is
equivalent to $\varphi\comp\bar{b}(f)\in C=\bar{I'}\subseteq\bar{R}$ $\all f\in
\bar{F}$ and $\all\varphi\in\bar{G}(m-n)\dual$, which is clearly satisfied if
and only if every entry of $\bar{b}$ is in $\bar{I'}$, i.e. if and only if $I'=
I$.
\begin{description}

\item[$1\implies2$] By what we have just seen, we have to prove that
$\gsHom_Y(\s{C},\s{C})=\gsHom_Y(\s{C},\so_Y)$, i.e. that given $q\in Q\iso
\gHom_{\bar{R}}(C,\bar{R})$ and $c\in C$, $q(c)=cq\in C$. Now, in any case $cq
\in\bar{R}$, and (by definition of $C$) $cq\in C$ if and only if $(cq)q'\in
\bar{R}$ $\all q'\in Q$. As $Q$ is an $\bar{R}$--algebra, we have $(cq)q'=cqq'=
(qq')(c)\in\bar{R}$.

\item[$2\implies1$] As $\s{Q}\iso\gsHom_Y(\s{C},\so_Y)=\gsHom_Y(\s{C},\s{C})$,
$\s{Q}$ has a natural structure of $\so_Y$--algebra, with multiplication given
by composition in $\gsHom_Y(\s{C},\s{C})$. Given $q,q'\in Q$, we have
\[{\bar{d'}}^2qq'=\bar{d'}q\bar{d'}q'=\bar{d'}q'\bar{d'}q={\bar{d'}}^2q'q\]
(since $\bar{R}$ is obviously in the centre of $Q$, and $\bar{d'}q,\bar{d'}q'
\in\bar{R}$ because $\bar{d'}\in\bar{I'}=C$). As $\bar{d'}$ is a {\nzdiv} in
$Q$ (if $p\in Q\iso\gHom_{\bar{R}}(C,C)$ is such that $\bar{d'}p=0$, then $\all
c\in C$ we have $0=(\bar{d'}p)(c)=\bar{d'}p(c)$, whence $p(c)=0$ because
$\bar{d'}$ is a {\nzdiv} in $\bar{R}$), this implies that $qq'=q'q$, i.e.
multiplication is commutative. The statement about unit is immediate to check.
\end{description}
    \end{proof}

    \begin{rema}\label{RCrem}
It follows from the proof that $1\implies2$ is true also without the assumption
$\depth(\s{I}_r(\alpha'))\ge2$, provided one knows that $\so_Y\mono\s{Q}$ and
$\gsHom_Y(\s{Q}/\so_Y,\so_Y)=0$.
    \end{rema}

    \begin{rema}\label{birpts}
If the hypotheses of \ref{RC} are satisfied, it is clear from the proof that
the natural map $\so_{\gY}\mono\s{Q}$ is an isomorphism precisely at those
points $y\in\gY$ where $\s{I}_r(\alpha')_y=\so_{Z,y}$.
    \end{rema}

                \section{Main theorem}

        \begin{theo}\label{mthm}
Let $\S$ be a minimal surface of general type and $\gS$ the {\std} graded
scheme $\stdsch{\S}{\ds_{\S}}$ (see  \ref{stdchar}), let
$\gph:\gS\to\gY\subset\gP=\gP(\w)$ (with $\dim\P=3$) be a \gbwcp,
and let $\s{E}(\gph)$ be the vector bundle on $\gP$ defined in
\ref{Edef}. Then there exists an exact sequence on $\gP$
\[0\to(\so\oplus\s{E}(\gph))\dual(-1-\sw{\w})\mor{\alpha=\begin{pmatrix}
\alpha^{(1)} \\
\alpha'
\end{pmatrix}}
\so\oplus\s{E}(\gph)\to\gph_*\so_{\gS}\to0,\]
such that $\alpha$ is a minimal morphism which satisfies the following
properties:
\begin{enumerate}
\item $\alpha$ is symmetric, i.e. $\alpha=\alpha\dual(-1-\sw{\w})$;
\item $\f=\det(\alpha)$ is an irreducible polynomial (defining $\gY=
\gproj\p/(\f)$);
\item $\s{I}_r(\alpha)=\s{I}_r(\alpha')$ (where $r:=\rk{\s{E}(\gph)}$);
\item $\gX:=\gSpec(\cok{\alpha})$ is a {\std} graded scheme and $\X$ is
a surface with only rational double points as singularities.
\end{enumerate}

Conversely, given a vector bundle $\s{E}$ on $\gP$ with $\s{E}(2)\in
\com{]3-\sw{\w},0[}$ and a minimal morphism
$\alpha:(\so\oplus\s{E})\dual(-1-\sw{\w})\to\so\oplus \s{E}$ satisfying
properties $1$, $2$, $3$, then $\cok\alpha$ is in a natural way a graded
coherent sheaf of commutative $\so_{\gY}$--algebras, and if $4$ also holds,
then, denoting by $\pr:\S\to\X$ a minimal resolution of singularities of $\X$,
$\S$ is a minimal surface of general type, $\X$ is its canonical model, $\pr$
extends to a morphism of {\std} graded schemes
$\gpr:\gS:=\stdsch{\S}{\ds_{\S}}\to\gX$, $\gph:\gS\mor{\gpr} \gX\to\gY$ is a
{\gbwcp} and $\s{E}=\s{E}(\gph)$.
        \end{theo}

        \begin{proof}
The first implication is now very easy to prove: by \ref{symmres} there exists
$\alpha$ minimal satisfying $1$; moreover, $2$ follows from \ref{det=Y}, $3$
from \ref{RC} and \ref{RCrem} (clearly $\so_{\gY}\mono\gph_*\so_{\gS}$ and
$\gsHom_{\gY}((\gph_*\so_{\gS})/\so_{\gY},\so_{\gY})=0$, because the support of
$(\gph_*\so_{\gS})/\so_{\gY}$ is a proper closed subset of $\gY$, as $\gph$ is
birational), and $4$ from \ref{gSpecA} and \ref{canmodstd}.

Conversely, denoting $\cok\alpha$ by $\s{Q}$, $\dual(-1-\sw{\w})$ by $\tdual$
and $\gsExt^1_{\gP}(-,\so(-1-\sw{\w}))$ by $E$, since $\alpha$ is injective
(because $\det(\alpha)\ne0$), there is an exact sequence
\[0\to(\so\oplus\s{E})\tdual\mor{\alpha}\so\oplus\s{E}\to\s{Q}\to0.\]
Applying $\gsHom_{\gP}(-,\so(-1-\sw{\w}))$ we obtain the exact sequence
\[0\to(\so\oplus\s{E})\tdual\mor{\alpha\tdual}\so\oplus\s{E}\to E(\s{Q})\to0\]
($\gsHom_{\gP}(\s{Q},\so(-1-\sw{\w}))=0$ because $\s{Q}\in\gsmo{\gY}$). As
$\alpha\tdual=\alpha$, there is then a natural isomorphism $E(\s{Q})\iso\s{Q}$.
Notice that $\depth(\s{I}_r(\alpha'))\ge2$: clearly it is enough to show that
if $y\in\gY$, then $\depth(\s{I}_r(\alpha')_y)\ge2$. As $\f\in
\s{I}_r(\alpha')_y$ is irreducible, this is true if and only if $(\f)=
\s{I}_{r+1}(\alpha)_y\subsetneq\s{I}_r(\alpha')_y=\s{I}_r(\alpha)_y$. If on the
contrary $(\f)=\s{I}_r(\alpha)_y$, then $\f$ would divide every entry of
$\beta_y$ (the matrix of cofactors of $\alpha_y$), whence $\f^{r+1}$ would
divide $\det(\beta_y)$, but $\det(\beta_y)=\f^r$ (this follows from the fact
that $\alpha_y\beta_y=\beta_y\alpha_y=\f\id$). Therefore $\s{Q}$ is a graded
coherent sheaf of commutative $\so_{\gY}$--algebras by \ref{RC}, and, denoting
by $\gps:\gX:=\gSpec\s{Q}\to\gY\subset\gP$ the natural map, by \ref{gSpecA}
$\s{Q}\iso\gps_*\so_{\gX}$ and $\gps$ is finite. Moreover, $\gps$ is birational
onto $\gY$ because there exists an open subset $\emptyset\ne U\subseteq\gY$
such that $(\s{I}_r(\alpha)/\s{I}_{r+1}(\alpha))\rest{U}=\so_U$ (this follows
from the fact that $(\f)=\s{I}_{r+1}(\alpha)\subsetneq\s{I}_r(\alpha)$ and $\f$
is irreducible), whence $\s{Q}\rest{U}=\so_U$ by \ref{birpts}. By
\ref{gsExtdual} there are natural isomorphisms in
$\gcoh{\gP}\cap\gsmo{\gps_*\so_{\gX}}$
\[\gps_*\ds_{\gX}\iso\gsExt^1_{\gP}(\s{Q},\so(-\sw{\w}))\iso E(\s{Q})(1)\iso
\s{Q}(1).\] Since also
$\gps_*(\gps^*\so_{\gP}(1))\iso\so_{\gP}(1)\otimes\gps_*\so_{\gX} \iso\s{Q}(1)$
by projection formula, it follows from \ref{gSpecA} that
$\ds_{\gX}\iso\gps^*\so_{\gP}(1)\iso\so_{\gX}(1)$. Therefore $\ds_{\X}\iso
(\so_{\gX})_1$ is an invertible sheaf (because $\gX$ is \std). As $\X$ has only
rational double points as singularities, everything will follow if we prove
that $\ds_{\X}$ is ample: indeed, in this case $\X$ is the canonical model of a
minimal surface of general type $\S$, the map $\pr:\S\to\X$ extends naturally
to a morphism $\gpr:\gS\to\gX$ (where $\gS:=\stdsch{\S}{\ds_{\S}}$), $\gph:=
\gps\comp\gpr$ is a {\gbwcp} and $\s{E}=\s{E}(\gph)$ by the part already proved
and by the uniqueness of the minimal resolution (see \ref{unisym}).

So it remains to prove that $\ds_{\X}$ is ample, and it is clearly enough to
show that $\ds_{\X}^m$ is ample, where $m$ is the least common multiple of
$\{\w_0,\w_1,\w_2,\w_3\}$. By definition this is true if and only if, given
$\s{G}\in\coh{\X}$, $\exi l\in\N$ such that the natural map
\[g_n:\so_{\X}\otimes_{\K}H^0(\X,\s{G}\otimes\ds_{\X}^{nm})\to
\s{G}\otimes\ds_{\X}^{nm}\]
is surjective for $n\ge l$. As $\so_{\P}(m)$ is an ample invertible sheaf (by
\cite[prop. 2.3]{D}), $\exi l\in\N$ such that the natural map
\[h_n:\so_{\P}\otimes_{\K}H^0(\P,{\ps}_*\s{G}\otimes\so_{\P}(nm))\to
{\ps}_*\s{G}\otimes\so_{\P}(nm)\]
is surjective for $n\ge l$. Since $\ps^*\so_{\P}(nm)\iso\ds_{\X}^{nm}$ (this is
an easy consequence of the fact that every stalk of $\so_{\gP}$ contains an
invertible element of degree $m$), by projection formula
${\ps}_*(\s{G}\otimes\ds_{\X}^{nm})\iso{\ps}_*\s{G}\otimes\so_{\P}(nm)$.
Moreover, there is a natural isomorphism $H^0(\X,\s{G}\otimes\ds_{\X}^{nm})\iso
H^0(\P,{\ps}_*(\s{G}\otimes\ds_{\X}^{nm}))$ (because $\ps$ is finite), and then
it is clear that $h_n$ factors through ${\ps}_*g_n$, which is therefore
surjective if $n\ge l$. By \ref{gSpecA}, then, also $g_n$ is surjective for
$n\ge l$.
        \end{proof}

    \begin{rema}\label{coeffdiff}
As we will see, the coefficients $\cre_j(\gph)$ of the sheaves $\so(j-2)$ (for
$3-\sw{\w}<j<0$) in $\s{E}(\gph)$ do depend on $\gph$. On the other hand, the
difference $\cre_j(\gph)-\cre_{3-\sw{\w}-j}(\gph)$ is uniquely determined by
the numerical invariants of $\S$ (and by $\w$, of course). Indeed, by
\ref{Fres}, \ref{Fbdres} and \ref{O2res} we have
\begin{multline*}
\cre_j(\gph)-\cre_{3-\sw{\w}-j}(\gph)=\cbrd_j-z^0_j-z^{-1}_{3-\sw{\w}-j}-
\cbrd_{3-\sw{\w}-j}+z^0_{3-\sw{\w}-j}+z^{-1}_j\\
=\chi(\gP,\gph_*\so_{\gS}(2)\otimes\cp{\bd{j}})-\chi(\gP,\cp{\bd{j}}(2))+
\chi(\gP,\cp{\bd{3-\sw{\w}-j}}(2)),
\end{multline*}
and $\chi(\gP,\gph_*\so_{\gS}(2)\otimes\cp{\bd{j}})=\sum_i(-1)^i
\chi(\gP,\gph_*\so_{\gS}(2)\otimes\bd{j}^i)$ depends only on $K^2_{\S}$,
$\chi(\so_{\S})$ and $\w$, whereas obviously
$\chi(\gP,\cp{\bd{j}}(2))$ and $\chi(\gP,\cp{\bd{3-\sw{\w}-j}}(2))$
depend only on $\w$.
    \end{rema}

    \begin{rema}\label{natmor}
The image in $\Hom_{\coq{}}$ of the component of $\alpha$ from $(\sd^2)^q$ to
$\sd^1(-1)^q$ can be identified with the natural map
\[\xymatrix{\sd^2\otimes_{\K}H^1(\gP,\gph_*\so_{\gS})
\ar[r] \ar[d]_{\upsilon\otimes\id} &
\sd^1(-1)\otimes_{\K}H^1(\gP,\gph_*\so_{\gS}(1)) \\
\sd^1(-1)\otimes_{\K}\p_1\otimes_{\K}H^1(\gP,\gph_*\so_{\gS})
\ar[ur]_{\id\otimes m}}\]
(this follows from \ref{beildiff}, taking into account that this map is not
changed by passing from the Beilinson complex to the present complex, see the
proof of \ref{Fres}).
    \end{rema}

            \chapter{Applications to surfaces with $p_g=q=2$, $K^2=
                4$}

In this chapter we want to determine explicitly the resolution, whose existence
is assured by \ref{mthm}, for {\gbwcp}s of a particular class of (minimal)
irregular surfaces of general type, namely those obtained as double covers of a
principally polarized abelian surface $(A,\Theta$), branched along a (smooth)
divisor of $\linsys{2\Theta}$. These surfaces (whose numerical invariants are
those indicated in the title) are interesting because their bicanonical map is
not birational and, if $\Theta$ is irreducible, they are examples of the
so--called non--standard case for the non birationality of the bicanonical map
(see \cite{Ci} and \cite{CM}). The standard case is, by definition, that of a
surface admitting a rational map to a curve whose general fibre is an
irreducible curve of geometric genus $2$ (this property easily implies that the
bicanonical map is not birational). In particular, in \cite{CM} it is proved
that a surface with $p_g=q=2$ presents the non--standard case if and only if it
is a double cover as above (with $\Theta$ irreducible, since if on the contrary
$\Theta$ is reducible, i.e. if $A$ is isomorphic to the product of two elliptic
curves, then the surface presents the standard case). As for the weights to use
for the projection (which we would like to be as small as possible), only two
of them can be equal to $1$ (since $p_g=2$) and, by what we have just said it
is clear that at least one must be $>2$, so that the most natural choice (which
we make) is $\w=(1,1,2,3)$.

We consider only the (more interesting) case of an irreducible polarization.
It turns out that even the problem of determining exactly the vector bundle
$\s{E}(\gph)$ (i.e., of finding the coefficients $\cre_j(\gph)$ of the sheaves
$\so(j-2)$) requires a deep knowledge of the canonical ring $\R$ of the
surface $\S$. So we first study it: the hard part is to determine the ring of
theta functions on $A$ $\R(A,\Theta):=\bigoplus_{n\ge0}H^0(A,\so_A(n\Theta))$
(because $\R$ is simply an integral extension of it, generated by a single
element of degree $1$). We manage to do it, finding explicit generators and
relations, in a purely algebraic way. Once we have the canonical ring, we can
compute the $\cre_j(\gph)$ using ${\rm \check{C}}$ech cohomology (and we see
that they actually depend on the choice of $\gph$). Then we choose $\gph$ as
simple as possible and briefly describe the computations (which require the use
of a computer) which allow us to determine a symmetric resolution $\alpha$ of
$\gph_*\so_{\gS}$ as in \ref{mthm}.

We will work over the field $\C$ of complex numbers.

        \section{Double covers of an abelian surface}

Let $A$\index{A@$A$} be an abelian surface and let
$\Theta$\index{Theta@$\Theta$} be a symmetric principal polarization on $A$.
Given a smooth divisor $B\in\linsys{2\Theta}$,\index{B@$B$} let $\S$ be the
surface obtained as double cover of $A$ branched along $B$, and denote by
$p:\S\to A$ the projection. Then $\S$ is a minimal surface of general type and
$p_g(\S)=q(\S)=2$, $K_{\S}^2=4$. In fact, as $\ds_{\S}\iso p^*\so_A(\Theta)$
and $p_*\so_{\S}\iso\so_A\oplus\so_A(-\Theta)$, we have
\[p_*\ds_{\S}^n\iso\so_A(n\Theta)\oplus\so_A((n-1)\Theta)\]
$\all n\in\Z$ ($\so_A(n\Theta)$ is the invariant part and $\so_A((n-1)\Theta)$
the antiinvariant if $n$ is even, and viceversa if $n$ is odd). $p$ being
finite it follows that
\[H^i(\S,\ds_{\S}^n)\iso H^i(A,p_*\ds_{\S}^n)\iso H^i(A,\so_A(n\Theta))\oplus
H^i(A,\so_A((n-1)\Theta))\]
$\all i\in\N$. Therefore, since $h^0(A,\so_A(n\Theta))=n^2$ for $n>0$,
$h^i(A,\so_A(n\Theta))=0$ for $n\ne0$ and $h^1(A,\so_A)=2$, we see that
$p_g(\S)=h^0(A,\so_A(\Theta))+h^0(A,\so_A)=2$, $q(\S)=h^1(A,\so_A)+
h^1(A,\so_A(-\Theta))=2$ (so that $\chi(\so_{\S})=1$) and (for $n>1$)
\[h^0(\S,\ds_{\S}^n)=h^0(A,\so_A(n\Theta))+h^0(A,\so_A((n-1)\Theta))=
n^2+(n-1)^2=1+2n(n-1),\] whence $\S$ is of general type and $K_{\S}^2=4$
(because $h^0(\S,\ds_{\S}^n)= \chi(\so_{\S})+\frac{n(n-1)}{2}K_{\S}^2$).
Moreover, $\S$ contains no rational curve (because the same is true for $A$),
which implies that $\S$ is minimal and coincides with its canonical model, i.e.
$\S\iso\proj\R$, where $\R:= \R(\S,\ds_{\S})$ is the canonical ring of $\S$.
Denoting by $\R'$ the (positively graded) ring\index{R'@$\R'$}
\[\R'=\R(A,\Theta):=\bigoplus_{n\ge0}H^0(A,\so_A(n\Theta))\]
and by $\xi$\index{xi@$\xi$} a generator of $H^0(A,\so_A)\subset
H^0(\S,\ds_{\S})=\R_1$, $\all n\in\N$ we can write
$\R_n=\R'_n\oplus\xi\R'_{n-1}$, so that $\R=\R'\oplus\xi\R'$. Then
multiplication in $\R$ is given by
\[(x+\xi y)(x'+\xi y')=xx'+syy'+(xy'+x'y)\xi\]
$\all x,x',y,y'\in\R'$, where $s\in H^0(A,\so_A(2\Theta))=\R'_2$\index{s@$s$}
is a section whose zero divisor is $B$.

        \section{The ring of theta functions of an abelian surface}

We will assume that the polarization $\Theta$ is irreducible, so that $A$ is
isomorphic to the Jacobian of a (smooth projective) curve $C$\index{C@$C$} of
genus $2$. Being hyperelliptic, $C$ is a double cover of $\P^1$ branched on six
(distinct) points. Denoting by $[y_0,y_1]$ the homogeneous coordinates on
$\P^1$, we can write $C=\proj R$, where $R$\index{R@$R$} is the graded ring
\begin{align*}
& R:=\C[y_0,y_1,z]/(z^2-F) & \deg(y_0)=\deg(y_1)=1,\deg(z)=3
\end{align*}
and $F$\index{F@$F$} is a homogeneous polynomial of degree $6$ vanishing
exactly at the branch points (notice that $R$ can be identified with the
canonical ring $\R(C,\ds_C)$ of $C$). Up to composing with an automorphism of
$\P^1$, we can assume that three of the branch points are $0$, $1$, $\infty$,
whence we can write $F$ as
\[F=y_0y_1 (y_1^4+\lambda y_0y_1^3+\mu y_0^2y_1^2+\nu y_0^3y_1+\ep y_0^4).\]
for some $\lambda,\mu,\nu,\ep\in\C$\index{lambda@$\lambda$}\index{mu@$\mu$}
\index{nu@$\nu$}\index{epsilon@$\ep$} with $1+\lambda+\mu+\nu+\ep=0$.

We will denote by $h:C\to\P^1$ the double cover and by $\iota:C\to C$ the
hyperelliptic involution (given, in weighted homogeneous coordinates
$[y_0,y_1,z]$, by $[a,b,c]\mapsto[a,b,-c]$). We will identify the abelian
surface $A=\Pic^0(C)$\index{A@$A$} with $\Pic^2(C)$ using the isomorphism given
by $\s{L}\mapsto\s{L}\otimes \so_C(2c_{\infty})$, where $c_{\infty}:=[0,1,0]\in
C$. We recall that $\Pic^2(C)$ is naturally isomorphic to a blow--down of
$C^{(2)}$ (the twofold symmetric product of $C$, whose closed points are in
canonical bijection with the set of effective divisors of degree $2$ on $C$):
the natural map $q:C^{(2)} \to\Pic^2(C)$ is an isomorphism everywhere, except
that it contracts to a point the exceptional divisor $E$, corresponding to the
linear system $g^1_2$ of $C$. Letting $r:C\times C\to C^{(2)}$ be the quotient
map and setting $D_1:= \{c_{\infty}\}\times C$ and
$D_2:=C\times\{c_{\infty}\}$, we can take as theta divisor $\Theta:=q\comp
r(D_1)=q\comp r(D_2)\subset A$.\index{Theta@$\Theta$}

Denoting by $\theta$\index{theta@$\theta$} a generator of
$H^0(A,\so_A(\Theta))$, the exact sequences
\[0\to\so_A((n-1)\Theta)\mor{\cdot\theta}\so_A(n\Theta)\to\so_{\Theta}(n\Theta)
\to0\]
induce exact sequences for $n>1$ (as $H^1(A,\so_A(m\Theta))=0$ for $m>0$)
\[0\to H^0(A,\so_A((n-1)\Theta))\mor{\cdot\theta}H^0(A,\so_A(n\Theta))\to
H^0(\Theta,\so_{\Theta}(n\Theta))\to0,\] whereas for $n=1$ we have
$\theta:H^0(A,\so_A)\isomor H^0(A,\so_A(\Theta))$ and
$H^0(\Theta,\so_{\Theta}(\Theta))\isomor H^1(A,\so_A)\iso\C^2$. As $\Theta \iso
C$ and (by adjunction) $\so_{\Theta}(\Theta)\iso\ds_{\Theta}$,
$H^0(\Theta,\so_{\Theta}(n\Theta))\iso H^0(C,\ds_C^n)\iso R_n$. Letting $R'$ be
the (graded) subring of $R$ given by $R':=\bigoplus_{d\ne1}R_d$\index{R'@$R'$}
we have therefore that $R'\iso\R'/(\theta)$. Since $\theta$ is a {\nzdiv} in
$\R'$, the ``hyperplane section principle'' (see \cite[p. 218]{R}) implies that
$\R'$ has the same generators as $R'$ plus $\theta$, and that the relations of
$\R'$ are given by lifts of those of $R'$ (in particular, the two rings have
the same number of relations and in the same degrees).

So, as a first step to compute $\R'$, we will find generators and
relations of $R'$. Setting $r_i:=y_0^{2-i}y_1^i\in R'_2$ ($i=0,1,2$), $s_i:=
y_0^{3-i}y_1^i\in R'_3$ ($i=0,1,2,3$) and $t_i:=y_iz\in R'_4$ ($i=0,1$), it is
clear that $R'$ is minimally generated (as a graded $\C$--algebra), by the
$r_i$, $s_i$, $t_i$ and $z$. As for relations, we have the following result.

    \begin{lemm}\label{curvering}
The kernel $J$ of the natural epimorphism of graded rings
\[P:=\C[r_0,r_1,r_2,s_0,s_1,s_2,s_3,z,t_0,t_1]\epi R'\]
is minimally generated by the following $37$ homogeneous polynomials (of which
$1$ has degree $4$, $6$ have degree $5$, $17$ have degree $6$, $10$ have degree
$7$ and $3$ have degree $8$):
\begin{align*}
& r_1^2-r_0r_2 & \text{(of degree $4$);}\\
& r_is_j-r_{i+1}s_{j-1} \text{ for $i=0,1$ and $j=1,2,3$} &
\text{(of degree $5$);}\\
& r_it_1-r_{i+1}t_0 \text{ for $i=0,1$} & \text{(of degree $6$);}\\
& s_is_j-r_0^{3-\lfloor\frac{i+j+1}{2}\rfloor}
r_1^{i+j-2\lfloor\frac{i+j}{2}\rfloor}r_2^{\lfloor\frac{i+j}{2}\rfloor}
\text{ for $0\le i\le j\le3$} & \text{(of degree $6$);}\\
& s_iz-r_it_0 \text{ for $i=0,1,2$ and } s_3z-r_2t_1 &
\text{(of degree $6$);}\\
& z^2-r_1(r_2^2+\lambda r_1r_2+\mu r_0r_2+\nu r_0r_1+\ep r_0^2) & \text{(of
degree $6$);}\\
& s_it_j-r_0^{2-\lfloor\frac{i+j+1}{2}\rfloor}
r_1^{i+j-2\lfloor\frac{i+j}{2}\rfloor}r_2^{\lfloor\frac{i+j}{2}\rfloor}z
\text{ for $0\le i\le3$ and $j=0,1$} & \text{(of degree $7$);}\\
& zt_i-s_{i+1}(r_2^2+\lambda r_1r_2+\mu r_0r_2+\nu r_0r_1+\ep r_0^2)
\text{ for $i=0,1$} & \text{(of degree $7$);}\\
& t_it_j-r_1r_{i+j}(r_2^2+\lambda r_1r_2+\mu r_0r_2+\nu r_0r_1+\ep r_0^2)
\text{ for $0\le i\le j\le1$} & \text{(of degree $8$).}
\end{align*}
    \end{lemm}

    \begin{proof}
Denoting by $J'$ the (homogeneous) ideal generated by the polynomials in the
above list, it is clear that $J'\subseteq J$. $\all n\in\N$ let $\tilde{R}_n$
be the subspace of $P_n$ defined in the following way: $\tilde{R}_n:=P_n$ if
$n\le3$, whereas for $m\ge2$
\begin{gather*}
\tilde{R}_{2m}:=\langle r_0^{m-i-j-2k}r_1^jr_2^it_0^k\st
0\le j,k\le 1;0\le i\le m-j-2k\rangle\oplus\langle r_2^{m-2}t_1\rangle,\\
\tilde{R}_{2m+1}:=\langle r_0^{m-1-i-j}r_1^jr_2^is_0^{1-k}z^k\st
0\le j,k\le 1;0\le i<m-j\rangle
\oplus\langle r_2^{m-1}s_i\st1\le i\le3\rangle.
\end{gather*}
Notice that, for $n>1$, $\dim_{\C}\tilde{R}_n=2n-1$, so that
$\dim_{\C}\tilde{R}_n=\dim_{\C}R'_n$ $\all n\in\N$. Moreover, it follows easily
from the definition that $\tilde{R}_n+J'_n=P_n$ $\all n\in\N$. Therefore
\[\dim_{\C}J'_n\ge\dim_{\C}P_n-\dim_{\C}\tilde{R}_n=\dim_{\C}P_n-\dim_{\C}R'_n=
\dim_{\C}J_n,\]
which implies (since $J'_n\subseteq J_n$) that $J'_n=J_n$ $\all n\in\N$, i.e.
that $J'=J$. The given list of polynomials is also a minimal set of generators
of $J$: indeed, it is not difficult to check that, if $J''$ denotes the ideal
generated by the whole list except one polynomial (say of degree $n$), then
$J''_n+\tilde{R}_n\subsetneq P_n$.
    \end{proof}

By what we said before, it follows that $\R'$ has $11$ generators ($\theta$ of
degree $1$, $3$ of degree $2$, $5$ of degree $3$ and $2$ of degree $4$) and
$37$ relations (in the same degrees as those of $R'$) and that, for instance,
denoting by $\rho_i$ lifts of $r_i$, since the only relation in degree $4$ in
$R'$ is $r_1^2=r_0r_2$, $\R'$ has also a single relation in degree $4$, and it
is of the form $\rho_1^2=\rho_0\rho_2+\theta\varphi$ for some $\varphi\in
\R'_3$. In order to determine this $\varphi$ and the analogous terms in the
other relations, we need to study the ring $\R'$ more directly.

Consider the divisor $E':=r^*(E)$ of $C\times C$ and observe that $E'=
\{(c,\iota(c))\st c\in C\}$.

    \begin{lemm}\label{Thetacohom}
Let $D:=D_1+D_2+E'$. Then $\all n\in\Z$ there are natural isomorphisms
\[H^0(A,\so_A(n\Theta))\iso H^0(C\times C,\so_{C\times C}(nD))^+,\]
where the superscript $^+$ denotes the subspace of invariants under the
natural action of $\Z/(2)$ induced by the map which exchanges the two factors
of $C\times C$.
    \end{lemm}

    \begin{proof}
Since $q_*\so_{C^{(2)}}\iso\so_A$, by projection formula $\so_A(n\Theta)\iso
q_*q^*\so_A(n\Theta)$, whence $H^0(A,\so_A(n\Theta))\iso
H^0(C^{(2)},q^*\so_A(n\Theta))$. Observe moreover that
\[H^0(C\times C,r^*q^*\so_A(n\Theta))\iso
H^0(C^{(2)},q^*\so_A(n\Theta)\otimes r_*\so_{C\times C}).\]
As the invariant part of $r_*\so_{C\times C}$ is $\so_{C^{(2)}}$, it follows
that
\[H^0(A,\so_A(n\Theta))\iso H^0(C^{(2)},q^*\so_A(n\Theta))\iso
H^0(C\times C,r^*q^*\so_A(n\Theta))^+.\]
To conclude, just notice that $q^*\so_A(\Theta)\iso\so_{C^{(2)}}(\Theta'+E)$
(where $\Theta':=r(D_1)=r(D_2)$), $r^*\so_{C^{(2)}}(\Theta')\iso
\so_{C\times C}(D_1+D_2)$ and $r^*\so_{C^{(2)}}(E)\iso\so_{C\times C}(E')$,
so that $r^*q^*\so_A(\Theta)\iso\so_{C\times C}(D)$.
    \end{proof}

We can therefore identify the elements of $H^0(A,\so_A(n\Theta))$ with
symmetric rational functions on $C\times C$ having as poles at most $nD$. We
will denote by $[y'_0,y'_1,z']$ the coordinates on the second factor of
$C\times C\iso\proj(R\times_{\C}R)$ (where $(R\times_{\C}R)_n:=
R_n\otimes_{\C}R_n$). Our aim is to find functions corresponding to generators
of $\R'$ (then it will be not difficult to determine the relations, because in
this way multiplication in the ring is simply given by multiplication of
functions). Now, some of them are easy to find: it is clear that $\theta:=1$ is
a generator of $H^0(A,\so_A(\Theta))$. As for $H^0(A,\so_A(2\Theta))$, besides
$1=\theta^2$ we have for instance $\frac{y_1\otimes y'_1}{y_0\otimes y'_0}$ and
$\frac{y_1\otimes y'_0+y_0\otimes y'_1}{y_0\otimes y'_0}$ (which have as poles
$2D_1+2D_2$). On the other hand, as we will see, a fourth generator of
$H^0(A,\so_A(2\Theta))$ is not so simple (it has poles also at $2E'$), and the
same is true for two generators of degree $3$ and one of degree $4$. Therefore
we need a systematic procedure which allows us to determine all sections of $n
\Theta$ (at least for $n\le 4$). Before we describe it, we need to know
explicitly the morphism $\R'\epi R'$: it is not difficult to prove the
following result.

    \begin{lemm}
The map $\R'_n=H^0(A,\so_A(n\Theta))\epi R'_n$ can be identified with
\[\gamma\mapsto y_0^n((\frac{y'_0}{z'})^n\gamma\rest{y'_0=z'=0,y'_1=1})\]
(notice that $\frac{y'_0}{{z'}^2}\rest{y'_0=z'=0,y'_1=1}=1$).
    \end{lemm}

For instance, it follows that (as expected) $\theta$ and $\theta^2$ are mapped
to $0$, whereas $\frac{y_1\otimes y'_1}{y_0\otimes y'_0}$ goes to $y_0y_1=r_1$
and $\frac{y_1\otimes y'_0+y_0\otimes y'_1}{y_0\otimes y'_0}$ to $y_0^2=r_0$
(so the still missing generator of degree $2$ will be a lift of $r_2$).

In order to shorten notation in the following we will write $a\odot
b$\index{$\odot$} for $a(y_0,y_1,z)\otimes b(y'_0,y'_1,z')+b(y_0,y_1,z)\otimes
a(y'_0,y'_1,z')$ and $a\wedge b$ for $a(y_0,y_1,z)\otimes
b(y'_0,y'_1,z')-b(y_0,y_1,z)\otimes a(y'_0,y'_1,z')$. Consider the
$\Z/(2)\times\Z/(2)$ Galois cover $H:=h\times h: C\times C\to\P^1\times\P^1$:
setting $D'_1:=\{[0,1]\}\times\P^1$ and $D'_2:= \P^1\times\{[0,1]\}$, we have
$H^*(D'_i)=2D_i$ and $H^*(\Delta_{\P^1})=\Delta_C +E'$ ($\Delta_X$ denotes the
diagonal in $X\times X$), whence
\[nD=H^*(\tilde{D}_n)-n\Delta_C-(n-2\lfloor\frac{n}{2}\rfloor)(D_1+D_2),\]
where $\tilde{D}_n:=n\Delta_{\P^1}+\lfloor\frac{n+1}{2}\rfloor(D'_1+D'_2)$.
Thus we see that $H^0(C\times C,\so_{C\times C}(nD))$ can  be identified with a
subspace of $H^0(C\times C,\so_{C\times C}(H^*(\tilde{D}_n)))$. More precisely,
every global section $\gamma$ of $\so_{C\times C}(H^*(\tilde{D}_n))$ can be
written in the form
\[\gamma=\frac{G}
{(y_0\otimes y'_0)^{\lfloor\frac{n+1}{2}\rfloor}(y_1\wedge y'_0)^n}\]
with $G\in(R\times_{\C}R)_{n+\lfloor\frac{n+1}{2}\rfloor}$, and then the
condition $\gamma\in H^0(C\times C,\so_{C\times C}(nD))^+$ is equivalent to the
following: $G$ vanishes on $\Delta_C$ of order $n$, $G$ vanishes on $D_1+D_2$
if $n$ is odd, and $G$ is symmetric or antisymmetric according to whether $n$
is even or odd. The problem of determining the $G$ satisfying the above
condition eventually translates into a problem of linear algebra: we are going
to see how only for $n=2$ (the cases $n=3$ and $n=4$ are similar, only slightly
more complicated). So we are looking for symmetric $G\in R_3\otimes_{\C}R_3$
vanishing on $\Delta_C$ of order $2$. Passing to affine coordinates
$(y=y_1/y_0,\tilde{z}=z/y_0)$, $(y'=y'_1/y'_0,\tilde{z}'=z'/y'_0)$ and dropping
tensor products, we want the $g$ vanishing on $\Delta_C$ of order $2$ and of
the form
\[g=u(y,y')+v(y)\tilde{z}'+\tilde{z}v(y')+2k\tilde{z}\tilde{z}'\]
with $u(y,y')=u(y',y)$ of degree $\le3$ in each variable, $\deg(v)\le3$ and $k
\in\C$. Now, the condition $g\rest{\Delta}=0$ immediately implies $v=0$, so
that $g=u(y,y')+2k\tilde{z}\tilde{z}'$. Setting $f(y):=F(1,y)$, we have
$(\tilde{z}-\tilde{z}')^2=f(y)+f(y')-2\tilde{z}\tilde{z}'$, so that $g$
vanishes on $\Delta_C$ of order $2$ if and only if the same is true for
\[g':=g+k(\tilde{z}-\tilde{z}')^2=u(y,y')+k[f(y)+f(y')].\]
This last condition is obviously satisfied if and only if $g'(y,y')$ is
divisible by $(y-y')^2$. Actually, as $g'(y,y')=g(y',y)$, it is enough to
require that $g'$ is divisible by $(y-y')$, i.e. that $g'(y,y)=u(y,y)+2kf(y)=
0$. Now it is clear that this gives a linear system of $7$ equations (since
$\deg(g'(y,y))=6$) with $11$ unknowns ($k$ and the coefficients of $u$, which
are $10$ because $u$ is symmetric). It turns out that the space of solutions
has the expected dimension $4=h^0(A,\so_A(2\Theta))$, and a base is given by
the following values of $k$ and $u(y,y')$ (the first three solutions correspond
to the already known functions $1$, $\frac{y_1\odot y'_0}{y_0\otimes y'_0}$ and
$\frac{y_1\otimes y'_1}{y_0\otimes y'_0}$):

$k=0$, $u=(y-y')^2$;

$k=0$, $u=(y+y')(y-y')^2$;

$k=0$, $u=yy'(y-y')^2$;

$k=-1$, $u=y^2{y'}^2(y+y')+2\lambda y^2{y'}^2+\mu yy'(y+y')+2\nu yy'+
\ep(y+y')$.

Carrying out a similar analysis also for $n=3$ and $n=4$, one eventually
obtains the following result.

    \begin{prop}
Identifying, as above, global sections of $\so_A(n\Theta)$ with rational
functions on $C\times C$, $\R'$ is generated (as a graded $\C$--algebra)
by the following (homogeneous) elements:
\[\theta:=1\text{ of degree $1$;}\]
\begin{align*}
& \rho_0:=\frac{y_1\odot y'_0}{y_0\otimes y'_0}, &
\rho_1:=\frac{y_1\otimes y'_1}{y_0\otimes y'_0},
\end{align*}
\begin{multline*}
\rho_2:=(y_1^3\odot y'_0{y'_1}^2+\lambda y_0y_1^2\odot y'_0{y'_1}^2+
\mu y_0y_1^2\odot {y'_0}^2y'_1+\nu y_0^2y_1\odot{y'_0}^2y'_1\\
+\ep y_0^2y_1\odot{y'_0}^3-z\odot z')/[(y_0\otimes y'_0)(y_1\wedge y'_0)^2]
\end{multline*}
of degree $2$;\index{rho_i@$\rho_i$}
\begin{align*}
& \sigma_0:=\frac{z\wedge{y'_0}^3}{(y_0^2\otimes{y'_0}^2)(y_1\wedge y'_0)}, &
& \sigma_1:=\frac{z\wedge{y'_0}^2y'_1}{(y_0^2\otimes{y'_0}^2)(y_1\wedge y'_0)},
&
\sigma_2:=\frac{z\wedge y'_0{y'_1}^2}{(y_0^2\otimes{y'_0}^2)(y_1\wedge y'_0)},
\end{align*}
\begin{multline*}
\sigma_3:=[y_1^2z\wedge{y'_0}^2{y'_1}^3+y_0y_1z\wedge(3y'_0{y'_1}^4+
4\lambda{y'_0}^2{y'_1}^3+3\mu{y'_0}^3{y'_1}^2+2\nu{y'_0}^4y'_1+\ep{y'_0}^5)\\
+y_0^2z\wedge(\mu{y'_0}^2{y'_1}^3+2\nu{y'_0}^3{y'_1}^2+3\ep{y'_0}^4y'_1)]/
[(y_0^2\otimes{y'_0}^2)(y_1\wedge y'_0)^3],
\end{multline*}
\begin{multline*}
\zeta:=(y_1z\odot y'_0z'-y_0y_1^3\odot y'_0{y'_1}^3-\lambda y_0y_1^3\odot
{y'_0}^2{y'_1}^2-\mu y_0y_1^3\odot{y'_0}^3y'_1-\nu y_0^2y_1^2\odot{y'_0}^3y'_1
\\
-\ep y_0^3y_1\odot{y'_0}^3y'_1)/[(y_0^2\otimes{y'_0}^2)(y_1\wedge y'_0)^2]
\end{multline*}
of degree $3$;\index{sigma_i@$\sigma_i$}\index{zeta@$\zeta$}
\[\tau_0:=\frac{y_1^3\wedge z'}{(y_0^2\otimes{y'_0}^2)(y_1\wedge y'_0)},\]
\begin{multline*}
\tau_1:=[\ep y_0^3y_1^2\wedge{y'_0}^2z'+
(y_1^5+2\lambda y_0y_1^4+3\mu y_0^2y_1^3+4\nu y_0^3y_1^2+3\ep y_0^4y_1)\wedge
y'_0y'_1z'\\
+(3y_0y_1^4+2\lambda y_0^2y_1^3+\mu y_0^3y_1^2)\wedge{y'_1}^2z']/
[(y_0^2\otimes{y'_0}^2)(y_1\wedge y'_0)^3]
\end{multline*}
of degree $4$.\index{tau_i@$\tau_i$}

Moreover, the images in $R'\iso\R'/(\theta)$ of $\rho_i$, $\sigma_i$, $\zeta$
and $\tau_i$ are, respectively, $r_i$, $s_i$, $z$ and $t_i$.
    \end{prop}

    \begin{rema}\label{invar}
$\theta$, $\rho_i$ and $\zeta$ are invariant with respect to the involution of
$A$ given by multiplication by $-1$, whereas $\sigma_i$ and $\tau_i$ are
antiinvariant.
    \end{rema}

    \begin{rema}
The above functions are not uniquely determined by the condition of being lifts
of the corresponding generators of $R'$: for instance, the $\rho_i$ could be
changed by adding multiples of $\theta^2=1$. On the other hand, the $\sigma_i$
are unique if we require them to be antiinvariant, since the elements of
$\R'_{<3}$ are all invariant.
    \end{rema}

Now we have to determine the relations of $\R'$, and it is clear that this
becomes a problem of linear algebra (which can be easily solved by a computer).
For instance, the unique relation of degree $4$ must be of the form $\rho_1^2=
\rho_0\rho_2+\theta\varphi$ with (since $\varphi$ is invariant)
\[\varphi=a_1\theta^3+a_2\theta\rho_0+a_3\theta\rho_1+a_4\theta\rho_2+a_5\zeta
\]
for some $a_i\in\C$. The condition to impose on the $a_i$ is then simply the
vanishing of $\rho_1^2=\rho_0\rho_2+\theta\varphi$, i.e. of its numerator,
which is a polynomial in $y_0$, $y_1$, $y'_0$, $y'_1$, $z$, $z'$ (not in this
case, but in higher degrees one has of course to use the fact that $z^2=
F(y_0,y_1)$ and ${z'}^2=F(y'_0,y'_1)$), whose coefficients depend linearly on
the $a_i$. The result of these computations is the following.

    \begin{prop}\label{thetaring}
$\R':=\R(A,\Theta)\iso Q'/I'$, where $Q'$ is the polynomial algebra
\[Q':=\C[\theta,\rho_0,\rho_1,\rho_2,\sigma_0,\sigma_1,\sigma_2,\sigma_3,\zeta,
\tau_0,\tau_1]\] with $\deg(\theta)=1$, $\deg(\rho_i)=2$,
$\deg(\sigma_i)=\deg(\zeta)=3$, $\deg(\tau_i)=4$, and $I'$ is the homogeneous
ideal which is minimally generated by the following polynomials (such that
$\deg(\rel{1})=4$, $\deg(\rel{i})=5$ for $2\le i\le7$, $\deg(\rel{i})=6$ for
$8\le i\le24$, $\deg(\rel{i})=7$ for $25\le i\le34$ and $\deg(\rel{i})=8$ for
$35\le i\le37$):\index{rel'_i@$\rel{i}$|(}
\begin{itemize}

\item [] $\rel{1}:=\rho_1^2-\rho_0\rho_2+\theta(-2\zeta-\mu\theta\rho_1+
\ep\theta^3)$;

\item [] $\rel{2}:=\rho_0\sigma_1-\rho_1\sigma_0-\theta^2\sigma_2$;

\item [] $\rel{3}:=\rho_0\sigma_2-\rho_1\sigma_1+\theta\tau_0$;

\item [] $\rel{4}:=\rho_0\sigma_3-\rho_1\sigma_2+
\theta(2\tau_1-\ep\theta\sigma_0-2\nu\theta\sigma_1-\mu\theta\sigma_2)$;

\item [] $\rel{5}:=\rho_1\sigma_1-\rho_2\sigma_0+
\theta(-2\tau_0+\mu\theta\sigma_1+2\lambda\theta\sigma_2+\theta\sigma_3)$;

\item [] $\rel{6}:=\rho_1\sigma_2-\rho_2\sigma_1-\theta\tau_1$;

\item [] $\rel{7}:=\rho_1\sigma_3-\rho_2\sigma_2+\ep\theta^2\sigma_1$;

\item [] $\rel{8}:=\rho_0\tau_1-\rho_1\tau_0+\theta(\mu\rho_1\sigma_1+
2\lambda\rho_1\sigma_2+2\rho_1\sigma_3+\ep\theta^2\sigma_1)$;

\item [] $\rel{9}:=\rho_1\tau_1-\rho_2\tau_0+\theta(-2\ep\rho_0\sigma_1-
2\nu\rho_1\sigma_1-\mu\rho_1\sigma_2+\ep\theta^2\sigma_2)$;

\item [] $\rel{10}:=\sigma_0^2-\rho_0^3+\theta^2(-\lambda\rho_0^2+\rho_0\rho_1-
\mu\theta^2\rho_0-\theta^2\rho_2-\nu\theta^4)$;

\item [] $\rel{11}:=\sigma_0\sigma_1-\rho_0^2\rho_1+\theta^2(\theta\zeta-
\lambda\rho_0\rho_1)$;

\item [] $\rel{12}:=\sigma_0\sigma_2-\rho_0^2\rho_2+\theta(-\rho_0\zeta+
\theta\rho_1\rho_2+\ep\theta^3\rho_0+\nu\theta^3\rho_1)$;

\item [] $\rel{13}:=\sigma_0\sigma_3-\rho_0\rho_1\rho_2+\theta(\rho_1\zeta-
\mu\theta^2\zeta-\ep\theta\rho_0^2-2\nu\theta\rho_0\rho_1-\mu\theta\rho_0\rho_2
-2\lambda\theta\rho_1\rho_2-\theta\rho_2^2-\lambda\ep\theta^3\rho_0-
2\lambda\nu\theta^3\rho_1-\nu\theta^3\rho_2)$;

\item [] $\rel{14}:=\sigma_1^2-\rho_0\rho_1^2-
\theta^2\rho_1(\lambda\rho_1+\rho_2)$;

\item [] $\rel{15}:=\sigma_1\sigma_2-\rho_0\rho_1\rho_2-\theta\rho_1\zeta$;

\item [] $\rel{16}:=\sigma_1\sigma_3-\rho_0\rho_2^2+\theta(-\rho_2\zeta+
\ep\theta\rho_0\rho_1+\lambda\ep\theta^3\rho_1+\ep\theta^3\rho_2)$;

\item [] $\rel{17}:=\sigma_2^2-\rho_1^2\rho_2-
\theta^2\rho_1(\ep\rho_0+\nu\rho_1)$;

\item [] $\rel{18}:=\sigma_2\sigma_3-\rho_1\rho_2^2+\theta^2(\ep\theta\zeta-
\nu\rho_1\rho_2)$;

\item [] $\rel{19}:=\sigma_3^2-\rho_2^3+\theta^2(\ep\rho_1\rho_2-\nu\rho_2^2-
\ep^2\theta^2\rho_0-\mu\ep\theta^2\rho_2-\lambda\ep^2\theta^4)$;

\item [] $\rel{20}:=\sigma_0\zeta-\rho_0\tau_0+
\theta(\theta\tau_1+\mu\rho_1\sigma_0+\lambda\rho_0\sigma_2+\rho_0\sigma_3-
\ep\theta^2\sigma_0-\nu\theta^2\sigma_1)$;

\item [] $\rel{21}:=\sigma_1\zeta-\rho_1\tau_0+
\theta\rho_1(\mu\sigma_1+\lambda\sigma_2+\sigma_3)$;

\item [] $\rel{22}:=\sigma_2\zeta-\rho_2\tau_0-
\theta\sigma_1(\ep\rho_0+\nu\rho_1)$;

\item [] $\rel{23}:=\sigma_3\zeta-\rho_2\tau_1+
\theta(-\nu\theta\tau_1+\ep\rho_0\sigma_2+\nu\rho_1\sigma_2+\mu\rho_2\sigma_2+
\lambda\ep\theta^2\sigma_2)$;

\item [] $\rel{24}:=\zeta^2-\rho_1(\rho_2^2+\lambda\rho_1\rho_2+\mu\rho_0\rho_2
+\nu\rho_0\rho_1+\ep\rho_0^2)-\theta^2\rho_1(\lambda\ep\rho_0+\lambda\nu\rho_1+
\nu\rho_2)$;

\item [] $\rel{25}:=\sigma_0\tau_0-\rho_0^2\zeta+\theta\rho_1(\theta\zeta-
\mu\rho_0^2-\lambda\rho_0\rho_1-\rho_0\rho_2-\nu\theta^2\rho_0)$;

\item [] $\rel{26}:=\sigma_1\tau_0-\rho_0\rho_1\zeta+
\theta\rho_1(-\mu\rho_0\rho_1-\lambda\rho_1^2-\rho_1\rho_2+\ep\theta^2\rho_0)$;

\item [] $\rel{27}:=\sigma_2\tau_0-\rho_1^2\zeta+
\theta\rho_0\rho_1(\ep\rho_0+\nu\rho_1)$;

\item [] $\rel{28}:=\sigma_3\tau_0-\rho_1\rho_2\zeta+
\theta(-\ep\theta\rho_0\zeta+\ep\rho_0\rho_1^2+\nu\rho_0\rho_1\rho_2+
\lambda\ep\theta^2\rho_1^2+\ep\theta^2\rho_1\rho_2)$;

\item [] $\rel{29}:=\sigma_0\tau_1-\rho_0\rho_1\zeta+\theta(-\theta\rho_2\zeta+
\lambda\rho_0\rho_1\rho_2+\rho_1^2\rho_2+\ep\theta^2\rho_0\rho_1+
\nu\theta^2\rho_1^2)$;

\item [] $\rel{30}:=\sigma_1\tau_1-\rho_1^2\zeta+
\theta\rho_1\rho_2(\lambda\rho_1+\rho_2)$;

\item [] $\rel{31}:=\sigma_2\tau_1-\rho_1\rho_2\zeta+
\theta\rho_1(-\ep\rho_0\rho_1-\nu\rho_1^2-\mu\rho_1\rho_2+\ep\theta^2\rho_2)$;

\item [] $\rel{32}:=\sigma_3\tau_1-\rho_2^2\zeta+\theta\rho_1(\ep\theta\zeta-
\ep\rho_0\rho_2-\nu\rho_1\rho_2-\mu\rho_2^2-\lambda\ep\theta^2\rho_2)$;

\item [] $\rel{33}:=\zeta\tau_0-\rho_1^2\sigma_3-\lambda\rho_1^2\sigma_2-
\mu\rho_0\rho_1\sigma_2-\nu\rho_0\rho_1\sigma_1-\ep\rho_0^2\sigma_1-
\ep\theta^3\tau_0$;

\item [] $\rel{34}:=\zeta\tau_1-\rho_1\rho_2\sigma_3-\lambda\rho_1^2\sigma_3-
\mu\rho_1^2\sigma_2-\nu\rho_1^2\sigma_1-\ep\rho_1^2\sigma_0+
\theta(\mu\rho_1\tau_1-\lambda\ep\theta^2\tau_0-\lambda\ep\theta\rho_0\sigma_2-
\ep\theta\rho_1\sigma_2)$;

\item [] $\rel{35}:=\tau_0^2-\rho_1(\rho_1^2\rho_2+\lambda\rho_1^3+
\mu\rho_0\rho_1^2+\nu\rho_0^2\rho_1+\ep\rho_0^3)+\ep\theta^2\rho_0\rho_1^2$;

\item [] $\rel{36}:=\tau_0\tau_1-\rho_1(\rho_1\rho_2^2+\lambda\rho_1^2\rho_2+
\mu\rho_0\rho_1\rho_2+\nu\rho_0\rho_1^2+\ep\rho_0^2\rho_1)+
\theta\rho_1(-\mu\rho_1\zeta+\ep\theta^2\zeta-2\lambda\ep\theta\rho_0\rho_1
-\ep\theta\rho_0\rho_2-2\lambda\nu\theta\rho_1^2-2\nu\theta\rho_1\rho_2)$;

\item [] $\rel{37}:=\tau_1^2-\rho_1(\rho_2^3+\lambda\rho_1\rho_2^2+
\mu\rho_1^2\rho_2+\nu\rho_1^3+\ep\rho_0\rho_1^2)+\ep\theta^2\rho_1^2\rho_2$.
\end{itemize}
    \end{prop}\index{rel'_i@$\rel{i}$|)}

    \begin{rema}
The above polynomials are not all lifts of those listed in \ref{curvering}: for
instance, $\rel{17}$ is a lift of $s_2^2-r_1^2r_2$ and not of $s_2^2- r_0r_2^2$
(but, of course, \ref{curvering} remains true if $s_2^2-r_1^2r_2$ is
substituted by $s_2^2-r_0r_2^2$). The reason for this and other changes was (to
try) to keep the $\rel{i}$ as simple as possible.
    \end{rema}

        \section{Computation of $\s{E}(\gph)$}

If $s\in\langle\theta^2,\rho_0,\rho_1\rangle$, it is not difficult to see that
$0\in B\subset A$ is a singular point. So we can assume that $s=\rho_2+a\rho_0+
b\rho_1+c\theta^2\in\R'_2$ for some
$a,b,c\in\C$.\index{a@$a$}\index{b@$b$}\index{c@$c$} Notice that, since
$\xi^2=s$, in $\R$ there is the relation
\[\rho_2=\xi^2-a\rho_0-b\rho_1-c\theta^2.\]
The structure of $\R$ then follows immediately from \ref{thetaring}.

    \begin{prop}\label{canring}
The canonical ring $\R=\R(\S,\ds_{\S})$ is isomorphic to $Q/I$, where $Q$ is
the polynomial algebra
\[Q:=\C[\xi,\theta,\rho_0,\rho_1,\sigma_0,\sigma_1,\sigma_2,\sigma_3,\zeta,
\tau_0,\tau_1]\] with $\deg(\xi)=\deg(\theta)=1$, $\deg(\rho_i)=2$,
$\deg(\sigma_i)=\deg(\zeta)= 3$, $\deg(\tau_i)=4$, and $I$ is the homogeneous
ideal which is minimally generated by the polynomials
$\crel{i}$\index{rel_i@$\crel{i}$|(} ($1\le i\le37$) defined as $\rel{i}$ in
\ref{thetaring}, except that $\rho_2$ is substituted by $\xi^2-a\rho_0-b\rho_1-
c\theta^2$.
    \end{prop}

We will write here only the first $\crel{i}$, which we will use directly later.
\begin{itemize}

\item [] $\crel{1}=a\rho_0^2+b\rho_0\rho_1+\rho_1^2-\xi^2\rho_0+
\theta(-2\zeta+c\theta\rho_0-\mu\theta\rho_1+\ep\theta^3)$;

\item [] $\crel{2}=\rho_0\sigma_1-\rho_1\sigma_0-\theta^2\sigma_2$;

\item [] $\crel{3}=\rho_0\sigma_2-\rho_1\sigma_1+\theta\tau_0$;

\item [] $\crel{4}=\rho_0\sigma_3-\rho_1\sigma_2+\theta(2\tau_1-
\ep\theta\sigma_0-2\nu\theta\sigma_1-\mu\theta\sigma_2)$;

\item [] $\crel{5}=a\rho_0\sigma_0+b\rho_1\sigma_0+\rho_1\sigma_1-\xi^2\sigma_0
+\theta(-2\tau_0+c\theta\sigma_0+\mu\theta\sigma_1+2\lambda\theta\sigma_2+
\theta\sigma_3)$;

\item [] $\crel{6}=a\rho_0\sigma_1+b\rho_1\sigma_1+\rho_1\sigma_2-\xi^2\sigma_1
+\theta(-\tau_1+c\theta\sigma_1)$;

\item [] $\crel{7}=a\rho_0\sigma_2+b\rho_1\sigma_2+\rho_1\sigma_3-\xi^2\sigma_2
+\theta^2(\ep\sigma_1+c\sigma_2)$;

\item [] $\crel{8}=\rho_0\tau_1-\rho_1\tau_0+\theta(\mu\rho_1\sigma_1+
2\lambda\rho_1\sigma_2+2\rho_1\sigma_3+\ep\theta^2\sigma_1)$;

\item [] $\crel{9}=\rho_1\tau_1+(a\rho_0+b\rho_1-\xi^2+c\theta^2)\tau_0+
\theta(-2\ep\rho_0\sigma_1-2\nu\rho_1\sigma_1-\mu\rho_1\sigma_2+
\ep\theta^2\sigma_2)$;

\item [] $\crel{10}=\sigma_0^2-\rho_0^3+\theta^2[-\lambda\rho_0^2+\rho_0\rho_1+
(a-\mu)\theta^2\rho_0+b\theta^2\rho_1-\theta^2\xi^2+(c-\nu)\theta^4]$.

\end{itemize}\index{rel_i@$\crel{i}$|)}

    \begin{rema}\label{regseq}
It is clear that $\R/(\theta)$ is isomorphic to the ring $R'\oplus\xi' R'$ with
multiplication defined by
\[(x+\xi'y)(x'+\xi'y')=xx'+(r_2+ar_0+br_1)yy'+(xy'+x'y)\xi'.\]
Since $\xi'$ is a {\nzdiv} in $R'\oplus\xi' R'$, we see that $(\theta,\xi)$ is
a regular sequence in $\R$.
    \end{rema}

From now on we set $\w:=(1,1,2,3)$, and we will write (as usual) $\p$ for
$\p(\w)$, $\gP$ for $\gP(\w)$ and $\pd_i$ for $\pd_i(\w)$.

Assume that $\gph:\gS\to\gP$ is a \gbwcp: $\gph$ is determined by a morphism of
graded rings $\rh:\p\to\R$, and (up to composing $\gph$ with an automorphism of
$\gP$) we can always assume that
\begin{align*}
& \rh(x_0)=\theta, & & \rh(x_1)=\xi, & & \rh(x_2)\in\langle\rho_0,\rho_1
\rangle\minus\{0\}, & \rh(x_3)\in\R_3\minus\C[\theta,\xi,\rh(x_2)]_3.
\end{align*}

    \begin{rema}
Not every choice of $\rh$ as above determines a \gbwcp. For instance, we must
have $\rh(x_3)\notin\C[\theta,\xi,\rho_0,\rho_1,\zeta]_3$, otherwise $\ph$
would factor through the involution of $\S$ which lifts multiplication by $-1$
on $A$ (see \ref{invar}). On the other hand, one can prove that $\gph$ is good
birational for a general choice of $\rh$ (see also \ref{excexi}). Later we will
prove that $\gph$ is good birational in a particular case.
    \end{rema}

Since $\sw{\w}=7$, $\pd_1=2$, $\pd_2=4$ the vector bundle $\s{E}(\gph)$ defined
in \ref{Edef} is in our case given by
\[\s{E}(\gph)=\so(-2)\oplus\sd^1(-1)^2\oplus\so(-3)^{\cre_{-1}(\gph)}\oplus
\so(-4)^{\cre_{-2}(\gph)}\oplus\so(-5)^{\cre_{-3}(\gph)}.\]
In order to compute the $\cre_j(\gph)$ we need some preliminary results.

First we study morphisms from $\sd_{\gP}^1$ to $\gph_*\so_{\gS}(1)$: the exact
sequence
\[\sko^{-3}\mor{\cdiff{\sko}^{-3}}\sko^{-2}\to\sd^1\to0\]
shows that $\Hom_{\gP}(\sd^1,\gph_*\so_{\gS}(1))\iso
\ker\Hom_{\gP}(\cdiff{\sko}^{-3},\gph_*\so_{\gS}(1))$. Remembering the
definition of $\cdiff{\sko}^{-3}$, this means that $\varphi\in
\Hom_{\gP}(\sd^1,\gph_*\so_{\gS}(1))$ corresponds to $6$ elements
\[\varphi_{i,j}\in\Hom_{\gP}(\so(-\w_i-\w_j),\gph_*\so_{\gS}(1))\iso
\R_{\w_i+\w_j+1}\]
for $0\le i<j\le3$, satisfying the $4$ relations
\begin{equation}\label{drel}
\rho(x_i)\varphi_{j,k}-\rh(x_j)\varphi_{i,k}+\rh(x_k)\varphi_{i,j}=0
\end{equation}
for $0\le i<j<k\le3$. On the other hand, applying
$\Hom_{\gP}(-,\gph_*\so_{\gS}(1))$ to the short exact sequence $\esd{1}$, we
obtain the exact sequence
\[\Hom_{\gP}(\sko^{-1},\gph_*\so_{\gS}(1))\to
\Hom_{\gP}(\sd^1,\gph_*\so_{\gS}(1))\mor{\delta}
\Ext_{\gP}^1(\so,\gph_*\so_{\gS}(1))\to0.\] Recall that
$\Ext_{\gP}^1(\so,\gph_*\so_{\gS}(1))\iso H^1(\gP,\gph_*\so_{\gS}(1))$ and
$h^1(\gP,\gph_*\so_{\gS}(1))=q(\S)=2$.

    \begin{lemm}\label{omegamor}
Let $\eta,\eta'\in\Hom_{\gP}(\sd^1,\gph_*\so_{\gS}(1))$ be such that
$\delta\rest{\langle\eta,\eta'\rangle}$ is an isomorphism. If $\rh(x_2)=
k_0\rho_0+k_1\rho_1$, then
\begin{align*}
\langle\eta_{0,1},\eta'_{0,1}\rangle & =\langle a\sigma_0+b\sigma_1+\sigma_2,
a\sigma_1+b\sigma_2+\sigma_3\rangle\subset\R_3/\R_1\R_2=
\langle\sigma_0,\sigma_1,\sigma_2,\sigma_3,\zeta\rangle,\\
\langle\eta_{1,2},\eta'_{1,2}\rangle & =\langle-k_0\tau_0-k_1\tau_1,(k_0b-k_1a)
\tau_0+k_0\tau_1\rangle\subset\R_4/(\R_1\R_3+\R_2^2)=
\langle\tau_0,\tau_1\rangle.
\end{align*}
    \end{lemm}

    \begin{proof}
First we show that if $\varphi_{0,1}\in\R_3$ is such that $x_2\varphi_{0,1},
x_3\varphi_{0,1}\in(\theta,\xi)$, then it can be extended to a $\varphi\in
\Hom_{\gP}(\sd^1,\gph_*\so_{\gS}(1))$ (of course, the converse is also true by
\eqref{drel}). Indeed, by hypothesis we can choose $\varphi_{0,2},\varphi_{1,2}
\in\R_4$ and $\varphi_{0,3},\varphi_{1,3}\in\R_5$ such that
\[x_0\varphi_{1,2}-x_1\varphi_{0,2}+x_2\varphi_{0,1}=
x_0\varphi_{1,3}-x_1\varphi_{0,3}+x_3\varphi_{0,1}=0,\]
and it remains to find $\varphi_{2,3}\in\R_6$ such that
\[x_0\varphi_{2,3}-x_2\varphi_{0,3}+x_3\varphi_{0,2}=
x_1\varphi_{2,3}-x_2\varphi_{1,3}+x_3\varphi_{1,2}=0.\]
By the above relations we have
\[x_1(x_2\varphi_{0,3}-x_3\varphi_{0,2})=x_2(x_0\varphi_{1,3}+x_3\varphi_{0,1})
-x_3(x_0\varphi_{1,2}+x_2\varphi_{0,1})=x_0(x_2\varphi_{1,3}-x_3\varphi_{1,2})
,\]
which implies (since $(x_0,x_1)$ is a regular sequence in $\R$ by
\ref{regseq}) that $\exi\varphi_{2,3}\in\R_6$ such that $x_0\varphi_{2,3}=
x_2\varphi_{0,3}-x_3\varphi_{0,2}$ and $x_1\varphi_{2,3}=x_2\varphi_{1,3}-
x_3\varphi_{1,2}$.

Now, it is clear that if $\varphi\in\Hom_{\gP}(\sd^1,\gph_*\so_{\gS}(1))$
factors through $\sko^{-1}$, then $\varphi_{0,1}\in(\theta,\xi)$. Therefore the
image of $\varphi_{0,1}$ in $\R_3/(\theta,\xi)=\R_3/\R_1\R_2$ doesn't depend on
$\varphi$ up to morphisms factoring through $\sko^{-1}$, and the first equality
will follow if we prove that $\eta_{0,1}:=a\sigma_0+b\sigma_1+\sigma_2$ and
$\eta'_{0,1}:=a\sigma_1+b\sigma_2+\sigma_3$ satisfy $\rho\eta_{0,1},
\rho\eta'_{0,1}\in(\theta,\xi)$ $\all\rho\in\R_+$. Passing to the quotient ring
$\R/(\theta)\iso R'\oplus\xi'R'$ (see \ref{regseq}), this is the same as
proving that $\bar{\eta}_{0,1}:=as_0+bs_1+s_2$ and $\bar{\eta}'_{0,1}:=
as_1+bs_2+s_3$ are such that $r\bar{\eta}_{0,1},r\bar{\eta}'_{0,1}\in(\xi')\
\all r\in R'_+$. Since ${\xi'}^2=ar_0+br_1+r_2$, in the overring $R=
\C[y_0,y_1,z]/(z^2-F)$ of $R'$ we have $\bar{\eta}_{0,1}=y_0{\xi'}^2$ and
$\bar{\eta}'_{0,1}=y_1{\xi'}^2$, whence $r\bar{\eta}_{0,1}=(y_0r){\xi'}^2$ and
$r\bar{\eta}'_{0,1}=(y_1r){\xi'}^2$, and $y_0r,y_1r\in R'$.

As for the second equality, observe again that if $\varphi\in
\Hom_{\gP}(\sd^1,\gph_*\so_{\gS}(1))$ factors through $\sko^{-1}$, then
$\varphi_{1,2}\in(\xi,k_0\rho_0+k_1\rho_1)$, so that it is enough to check that
if $\eta$ (respectively $\eta'$) is a morphism such that $\eta_{0,1}=a\sigma_0+
b\sigma_1+\sigma_2$ (respectively $\eta'_{0,1}=a\sigma_1+b\sigma_2+\sigma_3$),
then $\eta_{1,2}=-k_0\tau_0-k_1\tau_1$ (respectively $\eta'_{1,2}=
(k_0b-k_1a)\tau_0+k_0\tau_1$) in $\R_4/(\R_1\R_3+\R_2^2)$. The other cases
being similar, we will do it only for $\eta$ assuming $k_0=1$, $k_1=0$ (so that
$\rh(x_2)=\rho_0$ and we expect $\eta_{1,2}=-\tau_0$). By $\crel{2}$ and
$\crel{3}$ we have
\[\rho_0(a\sigma_0+b\sigma_1+\sigma_2)=a\rho_0\sigma_0+b\rho_1\sigma_0+
b\theta^2\sigma_2+\rho_1\sigma_1-\theta\tau_0,\]
and then it follows from $\crel{5}$ that
\[\rho_0(a\sigma_0+b\sigma_1+\sigma_2)=\xi^2\sigma_0-
\theta[-\tau_0+c\theta\sigma_0+\mu\theta\sigma_1+(2\lambda-b)\theta\sigma_2+
\theta\sigma_3],\]
whence we can take $\eta_{0,2}=\xi\sigma_0$ and $\eta_{1,2}=-\tau_0+
c\theta\sigma_0+\mu\theta\sigma_1+(2\lambda-b)\theta\sigma_2+\theta\sigma_3$.
It is also clear that a different choice would change $\eta_{1,2}$ by a
multiple of $\xi$.
    \end{proof}

We will use ${\rm\check{C}}$ech resolutions over the affine open covering
$\cov=\{\{x_i\ne0\}\st i=0,\dots,3\}$ of $\gP$ to compute cohomology.

    \begin{lemm}\label{cohomgen}
Let $\eta,\eta'$ be as in \ref{omegamor}. Then
\[\{\frac{\eta_{i,j}}{x_ix_j}\},
\{\frac{\eta'_{i,j}}{x_ix_j}\}\in Z^1(\cov,\gph_*\so_{\gS}(1))\]
and they correspond to generators of $H^1(\cov,\gph_*\so_{\gS}(1))\iso
H^1(\gP,\gph_*\so_{\gS}(1))$.
    \end{lemm}

    \begin{proof}
It is clear by definition that if $\varphi\in
\Hom_{\gP}(\sd^1,\gph_*\so_{\gS}(1))$ then $\{\frac{\varphi_{i,j}}{x_ix_j}\}$
is a cocycle, and that if $\varphi$ factors through $\sko^{-1}$ then
$\{\frac{\varphi_{i,j}}{x_ix_j}\}\in B^1(\cov,\gph_*\so_{\gS}(1))$.
    \end{proof}

We are now ready to compute the coefficients $\cre_j(\gph)$ of $\so(j-2)$ (for
$j=-1,-2,-3$) in $\s{E}(\gph)$. First observe that (by \ref{coeffdiff})
\[\cre_{-1}(\gph)-\cre_{-3}(\gph)=
\chi(\gP,\gph_*\so_{\gS}(2)\otimes\cp{\bd{-1}})-\chi(\gP,\cp{\bd{-1}}(2))+
\chi(\gP,\cp{\bd{-3}}(2))\] and that
$\chi(\gP,\cp{\bd{-1}}(2))=-\chi(\gP,\cp{\bd{-3}}(2))=1$ by \ref{O2res}
(indeed, $h^0(\gP,\cp{\bd{-1}}(2))=h^{-1}(\gP,\cp{\bd{-3}}(2))=1$ and
$h^{-1}(\gP,\cp{\bd{-1}}(2))=h^0(\gP,\cp{\bd{-3}}(2))=0$). By definition
$\cp{\bd{-1}}$ is the complex
\[0\to\so(-1)\mor{\begin{pmatrix}
-x_1 \\
x_0
\end{pmatrix}}\begin{matrix}
\so \\
\oplus \\
\so
\end{matrix}
\mor{\begin{pmatrix}
x_0 & x_1
\end{pmatrix}}\so(1)\to0\]
(with $\bd{-1}^0=\so(1)$). Therefore
\begin{multline*}
\chi(\gP,\gph_*\so_{\gS}(2)\otimes\cp{\bd{-1}})=\sum_{i=-2}^0(-1)^i
\chi(\gP,\gph_*\so_{\gS}(2)\otimes\bd{-1}^i)\\
=\chi(\gP,\gph_*\so_{\gS}(3))-2\chi(\gP,\gph_*\so_{\gS}(2))+
\chi(\gP,\gph_*\so_{\gS}(1))\\
=\chi(\so_{\S})+3K_{\S}^2-
2(\chi(\so_{\S})+K_{\S}^2)+\chi(\so_{\S})=K_{\S}^2=4,
\end{multline*}
so that $\cre_{-1}(\gph)-\cre_{-3}(\gph)=2$.

    \begin{lemm}\label{cpcohom}
Let $\rh(x_2)=k_0\rho_0+k_1\rho_1$. Then
\[H^0(\gP,\gph_*\so_{\gS}(2)\otimes\cp{\bd{-1}})\iso\C\oplus
\R_3/(\R_1\R_2+\langle a\sigma_0+b\sigma_1+\sigma_2,a\sigma_1+b\sigma_2+
\sigma_3\rangle)\]
(in particular, $h^0(\gP,\gph_*\so_{\gS}(2)\otimes\cp{\bd{-1}})=4$) and
\[H^0(\gP,\gph_*\so_{\gS}(2)\otimes\cp{\bd{-2}})\iso\langle\tau_0,\tau_1\rangle
/\langle-k_0\tau_0-k_1\tau_1,(k_0b-k_1a)\tau_0+k_0\tau_1\rangle.\]
Moreover, with this identification the image of
$H^0(\gP,\gps\mrs(2)\otimes\cp{\bd{-1}})$ is $\langle(0,\rh(x_3))\rangle\in
H^0(\gP,\gph_*\so_{\gS}(2)\otimes\cp{\bd{-1}})$.
    \end{lemm}

    \begin{proof}
$\gph_*\so_{\gS}(2)\otimes\cp{\bd{-1}}$ is the complex
\[0\to\gph_*\so_{\gS}(1)\mor{\begin{pmatrix}
-x_1 \\
x_0
\end{pmatrix}}\begin{matrix}
\gph_*\so_{\gS}(2) \\
\oplus \\
\gph_*\so_{\gS}(2)
\end{matrix}
\mor{\begin{pmatrix}
x_0 & x_1
\end{pmatrix}}\gph_*\so_{\gS}(3)\to0\]
(with $\bd{-1}^0=\gph_*\so_{\gS}(3)$). Its ${\rm\check{C}}$ech resolution is a
double complex whose total complex $\cp{C}$ is such that $H^i(\cp{C})\iso
H^i(\gP,\gph_*\so_{\gS}(2)\otimes\cp{\bd{-1}})$. $\cp{C}$ is explicitly given
by
\[C^i:=C^i(\cov,\gph_*\so_{\gS}(3))\oplus C^{i+1}(\cov,\gph_*\so_{\gS}(2)^2)
\oplus C^{i+2}(\cov,\gph_*\so_{\gS}(1))\]
and (denoting always by $\cp{d}$ the differentials of the various
$\cp{C}(\cov,\gph_*\so_{\gS}(j))$)
\[\cdiff{C}^i(u,(v^0,v^1),w):=
(x_0v^0+x_1v^1+d^i(u),
(-x_1w-d^{i+1}(v^0),x_0w-d^{i+1}(v^1)),
d^{i+2}(w))\]
$\all u\in C^i(\cov,\gph_*\so_{\gS}(3))$, $\all v^0,v^1\in
C^{i+1}(\cov,\gph_*\so_{\gS}(2))$ and $\all w\in
C^{i+2}(\cov,\gph_*\so_{\gS}(1))$.

First we fix
\[\bar{w}\in Z^2(\cov,\gph_*\so_{\gS}(1))\minus B^2(\cov,\gph_*\so_{\gS}(1))\]
(observe that $H^2(\cov,\gph_*\so_{\gS}(1))\iso H^2(\gP,\gph_*\so_{\gS}(1))\iso
\C$),
\[\bar{v}^0,\bar{v}^1\in C^1(\cov,\gph_*\so_{\gS}(2))\text{ such that
$d^1(\bar{v}^1)=x_0\bar{w}$, $d^1(\bar{v}^0)=-x_1\bar{w}$}\]
(this is possible because $d^2(x_0\bar{w})=d^2(x_1\bar{w})=0$ and
$H^2(\cov,\gph_*\so_{\gS}(2))=0$) and
\[\bar{u}\in C^0(\cov,\gph_*\so_{\gS}(3))\text{ such that
$d^0(\bar{u})=x_0\bar{v}^0+x_1\bar{v}^1$}\]
(again, this can be done because $d^1(x_0\bar{v}^0+x_1\bar{v}^1)=0$ and
$H^1(\cov,\gph_*\so_{\gS}(3))=0$). Now, let $(u,(v^0,v^1),w)\in Z^0(\cp{C})$:
by definition of $\cdiff{C}^0$, this is equivalent to the following conditions:
\begin{align*}
& d^0(u)=-x_0v^0-x_1v^1, & & d^1(v^0)=-x_1w, & & d^1(v^1)=x_0w, & d^2(w)=0.
\end{align*}
As $d^2(w)=0$, $\exiun k\in\C$ and $\exi w'\in
C^1(\cov,\gph_*\so_{\gS}(1))$ such that $w=k\bar{w}-d^1(w')$. Since
\[d^1(v^0-x_1w'-k\bar{v}^0)=-x_1w-x_1(k\bar{w}-w)+kx_1\bar{w}=0\]
and similarly $d^1(v^1+x_0w'-k\bar{v}^1)=0$, and since
$H^1(\cov,\gph_*\so_{\gS}(2))=0$, $\exi{v'}^0,{v'}^1\in
C^0(\cov,\gph_*\so_{\gS}(2))$ such that $d^0({v'}^0)=v^0-x_1w'-k\bar{v}^0$ and
$d^0({v'}^1)=v^1+x_0w'-k\bar{v}^1$. It follows that
\[(u,(v^0,v^1),w)+\cdiff{C}^{-1}(0,({v'}^0,{v'}^1),w')=
(u+x_0{v'}^0+x_1{v'}^1,(k\bar{v}^0,k\bar{v}^1),k\bar{w}).\]
As $d^0(u+x_0{v'}^0+x_1{v'}^1)=-kx_0\bar{v}^0-kx_1\bar{v}^1=-kd^0(\bar{u})$, we
see that $u+x_0{v'}^0+x_1{v'}^1+k\bar{u}$ defines an element $\rho\in
H^0(\cov,\gph_*\so_{\gS}(3))\iso\R_3$. Now we want to see how this $\rho$
changes for different choices of $w'$, ${v'}^0$ and ${v'}^1$. Instead of them
we could choose $w'+w''$, ${v'}^0+{v''}^0$ and ${v'}^1+{v''}^1$, where $w''$ is
an arbitrary element of $Z^1(\cov,\gph_*\so_{\gS}(1))$, ${v''}^i\in
C^0(\cov,\gph_*\so_{\gS}(2))$ are such that $d^0({v''}^0)=-x_1w''$ and
$d^0({v''}^1)=x_0w''$, and then $\rho$ would be changed by $\tilde{\rho}=
x_0{v''}^0+x_1{v''}^1\in Z^0(\cov,\gph_*\so_{\gS}(3))=
H^0(\cov,\gph_*\so_{\gS}(3))\iso\R_3$. If $w''=d^0(\tilde{w})\in
B^1(\cov,\gph_*\so_{\gS}(1))$, then ${v''}^0=-x_1\tilde{w}+\tilde{v}^0$ and
${v''}^1=x_0\tilde{w}+\tilde{v}^1$ for arbitrary $\tilde{v}^i\in
Z^0(\cov,\gph_*\so_{\gS}(2))=H^0(\cov,\gph_*\so_{\gS}(2))\iso\R_2$, so that
in this case $\tilde{\rho}=x_0\tilde{v}^0+x_1\tilde{v}^1\in\R_1\R_2$.
Otherwise, by \ref{cohomgen} we can assume $w''_{i,j}=
\frac{\eta_{i,j}}{x_ix_j}$ or $w''_{i,j}=\frac{\eta'_{i,j}}{x_ix_j}$. In the
first case (the other one is completely analogous) the ${v''}^i$ are given (up
to elements of $Z^0(\cov,\gph_*\so_{\gS}(2))$) by
\begin{align*}
{v''}^0_0 & =\frac{\eta_{0,1}}{x_0}, & {v''}^0_1 & =0, &
{v''}^0_2 & =-\frac{\eta_{1,2}}{x_2}, & {v''}^0_3 & =-\frac{\eta_{1,3}}{x_3},\\
{v''}^1_0 & =0, & {v''}^1_1 & =\frac{\eta_{0,1}}{x_1}, &
{v''}^1_2 & =\frac{\eta_{0,2}}{x_2}, & {v''}^1_3 & =\frac{\eta_{0,3}}{x_3},
\end{align*}
so that $\tilde{\rho}=\eta_{0,1}$ (which we can assume to be $a\sigma_0+
b\sigma_1+\sigma_2$ by \ref{cohomgen}). Now it is clear that the map
\[\begin{split}
H^0(\cp{C}) & \to\C\oplus\R_3/
(\R_1\R_2+\langle a\sigma_0+b\sigma_1+\sigma_2,a\sigma_1+b\sigma_2+\sigma_3
\rangle) \\
{[(u,(v^0,v^1),w)]} & \mapsto(k,[\rho])
\end{split}\]
is well defined, and it is easily seen to be an isomorphism (non canonical,
since it depends on the choices of $\bar{w}$, $\bar{v}^0$, $\bar{v}^1$ and
$\bar{u}$): this proves the first statement. The computation of
$H^0(\gP,\gph_*\so_{\gS}(2)\otimes\cp{\bd{-2}})$ is similar and we shall omit
it; we just want to underline the fact that (by \ref{omegamor}) the terms
$-k_0\tau_0-k_1\tau_1$ and $(k_0b-k_1a)\tau_0+k_0\tau_1$ are just $\eta_{1,2}$
and $\eta'_{1,2}$. It is straightforward to check that the map
\[H^0(\gP,\gps\mrs(2)\otimes\cp{\bd{-1}}):H^0(\gP,\cp{\bd{-1}}(2))\iso\C\to
H^0(\gP,\gph_*\so_{\gS}(2)\otimes\cp{\bd{-1}})\]
is as stated.
    \end{proof}

    \begin{coro}\label{cj}
If $\rh(x_3)\in\langle a\sigma_0+b\sigma_1+\sigma_2,
a\sigma_1+b\sigma_2+\sigma_3\rangle+\R_1\R_2$ then $\cre_{-1}(\gph)=3$ and
$\cre_{-3}(\gph)=1$, otherwise $\cre_{-1}(\gph)=2$ and $\cre_{-3}(\gph)=0$.

Let $\rh(x_2)=k_0\rho_0+k_1\rho_1$. If $k_0^2-bk_0k_1+ak_1^2=0$ then
$\cre_{-2}(\gph)=1$, otherwise $\cre_{-2}(\gph)=0$.
    \end{coro}

    \begin{proof}
By \ref{Fres} we have
\begin{gather*}
\cre_{-1}(\gph)=h^0(\gP,\gph_*\so_{\gS}(2)\otimes\cp{\bd{-1}})-z^0_{-1}-
z^{-1}_{-3}+k_{-1}+k_{-3},\\
\cre_{-2}(\gph)=h^0(\gP,\gph_*\so_{\gS}(2)\otimes\cp{\bd{-2}})-z^0_{-2}-
z^{-1}_{-2}+2k_{-2}.
\end{gather*}
By \ref{O2res} $z^0_{-1}=z^{-1}_{-3}=1$, $z^0_{-2}=z^{-1}_{-2}=0$ (whence
$k_{-2}=0$) and $k_{-3}=0$ (because also $z^0_{-3}=0$). Therefore
$\cre_{-1}(\gph)=2+k_{-1}$ and $\cre_{-2}(\gph)=
h^0(\gP,\gph_*\so_{\gS}(2)\otimes\cp{\bd{-2}})$, and then everything follows
from \ref{cpcohom}.
    \end{proof}

    \begin{rema}\label{excexi}
It is possible to prove that a {\gbwcp} $\gph$ with $\cre_{-1}(\gph)=3$ or
$\cre_{-2}(\gph)=1$ actually exists. For instance, it is easy to see that the
projection into $\gP(1,1,2,2,3,3)$ determined (in the obvious way) by
$a\sigma_0+b\sigma_1+\sigma_2$ and $a\sigma_1+b\sigma_2+\sigma_3$ is good
birational, and then it follows that there exists $\gph$ good birational into
$\gP(1,1,2,3)$ with $\rh(x_3)\in\langle
a\sigma_0+b\sigma_1+\sigma_2,a\sigma_1+b\sigma_2+\sigma_3\rangle$.
    \end{rema}

        \section{Explicit resolution}

We will assume that $\rh(x_2)=\rho_0$, $\rh(x_3)=\sigma_0$ (and, as usual,
$\rh(x_0)=\theta$, $\rh(x_1)=\xi$). First of all, we must check that this is an
admissible choice.

    \begin{lemm}
With $\rh$ as above, $\gph$ is a \gbwcp.
    \end{lemm}

    \begin{proof}
$\gph$ is a morphism because every generator of $\R$ is in
$\sqrt{(\theta,\xi,\rho_0)}$ ($\rho_1$ by $\crel{1}$, $\zeta$, $\sigma_i$ and
$\tau_i$ by the $\crel{j}$ containing their squares). As for birationality,
using $\crel{10}$ we see that $x_0^2(x_2+bx_0^2)\rho_1$ belongs to the subring
of $\R$ generated by $\theta$, $\xi$, $\rho_0$ and $\sigma_0$. This means that
(generically) two points of $\gS$ mapped by $\gph$ to the same point of $\gP$
have the same coordinate $\rho_1$, and then it is immediate to see that the
same is true also for $\zeta$, $\sigma_i$ ($i>0$) and $\tau_i$.
    \end{proof}

    \begin{rema}
$\crel{10}$ also shows that $\gY=\im\gph$ does not contain the points
$[0,0,1,0]$ and $[0,0,0,1]$, which are (by \ref{stdmori}) the only non {\std}
points of $\gP$ (and also the only singular points of $\P$).
    \end{rema}

By \ref{cj} we have $\s{E}:=\s{E}(\gph)=\so(-2)\oplus\so(-3)^2\oplus
\sd^1(-1)^2$, so that the resolution given by \ref{mthm} is of the form
\[0\to\so(-8)\oplus\so(-6)\oplus\so(-5)^2\oplus(\sd^2)^2\mor{\alpha}
\so\oplus\so(-2)\oplus\so(-3)^2\oplus\sd^1(-1)^2\mor{\gamma}\gph_*\so_{\gS}\to0
\]
with $\alpha=\alpha\dual(-8)$. As a first step to construct $\alpha$, we choose
$\gamma$; $\Hom_{\gP}(\so(-j),\gph_*\so_{\gS})$ will be identified, as usual,
with $\R_j$. Then we take $\gamma\rest{\so}:=1$, $\gamma\rest{\so(-2)}:=
\rho_1$, $\gamma\rest{\so(-3)^2}:=(\sigma_1,\zeta)$ and
$\gamma\rest{\sd^1(-1)^2}:=(\eta,\eta')$, where $\eta,\eta'\in
\Hom_{\gP}(\sd^1(-1),\gph_*\so_{\gS})\iso\Hom_{\gP}(\sd^1,\gph_*\so_{\gS}(1))$
are as in \ref{omegamor}. Making the computations (using the relations in $\R$)
it turns out that they can be chosen as follows:

$\eta_{0,1}=a\sigma_0+b\sigma_1+\sigma_2$;

$\eta_{0,2}=\xi\sigma_0$;

$\eta_{1,2}=-\tau_0+c\theta\sigma_0+\mu\theta\sigma_1+
(2\lambda-b)\theta\sigma_2+\theta\sigma_3$;

$\eta_{0,3}=\xi\rho_0^2$;

$\eta_{1,3}=-\rho_0\zeta-b\theta^2\zeta+(c+ab-a\lambda)\theta\rho_0^2+
b(b-\lambda)\theta\rho_0\rho_1-b\theta\xi^2\rho_0+\theta\xi^2\rho_1+
(a^2+bc-a\mu+\ep)\theta^3\rho_0+(ab-b\mu+\nu-c)\theta^3\rho_1-a\theta^3\xi^2+
(ac-a\nu+b\ep)\theta^5$;

$\eta_{2,3}=\theta\xi(\rho_0\rho_1-\lambda\rho_0^2+(a-\mu)\theta^2\rho_0+
b\theta^2\rho_1-\theta^2\xi^2+(c-\nu)\theta^4)$;

$\eta'_{0,1}=a\sigma_1+b\sigma_2+\sigma_3$;

$\eta'_{0,2}=\xi\sigma_1$;

$\eta'_{1,2}=b\tau_0+\tau_1-\ep\theta\sigma_0+(c-2\nu)\theta\sigma_1-
\mu\theta\sigma_2$;

$\eta'_{0,3}=\xi\rho_0\rho_1$;

$\eta'_{1,3}=b\rho_0\zeta+\rho_1\zeta+(a-4b^2+4b\lambda-\mu)\theta^2\zeta+
(a\mu-a^2+2ab^2-2ab\lambda-\ep)\theta\rho_0^2+
(c+a\lambda+2b\mu-2\nu-3ab+2b^3-2b^2\lambda)\theta\rho_0\rho_1+
(2a-2b^2+2b\lambda-\mu)\theta\xi^2\rho_0+(3b-2\lambda)\theta\xi^2\rho_1+
(a\nu+c\mu-2ac-\lambda\ep+2b^2c-2bc\lambda)\theta^3\rho_0+(2c\lambda+2b\nu-3bc-
2\lambda\nu-2b^2\mu+2b\lambda\mu)\theta^3\rho_1-
\theta\xi^4+(2c-\nu)\theta^3\xi^2+(2b^2\ep-2b\lambda\ep+c\nu-c^2)\theta^5$;

$\eta'_{2,3}=\theta\xi(\theta\zeta-\lambda\rho_0\rho_1)$.

    \begin{lemm}
$\gamma:=(1,\rho_1,\sigma_1,\zeta,\eta,\eta'):\so\oplus\so(-2)\oplus\so(-3)^2
\oplus\sd^1(-1)^2\to\gph_*\so_{\gS}$ is such that $H^1(\gP,\gamma(1))$ and
$H^0(\gP,\gamma(i))$ $\all i\in\Z$ are surjective.
    \end{lemm}

    \begin{proof}
By \ref{omegamor} and \ref{cohomgen} it is clear that
\[H^1(\gP,(\eta,\eta')(1)):H^1(\gP,(\sd^1)^2)\iso\C^2\to
H^1(\gP,\gph_*\so_{\gS}(1))\iso\C^2\] is an isomorphism. As for $H^0$, since
$H^0(\gP,\gph_*\so_{\gS}(i))\iso\R_i\ \all i\in\Z$, it is enough to prove that
$\R$ is generated (as graded $\p$--module) by
$\{1,\rho_1,\sigma_1,\zeta,\eta_{i,j},\eta'_{i,j}\}$, and this will follow if
we show that $\R$ is generated by $\{1,\rho_1,\sigma_1,\sigma_2,
\sigma_3,\zeta,\tau_0,\tau_1,\rho_1\zeta\}$ ($\sigma_2$ comes from
$\eta_{0,1}$, $\sigma_3$ from $\eta'_{0,1}$, $\tau_0$ from $\eta_{1,2}$,
$\tau_1$ from $\eta'_{1,2}$ and $\rho_1\zeta$ from $\eta'_{1,3}$). Now, it is
easy to see that we can always assume that a monomial in $\R$ contains at most
one of the variables $\sigma_i$, $\zeta$ and $\tau_i$, and with exponent at
most $1$. Moreover, the relation $\crel{1}$ implies that we can assume that
also $\rho_1$ appears with exponent at most $1$. It follows that $\R$ is
generated by
\[\{1,\rho_1,\sigma_1,\sigma_2,\sigma_3,\zeta,\tau_0,\tau_1,\rho_1\sigma_1,
\rho_1\sigma_2,\rho_1\sigma_3,\rho_1\zeta,\rho_1\tau_0,\rho_1\tau_1\}.\]
Moreover, using $\crel{3}$, $\crel{4}$ and $\crel{7}$ it is immediate to see
that $\rho_1\sigma_1$, $\rho_1\sigma_2$ and $\rho_1\sigma_3$ are actually not
needed, and similarly for $\rho_1\tau_0$ and $\rho_1\tau_1$ (using $\crel{8}$
and $\crel{9}$).
    \end{proof}

Then it follows from \ref{resdet} that $\gamma$ is surjective and $\ker\gamma
\iso(\so\oplus\s{E})\dual(-8)$. Now we want to construct a morphism $\beta:
(\so\oplus\s{E})\dual(-8)\to\so\oplus\s{E}$ such that
\[0\to(\so\oplus\s{E})\dual(-8)\mor{\beta}\so\oplus\s{E}\mor{\gamma}
\gph_*\so_{\gS}\to0\]
is exact, and (again by \ref{resdet}) it is enough to require that $\gamma\comp
\beta=0$ and $H^0(\gP,\beta(i))$ is injective $\all i\in\Z$. First we need to
fix some notation. The exact sequence
\[0\to\sd^1(-1)\to\sko^{-1}(-1)=\so(-2)^2\oplus\so(-3)\oplus\so(-4)
\mor{\cdiff{\sko}^{-1}(-1)}\so(-1)\to0\]
induces naturally an exact sequence
\[0\to\so\oplus\s{E}\mor{j}\tilde{\s{E}}\mor{\tilde{d}}\so(-1)^2\to0,\]
where $\tilde{\s{E}}:=\so\oplus\so(-2)\oplus\so(-3)^2\oplus
(\so(-2)^2\oplus\so(-3)\oplus\so(-4))^2$. Denoting $\dual(-8)$ by $\tdual$, we
have also an exact sequence
\[0\to\so(-7)^2\mor{\tilde{d}\tdual}\tilde{\s{E}}\tdual\mor{j\tdual}
(\so\oplus\s{E})\tdual\to0.\] It is then clear that giving
$\beta:(\so\oplus\s{E})\tdual\to\so\oplus\s{E}$ is equivalent to giving
$\tilde{\beta}:\tilde{\s{E}}\tdual\to\tilde{\s{E}}$ such that
$\tilde{d}\comp\tilde{\beta}=0$ and $\tilde{\beta}\comp \tilde{d}\tdual=0$ (of
course, $\tilde{\beta}=j\comp\beta\comp j\tdual$), and that $\beta$ is
symmetric if and only if $\tilde{\beta}$ is. Observe that $\tilde{\beta}$ can
be identified with a $12\times12$ matrix of homogeneous polynomials (in $\p$).
Moreover, the exact sequence
\[0\to\sd^2(-1)\to\sko^{-2}(-1)=\so(-3)\oplus\so(-4)^2\oplus\so(-5)^2\oplus
\so(-6)\to\sd^1(-1)\to0\]
induces naturally an exact sequence
\[0\to\sd^2(-1)^2\to\s{E}'\mor{p'}\so\oplus\s{E}\to0,\]
where $\s{E}':=\so\oplus\so(-2)\oplus\so(-3)^2\oplus\sko^{-2}(-1)^2$. While
$\beta$ need not factor through $p'$ (because $\Ext_{\gP}^1(\sd^2,\sd^2(-1))\ne
0$, and actually the morphism we are looking for does not factor), $\beta\comp
j\tdual$ factors through $p'$ (since
$\Ext_{\gP}^1(\tilde{\s{E}}\tdual,\sd^2(-1))=0$), so that we must have
$\tilde{\beta}=j\comp p'\comp\beta'$ for some $\beta':\tilde{\s{E}}\tdual\to
\s{E'}$, as we can see from the following diagram
\[\begin{CD}
0 @>>> \so(-7)^2 @>\tilde{d}\tdual>> \tilde{\s{E}}\tdual @>j\tdual>>
(\so\oplus\s{E})\tdual @>>> 0 \\
@. @. @VV{\beta'}V @VV{\beta}V \\
0 @>>> \sd^2(-1)^2 @>>> \s{E}' @>p'>> \so\oplus\s{E} @>j>> \tilde{\s{E}}
@>\tilde{d}>> \so(-1)^2
\end{CD}\]
which has exact rows, except at $\so\oplus\s{E}$. Then the condition
$\tilde{d}\comp\tilde{\beta}=0$ is automatically satisfied, whereas
$\tilde{\beta}\comp \tilde{d}\tdual=0$ is equivalent to
$p'\comp\beta'\comp\tilde{d}\tdual=0$ (because $j$ is injective). Setting
$\gamma':=\gamma\comp p':\s{E}'\to \gph_*\so_{\gS}$ (notice that, by
definition, $\gamma'=(1,\rho_1,\sigma_1,
\zeta,\{\eta_{i,j}\},\{\eta'_{i,j}\})$), $\gamma\comp\beta=0$ is clearly
equivalent to $\gamma'\comp\beta'=0$. Finally, $H^0(\gP,\beta(i))$ injective
$\all i\in\Z$ translates into the condition that
\begin{equation}\label{gsecseq}
0\to H^0(\gP,\so(i-7)^2)\mor{H^0(\gP,\tilde{d}\tdual(i))}
H^0(\gP,\tilde{\s{E}}\tdual(i))\mor{H^0(\gP,(j\comp p'\comp\beta')(i))}
H^0(\gP,\tilde{\s{E}}(i))
\end{equation}
is an exact sequence $\all i\in\Z$ (this is a consequence of the fact that,
$\all i\in\Z$, $H^0(\gP,j(i))$ is injective and $H^0(\gP,j\tdual(i))$ is
surjective, as $H^0(\gP,\so(i-7))=0$).

We are going to sketch only how it is possible to find such a $\beta'$. Since
\[\tilde{\s{E}}\tdual=\so(-4)^2\oplus\so(-5)^4\oplus\so(-6)^5\oplus\so(-8),\]
$\beta'$ is given by sections $a_1,a_2\in H^0(\gP,\s{E}'(4))$,
$b_1,\dots,b_4\in H^0(\gP,\s{E}'(5))$, $c_1,\dots,c_5\in H^0(\gP,\s{E}'(6))$,
$e\in H^0(\gP,\s{E}'(8))$. As we want $\gamma'\comp\beta'=0$, it is clear that
each of them must be in the kernel of the corresponding $H^0(\gP,\gamma'(i))$
(and must not be mapped to $0$ by $H^0(\gP,(j\comp p')(i))$ by
\eqref{gsecseq}). So, for instance, denoting by $g(\rho)$ (where
$\rho:\s{L}\iso\so(-i)\to \gph_*\so_{\gS}$ is one of the components of
$\gamma'$) a generator of $H^0(\gP,\s{L}(i))\iso\C$, it is easy to see that
$\ker H^0(\gP,\gamma'(4))$ is generated by $-x_1x_3g(1)+g(\eta_{0,2})$ and
$x_1g(\sigma_1)-g(\eta'_{0,2})$, so that we can choose them as $a_1$ and $a_2$.
Continuing in this way, it is not difficult to find candidates also for $b_i$,
$c_i$ and $e$: the $b_i$ must be chosen in such a way that, together with the
$x_la_i$ ($l=0,1$), they generate $\ker
H^0(\gP,\gamma'(5))/H^0(\gP,\sd^2(4)^2)$ (notice that before we didn't consider
the analogous quotient because $H^0(\gP,\sd^2(3))=0$), and similarly for the
$c_i$ and $e$. It is clear that then $\gamma'\comp\beta' =0$ and the sequences
\eqref{gsecseq} are exact, and it is also not hard to get that
$p'\comp\beta'\comp\tilde{d}\tdual=0$ and even that the resulting morphism from
$(\sd^2)^2$ to $\sd^1(-1)^2$ is already symmetric.

Once we have thus obtained the required $\beta$ (corresponding to
$\tilde{\beta}=j\comp p'\comp\beta'$), it remains only to symmetrize it. By
\ref{resdet} there exists an automorphism $\delta$ of $(\so\oplus\s{E})\tdual$
such that $\alpha:=\beta\comp\delta=\alpha\tdual$ and
\[0\to(\so\oplus\s{E})\tdual\mor{\alpha}\so\oplus\s{E}\mor{\gamma}
\gph_*\so_{\gS}\to0\] is exact. It is immediate to see that finding such a
$\delta$ is equivalent to finding an automorphism $\tilde{\delta}$ of
$\tilde{\s{E}}\tdual$ such that
$\tilde{\alpha}:=\tilde{\beta}\comp\tilde{\delta}=\tilde{\alpha}\tdual$ and
$\tilde{\delta}\comp\tilde{d}\tdual=\tilde{d}\tdual\comp s$ for some
automorphism $s$ of $\so(-7)^2$. Regarding $\tilde{\beta}$ as a matrix,
composing $\tilde{\beta}$ with an automorphism of $\tilde{\s{E}}\tdual$
corresponds to making ``homogeneous'' elementary operations on the columns: one
can multiply a column by a non zero constant, can exchange two columns
corresponding to morphisms from $\so(i_1)$ and $\so(i_2)$ with $i_1=i_2$, and
can add to a column a multiple of the ``right'' degree of another one. In our
case, since the component of $\tilde{\beta}$ from $(\sko^{-1}(-1)\tdual)^2\iso
(\sko^{-3})^2$ to $\sko^{-1}(-1)^2$ is already symmetric, we need only to
modify the first $4$ columns (those giving
$\tilde{\beta}\rest{\so(-8)\oplus\so(-6)\oplus\so(5)^2}$), which assures that
$\tilde{\delta}\comp\tilde{d}\tdual=\tilde{d}\tdual$. After these computations,
we eventually obtain the following symmetric matrix
$\tilde{\alpha}=j\comp\alpha\comp j\tdual$. Denoting by $p_2$, $p_4$, $q_2$,
$q_4$ the following homogeneous polynomials (subscripts denote degrees)
\begin{align*}
p_2 & :=(ab-a\lambda+c)x_2-bx_1^2+bcx_0^2,\\
p_4 & :=(a\lambda-ab-c)x_2^2+bx_1^2x_2+(a\mu-a^2-bc-\ep)x_0^2x_2+
ax_0^2x_1^2+(a\nu-ac-b\ep)x_0^4,\\
q_2 & :=b(\lambda-b)x_2-x_1^2+(b\mu-ab+c-\nu)x_0^2,\\
q_4 & :=\lambda x_2^2+(\mu-a)x_0^2x_2+x_0^2x_1^2+(\nu-c)x_0^4,
\end{align*}
the first $4$ columns of $\tilde{\alpha}$ are

\begin{tabular}{r|cccc}
\multicolumn{2}{c}{} & \\
& $\so(-8)$ & $\so(-6)$ & $\so(-5)$ & $\so(-5)$ \\
\hline
$\so$ & $0$ & $x_0x_3p_2$ & $x_0p_4$ & $ax_0^2x_3$ \\
$\so(-2)$ & $x_0x_3p_2$ & $2b(b-\lambda)x_0x_3$ & $x_0q_2$ & $-x_3$ \\
$\so(-3)$ & $x_0p_4$ & $x_0q_2$ & $0$ & $x_2+bx_0^2$ \\
$\so(-3)$ & $ax_0^2x_3$ & $-x_3$ & $x_2+bx_0^2$ & $0$ \\
$\so(-2)$ & $0$ & $-\lambda x_0x_1x_2$ & $0$ & $x_0^2x_1$ \\
$\so(-2)$ & $0$ & $x_2^2+\lambda x_0^2x_2$ & $-x_3$ & $-x_0^3$ \\
$\so(-3)$ & $0$ & $-x_1x_2$ & $0$ & $0$ \\
$\so(-4)$ & $0$ & $0$ & $x_1$ & $0$ \\
$\so(-2)$ & $x_0x_1q_4$ & $-x_0x_1(x_2+bx_0^2)$ & $0$ & $0$ \\
$\so(-2)$ & $x_3^2-x_2^3-x_0^2q_4$ & $x_0^2(x_2+bx_0^2)$ & $0$ & $0$ \\
$\so(-3)$ & $x_1x_2^2$ & $0$ & $0$ & $0$ \\
$\so(-4)$ & $-x_1x_3$ & $0$ & $0$ & $0$ \\
\multicolumn{2}{c}{}
\end{tabular}

and the last $8$ rows and columns (those giving the component of
$\tilde{\alpha}$ from $(\sko^{-3})^2$ to $\sko^{-1}(-1)^2$, which corresponds
to the component of $\alpha$ from $(\sd^2)^2$ to $\sd^1(-1)^2$) are

\begin{tabular}{r|cccccccc}
\multicolumn{2}{c}{} \\
& $\so(-6)$ & $\so(-6)$ & $\so(-5)$ & $\so(-4)$ & $\so(-6)$ & $\so(-6)$ &
$\so(-5)$ & $\so(-4)$ \\
\hline
$\so(-2)$ & $0$ & $0$ & $0$ & $0$ & $0$ & $0$ & $x_3$ & $-x_2$ \\
$\so(-2)$ & $0$ & $0$ & $0$ & $0$ & $0$ & $0$ & $0$ & $0$ \\
$\so(-3)$ & $0$ & $0$ & $0$ & $0$ & $-x_3$ & $0$ & $0$ & $x_0$ \\
$\so(-4)$ & $0$ & $0$ & $0$ & $0$ & $x_2$ & $0$ & $-x_0$ & $0$ \\
$\so(-2)$ & $0$ & $0$ & $-x_3$ & $x_2$ & $0$ & $0$ & $0$ & $0$ \\
$\so(-2)$ & $0$ & $0$ & $0$ & $0$ & $0$ & $0$ & $0$ & $0$ \\
$\so(-3)$ & $x_3$ & $0$ & $0$ & $-x_0$ & $0$ & $0$ & $0$ & $0$ \\
$\so(-4)$ & $-x_2$ & $0$ & $x_0$ & $0$ & $0$ & $0$ & $0$ & $0$ \\
\multicolumn{2}{c}{}
\end{tabular}
(of course, by symmetry we need not write the first $4$ rows entirely).

    \begin{rema}
It is easy to see that the component of $\alpha$ from $(\sd^2)^2$ to
$\sd^1(-1)^2$ is as expected from \ref{natmor} (observe that $\p_1=\langle
x_0,x_1\rangle$ and that the multiplication map $x_i:H^1(\gP,\gph_*\so_{\gS})
\to H^1(\gP,\gph_*\so_{\gS}(1))$ is $0$ for $i=0$ and is an isomorphism for
$i=1$).
    \end{rema}

    \begin{rema}
The fact that $\tilde{\alpha}=j\comp\alpha\comp j\tdual$ for some (uniquely
determined) $\alpha$ is equivalent to the fact that, denoting by
$\tilde{\alpha}^{(i)}$ the i$^{th}$ row of $\tilde{\alpha}$, we have (as it
is immediate to check) $\sum_{i=0}^3x_i\tilde{\alpha}^{(5+i)}=0$ and
$\sum_{i=0}^3x_i\tilde{\alpha}^{(9+i)}=0$.
    \end{rema}

    \begin{rema}
We checked that, with $\tilde{\alpha}$ defined as above, $\alpha$ actually
satisfies the rank condition and its determinant is an irreducible polynomial
(of degree $24$).
    \end{rema}

                        \begin{appendix}

                        \chapter{Abelian categories and derived categories}

In \ref{abelquot} we recall the definition and some basic properties about the
quotient of an abelian category by a localizing subcategory (see \cite{G} for
more details and proofs). In the rest of the appendix, besides proving
a couple of results for which we could not find a reference, we give (without
proofs) the essential elements of the theory of derived categories. A good
reference for this subject is \cite[chapter I]{H2} (where almost all proofs can
be found); see also \cite{V} for more details on triangulated categories and
derived categories. We refer to \cite{Bon} for the part on exceptional
sequences and mutations in a triangulated category.

                \section{Quotient of an abelian category}\label{abelquot}

        \begin{defi}
Let $\cat{A}$ be an abelian category. A full subcategory $\cat{C}$ of $\cat{A}$
is called {\em thick} if for every short exact sequence $0\to M'\to M\to M''\to
0$ in $\cat{A}$, $M$ is an object of $\cat{C}$ if and only if $M'$ and $M''$
are objects of $\cat{C}$.
        \end{defi}

Given an abelian category $\cat{A}$ and a thick subcategory $\cat{C}$ of
$\cat{A}$, it is possible to define the quotient category $\cat{A}/\cat{C}$,
whose objects are the same as those of $\cat{A}$, and whose morphisms are
defined in the following way. Given $M,N$ objects of $\cat{A}$, the set
\[I(M,N):=\{(M',N')\st M'\subseteq M;N'\subseteq N;M/M'\in\cat{C};N'\in
\cat{C}\}\]
is partially ordered by $(M',N')\le(M'',N'')$ if and only if $M''\subseteq M'$
and $N'\subseteq N''$. It is easy to prove that $(I(M,N),\le)$ is a direct set
(since $\cat{C}$ is thick), and so we can define
\[\Hom_{\cat{A}/\cat{C}}(M,N):=\lim_{(M',N')\in I(M,N)}\Hom_{\cat{A}}(M',N/N')
.\]
Again using the fact that $\cat{C}$ is thick, it is not difficult to see that
there is a natural way of defining composition of morphisms, so that
$\cat{A}/\cat{C}$ is indeed a category.

        \begin{prop}(\cite[prop. 1, p. 367]{G})
If $\cat{C}$ is a thick subcategory of an abelian category $\cat{A}$, then
$\cat{A}/\cat{C}$ is an abelian category and the natural functor $T:\cat{A}\to
\cat{A}/\cat{C}$ is exact.
        \end{prop}

        \begin{defi}
A thick subcategory $\cat{C}$ of an abelian category $\cat{A}$ is called a {\em
right localizing} subcategory\index{localizing subcategory} of $\cat{A}$ if
there exists a functor $S: \cat{A}/\cat{C}\to\cat{A}$ which is right adjoint of
$T:\cat{A}\to \cat{A}/\cat{C}$.
        \end{defi}

        \begin{rema}
Of course, one can define (in the obvious way) also the notion of left
localizing subcategory.
        \end{rema}

Clearly if $S$ exists it is unique up to isomorphism of functors, and it is
left exact (since it is a right adjoint). Moreover, we have the following
results (\cite[prop. 3, p. 371 and prop. 5, p. 374]{G}).

        \begin{prop}
The natural transformation (induced by the properties of adjoint functors)
$T\comp S\to\id_{\cat{A}/\cat{C}}$ is an isomorphism of functors.
        \end{prop}

        \begin{prop}\label{locsubcat}
Let $F:\cat{A}\to\cat{B}$ be an exact functor between abelian categories, let
$G:\cat{B}\to\cat{A}$ be a right adjoint of $F$, and assume that the associated
natural transformation $F\comp G\to\id_{\cat{B}}$ is an isomorphism of
functors. Then, denoting by $\ker F$ the full subcategory of $\cat{A}$ whose
objects are the $A$ such that $F(A)\iso0$ and by $T:\cat{A}\to\cat{A}/\ker F$
the natural functor, $\ker F$ is a right localizing subcategory of $\cat{A}$
and $F$ induces an exact equivalence of categories $F':\cat{A}/\ker F\to
\cat{B}$, whose quasi--inverse is $T\comp G$.
        \end{prop}

        \section{Triangulated categories}\label{trcat}

A triangulated category\index{triangulated category|(} is given by an additive
category $\cat{C}$ with some additional structures. First there is an additive
automorphism $T:\cat{C}\to \cat{C}$, called {\em translation} (or {\em shift})
functor; $\all n\in\Z$ the functor $T^n$ will be usually denoted by $[n]$. A
{\em triangle}\index{triangle} in $(\cat{C},T)$ is given by three objects $X$,
$Y$, $Z$ and by three morphisms $u:X\to Y$, $v:Y\to Z$, $w:Z\to X[1]$ of
$\cat{C}$, and will be usually denoted by $\tri{X}{Y}{Z}{u}{v}{w}$. A morphism
between two triangles\index{morphism of!triangles} $\tri{X}{Y}{Z}{u}{v}{w}$ and
$\tri{X'}{Y'}{Z'}{u'}{v'}{w'}$ is given by three morphisms $f:X\to X'$, $g:Y\to
Y'$ and $h:Z\to Z'$ such that
\[\begin{CD}
X @>u>> Y @>v>> Z @>w>> X[1] \\
@VVfV @VVgV @VVhV @VVf[1]V \\
X' @>u'>> Y' @>v'>> Z' @>w'>> X'[1]
\end{CD}\]
is a commutative diagram. Then a structure of triangulated category on
$(\cat{C},T)$ is given by a family of triangles (called {\em distinguished}
triangles)\index{distinguished triangle} which have to satisfy the following
four axioms.
\begin{description}

\item [(TR1)] A triangle isomorphic to a distinguished triangle is
distinguished.

Every morphism $u:X\to Y$ of $\cat{C}$ is contained in a distinguished triangle
$\tri{X}{Y}{Z}{u}{v}{w}$.

For every object $X$ of $\cat{C}$ the triangle $\tri{X}{X}{0}{\id_X}{}{}$ is
distinguished.

\item [(TR2)] The triangle $\tri{X}{Y}{Z}{u}{v}{w}$ is distinguished if and
only if the triangle $\tri{Y}{Z}{X[1]}{v}{w}{-u[1]}$ is distinguished.

\item [(TR3)] Given two distinguished triangles $\tri{X}{Y}{Z}{u}{v}{w}$,
$\tri{X'}{Y'}{Z'}{u'}{v'}{w'}$ and two morphisms $f:X\to X'$, $g:Y\to Y'$ such
that $u'\comp f=g\comp u$, there exists a morphism $h:Z\to Z'$ such that
$(f,g,h)$ is a morphism of triangles.

\item [(TR4)] Given three distinguished triangles
\begin{gather*}
\tri{X}{Y}{Z'}{u}{u'}{u''},\\
\tri{Y}{Z}{X'}{v}{v'}{v''},\\
\tri{X}{Z}{Y'}{v\comp u}{w'}{w''},
\end{gather*}
there exist two morphisms $f:Z'\to Y'$ and $g:Y'\to X'$, such that
$(\id_X,v,f)$ and $(u,\id_Z,g)$ are morphisms of triangles
\[\begin{CD}
X @>u>> Y @>u'>> Z' @>u''>> X[1] \\
@VV\id_XV @VVvV @VVfV @VV\id_{X[1]}V \\
X @>v\comp u>> Z @>w'>> Y' @>w''>> X[1]\\
@VVuV @VV\id_ZV @VVgV @VVu[1]V \\
Y @>v>> Z @>v'>> X' @>v''>> Y[1]
\end{CD}\]
and  $\tri{Z'}{Y'}{X'}{f}{g}{u'[1]\comp v''}$ is a distinguished triangle.
\index{triangulated category|)}
\end{description}

A functor between two triangulated categories is said to be {\em
exact}\index{exact functor} if it commutes (up to isomorphism) with translation
functors and if it carries distinguished triangles to distinguished triangles.
A functor $H:\cat{C}\to \cat{A}$ between a triangulated category and an abelian
category is called {\em cohomological} if for
every distinguished triangle $\tri{X}{Y}{Z}{u}{v}{w}$ of $\cat{C}$ the sequence
\[\cdots\to H(X)\mor{H(u)}H(Y)\mor{H(v)}H(Z)\mor{H(w)}H(X[1])\mor{H(u[1])}
H(Y[1])\to\cdots\]
is exact in $\cat{A}$.

    \begin{prop}\label{trcatprop}
A triangulated category $\cat{C}$ satisfies the  following properties:
\begin{enumerate}

\item the composition of two consecutive morphisms in a distinguished triangle
is zero;

\item for every object $M$ of $\cat{C}$, $\Hom_{\cat{C}}(M,-):\cat{C}\to
\cat{Ab}$ and $\Hom_{\cat{C}}(-,M):\cat{C}\opp\to\cat{Ab}$ are cohomological
functors;\footnote{$\cat{C}\opp$ has a natural structure of triangulated
category. $\cat{Ab}$ is the category of abelian groups.}

\item for every morphism $u:X\to Y$ in $\cat{C}$ a distinguished triangle of
the form $\tri{X}{Y}{Z}{u}{v}{w}$ (which exists by $\mathbf{(TR1)}$) is not
unique, but the object $Z$ is uniquely determined up to isomorphism; moreover,
$u$ is an isomorphism if and only if $Z\iso0$.
\end{enumerate}
    \end{prop}

For the rest of this section we will assume that $\K$ is a field and
that $\cat{C}$ is a $\K$--linear triangulated category (namely, a
triangulated category in which $\Hom_{\cat{C}}(X,Y)$ has a structure
of $\K$--vector space for every objects $X$ and $Y$ of $\cat{C}$
such that composition of morphisms is $\K$--bilinear and the shift
functor is $\K$--linear).

    \begin{defi}\label{exc}
An object $E$ of $\cat{C}$ is called {\em exceptional}\index{exceptional} if
$\Hom_{\cat{C}}(E,E)\iso\K$ and $\Hom_{\cat{C}}(E,E[k])=0$ for
$k\ne0$.

A sequence $(E_0,\dots,E_m)$ of objects of $\cat{C}$ is called an
{\em exceptional sequence}\index{exceptional!sequence} (or {\em
collection}) if each $E_i$ is exceptional and
$\Hom_{\cat{C}}(E_i,E_j[k])=0$ for $i>j$ and $\all k\in\Z$.

An exceptional sequence  $(E_0,\dots,E_m)$ is {\em
strong}\index{strong exceptional sequence} if
$\Hom_{\cat{C}}(E_i,E_j[k])=0$ for $k\ne0$ and $\all i,j$.

An exceptional sequence  $(E_0,\dots,E_m)$ is {\em
full}\index{full exceptional sequence} if $\{E_0,\dots,E_m\}$
generates $\cat{C}$ as a triangulated category (i.e., if every
strictly full triangulated subcategory of $\cat{C}$ containing the
$E_i$ coincides with $\cat{C}$).
    \end{defi}

    \begin{rema}\label{shiftexc}
Given $k_0,\dots,k_m\in\Z$, a sequence $(E_0,\dots,E_m)$ is
exceptional (and full) if and only if $(E_0[k_0],\dots,E_m[k_m])$ is
exceptional (and full). On the other hand, the property of being
strong is not preserved in general, unless $k_0=\cdots=k_m$.
    \end{rema}

We will assume also that
$\dim_{\K}\bigoplus_{k\in\Z}\Hom_{\cat{C}}(X,Y[k])<\infty$ for every
objects $X,Y$ of $\cat{C}$ and we will denote by $L_XY$ the object of
$\cat{C}$ defined (up to isomorphism) by the distinguished triangle
\[L_XY\to\bigoplus_{k\in\Z}(\Hom_{\cat{C}}(X[k],Y)\otimes_{\K}X[k])\mor{u}Y
\to(L_XY)[1]\]
(where $u$ is the natural morphism).

    \begin{prop}
If $\sigma=(E_0,\dots,E_m)$ is a (full) exceptional sequence of $\cat{C}$, then
for $i=0,\dots m-1$ the sequence
$(E_0,\dots,E_{i-1},L_{E_i}E_{i+1},E_i,E_{i+2},\dots,E_m)$ (called a {\em left
mutation}\index{left!mutation} of $\sigma$)\footnote{There is also a similar
notion of right mutation, which we will not consider.} is (full) exceptional,
too.
    \end{prop}

Given a sequence $(E_0,\dots,E_m)$, for $0\le j\le i\le m$ we define
inductively objects $L^{(j)}E_i$ by $L^{(0)}E_i:=E_i$ and
$L^{(j)}E_i:=L_{E_{i-j}}(L^{(j-1)}E_i)$ for $j>0$.

    \begin{coro}\label{dualmut}
If $\sigma=(E_0,\dots,E_m)$ is a (full) exceptional sequence of $\cat{C}$, then
the sequence $(L^{(m)}E_m,L^{(m-1)}E_{m-1},\dots,L^{(1)}E_1,L^{(0)}E_0=E_0)$
(called the {\em left dual}\index{left!dual} of $\sigma$) is (full)
exceptional, too.
    \end{coro}

                \section{Derived categories}\label{dercat}

If $\cat{A}$ is an abelian category or, more generally, an additive subcategory
of an abelian category, $C(\cat{A})$,\index{C(A)@$C(\cat{A})$|(}
$K(\cat{A})$\index{K(A)@$K(\cat{A})$|(} and
$D(\cat{A})$\index{D(A)@$D(\cat{A})$|(} will denote the category of complexes,
the homotopy category and the derived category of $\cat{A}$, respectively: we
recall the definitions.

The objects of $C(\cat{A})$ are just the complexes of $\cat{A}$ (by a complex
we always mean a cohomological complex, that is differentials have degree $+1$)
and the morphisms are just usual morphisms of
complexes.\index{C(A)@$C(\cat{A})$|)} A complex
\[\cdots\to X^{-1}\mor{\cdiff{X}^{-1}}X^{0}\mor{\cdiff{X}^{0}}X^{1}\to\cdots\]
will be denoted by $\cp{X}$; similarly, a morphism of complexes from $\cp{X}$
to $\cp{Y}$ given by morphisms $f^i:X^i\to Y^i$ (with $\cdiff{Y}^{i+1}\comp
f^i=f^{i+1}\comp\cdiff{X}^i$) will be denoted by $\cp{f}$.\index{$\cp{}$}

The objects of $K(\cat{A})$ are the same as those of $C(\cat{A})$, whereas the
morphisms are the homotopy equivalence classes of morphisms of complexes (two
morphisms $\cp{f},\cp{g}:\cp{X}\to\cp{Y}$ in $C(\cat{A})$ are homotopic if
there exist morphisms $k^{i}:X^i\to Y^{i-1}$ such that $f^i-g^i=\cdiff{Y}^{i-1}
\comp k^i+k^{i+1}\comp\cdiff{X}^i$).\index{K(A)@$K(\cat{A})$|)} By abuse of
notation we will denote in the same way a morphism of $C(\cat{A})$ and its
equivalence class in $K(\cat{A})$. Note that if $\cp{f}:\cp{X}\to\cp{Y}$ is a
morphism of $K(\cat{A})$, then the morphisms in cohomology
$H^i(\cp{f}):H^i(\cp{X})\to H^i(\cp{Y})$ are well defined: in particular,
$\cp{f}$ is said to be a {\em quasi--isomorphism}\index{quasi--isomorphism} if
$H^i(\cp{f})$ is an isomorphism $\all i\in\Z$.

The objects of $D(\cat{A})$ are always the complexes of $\cat{A}$, but the
morphisms are obtained from those of $K(\cat{A})$ by inverting formally the
quasi--isomorphisms: more precisely, a morphism in $D(\cat{A})$ from $\cp{X}$
to $\cp{Y}$ is represented by a couple of morphisms (of $K(\cat{A})$)
$(\cp{s}:\cp{Z}\to\cp{X},\cp{f}:\cp{Z}\to\cp{Y})$ where $\cp{s}$ is a
quasi--isomorphism and $\cp{Z}$ is another object of $D(\cat{A})$; on these
couples one introduces the equivalence relation generated by the direct
relation which identifies two couples $(\cp{s}:\cp{Z}\to\cp{X},\cp{f}:\cp{Z}\to
\cp{Y})$ and $(\cp{t}:\cp{W}\to\cp{X},\cp{g}:\cp{W}\to\cp{Y})$ if there exists
a morphism $\cp{a}:\cp{Z}\to\cp{W}$ such that $\cp{f}=\cp{g}\comp\cp{a}$ and
$\cp{s}=\cp{t}\comp\cp{a}$.\index{D(A)@$D(\cat{A})$|)} The equivalence class of
a couple $(\cp{s},\cp{f})$ will be denoted by $[(\cp{s},\cp{f})]$ (or by
$[(\cp{s},\cp{f})]_{\cat{A}}$ if there can be doubt about $\cat{A}$); if
$\cp{s}=\id_{\cp{X}}$, we will simply write $[\cp{f}]$ instead of
$[(\id_{\cp{X}},\cp{f})]$.

Clearly $C(\cat{A})$, $K(\cat{A})$ and $D(\cat{A})$ are all additive
categories, but, while $C(\cat{A})$ is abelian if $\cat{A}$ is abelian, this is
not the case for $K(\cat{A})$ and $D(\cat{A})$. On the other hand, $K(\cat{A})$
and $D(\cat{A})$ have a natural structure of triangulated categories, which we
are now going to describe briefly.

$\all n\in\Z$ the shift functor $[n]:C(\cat{A})\to C(\cat{A})$ is defined on
objects by
\begin{align*}
& X[n]^i:=X^{n+i}, & \diff{\cp{X}[n]}^i:=(-1)^n\cdiff{X}^{n+i}
\end{align*}
and on morphisms by $f[n]^i:=f^{n+i}$. Since the shift functors preserve
homotopic morphisms and quasi--isomorphisms, they induce shift functors
(denoted again by $[n]$) also on $K(\cat{A})$ and $D(\cat{A})$. Clearly in each
case $[n]$ is an additive automorphism and $[n]\comp[m]=[n+m]$ $\all n,m\in\Z$.
In order to be able to say what are the distinguished triangles of $K(\cat{A})$
and $D(\cat{A})$ we need to recall the definition of mapping cone of a morphism
of complexes.

            \begin{defi}
The {\em mapping cone}\index{mapping cone} of a morphism of complexes
$\cp{f}:\cp{X}\to\cp{Y}$ is the complex $\mc{f}$\index{MC(f)@$\mc{f}$} defined
by
\begin{align*}
& \mc{f}^n:=X^{n+1}\oplus Y^n, &
\diff{\mc{f}}^n:=\begin{pmatrix}
-\cdiff{X}^{n+1} & 0 \\
f^{n+1} & \cdiff{Y}^n
\end{pmatrix}.
\end{align*}
There are natural morphisms of complexes (inclusion and projection)
\begin{align*}
& \inc{f}:\cp{Y}\to\mc{f}, & \pro{f}:\mc{f}\to\cp{X}[1].
\end{align*}
            \end{defi}\index{j(f)@$\inc{f}$}\index{p(f)@$\pro{f}$}

Then a triangle in $K(\cat{A})$ is distinguished if it is isomorphic to a
triangle of the form $\tri{X}{Y}{\mc{f}}{\cp{f}}{\inc{f}}{\pro{f}}$ for some
morphism of complexes $\cp{f}$. A triangle in $D(\cat{A})$ is distinguished if
it is isomorphic to the image of a distinguished triangle of $K(\cat{A})$
under the natural functor $Q:K(\cat{A})\to D(\cat{A})$ (which is therefore
exact by definition). The distinguished triangles of $D(\cat{A})$ can be
characterized in the following way.

                \begin{prop}\label{derdist}
If $0\to\cp{X}\mor{\cp{f}}\cp{Y}\mor{\cp{g}}\cp{Z}\to0$ is a short exact
sequence of complexes, then there exists a morphism $\alpha:\cp{Z}\to\cp{X}[1]$
in $D(\cat{A})$ such that $\tri{\cp{X}}{\cp{Y}}{\cp{Z}}{[\cp{f}]}{[\cp{g}]}
{\alpha}$ is a distinguished triangle. Conversely, every distinguished triangle
in $D(\cat{A})$ is isomorphic to one of this form.
                \end{prop}

Given a complex $\cp{X}$ of $\cat{A}$ and $n\in\Z$, we will denote by
$X^{\le n}$ or $X^{<n+1}$ the complex defined by
\begin{align*}
& (X^{\le n})^i:=\begin{cases}
X^i & \text{if $i\le n$} \\
0 & \text{if $i>n$}
\end{cases}, &
\diff{X^{\le n}}^i:=\begin{cases}
\cdiff{X}^i & \text{if $i<n$} \\
0 & \text{if $i\ge n$}
\end{cases}.
\end{align*}
In a similar way we define the complex $X^{\ge n}=X^{>n-1}$. Then $\all n\in\Z$
there is a short exact sequence of complexes $0\to X^{\ge n}\to\cp{X}\to X^{<n}
\to0$, which gives a distinguished triangle not only in $D(\cat{A})$ (as we
know from \ref{derdist}) but even in $K(\cat{A})$. Indeed, defining
$\cp{\tbo{\cp{X}}{n}}:X^{<n}[-1]\to X^{\ge
n}$\index{delta_X,n@$\cp{\tbo{\cp{X}}{n}}$} by
\[\tbo{\cp{X}}{n}^i:=\begin{cases}
\cdiff{X}^{n-1} & \text{if $i=n$} \\
0 & \text{if $i\ne n$}
\end{cases}\]
it is very easy to prove the following result.

        \begin{lemm}\label{homotdist}
$\all\cp{X}\in K(\cat{A})$ and $\all n\in\Z$ the triangle $X^{<n}[-1]
\mor{\cp{\tbo{\cp{X}}{n}}}X^{\ge n}\mono\cp{X}\epi X^{<n}$ is distinguished in
$K(\cat{A})$.
        \end{lemm}

An object $X$ of $\cat{A}$ will be always considered as a complex concentrated
in position $0$. Obviously this defines functors $\cat{A}\to C(\cat{A})$,
$\cat{A}\to K(\cat{A})$ and $\cat{A}\to D(\cat{A})$, the first two of which are
clearly fully faithful. We have besides the following result.

        \begin{lemm}\label{AfullDA}
The natural functor $\cat{A}\to D(\cat{A})$ is fully faithful (i.e., $\cat{A}$
can be identified with a full subcategory of $D(\cat{A})$).
        \end{lemm}

$C^+(\cat{A})$,\index{C^+(A)@$C^+(\cat{A})$, $C^-(\cat{A})$, $C^b(\cat{A})$}
$K^+(\cat{A})$\index{K^+(A)@$K^+(\cat{A})$, $K^-(\cat{A})$, $K^b(\cat{A})$} and
$D^+(\cat{A})$\index{D^+(A)@$D^+(\cat{A})$, $D^-(\cat{A})$, $D^b(\cat{A})$}
will denote the full subcategories of $C(\cat{A})$, $K(\cat{A})$ and
$D(\cat{A})$ having as objects the complexes which are bounded from below (i.e.
complexes $\cp{X}$ such that $X^{i}=0$ for $i<<0$). $K^+(\cat{A})$ and
$D^+(\cat{A})$ are triangulated subcategories of $K(\cat{A})$ and $D(\cat{A})$,
and $D^+(\cat{A})$ is equivalent to the the full subcategory of
$D(\cat{A})$ having as objects the complexes $\cp{X}$ such that
$H^i(\cp{X})=0$ for $i<<0$. Similarly, $C^-(\cat{A})$, $K^-(\cat{A})$,
$D^-(\cat{A})$ and $C^b(\cat{A})$, $K^b(\cat{A})$, $D^b(\cat{A})$ will denote
the full subcategories consisting, respectively, of complexes which are bounded
from above and bounded in both directions.

        \section{Derived functors}\label{derfun}

Let $\cat{A}$ and $\cat{B}$ be abelian categories, and let $K^*(\cat{A})$ be
one of the homotopy categories $K(\cat{A})$, $K^+(\cat{A})$, $K^-(\cat{A})$ or
$K^b(\cat{A})$. If $F:K^*(\cat{A})\to K(\cat{B})$ is an exact functor (the only
case we are actually interested in is when $F$ is induced by an additive
functor $\cat{A}\to \cat{B}$, also denoted by $F$ by abuse of notation), then
it is not true in general that $F$ induces a functor $D^*(\cat{A})\to
D(\cat{B})$ (this happens if $F$ preserves quasi--isomorphisms, e.g. if $F:
\cat{A}\to \cat{B}$ is exact). In any case, one can define the right derived
functor\index{derived functor} of $F$ in the following way. It consists of an
exact functor $RF: D^*(\cat{A})\to D(\cat{B})$\index{RF@$RF$} together with a
morphism of functors $\xi:Q\comp F \to RF\comp Q$, satisfying the following
universal property: if $G:D^*(\cat{A}) \to D(\cat{B})$ is an exact functor and
$\zeta:Q\comp F\to G\comp Q$ is a morphism of functors, then there exists a
unique morphism $\eta:RF\to G$ such that $\zeta=(\eta\comp Q)\comp\xi$. The
analogous definition of left derived functor $LF$\index{LF@$LF$} can be
obtained just inverting the arrows $\xi$, $\zeta$ and $\eta$. It is clear from
the definition that derived functors are unique up to isomorphism, if they
exist.

A sufficient condition for the existence of $RF:D^+(\cat{A})\to D(\cat{B})$
(respectively $LF:D^-(\cat{A})\to D(\cat{B})$) is that $\cat{A}$ has enough
injectives (respectively projectives). If this is the case, then $RF$ can be
computed as follows: given $\cp{X}\in K^+(\cat{A})$, there exists a
quasi--isomorphism $\cp{X}\to\cp{I}$ in $K^+(\cat{A})$ such that each $I^j$ is
injective, and then one can define $RF(\cp{X}):=F(\cp{I})$. Setting
$R^iF(\cp{X}):=H^i(RF(\cp{X}))$, we see that, in particular, the functors
$R^iF$ restricted to $\cat{A}\subset D(\cat{A})$ are the usual right derived
functors of $F:\cat{A}\to\cat{B}$ (in case $F$ is left exact).

If $\cat{A}$ is an abelian category, then $\Hom_{\cat{A}}(-,-):\cat{A}\opp
\times\cat{A}\to\cat{Ab}$ induces a functor, again denoted by
\[\Hom_{\cat{A}}(-,-):K(\cat{A})\opp\times K(\cat{A})\to K(\cat{Ab}).\]

    \begin{prop}\label{Yoneda}
If $\cat{A}$ has enough injectives, then
\[R\Hom_{\cat{A}}(-,-):D(\cat{A})\opp\times D^+(\cat{A})\to D(\cat{Ab})\]
exists and $R^i\Hom_{\cat{A}}(-,-)\iso\Hom_{D(\cat{A})}(-,-[i])$ $\all i\in\Z$.
In particular, if $X,Y\in\cat{A}$, then $\Ext^i_{\cat{A}}(X,Y)\iso
\Hom_{D(\cat{A})}(X,Y[i])$.
    \end{prop}

        \section{Some related results}

        \begin{lemm}\label{BqisoA}
Let $\cat{A}$ be an abelian category and $\cat{B}$ a full abelian subcategory
of $\cat{A}$. Assume that for every surjective morphism $f:A\to B$ in $\cat{A}$
with $B\in\cat{B}$ there exists a morphism $g:C\to A$ in $\cat{A}$ such that
$C\in\cat{B}$ and $f\comp g$ is surjective. Then for every $\cp{X}\in
K^-(\cat{A})$ such that $H^i(\cp{X})\in\cat{B}$ $\all i\in\Z$, there exists a
quasi--isomorphism $\cp{s}:\cp{Y}\to\cp{X}$ in $K^-(\cat{A})$ such that $\cp{Y}
\in K^-(\cat{B})$.
        \end{lemm}

        \begin{proof}
It is not difficult to adapt to this more general setting the proof of the
first part of \cite[III, lemma 12.3]{H1}, where the particular case
$\cat{A}=\mo{A}$ and $\cat{B}=\fmo{A}$ ($A$ a noetherian ring) is considered.
        \end{proof}

        \begin{coro}\label{HomDA=HomDB}
Let $\cat{A}$ and $\cat{B}$ be abelian categories satisfying the same
hypotheses as in \ref{BqisoA}. Then $\all\cp{X},\cp{Y}\in D^-(\cat{B})$ the
natural map
\[\varphi:\Hom_{D^-(\cat{B})}(\cp{X},\cp{Y})\to
\Hom_{D^-(\cat{A})}(\cp{X},\cp{Y})\]
is an isomorphism. In other words, $D^-(\cat{B})$ can be identified with a full
subcategory of $D^-(\cat{A})$.
        \end{coro}

        \begin{proof}
$\varphi$ is injective: assume that
\[[(\cp{s}:\cp{Z}\to\cp{X},
\cp{f}:\cp{Z}\to\cp{Y})]_{\cat{B}}\in\Hom_{D^-(\cat{B})}(\cp{X},\cp{Y})\]
(with $\cp{s}$ quasi--isomorphism in $K^-(\cat{B})$) is such that
$\varphi([(\cp{s},\cp{f})]_{\cat{B}})=[(\cp{s},\cp{f})]_{\cat{A}}=0$. By
definition this means that there exists a quasi--isomorphism $\cp{t}:\cp{W}\to
\cp{Z}$ in $K^-(\cat{A})$ such that $\cp{f}\comp\cp{t}=0$ in $K^-(\cat{A})$.
By \ref{BqisoA} there exists a quasi--isomorphism $\cp{u}:\cp{V}\to\cp{W}$ with
$\cp{V}\in K^-(\cat{B})$, and then $\cp{f}\comp(\cp{t}\comp\cp{u})=0$ in
$K^-(\cat{B})$, i.e. $[(\cp{s},\cp{f})]_{\cat{B}}=0$.

$\varphi$ is surjective: given
\[[(\cp{s}:\cp{Z}\to\cp{X},
\cp{f}:\cp{Z}\to\cp{Y})]_{\cat{A}}\in\Hom_{D^-(\cat{A})}(\cp{X},\cp{Y})\]
(with
$\cp{s}$ quasi--isomorphism in $K^-(\cat{A})$), by \ref{BqisoA} there exists a
quasi--isomorphism $\cp{t}:\cp{W}\to\cp{Z}$ with $\cp{W}\in K^-(\cat{B})$, and
then
\[[(\cp{s},\cp{f})]_{\cat{A}}=[(\cp{s}\comp\cp{t},\cp{f}\comp\cp{t})]_{\cat{A}}
=\varphi([(\cp{s}\comp\cp{t},\cp{f}\comp\cp{t})]_{\cat{B}}).\]
        \end{proof}

        \begin{prop}\label{HomKA=HomDA}
Let $\cat{A}$ be an abelian category with enough injectives and let
$\cp{X},\cp{Y}\in K^b(\cat{A})$ be such that $\Ext_{\cat{A}}^i(X^p,Y^q)=0$ for
$0<i\le p-q$ and $\all p,q\in\Z$. Then there is a natural isomorphism
\[\Hom_{K(\cat{A})}(\cp{X},\cp{Y})\isomor\Hom_{D(\cat{A})}(\cp{X},\cp{Y}).\]
        \end{prop}

        \begin{proof}
First we prove that $\all q\in\Z$ the natural map
\[\Hom_{K(\cat{A})}(\cp{X},Y^q[t-q])\to
\Hom_{D(\cat{A})}(\cp{X},Y^q[t-q])\]
is an isomorphism if $t\le0$ and is injective if $t=1$. Fixing $q\in\Z$, we
denote by $K'_t$ (respectively $D'_t$) the (contravariant) functor
$\Hom_{K(\cat{A})}(-,Y^q[t-q])$ (respectively $\Hom_{D(\cat{A})}(-,Y^q[t-q])$)
and by $\epsilon'_t$ the natural transformation $K'_t\to D'_t$. We have to
prove that $\epsilon'_t(\cp{X})$ is an isomorphism if $t\le0$ and is injective
if $t=1$: as usual (since $\cp{X}$ is bounded) it is enough to show that the
same is true for $\epsilon'_t(X^{\ge p})$ $\all p\in\Z$, and we can proceed by
descending induction on $p$, because this is trivially true for $p>>0$.
Applying $\epsilon'_t$ to the triangle $\tri{X^{>p}}{X^{\ge p}}{X^p[-p]}{}{}{}$
(which is distinguished in $K^b(\cat{A})$ by \ref{homotdist}), we obtain a
commutative diagram with exact rows
\[\begin{CD}
K'_t(X^{>p}[1]) @>>> K'_t(\cp{A}[1]) @>>> K'_t(X^{\ge p}) @>>> K'_t(X^{>p})
@>>> K'_t(\cp{A}) \\
@VV{\epsilon}'_t(X^{>p}[1])V @VV{\epsilon}'_t(\cp{A}[1])V
@VV{\epsilon}'_t(X^{\ge p})V @VV{\epsilon}'_t(X^{>p})V
@VV{\epsilon}'_t(\cp{A})V \\
D'_t(X^{>p}[1]) @>>> D'_t(\cp{A}[1]) @>>> D'_t(X^{\ge p}) @>>> D'_t(X^{>p})
@>>> D'_t(\cp{A})
\end{CD}\]
(where $\cp{A}:=X^p[-p-1]$). By the inductive hypothesis
$\epsilon'_t(X^{>p}[1])=\epsilon'_{t-1}(X^{>p})$ is an isomorphism for $t\le1$
and $\epsilon'_t(X^{>p})$ is an isomorphism for $t\le0$ and is injective for
$t=1$. Since $\all p'\in\Z$
\[K'_t(X^p[-p'])\iso\begin{cases}
\Hom_{\cat{A}}(X^p,Y^q) & \text{if $p'=q-t$} \\
0 & \text{if $p'\ne q-t$}
\end{cases}\]
and $D'_t(X^p[-p'])\iso\Ext_{\cat{A}}^{p'-q+t}(X^p,Y^q)$ (by \ref{Yoneda}),
the fact that $\Ext_{\cat{A}}^i(X^p,Y^q)=0$ for $0<i\le p-q$ implies that
$\epsilon'_t(\cp{A}[1])$ is an isomorphism for $t\le0$ and is injective for
$t=1$, and that $\epsilon'_t(\cp{A})$ is always injective. Then it follows
from the five--lemma that $\epsilon'_t(X^{\ge p})$ is an isomorphism if $t\le0$
and is injective if $t=1$.

Denoting by $K_t$ (respectively $D_t$) the functor
$\Hom_{K(\cat{A})}(\cp{X},-[t])$ (respectively
$\Hom_{D(\cat{A})}(\cp{X},-[t])$) and by $\epsilon_t$ the natural
transformation $K_t\to D_t$, we have to prove that $\epsilon_t(\cp{Y})$ is an
isomorphism for $t=0$, and we will prove more generally that it is an
isomorphism if $t\le0$. As before, it is enough to show by induction on $q$
that $\epsilon_t(Y^{\le q})$ is an isomorphism if $t\le0$. Applying
$\epsilon_t$ to the distinguished triangle $Y^{<q}[-1]\to Y^q[-q]\to Y^{\le q}
\to Y^{<q}$ we obtain a commutative diagram with exact rows
\[\begin{CD}
K_t(Y^{<q}[-1]) @>>> K_t(\cp{B}) @>>> K_t(Y^{\le q}) @>>> K_t(Y^{<q})
@>>> K_t(\cp{B}[1]) \\
@VV{\epsilon}_t(Y^{<q}[-1])V @VV{\epsilon}_t(\cp{B})V
@VV{\epsilon}_t(Y^{\le q})V @VV{\epsilon}_t(Y^{<q})V
@VV{\epsilon}_t(\cp{B}[1])V \\
D_t(Y^{<q}[-1]) @>>> D_t(\cp{B}) @>>> D_t(Y^{\le q}) @>>> D_t(Y^{<q})
@>>> D_t(\cp{B}[1])
\end{CD}\]
(where $\cp{B}:=Y^q[-q]$). Now, $\epsilon_t(Y^{<q})$ and
$\epsilon_t(Y^{<q}[-1])=\epsilon_{t-1}(Y^{<q})$ are isomorphisms (by induction)
for $t\le0$. Moreover, by what we have already proved, $\epsilon_t(\cp{B})=
\epsilon'_t(\cp{X})$ is an isomorphism for $t\le0$ and $\epsilon_t(\cp{B}[1-q])
=\epsilon'_{t+1}(\cp{X})$ is an isomorphism for $t<0$ and is injective for
$t=0$, whence $\epsilon_t(Y^{\le q})$ is an isomorphism if $t\le0$ by the
five--lemma.
        \end{proof}

    \begin{rema}
The above result is proved in \cite{AO}, under the stronger hypothesis
$\Ext_{\cat{A}}^i(X^p,Y^q)=0$ for $i>0$.
    \end{rema}

\end{appendix}

\backmatter

\printindex

                                \end{document}